\documentclass[aos]{imsart}
\makeatletter
\def\journal@name{}
\def\journal@url{}
\makeatother
\input{preamble.tex}
\startlocaldefs

\theoremstyle{plain}
\newtheorem{theorem}{Theorem}[subsection]
\renewcommand{\thetheorem}{%
  \ifnum\value{subsection}=0
    \thesection.\arabic{theorem}%
  \else
    \thesubsection.\arabic{theorem}%
  \fi}
\newtheorem{lemma}[theorem]{Lemma}
\newtheorem{proposition}[theorem]{Proposition}
\newtheorem{corollary}[theorem]{Corollary}

\theoremstyle{definition}
\newtheorem{definition}{Definition}[section]

\newtheorem{remark}{Remark}[section]

\newcommand{\cY}{\mathcal{Y}}
\newcommand{\cD}{\mathcal{D}}
\newcommand{\cB}{\mathcal{B}}
\newcommand{\cI}{\mathcal{I}}
\newcommand{\cH}{\mathcal{H}}

\newcommand{\cP}{\mathcal{P}}
\newcommand{\cZ}{\mathcal{Z}}
\newcommand{\cT}{\mathcal{T}}
\newcommand{\E}{\mathbb{E}}
\newcommand{\Var}{\operatorname{Var}}
\newcommand{\Cov}{\operatorname{Cov}}
\newcommand{\1}{\mathbf{1}}

\tikzset{
  pou/.style={circle,draw=black!60,fill=white,inner sep=0pt,minimum size=5pt},
  pourel/.style={circle,draw=black,fill=black!70,inner sep=0pt,minimum size=6pt},
  hedge/.style={line width=3.4pt,line cap=round,line join=round,opacity=0.30},
  hlab/.style={font=\scriptsize},
  axlab/.style={font=\scriptsize},
  netnode/.style={circle,draw=black!70,fill=black!5,inner sep=1pt,
                  minimum size=13pt,font=\scriptsize}
}

\newcommand{\cC}{\mathcal{C}}
\newcommand{\Rall}{R_{\mathrm{all}}}
\newcommand{\Runb}{R_{\mathrm{unb}}}
\newcommand{\Rallstar}{R^{*}_{\mathrm{all}}}
\newcommand{\Runbstar}{R^{*}_{\mathrm{unb}}}

\endlocaldefs

\newcommand{\restatedresultheading}{}
\newtheorem*{restatedresult}{\restatedresultheading}
\newenvironment{resultrestatement}[2]{%
  \renewcommand{\restatedresultheading}{#1}%
  \begin{restatedresult}[#2]}{\end{restatedresult}}

\graphicspath{{figures/}}
\begin{document}
\addtocontents{toc}{\protect\setcounter{tocdepth}{-1}}
\begin{frontmatter}

\title{Minimax designs and estimators for causal inference under network
interference}
\runtitle{Minimax designs and estimators under network interference}

\begin{aug}
\author[A]{\fnms{Vishesh}~\snm{Karwa}}
\author[A]{\fnms{Edoardo}~\snm{Airoldi}}
\address[A]{Temple University}
\runauthor{V. Karwa and E. Airoldi}
\end{aug}

\begin{abstract}
We develop a design-based minimax theory for linear causal estimands under interference. We study the minimax risks of three problems: optimizing an estimator under a fixed design, optimizing an estimator and a design from a restricted class, or optimizing both over all possible designs. In each problem, we study the minimax risk over the class of all estimators and design-unbiased estimators separately. %
We allow any correctly specified exposure mapping, including the identity.
Our minimax theory represents causal inference under interference as a  constrained sampling design on an expanded finite-population. For a given exposure mapping, the expanded finite population consists of $\sum_i d_i$ potential outcomes, where $d_i$ is the number of exposure levels attainable by unit $i$.  We describe the constrained sampling design by a weighted \emph{sampling hypergraph}:
its vertices are $\sum_i d_i$ potential outcomes, and its hyperedges are the maximal sets
of potential outcomes that can be jointly sampled under feasible treatment
assignments. The hyperedge weights correspond to the sampling design induced by the randomization design.
We show that the minimax risks of the fixed-design, restricted-design, and all-design problems depend, respectively, on the weighted hypergraph induced by the randomization design, the collection of weighted hypergraphs attainable within the restricted class of designs, and the unweighted sampling hypergraph.
Our Bayes lower bounds compare potential outcome tables that differ on jointly missing sets of potential outcomes. We show that joint-missingness probabilities govern the lower bounds over all estimators, whereas the corresponding odds govern the lower bounds over unbiased estimators. We illustrate the theory across a wide range of interference structures and estimands, including average treatment effects without interference, total, direct, and contagion effects under neighborhood interference, and bipartite dose-response contrasts. The upper and lower bounds match in several settings.

\end{abstract}



\end{frontmatter}

\section{Introduction}
\label{sec:introduction}

In finite-population design-based causal inference, we define an \emph{estimation strategy} as a pair consisting of a randomized design and an estimator.  This terminology follows the finite-population sampling literature %
\citep[see, e.g.,][]{Godambe1955,CasselSarndalWretman1977,NedyalkovaTille2008}.
Two estimation strategies can be compared through their design-based risks, defined as the expected mean squared error induced by a design--estimator pair. The expectation is taken over the randomization and the potential outcomes are held fixed.

Most existing works under interference evaluate a specified estimation strategy rather than determining the smallest risk attainable across strategies. For a given exposure mapping and estimand, one typically pairs a design with a Horvitz--Thompson or H\'ajek-type estimator, and bounds the design-based risk of the resulting estimation strategy, see for example, \citet{HudgensHalloran2008,AronowSamii2017,UganderYin2023}. In some cases, the design is also optimally chosen from a prescribed family to improve the precision of a fixed estimator \citep{BairdEtAl2018,JagadeesanEtAl2020,Leung2022,VivianoEtAl2023Causal,CandoganChenNiazadeh2024}.
Such an analysis establishes how the proposed strategy performs.  It cannot determine whether the design-based risk can be reduced further by choosing a different estimator under the same design, by jointly choosing a different design--estimator pair within some constraints on the design, or by choosing a design--estimator pair without any constraints. %
A general minimax theory covering these three problems across interference structures, exposure mappings and estimands is lacking. In this paper, we develop a design-based minimax analysis of these three problems for finite linear causal estimands under interference for an arbitrary, but correctly specified exposure mapping, including the identity.  We consider the worst-case risk over uniformly bounded potential outcome tables. In each problem, the minimax risk is defined as the smallest such worst-case design-based risk over the choices available to the statistician see Section \ref{sec:setup-strategy} for formal definitions.

The three problems correspond to situations that arise at the design and analysis stages of an experiment. At the analysis stage, the design is given, so the risk is minimized only over the choice of estimators. This gives the \emph{fixed-design} minimax problem. At the design stage, the risk can be minimized jointly over the design and estimator, but the design may be restricted to a prescribed family, for example, because of a treatment budget or because the interference graph may be unknown at the design stage. This gives the \emph{restricted-design} minimax problem. When the risk is minimized jointly over all designs and estimators, we obtain the \emph{all-design} minimax problem. Comparing the three problems shows whether a large risk is due to the given design, a design restriction, or a limit of the experiment.

For each problem, we study the minimax risk over two classes of estimators: the class of all estimators and the smaller class of design unbiased estimators.  %
We define the price of unbiasedness as the ratio of the minimax risk over unbiased estimators to the minimax risk over all estimators.  This ratio measures the multiplicative cost of requiring unbiasedness. The theory applies under no interference, partial interference, neighborhood
interference, and bipartite interference, and permits an arbitrary, but correctly
specified exposure mapping.  This includes the identity exposure mapping,
which retains the full treatment configuration relevant to a unit, as well as
coarser mappings based on treatment counts or thresholds.  The estimands are finite linear functionals of the resulting potential outcomes, including
one-arm means, two-arm contrasts, and policy means and policy contrasts.
Our examples include the average treatment effect without interference;
the total and direct treatment effects under partial and neighborhood
interference; a complex-contagion contrast under neighborhood interference;
and dose-response contrasts under bipartite interference.  The all-design problem was recently studied under neighborhood interference and identity exposure mapping for the all estimator class by \citet{KandirosHarshawSavje2026}; Section~\ref{subsec:related-khs} gives a detailed comparison.

\paragraph*{Overview of the results} At the center of our minimax theory is a constrained sampling design on an expanded finite population consisting of potential outcomes. This sampling design is induced by a randomization design, interference graph and an exposure mapping. In particular, the interference graph, exposure mapping, and assignment restrictions constrain which potential outcomes can be jointly observed, while the randomization design assigns probabilities to the feasible samples. We represent these constraints by a \emph{sampling hypergraph} and, for a given randomization design, the sampling probabilities by weights on its hyperedges, which we call the \emph{induced sampling design}. The nodes of the hypergraph correspond to potential outcomes. Each hyperedge records a maximal feasible sample, and its weight gives the probability of observing the corresponding sample under the randomization design. We describe these objects next.

Consider a finite population of $n$ units, indexed by $i\in[n]$ with a binary treatment assignment vector $z = (z_1, \ldots, z_n) \in \{0,1\}^n$. The unrestricted potential outcome table has $n$ rows and $2^n$ columns corresponding to the fixed but unknown potential outcomes $Y_i(z)$.  A correctly specified exposure mapping reduces this table to a \emph{ragged table}.  Row $i$ of this ragged table contains $d_i$ cells, each corresponding to a \emph{reduced} potential outcome $Y_i(e)$, one for each exposure level $e \in \mathcal D_i$, where $\mathcal D_i$ is the set of exposure levels attainable by unit $i$. Let $d_i = |\mathcal D_i| \geq 2$ be the number of exposure levels attainable by unit $i$.  Hence, the ragged table contains $|K| =\sum_{i=1}^nd_i$ potential outcome cells, each cell is indexed by the attainable unit--exposure pairs $(i,e)$. %
The expanded finite population is the set $K=\{(i,e):i\in[n],\ e\in\mathcal D_i\}$. Each element $(i,e)\in K$ carries the fixed potential-outcome value $Y_i(e)$. A single treatment assignment selects a sample of $n$ potential outcomes from the expanded finite population $K$ and thereby produces a missingness pattern on the ragged table of potential outcomes.  Within each row of the table, exactly one cell is observed and all remaining cells are missing. This is the usual within-unit counterfactual missingness, and corresponds to a sampling constraint that only one potential outcome per row can be observed. The missingness pattern also has an across-row component. In particular, the interference graph, exposure mapping, and restrictions, if any, on the set of treatment assignments together constrain which cells from different rows can be sampled together. %

\begin{figure}[t]
  \centering
  \includegraphics[width=\textwidth]{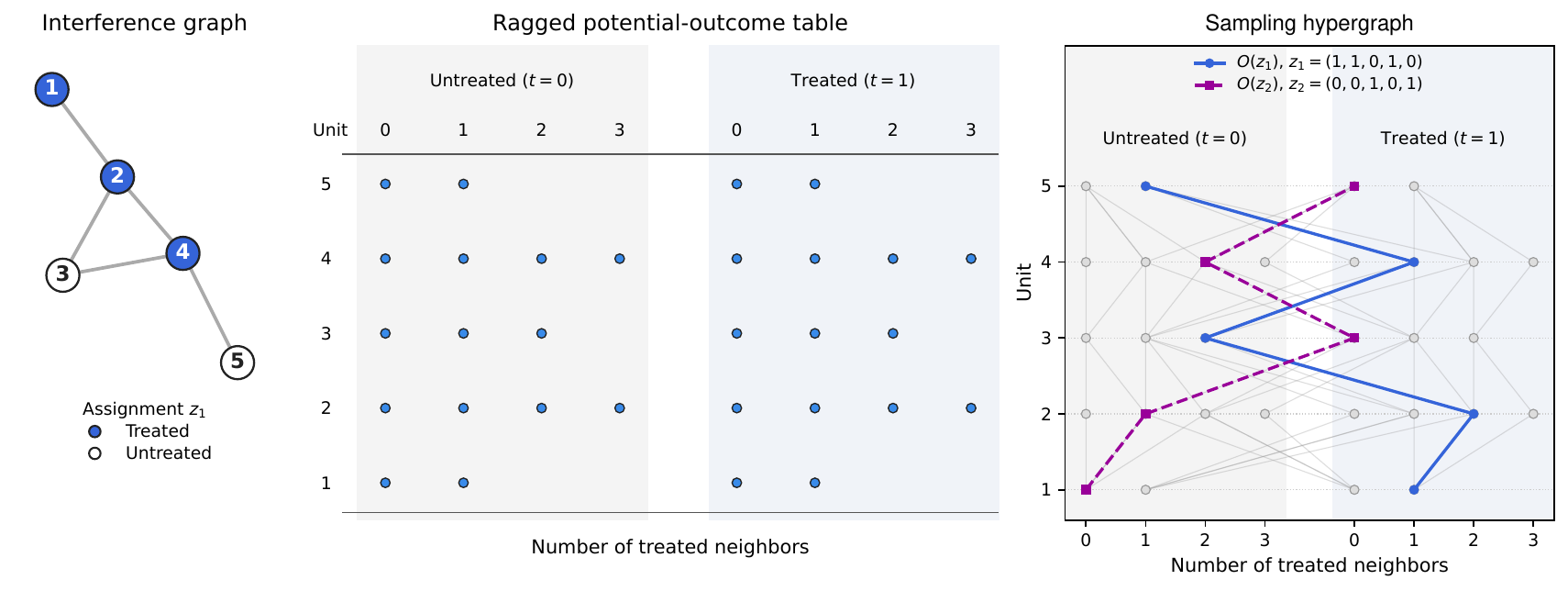}
  \caption{Construction of the (unweighted) sampling hypergraph for an interference graph on five units, with exposure $f_i(z)=(z_i,\sum_{j\in N_i}z_j)$.
  Left: filled nodes are treated under $z_1$. Middle: rows index units and
  columns are grouped by unit's treatment $t$, with headers giving the
  treated-neighbor count $c$ within each group. A unit's exposure is $(t,c)$. %
  Each dot represents the potential outcome $Y_i(t,c)$; blank entries indicate unattainable exposures. %
  Right: each hyperedge shows the cells of potential outcomes that can be sampled together, and selects one cell per row. The hyperedges $O(z_1)$ and $O(z_2)$ are highlighted; the remaining hyperedges are gray.}
  \label{fig:ragged-table-hypergraph}
\end{figure}

For each feasible assignment $z$, let $M(z)\subseteq K$ denote the set of potential outcome cells of the ragged table that are jointly missing under $z$, and let $O(z) $ denote its complement, the set of feasible potential outcome cells jointly observed under $z$.  The collection of sets $O(z)$ of potential outcomes, as $z$ ranges over all feasible treatment assignments, forms the hyperedge set of the \emph{sampling hypergraph}.
Its nodes are the elements of the expanded finite population $K$, representing the unit--exposure pairs $(i,e)$. Each hyperedge represents one feasible sample by recording the maximal set of potential outcome cells that can be jointly observed under a feasible assignment. Note that the hypergraph includes every feasible sample, including those assigned probability zero by a particular randomization design. A randomized design assigns each hyperedge a probability of jointly observing the potential outcomes in the set $O(z)$. For a given design, we call the sampling hypergraph equipped with these design-induced probabilities a \emph{weighted sampling hypergraph}, see Section \ref{sec:lower-bounds} for formal definitions. %
Figure~\ref{fig:ragged-table-hypergraph} illustrates the construction of a sampling hypergraph for an exposure mapping that records a unit's treatment along with the number of treated neighbors.

The minimax risks of the three problems depend on the induced sampling design and the corresponding (weighted) sampling hypergraph as follows. %
In the fixed-design problem, the randomization design is given and only the estimator can be optimized. The minimax risk depends on this randomization design only through its induced sampling design, represented by the weighted sampling hypergraph. In the restricted-design problem, we optimize over randomization designs in a prescribed class of designs and over estimators. The minimax risk therefore depends on the collection of induced sampling designs attainable within that class. In the all-design problem, every randomization design on the feasible assignment set is allowed. We can then optimize over all probability distributions on feasible samples and over estimators. Hence, when the estimand and potential outcome class is fixed, the minimax risk of the all-design problem depends only on the unweighted sampling hypergraph
(Proposition~\ref{thm:three-problem-characterization}).

The general lower-bound theory is developed in Section~\ref{sec:lower-bounds}. %
All of our lower bounds are obtained by constructing prior distributions over potential outcome tables. %
The key idea is a design-based indistinguishability principle: two potential-outcome tables with different estimand values cannot be distinguished using the observed data, whenever %
all potential outcomes on which the tables differ are missing under the randomization design. %
Specifically, consider two potential outcome tables that differ only on a set $B$ of estimand-relevant potential outcomes.  For a random assignment $Z$, consider the event $B\subseteq M(Z)$, meaning that every potential outcome in the set $B$ is missing under the treatment assignment $Z$; equivalently, $B\cap O(Z)=\varnothing$.  On this event, the two tables produce the same observed data even though the corresponding estimand values differ.  We place an equal-weight prior on these two tables. The resulting Bayes risk provides a lower bound on the worst-case design-based risk. %

We extend this two-table construction to a family of tables by varying potential outcomes within disjoint sets of estimand-relevant cells. We call this collection of disjoint sets a \emph{potential outcome packing} defined on the vertices of the sampling hypergraph. %
Formally, let $B_1,\ldots,B_m$ be disjoint sets of estimand-relevant potential outcomes. Let $c_b$ be the contribution to the estimand from a signed perturbation of the potential outcomes in $B_b$,
and let $\rho_P(B_b)$ be the probability under the randomization design $P$ that every cell in $B_b$ is missing.  Then the fixed-design minimax risks $\Rall(P)$ and $\Runb(P)$ over all estimators and over unbiased estimators respectively satisfy
\begin{equation}
 \Rall(P)\geq\sum_{b=1}^m c_b^2\,\rho_P(B_b),
 \qquad
 \Runb(P)\geq\sum_{b=1}^m c_b^2\,\frac{\rho_P(B_b)}{1-\rho_P(B_b)}.
 \label{eq:intro-main-bound-v9}
\end{equation}
 (Theorem~\ref{thm:fixed-joint-missingness}).  For all estimators, the lower bound depends on the probability that the relevant potential outcomes are jointly missing; for unbiased estimators, the lower bound depends on the odds of the
same event instead of its probability.  %
For the restricted-design problem, we may choose a different potential outcome packing for each feasible design within the restricted class to strengthen its lower bound, then take the infimum of the resulting bounds over those designs. When all randomization designs are allowed, we obtain lower bounds for general linear estimands that depend only on the sampling hypergraph and the estimand coefficients (Theorem~\ref{thm:global-hypergraph-lower}). We present two simpler bounds for estimands that contrasts two exposure levels, called \emph{arms}. 
The first uses the \emph{sampling rank} of one estimand arm, defined as the largest number of its potential outcomes that a single assignment can reveal. This rank plays the role of an effective sample size within the arm. For some two-arm contrasts, the one arm rank based lower bound can be weaker. In such a case, we derive a second lower bound that accounts for which potential outcomes in the two arms can be observed together
(Theorems~\ref{cor:rank-lower} and~\ref{thm:paired-sets-main}).

For a fixed design problem, we show that the minimax risk depends on the entire induced sampling design and not just the first-order inclusion probabilities and sampling hypergraph (Theorem~\ref{thm:cluster-with-noise}).  We give sufficient conditions under which exposure Horvitz--Thompson (H-T) is minimax rate optimal for a fixed design problem (Theorem~\ref{thm:fixed-ht-comparison}). We also illustrate a randomized-saturation counterexample where the exposure H-T is not minimax rate optimal
(Proposition~\ref{prop:fixed-ht-saturation-conditioning}). For the restricted- and all-design problems, our upper bounds pair exposure
HT, conditional HT, or a biased shrinkage estimator with designs adapted to the estimand. We provide two new constructions: an \emph{ordered-proposal} design, which thins focal units into an independent set, and a budget-constrained \emph{lit-ball}
design, which treats unions of selected neighborhoods.

We apply the general bounds to several causal estimands under four interference structures.
\begin{enumerate}[(i)]
\item \emph{No interference} (Section~\ref{sec:no-interference}).  For the
average treatment effect, we derive all-design minimax risks of $4/n$ over unbiased estimators and $(4+o(1))/n$ over all estimators
(Theorem~\ref{thm:ate-joint-benchmark}). We also show that if at most $m$
units can be treated, both risks are of order $1/m$
(Theorem~\ref{thm:ate-fixed-count}).
\item \emph{Partial interference} (Section~\ref{sec:partial-interference}).
We show that when the cluster memberships are known, the total and direct effects have
minimax rate $1/\xi$ over both estimator classes, where $\xi$ is the number of
clusters.  If the design must be chosen without knowing the memberships at the design stage, we show that the minimax
risk of the total effect increases by a factor that can be exponential in the
cluster size (Theorem~\ref{thm:invariant-tte}).
\item \emph{Neighborhood interference}
(Section~\ref{sec:neighborhood-interference}).  For the direct treatment
effect, we construct core--periphery graphs on which unbiasedness
multiplies the minimax risk by a factor of up to order $\sqrt n$, because a
biased strategy can omit a small, hard-to-observe part of the population
(Theorem~\ref{cor:dte-price-is-a-rate}), thereby reducing the minimax risk.  For the total
treatment effect on graphs with restricted growth, we show that graph cluster
randomization with the exposure Horvitz--Thompson estimator attains the
minimax rate $(\bar d+1)/n$ (Theorem~\ref{thm:restricted-growth-gcr-tte}). Under a sufficiently large budget of $m$ treated units,
the rate is $(d+1)/m$ on restricted-growth $d$-regular graphs (Theorem~\ref{thm:hard-budget-regular-tte}). 

We also study the minimax risk of a complex-contagion contrast that compares treated units above and below a threshold number of treated neighbors, motivated by threshold models of contagion \citep{Granovetter1978,CentolaMacy2007}. We show that lower bounds based only on which pairs of potential outcomes can be observed together (for example, through a conflict graph) understate the minimax risk by a growing factor (Theorem~\ref{thm:contagion-main-bracket} and Corollary~\ref{cor:contagion-conflict-gap}).
\item \emph{Bipartite interference} (Section~\ref{sec:bipartite-interference}). We study the minimax risk of dose--response contrasts between nearby doses when each outcome unit is linked to $d$ treatment units drawn uniformly at random from a common pool, and the exposure of an outcome unit is its number of treated neighbors. We show that lower bounds based only on which pairs of potential outcomes can be observed together again understate the minimax risk, by a factor that grows as the dose moves away from the endpoints (Theorem~\ref{thm:bipartite-central-main} and Corollary~\ref{cor:bipartite-full-projection-gap}).
\end{enumerate}

The rest of the paper is organized as follows. Section~\ref{sec:related-work} discusses related work;
Sections~\ref{sec:setup-strategy}--\ref{sec:fixed-design-ht} develop the
framework, lower bounds, and fixed-design comparison; and
Sections~\ref{sec:no-interference}--\ref{sec:bipartite-interference} give the
applications.  Section~\ref{sec:discussion} concludes, and the appendices
contain proofs and additional results.

\section{Related work}
\label{sec:related-work}

Design-based inference under interference was developed for partial
interference by \citet{Sobel2006What} and \citet{HudgensHalloran2008}.
Exposure mappings under general interference were formalized by
\citet{Manski2013} and \citet{AronowSamii2017}.
\citet{ZiglerPapadogeorgou2021} developed the framework for bipartite
interference. Most literature on causal inference under interference develops designs or estimators for particular exposure
mappings and estimands: two-stage designs under partial interference
\citep{HudgensHalloran2008,BasseFeller2017}, graph-cluster randomization for
the total effect \citep{ugander2023,UganderYin2023}, independent-set and
sequential designs for direct or peer effects
\citep{KarwaAiroldi2018,Toulis2013Estimation}, conflict-graph designs
\citep{KandirosEtAlConflict2026}, and estimators based on low-order outcome
models \citep{YuEtAl2022,CortezRodriguezEtAl2023}.  Their risk bounds are
upper bounds on the corresponding minimax risks, but most of this work does
not study lower bounds.  

A line of work fixes the estimator and chooses a design to
minimize its risk under a model for the potential outcomes. For example, under partial interference, \citet{BairdEtAl2018,CaiEtAl2022} choose treatment saturations
under random-effects and linear interference models, respectively.
Other designs aim to reduce bias for the direct effect \citep{JagadeesanEtAl2020} or compare the bias of total-effect estimators under two clusterings \citep{PougetAbadieEtAl2018}. Some works study the fixed-design problem in limited settings. \citet{SussmanAiroldi2017} derive minimum integrated-variance estimators within the class of linear unbiased estimators, \citet{KarwaAiroldi2023Admissibility} study the admissibility of the
Horvitz--Thompson estimator within the class of linear estimators, and
\citet{CortezRodriguezEtAl2023} bound the risk of every estimator under
Bernoulli assignment and bounded-order interactions.

Closer to our restricted-design problem is work that fixes the estimator and
chooses a design within a prescribed family to minimize a worst-case risk.
Without interference, \citet{Kallus2018Optimal,Kallus2021} optimize the
design for difference in means over a class of potential outcome models.
Under interference, related work optimizes the clustering for difference
in means \citep{VivianoEtAl2023Causal} or the dependence in treatment
assignment for Horvitz--Thompson estimation
\citep{CandoganChenNiazadeh2024,ThiyageswaranEtAl2026}.
Under spatially decaying interference,
\citet{Leung2022,Leung2025CrossCluster} choose the number of clusters for
estimators that vary within a small family. These works target the global
average treatment effect under interference. Our restricted-design problem
optimizes over both the designs in a class and all estimators, for general
exposure mappings and linear estimands. This formulation draws on
finite-population sampling, where an estimation strategy pairs a sampling
design with an estimator \citep{Godambe1955,FienbergTanur1987}.

Minimax analyses that optimize jointly over all designs and estimators
are recent, both with and without interference.  Without interference, in a
superpopulation model with covariates, \citet{Bai2023} shows that complete
randomization with the difference-in-means estimator is minimax optimal over all designs and linear estimators, when the class of outcome
distributions is permutation-invariant.  For the
average treatment effect with bounded outcomes, \citet{Hull2026} shows that
independent assignment, rather than complete randomization, is minimax
optimal over all designs and estimators, because the random number of treated units carries information when outcomes are bounded.
\citet{SudijonoDobribanTchetgen2026} study binary outcomes and characterize
the minimax risk up to its second-order term.  It is attained by Bernoulli
randomization with a nonlinear shrinkage estimator, whereas complete
randomization with the difference-in-means estimator is optimal only to
first order.  Under interference,
\citet{FaridaniNiehaus2026} derive the minimax rate for the global average
treatment effect over all designs and linear estimators under spatially
decaying interference.
\citet{BorusyakHullMunro2026} jointly choose a design and a regression
estimator under a linear exposure mapping.

\paragraph*{Comparison with Kandiros, Harshaw and S\"avje}
\label{subsec:related-khs}
The work closest to ours is \citet{KandirosHarshawSavje2026}, who study the all-design minimax problem under neighborhood interference and the identity exposure mapping, without any restrictions on the set of treatment assignments. They show that the minimax rate for this setting is governed by a \emph{conflict graph} that records which pairs of estimand-relevant potential outcomes cannot be observed together under any treatment assignment. Their general results cover equal-weight two-arm contrasts, with applications to the direct and total treatment effects. They consider potential outcome classes defined by moment restrictions of order $q\geq2$. Our uniformly bounded class corresponds to their setting $q=\infty$.

The two formulations use different combinatorial objects.  The conflict graph
of \citet{KandirosHarshawSavje2026} records pairwise conflicts: an edge in the conflict graph joins two potential outcomes that cannot be observed
together.  Our sampling hypergraph instead records the sets of potential
outcomes that can be jointly observed, and therefore retains higher-order
information. Apart from the higher order information, another distinction between the two objects is that the sampling hypergraph describes feasibility by listing maximal feasible samples, whereas the conflict graph describes infeasibility through pairs that cannot be observed together. %
Under neighborhood interference with the identity exposure mapping and unrestricted assignments, \citet[Proposition~2.2]{KandirosHarshawSavje2026} show that the labelled conflict graph determines the all-design minimax risk. %
In this setting, we show that the conflict graph and sampling hypergraph determine one another (Theorem~\ref{thm:conflict-hypergraph-equivalence} in
Appendix~\ref{subsec:pairwise-exact}). In particular, if every pair in a collection of potential outcomes can be jointly observed, possibly under different assignments, then the entire collection can be observed under one assignment. Our sampling-rank lower bound (Theorem~\ref{cor:rank-lower}) reduces to the independence-number lower bound of \citet[Theorem~3.2]{KandirosHarshawSavje2026} when $q = \infty$.

However, when the exposure mapping is coarsened, or there are constraints on the set of possible treatment assignments, the conflict graph, in general, is no longer sufficient to determine the all-design minimax risk. In such cases, the sampling hypergraph retains the joint sampling constraints needed for the characterization of the all-design minimax risk.  For example, under a count or threshold mapping, a pairwise-observable collection of potential outcomes need not be jointly observable.  In our complex-contagion and bipartite examples, this higher-order constraint causes a conflict-graph lower bound to understate the minimax risk by a growing factor.  Even under neighborhood interference, the conflict graph can be insufficient when there are  constraints on the set of feasible treatment assignments.  For example, under a budget on the number of treated units, the minimax rates for the direct and total treatment effects under neighborhood interference %
are governed by higher-order sampling constraints that the conflict graph does not record; see Appendix~\ref{subsec:app-budgeted-dte} and Section~\ref{subsec:hard-budget-regular-tte}.

Beyond this overlap in the all-design minimax problem, we also consider the fixed- and restricted-design minimax problems.  For these problems, neither the unweighted sampling hypergraph nor the conflict graph is sufficient. Instead, the lower bounds also depend
on the probabilities assigned to jointly observable sets, not just on which
sets are feasible.  Our framework also covers arbitrary linear estimands, including one-arm means, weighted two-arm and multi-arm contrasts, and policy contrasts, for the three minimax problems. For each minimax problem, we also compare the minimax risks over all estimators and over unbiased estimators.
\citet{KandirosHarshawSavje2026} discuss a bias--variance trade-off, dropping
some conflict-graph vertices to improve their upper bound, and illustrate the
improvement under a second-moment restriction.  Their lower bounds apply to all estimators and hence also to unbiased estimators, but need not be sufficient to determine the minimax rate over the unbiased class and do not establish a cost of unbiasedness. We derive lower bounds that exploit the unbiasedness restriction and show that this cost can be a factor of order $\sqrt n$
(Theorem~\ref{cor:dte-price-is-a-rate}). In turn, \citet{KandirosHarshawSavje2026} cover outcome classes with finite $q$ and
give a general spectral upper bound, neither of which we study.
Remarks~\ref{rem:dte-khs-regular}, \ref{rem:dte-khs-comparison},
and~\ref{rem:tte-khs-comparison} compare the direct- and total-effect results
with theirs, and
Appendix~\ref{subsec:conflict-comparison} gives a more detailed comparison.

\section{Potential outcomes, constrained sampling design, and the three minimax problems}
\label{sec:setup-strategy}

In this section, we define the potential outcome table under interference, its
missing-data representation using a constrained sampling problem induced by a given randomization design, and the three minimax problems.

\subsection{Potential outcomes under interference}
\label{subsec:population}

Let $[n]=\{1,\ldots,n\}$ index a finite population under binary treatment.
The feasible assignments form a set $\cZ_0\subseteq\{0,1\}^n$, and a randomization design is
a law $P\in\cP(\cZ_0)$ generating $Z\sim P$. Thus $P$ is a probability distribution on treatment assignments. We use the term ``design'' to mean randomization design unless explicitly referring to the induced sampling design.  Each unit has a fixed but unknown potential outcome $Y_i(z)$ for every $z\in\{0,1\}^n$; all randomness comes
from the design.  Without further structure the table
$\mathbb T=\{Y_i(z):i\in[n],z\in\{0,1\}^n\}$ has $n2^n$ entries.  We reduce
this table using an interference graph and an exposure mapping.

Let $G=(V,E)$ be an undirected graph with vertices given by the units in the finite population $[n]$. %
The lack of an edge $\{i,j\}$ in $G$ encodes the assumption that the potential outcome of unit $i$ is not affected by the treatment assigned to unit $j$, and vice versa. Formally, let \(N_i\) denote the neighbors of unit \(i\) in $G$ and \(N_i^+=\{i\}\cup N_i\) denote the closed neighborhood of $i$.  More generally, \(B_r(i)=\{j:\operatorname{dist}_G(i,j)\leq r\}\) denotes the units within distance \(r\) of \(i\), so that \(B_1(i)=N_i^+\). We assume
\begin{equation}
  z_{N_i^+}=z'_{N_i^+}\quad\Longrightarrow\quad Y_i(z)=Y_i(z')
  \qquad\text{for each }i.
  \label{eq:graph-interference}
\end{equation}
Different assumptions on interference can be specified as restrictions on  $G$. For example, an empty $G$ encodes the \emph{no interference} assumption, disjoint cliques give \emph{partial interference}, and arbitrary $G$ gives \emph{neighborhood interference} (extendable to $r$-hop structure by enlarging $N_i$ to $B_r(i)$).  A complete $G$ encodes global interference.

For each unit $i$, let
\(h_i(\cdot;\theta_i):\{0,1\}^{N_i^+}\longrightarrow\cD_i\)
be a function of the treatment assignments in $i$'s closed neighborhood.
The corresponding \emph{effective treatment} function or \emph{exposure mapping} \citep{Manski2013,AronowSamii2017}
is defined on the full assignment vector $z$ by
\[
 f_i:\{0,1\}^{n}\longrightarrow\cD_i,
 \qquad f_i(z;\theta_i)=h_i(z_{N_i^+};\theta_i).
\]
Thus $G$ specifies which treatment assignments may enter the exposure,
and $h_i$ specifies how they are summarized.  We call $e=f_i(z;\theta_i)$
the \emph{effective treatment level} or \emph{exposure level} of unit $i$.
 It is \emph{correctly specified} when, for every unit $i$ and assignments $z,z'$,
\begin{equation}
  f_i(z; \theta_i)=f_i(z';\theta_i)\quad\Longrightarrow\quad Y_i(z)=Y_i(z').
  \label{eq:exposure-sufficiency}
\end{equation}
We suppress $\theta_i$ when it is clear and write $f=\{f_i\}$.  The identity
mapping $f_i(z)=z_{N_i^+}$ retains the full neighborhood assignment.  Two
coarser mappings used later are the anonymous count
$f_i(z)=(z_i,\sum_{j\in N_i}z_j)$ and, for $\theta\geq1$, the threshold
$f_i(z;\theta)=(z_i,\1\{\sum_{j\in N_i}z_j\geq\theta\})$.  %
In general, we take $\cD_i=f_i(\{0,1\}^n)$ to contain exactly the exposure levels attainable by unit $i$.  These levels depend on its neighborhood and exposure
mapping.  For example, a unit with no neighbors cannot attain an exposure level
requiring at least one treated neighbor.

We assume that the exposure mapping is correctly specified. Hence,
all assignments giving unit $i$ the same exposure $e$ have the same
potential outcome.  With a slight abuse of notation, we denote this common
value by $Y_i(e)$, the \emph{reduced potential outcome} under exposure mapping $f$.  We collect these values in the table of reduced potential outcomes $\{Y_i(e):i\in[n],\ e\in\cD_i\}$. Since the attainable exposure levels can differ across units, the table of reduced potential outcomes is ragged.

\begin{definition}[Potential outcome indices]
Let $K=\{(i,e):i\in[n],\ e\in\cD_i\}$
be the index set of the reduced potential outcome table.  Each
$k=(i,e)\in K$ is called a \emph{potential outcome index}; the corresponding potential outcome is
$Y_k=Y_i(e)$.
\end{definition}

Thus $(Y_k)_{k\in K}$ is a vector of length
$|K|=\sum_i|\cD_i|$.  We use $(i,e)$ when the row and exposure matter and
$k\in K$ for the vector representation, and hereafter usually omit the word
``reduced.'' When referring to collections of potential outcomes, we
identify each potential outcome with its unit--exposure index, $Y_i(e)\leftrightarrow(i,e)$. We assume that the (reduced) potential outcomes are uniformly bounded, and, without loss of generality, we normalize the bound to one: for a given exposure mapping \(f\), the outcome class is \(\cY(f)=\{Y=(Y_k)_{k\in K}:|Y_k|\leq1\text{ for every }k\in K\}\). We call a potential-outcome table $Y$ \emph{admissible} if $Y\in\cY(f)$.

\subsection{Estimands}
\label{subsec:estimands}

The estimands of interest are known linear functionals of the reduced potential outcomes $Y_i(e)$. A generic estimand is
\begin{equation}
  \tau_a(Y)=\sum_{(i,e)\in K}a_{(i,e)}Y_i(e),
  \label{eq:general-estimand}
\end{equation}
with $K_a=\{(i,e)\in K:a_{(i,e)}\neq0\}$ the set of \emph{estimand relevant} potential outcome indices. This formulation includes finite-population means and contrasts.  A finite-population mean at a fixed exposure level $e$ attainable by every unit, $e\in\bigcap_i\cD_i$, is 
$$\mu_e(Y) = \frac1n\sum_{i=1}^n Y_i(e).$$ 
Linear combinations of fixed-exposure means are obtained by choosing $a_{(i,e)} = \frac{\alpha_e}{n}$ for each such level $e$, which gives
\begin{equation}
  \tau_\alpha(Y) = \sum_e \alpha_e \mu_e(Y) = \sum_e\alpha_e\left\{\frac1n\sum_{i=1}^n Y_i(e)\right\}.
  \label{eq:fixed-exposure-comparison}
\end{equation}
where both sums over $e$ range over the exposure levels attainable by every unit.
An average causal effect is a contrast $\tau_{\alpha}(Y)$ with $\sum_e\alpha_e=0$. A policy $\Gamma$ is a specified probability distribution on treatment assignments. Policy means, with coefficients $\omega_{i,\Gamma}(e)=\Pr_\Gamma\{f_i(Z)=e\}$, follow similarly:
\begin{equation}
  \mu_\Gamma(Y)=\frac1n\sum_{i=1}^n\sum_{e\in\cD_i}
    \omega_{i,\Gamma}(e)Y_i(e).
  \label{eq:policy-mean}
\end{equation}
A policy effect is a contrast of policy means.

\subsection{The sampling design on an expanded finite population}
\label{subsec:observed-data}

We now formulate the patterns of missing and observed potential outcomes as a constrained sampling problem from an expanded finite population. This population consists of the potential-outcome indices $K = \{(i,e): i \in [n], e \in \mathcal D_i\}$, with fixed values $(Y_k)_{k\in K}$ assigned to each element of $K$. Note that the finite population changes with the exposure mapping.  In general $|K| = \sum_i d_i $. For example, in the simplest case of no-interference, where $f_i(z) = z_i$, and $d_i = 2$, the expanded finite population consists of $2n$ elements. A random treatment assignment $Z\sim P$ selects a sample of size $n$ from this population. For $k=(i,e)\in K$, define the \emph{sampling} or the \emph{observation} indicator $I_k(Z)=\1\{f_i(Z)=e\}$, which records whether the potential outcome $Y_i(e)$ is sampled (or equivalently, observed). The observed and missing index sets of potential outcomes are defined as
\begin{equation}
  O_f(Z)=\{k\in K:I_k(Z)=1\},
  \qquad
  M_f(Z)=K\setminus O_f(Z).
  \label{eq:observed-missing-patterns}
\end{equation}
Equivalently, for a fixed treatment assignment $z$, $O_f(z)=\{(i,e)\in K:f_i(z)=e\}$. Thus, $O_f(Z) \subset K$ indexes the set of potential outcomes of the ragged table that are jointly observed under a treatment assignment $Z$.
The observed outcomes are \(Y_{O_f(Z)}=\{Y_k:k\in O_f(Z)\}\).
The \emph{observed data} is defined as the vector of the treatment assignment $Z$ and the observed outcomes \((Z,Y_{O_f(Z)})\). For any feasible treatment assignment $z \in \cZ_0$, we call $O_f(z)$ a \emph{feasible sample}. The design $P$ induces a law on the feasible samples.

\begin{definition}[Induced sampling design]
\label{def:observation-law}
The randomization design $P$, together with the exposure mapping $f$, induces
a sampling design on $K$, defined as %
\[
 Q_{P,f}(A)=P\{O_f(Z)=A\}
 =\sum_{z\in\cZ_0:\,O_f(z)=A}P(Z=z),\qquad A\subseteq K.
\]
The randomization design $P$ assigns probabilities to treatment assignments;
the induced sampling design $Q_{P,f}$ assigns probabilities to observed samples
of potential-outcome indices.
\end{definition}
For fixed $P$ and $f$, this sampling design does not depend on the
potential-outcome values. Changing the table $Y$ changes the values observed,
but leaves the probability of selecting each sample unchanged. Equivalently,
$Q_{P,f}$ is the joint distribution of the sample inclusion indicators
$(I_k(Z))_{k\in K}$. In missing-data terms, the fixed potential-outcome table is incompletely observed: $I_k(Z)=1$ indicates that $Y_k$ is observed, and $I_k(Z)=0$ indicates that it is missing. The induced sampling design $Q_{P,f}$ therefore specifies the distribution of the missingness pattern. For fixed $P$ and $f$, this distribution is known and is the same for every potential-outcome table. See \citet{Rubin1978} for the missing-data formulation of causal inference and \citet{DingLi2018} for a review. %

Two types of sampling constraints are imposed on which sets of potential outcomes can be sampled together in $O_f(Z)$. The first is within-unit counterfactual missingness. Each treatment assignment selects exactly one cell from each row of the ragged table; all remaining potential outcomes are missing. Thus $\sum_{e\in\cD_i}I_{(i,e)}(Z)=1$ for every $i\in[n]$, and the sample $O_f(Z)$ has size $n$. The second sampling constraint is across rows: here, the interference graph, exposure mapping, and restrictions on treatment assignments additionally constrain which potential outcomes from different units can be observed together. The sampling hypergraph described in Section~\ref{subsec:packings-and-the-three-problems} records these constraints implicitly by recording all feasible samples, and the hyperedge weights give the sample probabilities $Q_{P,f}(A)$.

The first-order inclusion probability that a relevant potential outcome is observed is
\begin{equation}
  \pi_{P,k}=P\{k\in O_f(Z)\}.
  \label{eq:first-order-pi}
\end{equation}
These are one-dimensional margins of the induced sampling design $Q_{P,f}$. For $k=(i,e)$, this is $\pi_{P,k}=P\{f_i(Z)=e\}$, also called the exposure probability of unit $i$ at exposure $e$.  For an exposure level $e$, define the \emph{observed exposure count} $n_e(Z)=\sum_{i:e\in\cD_i}\1\{f_i(Z)=e\}$.
This is the number of units revealing $Y_i(e)$ under assignment $Z$;
its expectation under randomization design $P$ is $\E_P n_e(Z)=\sum_{i:e\in\cD_i}\pi_{P,(i,e)}$. 

\subsection{The three minimax problems and the estimator classes}
\label{subsec:three-minimax}

Given \(P\), let $T_P$ denote any estimator and \(\cT(P,f)\) be the set of all estimators measurable with respect to
\((Z,Y_{O_f(Z)})\).

\begin{definition}[Estimation strategy]
\label{def:estimation-strategy}
For estimating $\tau_a$ under exposure mapping $f$, an \emph{estimation
strategy} is a pair $(P,T_P)$ consisting of a randomization design $P\in\cP(\cZ)$ and an
estimator $T_P\in\cT(P,f)$.
\end{definition}

The fixed-design minimax risk over all estimators is
\begin{equation}
  \Rall(P;f,a)
  =\inf_{T_P\in\cT(P,f)}\sup_{Y\in\cY(f)}
  \E_P\left[\{T_P-\tau_a(Y)\}^2\mid Y\right].
  \label{eq:minimax-risk-all}
\end{equation}
An estimator is design-unbiased when
\(\E_P(T_P\mid Y)=\tau_a(Y)\) for every \(Y\in\cY(f)\).
A design $P$ is positive for \(\tau_a\) when every relevant potential outcome can
be observed with positive probability:
\begin{equation}
  \cP_+(f,a)=\{P\in\cP(\cZ_0):
  \pi_{P,k}>0\text{ for every }k\in K_a\}.
  \label{eq:positive-design-class}
\end{equation}
Positivity is required for unbiasedness, but not for unrestricted
estimation, since a biased estimator remains well defined even when a relevant entry is always missing.

The corresponding fixed-design minimax risk is
\begin{equation}
  \Runb(P;f,a)
  =\inf_{\substack{T_P\in\cT(P,f)\\T_P\ \mathrm{exactly\ unbiased}}}
  \sup_{Y\in\cY(f)}\Var_P(T_P\mid Y),
  \label{eq:minimax-risk-unb}
\end{equation}
with value \(+\infty\) when the unbiased class is empty.  We always have
\(\Rall(P;f,a)\leq\Runb(P;f,a)\).  When
$f$ and $a$ are fixed, we write $\Rall(P)$ and $\Runb(P)$. For a design class \(\cC\subseteq\cP(\cZ)\), let
\(\cC_+=\cC\cap\cP_+(f,a)\) be the class of positive designs within $\cC$. Define
\begin{equation}
  \Rallstar(\cC;f,a)=\inf_{P\in\cC}\Rall(P;f,a),
  \qquad
  \Runbstar(\cC;f,a)=\inf_{P\in\cC_+}\Runb(P;f,a).
  \label{eq:class-minimax}
\end{equation}
We suppress $f$ and
$a$ when they are clear.  This gives three minimax problems:
\[
\begin{array}{lll}
\text{Problem 1:}&\Rall(P),\ \Runb(P),&\text{the design is fixed};\\
\text{Problem 2:}&\Rallstar(\cC),\ \Runbstar(\cC),&
  \text{the design class is restricted};\\
\text{Problem 3:}&\Rallstar(\cP(\cZ)),\ \Runbstar(\cP(\cZ)),&
  \text{all designs are allowed}.
\end{array}
\]
The estimand coefficients, exposure mapping, and outcome class are fixed in all
three problems; what changes is whether the design is fixed, optimized within a
restricted class, or optimized over all distributions on \(\cZ_0\).  In the next
section, we show that the minimax risks of these three problems are governed by three different objects: the law of \(O_f(Z)\) under a fixed design, the set of such laws that a design class can generate, and the collection of observation patterns possible over all feasible assignments. The notation \(\cP(\cZ_0)\) denotes all probability laws supported on an
assignment set \(\cZ_0\), for example an exact-count set \(\cZ_m\) or a budget
set \(\cZ_{\leq m}\).  We reserve \(\cC\) for design restrictions that are not
only support restrictions, such as product Bernoulli designs or the label-invariant
class in Section~\ref{sec:partial-interference}.

Some upper bound constructions use auxiliary randomization $U$, such as a
randomly selected set of units on which treatment is randomized.  An estimator
that uses $U$ in addition to the observed data need not belong to $\cT(P,f)$.
The following lemma (proof in Appendix~\ref{app:setup-proofs}) gives a \emph{Rao--Blackwellized} estimator in $\cT(P,f)$
with no greater mean squared error for every admissible table.  %

\begin{lemma}[Rao--Blackwellized estimator]
\label{lem:rao-blackwell}
Suppose $U$ takes values in a finite or countable set $\mathcal U$, and the joint law of $U$ and assignment $Z$
is known and does not depend on $Y$.  Let $D=(Z,Y_{O_f(Z)})$ and let
$T(U,D)$ have finite second moment for every admissible table.  Then
$\widetilde T(D)=\E[T(U,D)\mid D]$ depends only on the observed data $D$ and satisfies
\begin{equation}
 \E[(\widetilde T-\tau_a(Y))^2\mid Y]
 \leq \E[(T-\tau_a(Y))^2\mid Y]
 \qquad\text{for every }Y\in\cY(f).
 \label{eq:rb-domination}
\end{equation}
If $T$ is unbiased, so is $\widetilde T$, and its variance
is no larger.
\end{lemma}

\section{Minimax lower bounds from jointly missing potential outcomes}
\label{sec:lower-bounds}
The lower bounds in this section are based on the following principle.
Let $\tau_a$ be a linear estimand, with estimand coefficient vector $a$ on the relevant potential
outcomes.  Consider two potential outcome tables $Y,Y'\in\cY(f)$ that
agree outside a set $B\subseteq K_a$ and differ in the potential outcome at every index in $B$.
We choose the potential outcome values at $B$ so that $\tau_a(Y)\neq\tau_a(Y')$. When all potential outcomes in $B$ are missing, the two tables produce identical observed data. Hence, under this event, the design cannot distinguish the corresponding estimand values.  Therefore, an adversarial perturbation on $B$ gives a lower bound on the design-based mean squared error determined by the difference in estimand values and the probability that the potential outcomes in $B$ are jointly missing. Using this principle, Section~\ref{subsec:the-detection-principle} derives
lower bounds under a fixed design, for two tables that differ on a
single set $B$ and a disjoint collection of sets.
Section~\ref{subsec:three-problem-bounds} shows that the fixed-design minimax
risk depends on the design only through its induced sampling design, the
restricted-design risks only through the set of attainable induced sampling designs,
and the all-design risks only through the sampling hypergraph; the
fixed-design bounds therefore extend to the other two problems.  Section~\ref{subsec:all-designs-lower}
gives two lower-bound constructions for the all-design problem.

\subsection{Lower bounds under a fixed design}
\label{subsec:the-detection-principle}

Fix an exposure mapping $f$, estimand coefficients $a$, and a design $P$
on an allowed assignment set $\cZ_0$.  Recall that $K_a$ is the set of relevant potential outcomes with nonzero estimand coefficients.

\begin{definition}[Perturbation set]
\label{def:perturbation-set}
A \emph{perturbation set} is a nonempty subset $B\subseteq K_a$ of potential outcomes. Two admissible potential outcome tables $Y,Y'\in\cY(f)$ are said to
\emph{differ exactly on $B$} if $  \{k\in K:Y_k\neq Y'_k\}=B.$
\end{definition}

Consider two tables $Y,Y'\in\cY(f)$ that differ exactly on a set $B$ and
have different estimands $\tau_a(Y)\neq\tau_a(Y')$.  The set and the tables
may depend on the design $P$, but are chosen before the assignment $Z\sim P$
is realized.  Let $M=\{O_f(Z)\cap B=\varnothing\}$ be the event that all
potential outcomes in $B$ are jointly missing.  Under $M$, the tables
$Y$ and $Y'$ produce the same observed data $(Z,Y_{O_f(Z)})$, so their
estimands cannot be distinguished.  We place an equal-weight two-point prior
on the two tables, independently of $Z$.  This is the two-point construction
underlying Le Cam's method \citep{LeCam1973}, and its Bayes risk gives a lower
bound on the worst-case design mean squared error.

\begin{definition}[Joint missingness probability]
\label{def:joint-missingness}
For a nonempty set $B\subseteq K_a$ of potential outcomes, its \emph{joint missingness probability}
under design $P$ is
  $\rho_P(B)=P\{O_f(Z)\cap B=\varnothing\}.$
\end{definition}

The joint missingness probability can be read from the induced sampling
design $Q_{P,f}$ defined in Section~\ref{subsec:observed-data} as
$\rho_P(B)=\sum_{A:A\cap B=\varnothing}Q_{P,f}(A)$.

\begin{lemma}[Two-table lower bound]
\label{lem:two-table-general}
Let $Y,Y'\in\cY(f)$ differ exactly on a nonempty set $B\subseteq K_a$, and
let $\Delta=|\tau_a(Y)-\tau_a(Y')|>0$.  For every design $P$,
\[
 \Rall(P;f,a)\geq\frac{\Delta^2}{4}\rho_P(B), \mathrm{ and} \qquad
 \Runb(P;f,a)\geq\frac{\Delta^2}{4} \frac{\rho_P(B)}{1-\rho_P(B)}.
\]

If $\rho_P(B)=1$, no unbiased estimator exists and the second bound is $+\infty$.
\end{lemma}

The proof is in Appendix~\ref{app:fixed-design-refinements}.  Under the
two-point prior, the posterior variance of the estimand is $\Delta^2/4$ on the
event that all potential outcomes in $B$ are missing and zero otherwise, which gives
the first bound.  Unbiasedness gives the stronger second bound,
in which the joint missingness probability $\rho_P(B)$ is replaced by its
odds $\rho_P(B)/\{1-\rho_P(B)\}$.  The choice of $B$ involves a tradeoff: a
larger perturbation set $B$ can increase the separation $\Delta$ between the estimands, but it decreases the probability $\rho_P(B)$ that all the potential outcomes in $B$ are jointly missing. For bounded outcomes, we use perturbations based on  signs.  Fix perturbation signs
$s_k\in\{-1,1\}$ and set $Y_k=s_k$ and $Y'_k=-s_k$ for $k\in B$, with both
tables zero elsewhere.  Then $\tau_a(Y)=c_B$ and $\tau_a(Y')=-c_B$, where
$c_B=\sum_{k\in B}a_ks_k$, so $\Delta=2|c_B|$ and
Lemma~\ref{lem:two-table-general} gives the lower bounds $c_B^2\rho_P(B)$ and
$c_B^2\rho_P(B)/\{1-\rho_P(B)\}$.  Choosing $s_k$ to have the sign of $a_k$
makes $c_B=\sum_{k\in B}|a_k|$ positive.

\paragraph*{Several disjoint perturbation sets}

We now extend the two-table construction from a single perturbation set to independent perturbations of several disjoint sets, with one independent random sign for each set.

\begin{definition}[Perturbation packing]
\label{def:perturbation-packing}
A \emph{perturbation packing} is a finite collection
$\cB=\{B_1,\ldots,B_m\}$ of pairwise disjoint perturbation sets.
The sets need not cover all estimand-relevant potential outcomes in $K_a$.
\end{definition}

Fix the perturbation signs $s_{b,k}\in\{-1,1\}$ for $k\in B_b$.
Consider an $m$-dimensional sign vector $\theta=(\theta_1,\ldots,\theta_m)\in\{-1,1\}^m$,
and for each $\theta$ define a potential outcome table $Y^\theta=(Y_k^\theta)_{k\in K}$ by
\[
 Y_k^\theta=
 \begin{cases}
  s_{b,k}\theta_b,& k\in B_b,\\
  0,& k\notin\bigcup_{b=1}^m B_b.
 \end{cases}
\]
Each sign vector $\theta$ gives one potential outcome table, so the $2^m$
sign vectors give $2^m$ distinct tables $Y^{\theta}$.  We place a uniform
\emph{hypercube prior} on these tables: the components
$\theta_1,\ldots,\theta_m$ are independent Rademacher signs, taking the values
$\pm1$ with probability $1/2$ each, drawn independently of $Z\sim P$. The corresponding estimand values are
\[
 \tau_a(Y^\theta)=\sum_{b=1}^m c_b\theta_b,
\]
where $ c_b=\sum_{k\in B_b}a_ks_{b,k}.
 $
Observing any potential outcome in $B_b$ reveals $\theta_b$, because the perturbation
signs $s_{b,k}$ are fixed.  If every potential outcome in $B_b$ is missing,
$\theta_b$ remains equally likely to be $+1$ or $-1$ given the observed data.
The unresolved signs remain independent, so their contributions to the
posterior variance add.  This is the hypercube construction used in Assouad's
method \citep{Assouad1983}: changing one sign changes the table only on the
corresponding set $B_b$.  It gives the following lower bound.

\begin{theorem}[Lower bounds from disjoint perturbations]
\label{thm:fixed-joint-missingness}
\label{thm:blocklb-all}
\label{thm:blocklb}
Fix a design $P$, a perturbation packing $\cB=\{B_1,\ldots,B_m\}$, and signs
defining $c_b$.
\begin{enumerate}
\item[\textup{(i)}] Over all estimators,
\begin{equation}
  \Rall(P;f,a)\geq\sum_{b=1}^m c_b^2\rho_P(B_b).
  \label{eq:blocklb-all}
\end{equation}
\item[\textup{(ii)}] Over unbiased estimators,
\begin{equation}
  \Runb(P;f,a)
  \geq\sum_{b:c_b\neq0}c_b^2
       \frac{\rho_P(B_b)}{1-\rho_P(B_b)}.
  \label{eq:blocklb-unb}
\end{equation}
If $c_b\neq0$ and $\rho_P(B_b)=1$ for some set, no unbiased
estimator exists and the bound is $+\infty$.
\end{enumerate}
\end{theorem}

The proof is in Appendix~\ref{app:fixed-design-refinements}.  The bound holds even
though the observation events of different perturbation sets $B_b$ are dependent, since they are
determined by the common assignment.  When $m=1$, the theorem recovers the two-table lower bound of Lemma \ref{lem:two-table-general} above. Taking the singleton packing $\cB=\{\{k\}:k\in K_a\}$ gives a lower bound in terms of the inclusion probabilities $\pi_{P,k}=P\{k\in O_f(Z)\}$.

\begin{corollary}[Singleton lower bound]
\label{cor:singleton-lower}
For every design $P$,
\begin{equation}
 \Rall(P;f,a)\geq\sum_{k\in K_a}a_k^2(1-\pi_{P,k}).
 \label{eq:singleton-all}
\end{equation}
If $\pi_{P,k}>0$ for every $k\in K_a$, then
\begin{equation}
 \Runb(P;f,a)\geq\sum_{k\in K_a}a_k^2\bigl(\pi_{P,k}^{-1}-1\bigr).
 \label{eq:singleton-unb}
\end{equation}
If $\pi_{P,k}=0$ for some $k\in K_a$, no unbiased estimator
exists.
\end{corollary}

The proof is in Appendix~\ref{app:fixed-design-refinements}.  Rarely observed potential outcomes constrain unbiased estimation: an estimator that allows bias can ignore such a potential outcome and accept the bias, whereas an unbiased estimator must compensate for its small inclusion probability.  The unbiased lower bound matches the diagonal variance term of the corresponding Horvitz--Thompson estimator, which uses inverse inclusion probabilities as weights.

\subsection{The sampling hypergraph and the three minimax problems}
\label{subsec:packings-and-the-three-problems}
\label{subsec:three-problem-bounds}
We now characterize the objects that govern the minimax risks of the three
problems introduced in Section~\ref{subsec:three-minimax}.   The central objects are the sampling hypergraph and its weighted version, defined below. 
Recall that $O_f(Z)$ indexes the potential outcomes observed under assignment $Z$.
The sampling hypergraph records the maximal collections of potential outcomes that can be jointly observed.

\begin{definition}[Sampling hypergraph]
\label{def:co-observation-hypergraph}
For an assignment set $\cZ_0\subseteq\cZ$, the sampling hypergraph is defined as $\cH_{\mathrm{obs}}(f;\cZ_0) =\bigl(K, \mathcal E\bigr),$ where $K$ is the vertex set and $\mathcal E$ is the set of hyperedges. Each vertex $k=(i,e)\in K$ represents the reduced potential outcome $Y_i(e)$, and the hyperedge set is
$\mathcal E=\{O_f(z):z\in\cZ_0\}$.
\end{definition}
For a treatment assignment $z\in\cZ_0$, the corresponding hyperedge
$E=O_f(z)$ records exactly one potential outcome per unit and therefore has size $n$.
Thus every hyperedge is maximal.
The hyperedge set includes every feasible sample, including samples assigned probability zero by a particular randomization design.
When $\cZ_0=\cZ$, we write $\cH_{\mathrm{obs}}(f)$.
When $f$ and $\cZ_0$ are fixed, we often abbreviate the hypergraph by $\cH$.

We next define a \emph{weighted sampling hypergraph}. The induced sampling design $Q_{P,f}$, defined in Section~\ref{subsec:observed-data}, supplies its hyperedge weights.

\begin{definition}[Weighted sampling hypergraph]
\label{def:weighted-co-observation-hypergraph}

For a randomization design $P$ on $\cZ_0$, the weighted sampling hypergraph
is $(\cH_{\mathrm{obs}}(f;\cZ_0),w_P)$, with weights
$w_P(A)=Q_{P,f}(A)$ for every $A\in\mathcal E$.
\end{definition}

When the design varies within a class $\cC$, its induced sampling design varies
within the following set of laws, corresponding to the weighted hypergraphs
attainable within the class.

\begin{definition}[Attainable induced sampling designs]
\label{def:attainable-observation-laws}
For a design class $\cC$, the set of attainable induced sampling designs is
$\mathcal Q_f(\cC)=\{Q_{P,f}:P\in\cC\}$.

\end{definition}

Finally, suppose that every design on $\cZ_0$ is allowed. We can then optimize over all probability distributions on $\{O_f(z):z\in\cZ_0\}$. The minimax risk then depends only on the unweighted sampling hypergraph. %
The next proposition shows that these objects determine the minimax risks of the fixed design, constrained design and all designs problem.

\begin{proposition}[Sampling-design reduction of the three minimax problems]
\label{thm:three-problem-characterization}

Fix $f$, $a$, and $\cY(f)$.
\begin{enumerate}
\item[\textup{(i)}] For a fixed design $P$, both $\Rall(P)$ and
$\Runb(P)$ depend on $P$ only through $Q_{P,f}$.
\item[\textup{(ii)}]
For a restricted design class $\cC$, the joint minimax risk $\Rallstar(\cC)$ depends on the class only through $\mathcal Q_f(\cC)$. The corresponding risk $\Runbstar(\cC)$ depends only on $\mathcal Q_f(\cC_+)$, where $\cC_+=\cC\cap\cP_+(f,a)$.
\item[\textup{(iii)}]
When all designs on $\cZ_0$ are allowed, both $\Rallstar(\cP(\cZ_0))$ and $\Runbstar(\cP(\cZ_0))$ depend on the assignment and exposure model only through $\cH_{\mathrm{obs}}(f;\cZ_0)$.
\end{enumerate}
\end{proposition}

The proof is in Appendix~\ref{app:observation-law-proofs}. Proposition~\ref{thm:three-problem-characterization} lets us work with
induced sampling designs $Q_{P,f}$ in place of assignment laws $P$.  To obtain a lower bound at a given induced sampling design, we compute the joint missingness probabilities
of a perturbation packing and substitute them into
the lower bound from disjoint perturbations (Theorem~\ref{thm:fixed-joint-missingness}).  Taking the supremum over
packings gives a lower bound for that law; taking the infimum over the
laws attainable by a design class gives a lower bound for its joint
minimax risk.  Corollary~\ref{thm:restricted-joint-missingness} in
Appendix~\ref{app:fixed-design-refinements} states these bounds for both
estimator classes.

\begin{remark}
The minimax risk of the fixed design problem (part (i)) depends on the full induced sampling design $Q_{P,f}$ and cannot be reduced to just the supported hyperedges or the first-order inclusion probabilities.  Theorem~\ref{thm:cluster-with-noise} in
Appendix~\ref{subsec:app-cluster-with-noise} gives two designs for the
no-interference ATE with the same sampling hypergraph and the same
first-order inclusion probabilities, but the unbiased fixed-design risks of these two designs are of order $1/n$ and one, respectively.
\end{remark}

\paragraph*{Comparison with the conflict graph}
\label{par:conflict-graph-comparison}
The sampling hypergraph records which collections of potential outcomes can be observed together. To relate the sampling hypergraph to the conflict graph used by \citet{KandirosHarshawSavje2026}, we distinguish between observing pairs of potential outcomes and observing entire collections of potential outcomes. 
Formally, for an exposure mapping $f$ and a nonempty feasible assignment set $\cZ_0$, we say two distinct potential outcomes $Y_k,Y_\ell$ have a \emph{pairwise conflict} if $I_k(z)I_\ell(z)=0$ for every $z\in\cZ_0$: at least one potential outcome in the pair must remain missing under every feasible assignment. For each $k\in K$, define the \emph{revelation set} $\cZ_k=\{z\in\cZ_0:I_k(z)=1\}$ as the set of treatment assignments that reveal the potential outcome $Y_k$. Thus $Y_k$ and $Y_\ell$ have a pairwise conflict if and only if $\cZ_k\cap\cZ_\ell=\varnothing$. The conflict-graph of \citet{KandirosHarshawSavje2026}, is defined as $H_{\rm conf}$, with vertex set $K$, with an edge joining each conflicting pair of potential outcomes.

For a collection $S\subseteq K$ of potential outcomes attainable under $\cZ_0$, we call $S$ \emph{pairwise observable} if, for every pair of distinct $k,\ell\in S$, there exists an assignment $z\in\cZ_0$ with $I_k(z)I_\ell(z)=1$. The assignment may differ across pairs. Equivalently, $\cZ_k\cap\cZ_\ell\ne\varnothing$ for every distinct $k,\ell\in S$. Thus $S$ is pairwise observable if and only if it is an independent set in $H_{\rm conf}$. On the other hand, we call the set $S$ \emph{jointly observable} if there exists one assignment $z\in\cZ_0$ such that $I_k(z)=1$ for every $k\in S$. Equivalently, $\bigcap_{k\in S}\cZ_k\ne\varnothing$, or $S\subseteq O_f(z)$ for some $z\in\cZ_0$. In the sampling hypergraph, this means that $S$ is contained in a hyperedge.

Every jointly observable collection of potential outcomes is pairwise observable, but the converse need not hold. Thus the conflict graph implicitly records which pairs can be observed together by recording pairwise conflicts, whereas the sampling hypergraph records which entire collections of potential outcomes can be observed together. Appendix~\ref{sec:pairwise} gives a formal comparison of the conflict graph and the sampling hypergraph, establishes sufficient conditions under which pairwise observability implies joint observability, and provides examples where pairwise observability does not imply joint observability.

\subsection{Lower bounds valid for every design}
\label{subsec:all-designs-lower}
\label{subsec:simple-hypergraph-bounds}

Proposition~\ref{thm:three-problem-characterization} shows that the joint
minimax risks over $\cP(\cZ_0)$ depend only on the sampling hypergraph.
We give two lower-bound constructions based on this hypergraph, for estimands
that contrast two arms, and use them in the no-interference, neighborhood,
and bipartite examples in Sections~\ref{sec:no-interference},
\ref{sec:neighborhood-interference}, and~\ref{sec:bipartite-interference}.  In
some examples, perturbations within one arm give a matching lower bound; in
others, we need perturbations in both arms.  The one-arm result also applies
to one arm of a contrast with more than two arms.  Lower bounds for general
linear estimands are given in the general hypergraph lower bound (Theorem~\ref{thm:global-hypergraph-lower}) in
Appendix~\ref{app:hypergraph-certificates}. The next theorem uses the largest number of potential outcomes from one
estimand arm that a single feasible assignment can observe. We call this
number the \emph{sampling rank} of that arm.

\begin{theorem}[Reciprocal sampling-rank lower bound]
\label{cor:rank-lower}
Let $V\subseteq K_a$ contain $N_V=|V|$ relevant potential outcomes from one
estimand arm, each with estimand coefficient satisfying $|a_k|=A/N_V$. Write $\cH_{\mathrm{obs}}(f;\cZ_0)=(K,\mathcal E)$. The sampling rank
$r(V)=\max_{E\in\mathcal E}|E\cap V|$ is the largest number of outcomes
in $V$ observable under a single feasible assignment. If $r(V)\geq1$, then
\begin{equation}
  \Runbstar\bigl(\cP(\cZ_0)\bigr)
  \geq A^2\left\{\frac1{r(V)}-\frac1{N_V}\right\}.
  \label{eq:rank-unb}
\end{equation}
For all estimators, if $1\leq r(V)\leq(1-\epsilon)N_V$ for some
$0<\epsilon<1$, then
\begin{equation}
  \Rallstar\bigl(\cP(\cZ_0)\bigr)
  \geq \frac{A^2c_\epsilon}{r(V)}
  \label{eq:rank-all}
\end{equation}
for a constant $c_\epsilon>0$ depending only on $\epsilon$.
\end{theorem}

The rank $r(V)$ is the largest exposure count in the chosen arm over all feasible assignments.
For an arm corresponding to exposure $d$ across all eligible units,
this is $r(V)=\max_{z\in\cZ_0}n_d(z)$. The theorem is informative when the sampling rank $r(V)$ is smaller than the
arm size $N_V$.  If $r(V)=0$, the entire arm is unobservable: the all-estimator risk
is at least $A^2$, and no unbiased estimator exists.  If $1\leq r(V)\leq(1-\epsilon)N_V$, both lower bounds have order $A^2/r(V)$. At the other extreme, when $r(V)=N_V$, the one-arm bounds are uninformative. This can happen for each arm separately even when no assignment observes both arms in full.  We therefore need a lower bound that accounts for which parts
of the two arms can be observed together.
Specifically, for two arm contrasts, we count pairs of perturbation sets rather than individual potential outcomes. Choose $m$ pairs of perturbation sets $(B_b^+,B_b^-)$, where $B_b^+$ contains potential outcomes from the first arm and $B_b^-$ from the second, and all $2m$ sets are nonempty and pairwise disjoint. For this choice of pairs, define
\[
 \mathcal M=\max_{E\in\cH_{\mathrm{obs}}}
 \#\{b:E\cap B_b^+\neq\varnothing\ \text{and}\ E\cap B_b^-\neq\varnothing\}.
\]
Thus $\mathcal M$ is the largest number of pairs of perturbation sets for which a single assignment observes at least one potential outcome on both sides. Like the sampling rank $r(V)$ in the sampling-rank lower bound (Theorem~\ref{cor:rank-lower}), it is determined by the hypergraph once the perturbation sets are chosen.
\begin{theorem}[Lower bound from paired perturbation sets]
\label{thm:paired-sets-main}
For the pairs above, suppose that each of the $2m$ sets carries the same
estimand weight,
\[
 \sum_{k\in B_b^+}|a_k|=\sum_{k\in B_b^-}|a_k|=\frac{A}{m},
 \qquad b=1,\ldots,m.
\]
Then
\[
 \Rallstar\bigl(\cP(\cZ_0)\bigr)\geq A^2\frac{m-\mathcal M}{m^2}.
\]
The same bound holds for $\Runbstar\bigl(\cP(\cZ_0)\bigr)$.
\end{theorem}

Both theorems use the hypercube prior of
the lower bound from disjoint perturbations (Theorem~\ref{thm:fixed-joint-missingness}).  To prove
the sampling-rank lower bound (Theorem~\ref{cor:rank-lower}), we choose a perturbation packing
$\cB=\{B_1,\ldots,B_m\}$ and bound how many of its sets an allowed assignment
can intersect.  Observing one potential outcome in $B_b$ reveals its shared
Rademacher sign $\theta_b$; the other potential outcomes in $B_b$ need not be
observed.  If every set has coefficient magnitude $|c_b|=c$ and every
hyperedge intersects at most $r<m$ sets, then every assignment leaves at least
$m-r$ sets of potential outcomes entirely missing.  These missing
sets contribute at least $(m-r)c^2$ to the posterior variance.  Averaging over
any design gives the same lower bound on the Bayes risk and hence on the joint
minimax risk.  For the unbiased bound in Theorem~\ref{cor:rank-lower}, the
singleton packing instead uses
$\sum_{k\in V}\pi_{P,k}\leq r(V)$.  The full proof is in
Appendix~\ref{app:hypergraph-certificates}.

For Theorem~\ref{thm:paired-sets-main}'s proof, there is one independent Rademacher
sign for each of the $2m$ sets.  A pair contributes $(A/m)^2$ to the posterior
variance unless the assignment observes a potential outcome on both sides, and at
most $\mathcal M$ pairs can be resolved.  The bound for unbiased
estimators follows because they form a smaller class.  The full proof is in
Appendix~\ref{subsec:block-certificates}.
For a two-arm contrast with coefficients $\pm1/n$, the equal-weight condition
in Theorem~\ref{thm:paired-sets-main} requires all $2m$ sets to have the same
cardinality.  Theorem~\ref{thm:arm-paired-certificate} in
Appendix~\ref{subsec:block-certificates} allows the weights to vary across
pairs and gives the corresponding weighted lower bound.

\section{Minimax estimation under a fixed design}
\label{sec:fixed-design-ht}

In this section, we consider the fixed-design problem and give matching upper
and lower bounds for a class of designs with bounded positive same-exposure
dependence, under which the exposure Horvitz--Thompson estimator is
minimax-rate optimal \citep{HorvitzThompson1952,AronowSamii2017}.  We then
consider designs with stronger positive dependence, where the H-T estimator remains minimax optimal but a larger
joint-missingness lower bound is needed to match its variance. %
We conclude with an example in which the exposure HT estimator
is not minimax-rate optimal.  %

\subsection{Bounded positive dependence}
\label{subsec:fixed-ht-bounded-dependence}

Fix an exposure mapping \(f\), an estimand coefficient vector \(a\), and a
design \(P\) such that \(\pi_{P,k}>0\) for every \(k\in K_a\).  Write
\(I_{i,e}=\1\{f_i(Z)=e\}\) and
\(\pi_{P,i}(e)=P\{I_{i,e}=1\}\).  For distinct units \(i\) and \(j\), let
\[
 \phi_{P,ij}(e,e)
 =\frac{P\{I_{i,e}=1,I_{j,e}=1\}}
 {\pi_{P,i}(e)\pi_{P,j}(e)}-1.
\]
Thus \(\phi_{P,ij}(e,e)\) is the normalized covariance between the
indicators for observing \(Y_i(e)\) and \(Y_j(e)\).  Define
\begin{align}
 D_a(P)
 &=\sum_{(i,e)\in K_a}a_{i,e}^2
   \{\pi_{P,i}(e)^{-1}-1\}, \label{eq:fixed-diagonal}\\
 \mathrm{OD}_a^+(P)
 &=\sum_e\sum_{\substack{i\ne j\\(i,e),(j,e)\in K_a}}
   |a_{i,e}a_{j,e}|\{\phi_{P,ij}(e,e)\}^+.
 \label{eq:fixed-positive-off-diagonal}
\end{align}
$D_a(P)$ is the diagonal term in the variance of the exposure HT
estimator.  $OD_a^+(P)$ is the positive off-diagonal contribution from pairs
of units at the same exposure.  Let
\(L_a=|\{e:(i,e)\in K_a\text{ for some }i\}|\) be the number of exposure
levels represented in the estimand.

\begin{definition}[Bounded positive same-exposure dependence]
\label{def:bounded-positive-dependence}
For an exposure mapping \(f\) and estimand \(a\), a design \(P\) has bounded
positive same-exposure dependence, or bounded positive dependence for short,
if $\mathrm{OD}_a^+(P)\leq C_0D_a(P)$
for a fixed constant \(C_0\). %
\end{definition}
The condition restricts only positive dependence between observation indicators for the same exposure. It also holds for designs when same-exposure indicators are pairwise independent or nonpositively correlated. Positive dependence is also allowed, provided that \(\mathrm{OD}_a^+(P)\) is of the same order as \(D_a(P)\). No separate condition on negative covariance terms is needed: Lemma~\ref{supp:lem:same-exposure-sign} bounds their total magnitude by \(D_a(P)+\mathrm{OD}_a^+(P)\).

\begin{theorem}[Bounded positive dependence]
\label{thm:fixed-bounded-positive-dependence}
Suppose that \(P\) has bounded positive same-exposure dependence with
\(\mathrm{OD}_a^+(P)\leq C_0D_a(P)\), and that \(L_a\leq L_0\).  Then
\begin{equation}
 D_a(P)\leq\Runb(P;f,a)
 \leq2L_0(1+C_0)D_a(P).
 \label{eq:fixed-bpd-unbiased}
\end{equation}
If, in addition,
\(\pi_{P,k}\geq c>0\) for every \(k\in K_a\), with \(c\) independent of the
population size, then
\begin{equation}
 cD_a(P)\leq\Rall(P;f,a)\leq\Runb(P;f,a)
 \leq2L_0(1+C_0)D_a(P).
  \label{eq:fixed-bpd-all}
\end{equation}
\end{theorem}

Thus \(\Runb(P;f,a)\asymp D_a(P)\), and when the inclusion probabilities are
bounded away from zero, \(\Rall(P;f,a)\asymp\Runb(P;f,a)\asymp D_a(P)\), with
constants depending only on \(C_0\), \(L_0\) and \(c\). The upper bounds are attained by the exposure Horvitz--Thompson estimator
\[
 \widehat\tau_{\mathrm{HT},P}
 =\sum_{(i,e)\in K_a}a_{i,e}
   \frac{I_{i,e}Y_i(e)}{\pi_{P,i}(e)},
\]
which is unbiased. %
Write $V_{\mathrm{HT}}(P)=\sup_{Y\in\cY(f)}
\Var_P\{\widehat\tau_{\mathrm{HT},P}\mid Y\}$ for its worst-case variance. Proposition~\ref{supp:prop:exposure-wise-ht} gives
\begin{equation}
 V_{\mathrm{HT}}(P)
 \leq 2L_a\{D_a(P)+\mathrm{OD}_a^+(P)\}.
 \label{eq:fixed-ht-upper}
\end{equation}
The first inequality in equation~\eqref{eq:fixed-bpd-unbiased} is the unbiased
singleton lower bound of Corollary~\ref{cor:singleton-lower}.  %

We now show that commonly used designs, estimands, and exposure mappings
satisfy bounded positive dependence.  For these examples, exposure HT is
minimax-rate optimal over both all estimators and unbiased
estimators.  Detailed calculations are in
Appendix~\ref{subsec:app-fixed-ht-examples}.

\medskip\noindent\textit{Example 1 (No-interference ATE).}
Consider the Average Treatment Effect (ATE) $(1/n) \sum_i Y_i(1) - Y_i(0)$ with exposure mapping \(f_i(Z)=Z_i\).
For each \(e\in\{0,1\}\) and distinct units \(i\ne j\), the observation
indicators \(I_{i,e}\) and \(I_{j,e}\) are independent under Bernoulli
assignment.  Under complete randomization, their covariance is nonpositive.
Thus \(\mathrm{OD}_a^+(P)=0\), and both designs have bounded positive
dependence for the ATE.  For Bernoulli assignment,
\(D_a(P)=n^{-1}\{p^{-1}+(1-p)^{-1}-2\}\), and the same formula holds with \(p\)
replaced by \(n_1/n\) under complete randomization of \(n_1\) units.  If \(p\), or \(n_1/n\), is bounded away from zero and one, Theorem~
\ref{thm:fixed-bounded-positive-dependence} gives
\(\Rall(P;f,a)\asymp\Runb(P;f,a)\asymp n^{-1}\).

\medskip\noindent\textit{Example 2 (Bounded-degree neighborhood interference).}
Consider Bernoulli assignment with treatment probability bounded away from
zero and one, and let the maximum degree of the interference graph be
bounded.  Under the identity exposure mapping, the direct treatment effect
compares \((1,0_{N_i})\) with \((0,0_{N_i})\), and the total treatment effect
compares \((1,1_{N_i})\) with \((0,0_{N_i})\).  %
Each closed neighborhood \(N_i^+\) overlaps only a bounded
number of other closed neighborhoods.  Consequently,
\[
 D_a(P)\asymp n^{-1},
 \qquad
 \mathrm{OD}_a^+(P)=O(n^{-1}).
\]
The relevant exposure probabilities are bounded away from zero, so Theorem~
\ref{thm:fixed-bounded-positive-dependence} again gives
\(\Rall(P;f,a)\asymp\Runb(P;f,a)\asymp n^{-1}\).

The next example shows that uniformly positive exposure probabilities are sufficient, but not necessary, for fixed design minimax optimality of HT.  

\medskip\noindent\textit{Example 3 (A rare exposure on a regular graph).}
For the DTE on a \(d\)-regular graph, let the Bernoulli treatment probability
be \(p=1/(d+2)\).  Bounded degree is sufficient but not necessary
for bounded positive dependence: even when \(d\) grows, this design, exposure
mapping, and estimand have bounded positive dependence because
\(D_a(P)\asymp(d+1)/n\) and
\(\mathrm{OD}_a^+(P)=O\{(d+1)/n\}\).  However, The treated-exposure probability is of order \(1/(d+1)\), and is not bounded away from
zero.  Direct calculation gives
\(D_a(P)+\mathrm{OD}_a^+(P)\asymp(d+1)/n\), and
Proposition~\ref{prop:dte-regular-all-estimator} gives a matching
all-estimator lower bound when \(n\geq6(d+2)^3\).  Hence
\(\Rall(P;f,a)\asymp\Runb(P;f,a)\asymp(d+1)/n\).

\subsection{When the singleton lower bound is not enough}
\label{subsec:fixed-ht-stronger-dependence}

We now consider designs for which positive same-exposure dependence is not
bounded by \(D_a(P)\).  We give a sufficient condition under which exposure
HT remains minimax-rate optimal.  However, the singleton lower bound of the previous
subsection does not match the HT upper bound for these designs, so we use a
joint-missingness packing-based lower bound.  For a perturbation packing
\(\cB\), with signs chosen to agree with the signs of the estimand
coefficients, write
\[
 L_{\cB}(P)
 =\sum_{B\in\cB}\left(\sum_{k\in B}|a_k|\right)^2\rho_P(B).
\]
The lower bound from disjoint perturbations (Theorem~\ref{thm:fixed-joint-missingness}(i)) gives
\(\Rall(P;f,a)\geq L_{\cB}(P)\).  We say that positive dependence is
packing-matched when \(\mathrm{OD}_a^+(P)=O\{L_{\cB}(P)\}\) for such a
packing.

\begin{theorem}[HT optimality under packing-matched dependence]
\label{thm:fixed-ht-comparison}
Suppose that \(L_a\leq L_0\), \(\pi_{P,k}\geq c\) for every \(k\in K_a\),
and
\begin{equation}
 \mathrm{OD}_a^+(P)\leq C_0L_{\cB}(P),
 \label{eq:packing-matches-positive-dependence}
\end{equation}
where \(L_0\), \(c>0\), and \(C_0\) are independent of the population size.
Then
\begin{equation}
 \begin{aligned}
 \max\{cD_a(P),L_{\cB}(P)\}
 &\leq\Rall(P;f,a)\leq\Runb(P;f,a)\\
 &\leq V_{\mathrm{HT}}(P)
 \leq2L_0\{c^{-1}+C_0\}\Rall(P;f,a).
 \end{aligned}
 \label{eq:fixed-ht-comparison}
\end{equation}
Thus exposure HT remains rate optimal when joint missingness matches positive dependence.
\end{theorem}

The theorem replaces the diagonal lower bound by \(L_{\cB}(P)\).  Therefore,
even when \(\mathrm{OD}_a^+(P)/D_a(P)\) grows, HT remains rate optimal if a
packing of jointly missing potential outcomes has the same order as
\(\mathrm{OD}_a^+(P)\).  The next two examples have this property and do not
satisfy bounded positive dependence when the cluster size grows.

\medskip\noindent\textit{Example 4 (Cluster randomization).}
Partition the population into
\(\xi\) administrative clusters of size \(m\), and assign each cluster
entirely to treatment with probability \(p\), independently across clusters.
Assume that $p$ is bounded away from zero and one.
The interference graph is empty and the exposure is \(f_i(Z)=Z_i\); the
clusters only determine which units are assigned together.  For the ATE,
outcomes in the same cluster and treatment arm are observed together.  Hence
\(\mathrm{OD}_a^+(P)/D_a(P)\) grows with \(m\).  A packing with one set for
each cluster and arm satisfies
\(L_{\cB}(P)\asymp\mathrm{OD}_a^+(P)\), which gives
\(\Rall(P;f,a)\asymp\Runb(P;f,a)\asymp1/\xi\).

\medskip\noindent\textit{Example 5 (Randomized saturation).}
For the policy effect $\mu_{\Gamma_q}(Y)-\mu_{\Gamma_{q'}}(Y)$,
clusters independently draw a saturation $Q_g\in\{q,q'\}$ and treat
$mQ_g$ units uniformly. The exposure $f_i(Z)=(Z_i,Q_{g(i)})$ records
a unit's treatment and saturation. A cluster at saturation $q'$ reveals
no outcomes at $q$. A packing by cluster and saturation therefore
matches the positive-dependence term and gives both risks of order
$1/\xi$ under the conditions in Appendix~\ref{supp:ht-upper-bounds}.

\subsection{When exposure Horvitz--Thompson is not rate optimal}
\label{subsec:fixed-ht-saturation-conditioning}

The preceding results match the HT variance either with the singleton lower
bound or with a lower bound from jointly missing blocks.  In this subsection,
we give a two-stage design for which exposure HT is suboptimal but a
conditional HT estimator is minimax-rate optimal. Consider a finite population of \(n=\xi m\) units partitioned into \(\xi\)
administrative clusters of common even size \(m\).  Each administrative
cluster contains \(m/2\) two-person households.  Interference is confined to
households, so the interference graph \(G\) is the matching formed by the
household pairs.  Let \(h(i)\) denote the other member of unit \(i\)'s
household.  We use the identity exposure mapping on this two-person
neighborhood, $f_i(Z)=(Z_i,Z_{h(i)})\in\{0,1\}^2$.
The estimand compares the two target exposures \((1,1)\), in which both
household members are treated, and \((0,0)\), in which both are untreated:
\[
 \tau(Y)=\frac1n\sum_{i=1}^n\{Y_i(1,1)-Y_i(0,0)\}.
\]
The treatment is assigned in two stages.  First, each administrative cluster
independently draws a saturation \(Q_g\in\{q_L,q_H\}\), where
\(P(Q_g=q_H)=r\) and \(P(Q_g=q_L)=1-r\).  Second, exactly \(mQ_g\) units are
treated uniformly within the cluster.
Write $\pi_i(e\mid q)=P\{f_i(Z)=e\mid Q_{g(i)}=q\}$ for the conditional
exposure probability and $\pi_i(e)=P\{f_i(Z)=e\}$ for its marginal
counterpart. The cluster's realized saturation is observed from $Z$.
The conditional HT estimator $\widehat\tau_{\mathrm{cond}}$ uses
$\pi_i(e\mid Q_{g(i)})$ in place of the marginal probabilities in the HT denominators.
Appendix~\ref{subsec:app-fixed-ht-saturation-conditioning} gives the
estimator and the two sets of probabilities.

\begin{proposition}[Conditioning on randomized saturation]
\label{prop:fixed-ht-saturation-conditioning}
Under the design above, assume that $r,q_L,q_H$ are bounded away from zero and one, that $q_L<q_H$ with $q_H-q_L$ bounded away from zero, and that $mq$ is an integer satisfying $2\leq mq\leq m-2$ for each $q\in\{q_L,q_H\}$, with all bounds uniform in the population size. Then
\[
 \Rall(P;f,a)\asymp\Runb(P;f,a)\asymp\frac1n,
 \qquad
 V_{\mathrm{HT}}(P)\asymp\frac1\xi.
\]
\end{proposition}

The estimator \(\widehat\tau_{\mathrm{cond}}\) is unbiased and
has worst-case variance of order \(1/n\), so it attains the rates in
Proposition~\ref{prop:fixed-ht-saturation-conditioning}.  By contrast,
exposure HT uses the marginal exposure probabilities and retains the
cluster-level variation induced by the shared \(Q_g\).  Consequently, when
the cluster size \(m=n/\xi\) grows, exposure HT is not minimax-rate optimal.
The proof is in
Appendix~\ref{subsec:app-fixed-ht-saturation-conditioning}, and
Appendix~\ref{subsec:app-fixed-ht-no-interference-saturation} gives the same
separation under no interference. \(\widehat\tau_{\mathrm{cond}}\)  uses observed design information in the sense of \citet{PashleyBasseMiratrix2021}

\section{Experiments without interference}
\label{sec:no-interference}In this section, we illustrate the lower and upper bounds for the Average Treatment Effect (ATE) under no interference.  We consider two problems: the unrestricted
problem, where all designs on $\{0,1\}^n$ are allowed, and a constrained
problem, where the number of treated units is fixed or bounded above by $m$.  The
unrestricted minimax rate is $1/n$ and with exactly $m$ treated units, the rate is $1/m+1/(n-m)$.

\subsection{No-interference ATE benchmark}
\label{subsec:ate-benchmark}
Let $\cZ=\{0,1\}^n$, $f_i(z)=z_i$, and
\begin{equation}
  \tau_{\rm ATE}(Y)
  =\frac1n\sum_{i=1}^n\{Y_i(1)-Y_i(0)\}.
  \label{eq:ate-target}
\end{equation}

\begin{theorem}[Joint minimax risk of ATE without interference]
\label{thm:ate-joint-benchmark}
Let $\cZ = \{0,1\}^n$.  For $\tau_{\mathrm{ATE}}$, over all designs,
$\Rallstar(\cP(\cZ))=(4+o(1))/n$ over all estimators and
$\Runbstar(\cP(\cZ))=4/n$ over unbiased estimators.
\end{theorem}

Thus $\Runbstar(\cP(\cZ))/\Rallstar(\cP(\cZ))=1+o(1)$: unbiasedness
does not change the minimax rate, or even the first-order constant, for the
no-interference ATE. The unbiased minimax risk in Theorem~\ref{thm:ate-joint-benchmark} is attained by the estimation strategy  $(P_{1/2},\widehat\tau_{\rm HT})$ 
where  $P_{1/2}$ is the independent Bernoulli assignment with probability \(1/2\) and $\widehat\tau_{\rm HT}=2n^{-1}\sum_{i}\{Z_iY_i^{\rm obs}-(1-Z_i)Y_i^{\rm obs}\}$ is the usual HT estimator.
The all-estimator upper bound in Theorem~\ref{thm:ate-joint-benchmark} comes from
the strategy $(P_{1/2},\widehat\tau_{\alpha_n})$ with the shrunken estimator
\[
  \widehat\tau_{\alpha_n}
  =\frac{\alpha_n}{n}\sum_{i=1}^n
  \{Z_iY_i(1)-(1-Z_i)Y_i(0)\},
  \qquad
  \alpha_n=\frac{2\sqrt n}{\sqrt n+1}.
\]
Its worst-case risk is at most $4/(\sqrt n+1)^2$, a finite-sample
improvement over HT with the same first-order constant $4$.

For the all-design minimax risk over all estimators, \citet{SudijonoDobribanTchetgen2026}
give an exact characterization of the joint minimax risk and identify
its second-order term. On our bounded potential outcome class with scale $[-1,1]$, their expansion is
$4/n-4C_A n^{-4/3}+o(n^{-4/3})$, where $C_A>0$ is their Airy-function
constant. They show that independent Bernoulli assignment with treatment
probability $1/2$, paired with an optimal estimator, attains the
joint minimax risk at each finite $n$. Our proof in
Appendix~\ref{app:no-interference-proofs} uses the weighted sampling
hypergraph to obtain the sharp first-order lower bound and treats the
exact unbiased estimator class separately.

\subsection{No-interference ATE under a constraint on the number of treated units}
\label{subsec:ate-fixed-count}

We restrict the number of treated units, either exactly or through a budget.  For
$m\in\{1,\ldots,n-1\}$, let
\[
  \cZ_m=\left\{z\in\{0,1\}^n:\sum_i z_i=m\right\},
  \qquad
  \cZ_{\le m}=\left\{z\in\{0,1\}^n:\sum_i z_i\leq m\right\}
\]

be the exact-count and budget assignment spaces.

\begin{theorem}[Minimax risk of ATE under constraints on the number of treated units]
\label{thm:ate-fixed-count}
For the ATE in equation~\eqref{eq:ate-target} and every \(m\in\{1,\ldots,n-1\}\), the
following statements hold.
\begin{enumerate}[(i)]
\item The joint minimax risk over all designs with exactly $m$ treated units is
\begin{equation}
  \Rallstar(\cP(\cZ_m))
  \asymp
  \Runbstar(\cP(\cZ_m))
  \asymp
  \frac1m+\frac1{n-m}.
  \label{eq:ate-fixed-count-rate}
\end{equation}

\item Let $k_m=\min\{m,\lfloor n/2\rfloor\}$. The joint minimax risk over all designs with at most $m$ treated units is 
\begin{equation}
  \Rallstar(\cP(\cZ_{\le m}))
  \asymp
  \Runbstar(\cP(\cZ_{\le m}))
  \asymp
  \frac1{k_m}+\frac1{n-k_m}
  \asymp
  \frac1m.
  \label{eq:ate-budget-rate}
\end{equation}
\end{enumerate}
The implied constants are universal and do not depend on \(m\) or \(n\).
\end{theorem}

Thus, uniformly over $m$, unbiasedness does not change the minimax rate
under either constraint on the number of treated units.

The proof is in Appendix~\ref{subsec:app-fixed-count-ate}. Let $P_k$
be uniform on $\cZ_k$ and $S$ the treated set. Write
$\widehat\tau^{(k)}_{\rm DIM}=\bar Y_{S}(1)-\bar Y_{S^c}(0)$ for the
difference between the treated and control sample means.
The unbiased strategies $(P_m,\widehat\tau^{(m)}_{\rm DIM})$ and
$(P_{k_m},\widehat\tau^{(k_m)}_{\rm DIM})$ attain the exact-count and budget
rates, respectively, and are also rate optimal over all estimators.  Both
rates have the form $1/k+1/(n-k)$: the exact constraint fixes $k=m$, whereas the budget permits the balanced choice $k_m=\min\{m,\lfloor n/2\rfloor\}$.

For the lower bound, an exact-count assignment reveals $m$ treated and $n-m$
control potential outcomes, so the sampling ranks of the two arms
(Theorem~\ref{cor:rank-lower}) give the scale $1/\min\{m,n-m\}$.  Pairing the
two potential outcomes of each unit, which are never observed together
(Proposition~\ref{prop:crossblocklb}), sharpens this to $1/m+1/(n-m)$
(Remark~\ref{rem:fixed-count-fpc}).  Under the unrestricted assignment space
each arm has sampling rank $n$, but no assignment observes both
potential outcomes of a unit; the same pairing then gives the $1/n$ scale.

\begin{remark}[The conflict graph under a constraint on the number of treated units]
\label{rem:fixed-count-conflict-graph}
For \(2\leq m\leq n-2\), the pairwise conflict graph of both $\cZ_m$ and
$\cZ_{\le m}$ is the same as for the unrestricted ATE: the only pairwise
conflicts are between \(Y_i(1)\) and \(Y_i(0)\) of the same unit. A
lower bound that uses only this graph therefore gives only the no-interference scale \(1/n\).  The
sampling hypergraph also records the condition that every assignment in $\cZ_m$ must
reveal exactly \(m\) treated and \(n-m\) control potential outcomes, and this
restriction on hyperedge composition produces the rate
\(1/m+1/(n-m)\). The rate differs from \(1/n\)
when \(m=o(n)\) (Appendix~\ref{subsec:app-fixed-count-ate}).

\end{remark}

\section{Graph-aware and graph-unaware designs under partial interference}
\label{sec:partial-interference}
\label{sec:two-stage}

In this section, we consider partial interference, where the population is
partitioned into disjoint clusters. We compare the minimax risk over two
design classes. In the graph-aware class, every design is allowed, including
designs that use the cluster memberships to assign treatments.  In the
graph-unaware class, the design cannot use the memberships.  We show
that it suffices to consider permutation-invariant designs for the graph-unware risk. %
The graph-unaware problem arises when the
interference structure is unavailable at the design stage but is collected for
the analysis, as can happen when social contacts are measured during follow-up
\citep{AlexandriaHudgensAiello2023}.  It can also represent an situation
where individual-level randomization may be more practical even when the
memberships are known.  In both cases, we assume that the memberships are
available to the estimator.  Comparing the two minimax risks quantifies the price of not using the cluster memberships at the design stage.

\subsection{Partial interference, exposure mapping and estimands}
\label{subsec:two-stage-exposure}

The population is partitioned into $\xi$ clusters of common size $m$, so that
$n=\xi m$.  For a unit $i$, let $g(i)$ denote its cluster and let
$N_i=\{j:g(j)=g(i),\ j\neq i\}$ be its $m-1$ cluster peers.  Thus the
interference graph $G$ is a disjoint union of
cliques of size $m$.  We assume anonymous within-cluster interference:
$Y_i(z)$ depends on $z_{N_i}$ only through the number of treated peers
$C_i(z)=\sum_{j\in N_i}z_j$.  The exposure mapping is therefore
\begin{equation}
  f_i(z)=\bigl(z_i,C_i(z)\bigr)\in\cD=\{0,1\}\times\{0,\ldots,m-1\}.
  \label{eq:two-stage-exposure}
\end{equation}
Write an exposure level as $e=(t,a)$, where $t$ is own treatment and $a$
is the number of treated peers. Define the exposure mean as
$\mu_{t,a}=n^{-1}\sum_{i=1}^nY_i(t,a)$. The two estimands studied in this section are
\begin{equation}
  \tau_{\mathrm{TTE}}=\mu_{1,m-1}-\mu_{0,0},
  \qquad
  \tau_{\mathrm{DTE}}(c)=\mu_{1,c}-\mu_{0,c},
  \qquad 0\leq c\leq m-1.
  \label{eq:invariant-estimands}
\end{equation}
The total treatment effect compares a treated unit whose $m-1$ peers are all
treated with an untreated unit whose peers are all untreated.  The direct
treatment effect with $c$ treated peers holds the number of treated peers
fixed at $c$ and changes only the unit's own treatment.  In particular,
$\tau_{\mathrm{DTE}}(0)$ is the usual direct treatment effect with no treated
peers.
For an assignment $z$, let $k_g(z)$ be the number of treated units in cluster
$g$, and let $K_g=k_g(Z)$ be the random number of treated units in cluster
$g$.  We index exposures by the number of treated peers rather than by the
number of treated units in the cluster, $k_{g(i)}(z)=C_i(z)+z_i$: a contrast
at a fixed number of treated units in the cluster changes both the unit's own
treatment and its number of treated peers, so it is not a direct effect
(Remark~\ref{rem:peer-not-total} in Appendix~\ref{app:invariant-proofs}).

\subsection{Graph-unaware and permutation-invariant designs}
\label{subsec:invariant-designs}
\label{subsec:fixed-peer-count-arm}%

We now introduce the graph-unaware design class.  A graph-unaware design is a law
on $\{0,1\}^n$ that does not use the interference graph.  Such a law need not
be permutation invariant: for example, it could use unit labels or baseline
covariates.  We study the narrower permutation-invariant class because %
Proposition~\ref{prop:unknown-partition} below shows that it suffices to consider this class when the design must be chosen without the
cluster memberships.

\begin{definition}[Permutation-invariant design]
\label{def:label-invariant}
A design $P$ is \emph{permutation invariant}, or equivalently exchangeable, if
$P(\sigma z)=P(z)$ for every
permutation $\sigma$ of $[n]$ and every $z\in\{0,1\}^n$.  We write
$\cC_{\mathrm{inv}}$ for the class of permutation-invariant designs.
\end{definition}

Let $\Pi$ be the set of all partitions of $[n]$
into $\xi$ clusters of size $m$.  The true partition $\pi\in\Pi$ determines the
exposure mapping $f^\pi$, the outcome class $\cY(f^\pi)$ and the estimand
coefficients $a^\pi$.  In the graph-unaware problem, the experimenter first
chooses a single assignment law $P$ without knowing $\pi$.  The same law must
be used for every $\pi\in\Pi$, although the estimator may use the partition
after it is revealed at the analysis stage.

\begin{proposition}[Permutation invariance is sufficient when memberships are unknown]
\label{prop:unknown-partition}
For every fixed reference partition $\pi_0\in\Pi$,
\begin{equation}
  \inf_{P\in\cP(\cZ)}\ \sup_{\pi\in\Pi}\ \Rall(P;f^\pi,a^\pi)
  \;=\;
  \inf_{P\in\cC_{\mathrm{inv}}}\Rall(P;f^{\pi_0},a^{\pi_0}),
  \label{eq:unknown-partition}
\end{equation}
and the same identity holds with $\Runb$ in place of $\Rall$.
\end{proposition}

Proposition~\ref{prop:unknown-partition} says that permutation-invariant designs
are sufficient for the graph-unaware problem.  The left side of
equation~\eqref{eq:unknown-partition} also includes graph-unaware laws that are not
permutation invariant, such as heterogeneous Bernoulli assignment, but these
cannot improve the worst-partition minimax risk.

The problem is symmetric in the unit labels: renaming the units and their
cluster memberships together changes neither the outcome class nor the
estimand.  The proof therefore symmetrizes any fixed design by averaging it
over unit relabellings.
Because the true partition may be any element of $\Pi$, the worst-partition
risk dominates the average over relabelled partitions.  The random relabelling
is auxiliary randomization, so Lemma~\ref{lem:rao-blackwell} shows that the
symmetrized design has no larger risk.  The partition is
unknown only when the design is chosen and remains available to the estimator,
while the cluster size $m$ may be used at both stages.
Appendix~\ref{subsec:app-invariant-designs}
gives the proof. We next characterize permutation-invariant designs.
\begin{lemma}[Permutation-invariant designs randomize only the number of treated units]
\label{lem:invariant-representation}
A design is permutation invariant if and only if it draws the number $M$ of
treated units from a law on $\{0,1,\ldots,n\}$ and then treats a uniformly
chosen set of $M$ units.
\end{lemma}

Bernoulli assignment at a common rate and complete randomization of the
population are permutation invariant, as is any mixture of them.  Conditional
on $M$, the number $K_g$ of treated units in each cluster is hypergeometric.

\begin{theorem}[The price of graph-unaware assignment under partial interference]
\label{thm:invariant-tte}
\label{thm:invariant-dte}
Let $\xi\geq2$.
\begin{enumerate}
\item[\textup{(i)}] \emph{Graph-aware designs.}  The total treatment effect
and every direct treatment effect $\tau_{\mathrm{DTE}}(c)$,
$0\leq c\leq m-1$, satisfy
\begin{equation}
  \frac1\xi\leq\Rallstar(\cP(\cZ))\leq\frac4\xi,
  \qquad
  \frac2\xi\leq\Runbstar(\cP(\cZ))\leq\frac4\xi .
  \label{eq:invariant-all-designs}
\end{equation}
\item[\textup{(ii)}] \emph{Graph-unaware designs: total treatment effect.}
Over permutation-invariant designs, the minimax risks satisfy
\begin{equation}
  \frac18\min\Bigl\{1,\frac{2^m}\xi\Bigr\}
  \leq\Rallstar(\cC_{\mathrm{inv}})
  \leq\Runbstar(\cC_{\mathrm{inv}})
  \leq\min\Bigl\{4,\frac{2^{m+1}}\xi\Bigr\}.
  \label{eq:invariant-tte-invariant}
\end{equation}
\item[\textup{(iii)}] \emph{Graph-unaware designs: direct treatment effect
away from the extremes.}  Fix $\delta\in(0,1/2)$.  If
$\delta m\leq c\leq(1-\delta)m-1$ and $m\geq m_0(\delta)$, then
the minimax risks satisfy
\begin{equation}
  \Runbstar(\cC_{\mathrm{inv}})\asymp_\delta\frac{\sqrt m}{\xi},
  \qquad
  \Rallstar(\cC_{\mathrm{inv}})\asymp_\delta
  \min\Bigl\{1,\frac{\sqrt m}{\xi}\Bigr\}.
  \label{eq:invariant-dte-invariant}
\end{equation}
\item[\textup{(iv)}] \emph{Graph-unaware designs: direct treatment effect at
the two extremes.}  For $c=0$ and $c=m-1$, the minimax risks satisfy
\begin{equation}
  \Rallstar(\cC_{\mathrm{inv}})\asymp
  \Runbstar(\cC_{\mathrm{inv}})\asymp\frac1\xi,
  \qquad c\in\{0,m-1\}.
  \label{eq:invariant-dte-endpoints}
\end{equation}
\end{enumerate}
\end{theorem}

The theorem is driven by two properties of the hypergeometric distribution;
see Lemma~\ref{lem:count-bottleneck}.  For any fixed
$\delta\in(0,1/2)$ and $\delta m\leq\ell\leq(1-\delta)m$,
every permutation-invariant design satisfies
$P(K_g=\ell)\leq C_\delta m^{-1/2}$.  Also, conditional on $M$, at least one
of $K_g=0$ and $K_g=m$ has probability at most $2^{-m}$.  The first fact gives
the $\sqrt m$ cost in part~\textup{(iii)}, and the second gives the $2^m$ cost
in part~\textup{(ii)}.
The cluster counts required by the estimands explain these two costs.  The
total effect needs clusters with $K_g=0$ and $K_g=m$.  The direct effect with
$c$ treated peers needs $K_g=c$ for its untreated units and $K_g=c+1$ for its
treated units.  At $c=0$ or $c=m-1$, a sparse or dense Bernoulli design makes
both required counts common, so there is no rate cost. Consequently, for the total effect no permutation-invariant design is
consistent once $m$ is larger than a constant multiple of $\log\xi$, whereas
a membership-aware design retains rate $1/\xi$.  For a direct effect satisfying
$\delta m\leq c\leq(1-\delta)m-1$, the cost is $\sqrt m$ for unbiased
estimators and $\min\{\xi,\sqrt m\}$ for all estimators.

The upper bounds use the exposure HT estimator from
Section~\ref{sec:fixed-design-ht}.
When memberships are available, independently give each cluster one of the
two required treated counts with probability one half and then sample
uniformly within the cluster; this attains $1/\xi$.  Without memberships,
Bernoulli rates
$1/2$, $(c+1/2)/m$, and $1/(m+1)$ attain the displayed total-effect bound,
the direct-effect bound for $\delta m\leq c\leq(1-\delta)m-1$, and the bound
for $c=0$, respectively.  A fair coin that treats the whole
population or none of it gives the constant upper bound for the total effect.
Appendix~\ref{subsec:app-invariant-designs} gives the proof, including lower
bounds that apply uniformly to every permutation-invariant design.

\begin{remark}[The price of unbiasedness in the graph-unaware class]
\label{rem:invariant-price}
Unbiasedness has no rate cost for the total effect or for the direct
effects at the two extremes $c=0$ and $c=m-1$, with or without the
memberships.  For $\delta m\leq c\leq(1-\delta)m-1$, it has a cost under
permutation-invariant designs only when $\xi\lesssim\sqrt m$: the minimax risk
is then of constant order over all estimators but of order $\sqrt m/\xi$ over
unbiased estimators.  When $\xi\gtrsim\sqrt m$, both are of order
$\sqrt m/\xi$.
\end{remark}

\section{Neighborhood interference}
\label{sec:neighborhood-interference}
\label{subsec:dte-tte-neighborhood}

In this section, we study minimax risks under neighborhood interference for three different estimands. We consider the direct and total
treatment effects with the identity exposure mapping $f_i(z)=(z_i,z_{N_i})$ and a  spillover effect under the threshold exposure mapping
$  f_i(z)=\bigl(z_i,\1\{\sum_{j\in N_i}z_j\geq\theta\}\bigr),$
for an integer $\theta\geq1$. The lower bounds for these three estimands come
from different constraints on which sets of potential outcomes can be jointly observed, and the third estimand illustrates that pairwise observability alone can
miss such a constraint.

\subsection{The direct treatment effect}
\label{subsubsec:dte-independent-set-thinning}

The direct treatment effect is the effect of treating a unit vs not, when all its neighbors are untreated: %
\begin{equation}
  \tau_{\rm DTE}(Y)
  =\frac1n\sum_{i=1}^n
    \{Y_i(1,0_{N_i})-Y_i(0,0_{N_i})\}.
  \label{eq:dte-neighborhood-target}
\end{equation}
The treatment-arm potential outcome $Y_i(1,0_{N_i})$ is observed only when unit $i$ is treated and all its neighbors are untreated.  Hence, units whose treatment-arm potential outcomes are jointly observed form an
independent set of $G$. Conversely, treating an independent set while leaving all other units untreated reveals the treatment-arm potential outcome of every
unit in the independent set.  Therefore, the treatment arm's sampling rank is $\alpha(G)$, the independence number of $G$. The control arm  potential outcome $Y_i(0,0_{N_i})$ has sampling rank $n$, since the $z = 0$ reveals all $n$ control arm potential outcomes. We have the following result: %

\begin{theorem}[Minimax rates for the neighborhood DTE]
\label{thm:dte-rates}
There are universal constants $c,C>0$ such that the following hold.
\begin{enumerate}
\item[\textup{(i)}] For every graph $G$,
\[
  \frac{c}{\alpha(G)}
  \leq\Rallstar(\cP(\cZ);f,a_{\rm DTE})
  \leq\Runbstar(\cP(\cZ);f,a_{\rm DTE}).
\]
\item[\textup{(ii)}] If $G$ is $d$-regular, then
\begin{equation}
  \frac{c}{\alpha(G)}
  \leq
  \Rallstar(\cP(\cZ);f,a_{\rm DTE})
  \leq
  \Runbstar(\cP(\cZ);f,a_{\rm DTE})
  \leq
  C\,\frac{d+1}{n}.
    \label{eq:dte-regular-main-sandwich}
\end{equation}
\item[\textup{(iii)}] If $G$ is a uniformly random simple \(d\)-regular
graph, there are universal
constants \(d_0,c_0>0\) such that, whenever
\(d_0\leq d<n-c_0n/\log n\) and \(nd\) is even, with probability
\(1-o(1)\) as \(n\to\infty\),
\begin{equation}
  \frac{d}{n\log d}
  \lesssim
  \Rallstar(\cP(\cZ);f,a_{\rm DTE})
  \leq
  \Runbstar(\cP(\cZ);f,a_{\rm DTE})
  \lesssim
  \frac dn .
  \label{eq:dte-random-regular-bracket}
\end{equation}
\end{enumerate}
\end{theorem}

The following argument establishes the lower bound $c/\alpha(G)$.  If
$\alpha(G)\leq n/2$, applying the sampling-rank lower bound
(Theorem~\ref{cor:rank-lower}) to the treatment arm gives the stated lower bound.  If $\alpha(G)>n/2$, it is enough to obtain a lower bound of
order $1/n$. %
The two target potential
outcomes of each unit cannot be observed together, so their inclusion probabilities sum to at most one. Hence, summing the singleton lower bound
(Corollary~\ref{cor:singleton-lower}) over units gives the required
$1/n$ bound. Part (iii) combines the upper bound
in part (ii) with the estimate \(\alpha(G)\lesssim n\log d/d\) for a uniformly
random regular graph
\citep{FriezeLuczak1992,CooperEtAl2002,KrivelevichEtAl2001}.  The complete proof is given in Appendix~\ref{subsec:app-dte-ordered-thinning}.

The upper bounds in parts (ii) and (iii) are attained by a design that selects a random independent set \(\cI\), assigns its units to treatment or control by fair coins, and assigns all other units to control.  Since every neighbor of a unit in \(\cI\) is untreated, a unit in $\cI$ reveals the target potential outcomes $Y_i(1,0_{N_i})$ or $Y_i(0,0_{N_i})$.  The design is a variant of the
independent-set design of \citet{KarwaAiroldi2018}.  Their design builds a
maximal independent set by repeatedly selecting a random remaining unit and
removing its neighbors, and its inclusion probabilities generally have to be
approximated by Monte Carlo.  We instead use the following construction, whose
inclusion probabilities have a closed form.

\begin{definition}[Ordered proposal design]
\label{def:dte-ordered-proposal}
Fix an order \(\prec\) on $[n]$ and proposal probabilities
\(p_i\in(0,1]\), $i\in[n]$.  Each unit \emph{proposes} independently with
probability $p_i$.  A proposing unit $i$ is included in $\cI$ when no
neighbor $j\prec i$ proposes.  Units in \(\cI\) are independently assigned to
treatment or control with equal probability, and every other unit is left
untreated.  We write \(P_{\prec,p}\) for the resulting assignment law, where $p = \{p_i\}$ is the collection of proposal probabilities.  
\end{definition}

The set $\cI$ is an independent set because, whenever two neighboring units
propose, only the earlier one can be included.  Inclusion depends on earlier
\emph{proposals} rather than earlier inclusions, so, unlike the greedy
construction, $\cI$ need not be maximal.  For a fixed order $\prec$, unit
\(i\) is included in $\cI$ with probability
\[
  q_i=p_i\prod_{j\in N_i:j\prec i}(1-p_j),
\]
and then reveals each of its two target potential outcomes with probability
\(q_i/2\).  If the order $\prec$ is uniformly random and every unit proposes, then
\(q_i=1/(d_i+1)\), and unit \(i\) is included in $\cI$ exactly when it precedes all of
its neighbors in the ordering $\prec$. The inverse-probability estimator of $\tau_{\mathrm{DTE}}$ based on $\cI$ is
\[
 \widehat\tau_{\cI}
 =\frac2n\sum_{i=1}^n\frac{Y_i^{\rm obs}}{q_i}
 \left[
  \1\{i\in\cI,\ f_i(Z)=(1,0_{N_i})\}
  -\1\{i\in\cI,\ f_i(Z)=(0,0_{N_i})\}
 \right].
\]
For either target exposure $b\in\{(1,0_{N_i}),(0,0_{N_i})\}$, consistency gives
\[
 \1\{i\in\cI,\ f_i(Z)=b\}Y_i^{\rm obs}
 =\1\{i\in\cI,\ f_i(Z)=b\}Y_i(b).
\]
Because $\E\bigl[\1\{i\in\cI,\ f_i(Z)=b\}Y_i^{\rm obs}\bigr] =\frac{q_i}{2}Y_i(b),$ $\hat \tau_{\cI}$ is unbiased.
However, $\tau_{\cI}$ uses the auxiliary randomization $\cI$, it does not belong to the estimator class $\cT(P,f)$.  Therefore, we 
use its Rao--Blackwellized version
$\widehat\tau_{\rm RB}=\E[\widehat\tau_{\cI}\mid Z,Y_{O_f(Z)}]$.
Lemma~\ref{lem:rao-blackwell} shows that $\widehat\tau_{\rm RB}$ remains unbiased and has no larger variance. Proposition~\ref{prop:dte-ordered-proposal-thinning} in Appendix~\ref{subsec:app-dte-ordered-thinning} upper bounds the worst-case variance of $\widehat\tau_{\rm RB}$ for every order and choice of proposal probabilities.  %
Appendix~\ref{subsec:app-dte-ordered-thinning} also gives choices of order and
proposal probabilities for bounded-degree graphs, disjoint stars, and disjoint cliques of equal size.

\begin{remark}[Comparison with Kandiros, Harshaw and S\"avje]
\label{rem:dte-khs-regular}
At $q=\infty$, \citet{KandirosHarshawSavje2026} obtain the same
independence-number lower bound.  Under the identity exposure mapping and
unrestricted assignments, the induced conflict graph on the
estimand-relevant potential outcomes characterizes their joint
observability, so their conflict-graph bound agrees with our sampling-rank
bound.  For the upper bound, their conflict-graph design and our ordered
proposal design are related but distinct constructions: The conflict graph design operates on the exposure vertices of the conflict graph, whereas the ordered proposal design operates directly
on the units of $G$.  Both attain order $(d+1)/n$ on regular graphs and hence
give the same logarithmic bracket for random regular graphs.
\end{remark}

\paragraph*{A rate price of unbiasedness}
\label{subsubsec:dte-price-core-periphery}
In the graph families considered above, unbiased and all estimator classes have the same minimax rate. This need not hold on every graph.  Consider a \emph{core-periphery} family with a core of $S$ units that forms a
clique, $1\leq S<n$, and let the remaining $n-S$ units form the periphery.
Every periphery unit is adjacent to every core unit, and the subgraph induced by
the periphery has maximum degree at most a fixed constant $r$.

\begin{theorem}[Unbiasedness can cost a rate]
\label{cor:dte-price-is-a-rate}
For a core--periphery graph,
\[
  \Rallstar(\cP(\cZ);f,a_{\rm DTE})\asymp\frac1n+\frac{S^2}{n^2},
  \qquad
    \Runbstar(\cP(\cZ);f,a_{\rm DTE})\asymp\frac Sn,
\]
so that the ratio of the two minimax risks is of order $\min\{S,n/S\}$,
with constants depending only on $r$.  The price is largest, of order
$\sqrt n$, when $S\asymp\sqrt n$.
\end{theorem}

The gap follows from a bias--variance trade-off.  No assignment can simultaneously reveal
the treatment-arm potential outcomes of the core and the periphery.
A biased estimator can ignore the core, incurring bias of order $S/n$
to reduce the variance to order $1/n$. %
On the other hand, an unbiased estimator must estimate the core's contribution, which costs variance of order $S/n$.  The proof is given in
Appendix~\ref{subsec:app-core-periphery-dte}.

\begin{remark}[Comparison with the Kandiros--Harshaw--S\"avje bias--variance trade-off]
\label{rem:dte-khs-comparison}
\citet{KandirosHarshawSavje2026} also use a bias--variance trade-off in their upper bounds by dropping a set of conflict-graph vertices %
Under
their second-moment restriction ($q=2$), when the conflict graph is a star,
their upper bound for the unbiased conflict-graph design is of order
$1/\sqrt n$, while dropping the center gives an upper bound of order $1/n$.
Their lower bound shows that $1/n$ is the all-estimator minimax rate, but it
does not lower-bound the risk over unbiased estimators.  %
By contrast, Theorem~\ref{cor:dte-price-is-a-rate} lower-bounds every unbiased
stratergy and shows that unbiasedness can cost a factor of order
$\sqrt n$.
\end{remark}

\subsection{The total treatment effect}
\label{subsubsec:tte-packed-clusters}

We now study the minimax risk of the Total Treatment Effect (TTE) with and without a treatment budget. The TTE compares the potential outcome when treating the entire closed neighborhood
of unit $i$ with leaving it untreated:
\begin{equation}
  \tau_{\rm TTE}(Y)
  =\frac1n\sum_{i=1}^n
    \{Y_i(1_{N_i^+})-Y_i(0_{N_i^+})\}.
  \label{eq:tte-neighborhood-target}
\end{equation}
The TTE creates two distinct constraints on jointly observed potential outcomes. Hence, its minimax lower bound depends on two graph quantities. The first quantity identifies closed neighborhoods that are far enough to serve as
disjoint perturbation sets (Definition~\ref{def:perturbation-set}).  Define
\[
  \mathsf{Pack}_1(G)=
  \max_{C\subseteq V:\operatorname{dist}_G(c,c')\geq3}
  \sum_{c\in C}|B_1(c)|^2,
\]
where the distance condition applies to distinct units in $C$, which serve as neighborhood centers.  The quantity
$\mathsf{Pack}_1(G)$ is a weighted packing of disjoint closed neighborhoods.
On a $d$-regular graph, it equals $(d+1)^2$ times the maximum number of
centers at pairwise distance at least three. The second quantity measures how many potential outcomes from the two arms can be revealed simultaneously by a
single assignment.  For an assignment \(z\), let
\(T=\{i:z_i=1\}\) be the set of treated units and define
\[
 A_T=\{i:B_1(i)\subseteq T\},
 \qquad
 B_T=\{i:B_1(i)\subseteq V\setminus T\}.
\]
Thus $A_T$ and $B_T$ index the observed treatment- and control-arm outcomes.
The quantity $M(G)=\max_{T\subseteq V}\min\{|A_T|,|B_T|\}$ is the largest
number $r$ for which one assignment reveals at least $r$ outcomes from
each arm. For the upper bound, we consider graphs with \emph{restricted growth}
from radius one. With $B_r(i)$ the $r$-hop closed neighborhood, this
means $|B_{r+1}(i)|\leq\kappa|B_r(i)|$ for every unit $i$ and $r\geq1$,
where $\kappa$ is the restricted-growth coefficient. Let $d_i=|B_1(i)|-1$ and $\bar d=n^{-1}\sum_i d_i$.

\begin{theorem}[Minimax rates for the neighborhood TTE]
\label{thm:restricted-growth-gcr-tte}
\label{prop:tte-main-lower}
The following hold.
\begin{enumerate}
\item[\textup{(i)}] For every graph $G$ with $M(G)\geq1$,
\[
  \max\left\{\frac{\mathsf{Pack}_1(G)}{n^2},
             \frac{1}{16M(G)}\right\}
  \leq\Rallstar(\cP(\cZ);f,a_{\rm TTE})
  \leq\Runbstar(\cP(\cZ);f,a_{\rm TTE}).
\]
If $M(G)=0$, then $\Rallstar(\cP(\cZ);f,a_{\rm TTE})\geq1$.
\item[\textup{(ii)}] If $G$ has restricted-growth coefficient $\kappa$, then
\[
 \frac{\bar d+1}{\kappa n}
 \leq\Rallstar(\cP(\cZ);f,a_{\rm TTE})
 \leq\Runbstar(\cP(\cZ);f,a_{\rm TTE})
 \leq C(\kappa)\frac{\bar d+1}{n},
\]
where $C(\kappa)$ depends only on $\kappa$.  %
\end{enumerate}
\end{theorem}

The detailed proof is in
Appendix~\ref{subsec:app-tte-packed-clusters}.
 The proof uses Lemma~\ref{lem:restricted-growth-packing} which constructs a packing by repeatedly
choosing a largest remaining closed neighborhood centered at a unit and excluding other centers within
distance two.  The restricted growth condition then bounds the sum of neighborhood sizes over the excluded centers by $\kappa$ times the squared size of the chosen neighborhood, so $\mathsf{Pack}_1(G)\geq n(\bar d+1)/\kappa$, giving the lower bound in part (ii).

Graph-cluster randomization \citep{ugander2023} attains the upper bound in part (ii). Choose a maximal set $C$ of units at pairwise distance at
least three.  These units serve as cluster centers. Assign every unit to its
nearest center, and independently treat each cluster with probability $1/2$.  Restricted growth ensures that each closed neighborhood intersects only a bounded number of clusters.  Both TTE exposures therefore have probabilities bounded away from zero, and the exposure Horvitz--Thompson estimator has worst-case variance at most $C(\kappa)(\bar d+1)/n$.  Proposition~\ref{prop:restricted-growth-detail}
gives the details. Part (ii) holds for every restricted-growth graph, including graphs with unequal degrees.  For a restricted-growth regular graph, the rate is $(d+1)/n$, and $d$ may grow with $n$.

For graphs without restricted growth,
Appendix~\ref{app:tte-general-refinements} gives an upper bound using a distance-two
independent-set design with small proposal probabilities. It applies the ordered proposal construction to
$G^2$, so the selected units have disjoint closed neighborhoods, and assigns
each selected neighborhood entirely to treatment or control.  It gives an
unbiased upper bound on every graph.  On $d$-regular spectral
expanders, it gives the bracket
\[
 \frac dn\lesssim\Rallstar
 \leq\Runbstar
 \lesssim\min\left\{1,\frac{d^2+1}{n}\right\}.
\]

\begin{remark}[Comparison with conflict-graph bounds]
\label{rem:tte-khs-comparison}
For the TTE on a restricted-growth $d$-regular graph, the bounds of
\citet{KandirosHarshawSavje2026} give
\[
 \max\left\{\frac1n,\frac{d+1}{n^{3/2}}\right\}
 \ \lesssim\ R^*_{\rm all}
 \ \leq\ R^*_{\rm unb}
 \ \lesssim_{\kappa}\ \frac {d+1}n .
\]
Their stated lower and upper bounds can therefore differ by a factor of order
$\min\{d+1,\sqrt n\}$.  Our packing lower bound is of order $(d+1)/n$ and
closes this gap (see Appendix~\ref{subsec:tte-khs-comparison}). %
\end{remark}

\subsubsection{The total effect under a treatment budget}
\label{subsec:hard-budget-regular-tte}

Consider a budget constraint of at most $m$ treated units. With a treatment budget, the choice of which units to treat matters. To
reveal a treatment-arm potential outcome $Y_i(1_{N_i^+})$, the design must treat a unit along with
all its neighbors, so the budget should cover as many complete neighborhoods
as possible.  Fix a treatment budget $m$ and restrict the assignment set to  $\cZ_{\leq m}=\{z\in\{0,1\}^n:|T(z)|\leq m\},$
where $T(z)$ is the set of treated units.  The corresponding design class is
  $\cP(\cZ_{\leq m})$.  Under a given assignment, let $S$ be the set of units whose
treatment-arm potential outcomes are observed.  The design must then treat every unit
in the set
\[
N[S]=\bigcup_{i\in S}B_1(i),
\]
so any such set must satisfy $|N[S]|\leq m$.  The number of treatment-arm
potential outcomes that can be observed together therefore depends on how
the neighborhoods overlap.  When they overlap heavily, the same treated units
may reveal several outcomes; when they overlap less, the budget is used up
more quickly.

\begin{theorem}[Budgeted TTE under restricted growth]
\label{thm:hard-budget-regular-tte}
Let $G$ be $d$-regular with restricted-growth coefficient $\kappa$.  There is
a constant $c_{\rm rg}(\kappa)\in(0,1]$ such that, if $1\leq m\leq n/2$ and
$(d+1)/m\leq c_{\rm rg}(\kappa)$, then
\[
 c(\kappa)\frac{d+1}{m}
 \leq\Rallstar(\cP(\cZ_{\leq m});f,a_{\rm TTE})
 \leq\Runbstar(\cP(\cZ_{\leq m});f,a_{\rm TTE})
 \leq C(\kappa)\frac{d+1}{m}.
\]
\end{theorem}

The condition $(d+1)/m\leq c_{\rm rg}(\kappa)$ ensures that the budget is sufficiently large for positivity to hold in the
upper bound construction below.  For fixed $\kappa$, both minimax risks vanish
at rate $(d+1)/m$ when $d+1=o(m)$.
Without the budget, Theorem~\ref{thm:restricted-growth-gcr-tte} gives rate
$(d+1)/n$ on the same graph class.  When the budget is large enough relative
to the neighborhood size, it therefore replaces $n$ by $m$ in the minimax
rate. The lit-ball design attains the upper bound.

\begin{definition}[Lit-ball design]
\label{def:hard-budget-lit-ball}
Fix $1\leq\ell\leq\lfloor m/(d+1)\rfloor$.  Draw $C_{\rm lit}$ uniformly from
all $\ell$-element subsets of $[n]$.  Assign every unit in
$\bigcup_{i\in C_{\rm lit}}B_1(i)$ to treatment and the rest to control.
\end{definition}

Choosing $\ell\leq m/(d+1)$ respects the budget.  For the theorem, take
$\ell\asymp_\kappa m/(d+1)$, where the comparison constants depend only on
$\kappa$.  Units in $C_{\rm lit}$ reveal the treatment arm.  A unit reveals
the control arm when no selected unit lies within two hops.  Restricted growth
and the condition $(d+1)/m\leq c_{\rm rg}(\kappa)$ keep these control-exposure probabilities bounded away from zero.
Inverse-probability weighting over the sampled units gives an unbiased
estimator.  Its Rao--Blackwellized version depends only on the observed
assignment and outcomes and has variance of order $(d+1)/m$.

The lower bound uses a perturbation packing construction, which we describe
next.  Choose a maximal set $C$ of units at pairwise distance at least five.
The two-hop neighborhoods $B_2(c)$, $c\in C$, are then pairwise disjoint.
If $Y_i(1_{N_i^+})$ for a unit $i\in B_1(c)$ is observed, at least $d+1$
units inside $B_2(c)$ must be treated.  An assignment with at most $m$
treated units can therefore reveal treatment-arm outcomes in at most
$J=\lfloor m/(d+1)\rfloor$ of the closed neighborhoods $B_1(c)$, $c\in C$. When $|C|\geq2J$, partition $C$ into $2J$ groups $C_1,\ldots,C_{2J}$ of
nearly equal size.  For each group, define the perturbation set
\[
 \mathcal B_b
 =\left\{(i,1):i\in\bigcup_{c\in C_b}B_1(c)\right\},
 \qquad b=1,\ldots,2J.
\]
Thus $\mathcal B_b$ contains the treatment-arm potential outcomes
of all units in the closed neighborhoods centered at units in $C_b$.
These perturbation sets are disjoint because the neighborhoods $B_1(c)$,
$c\in C$, are pairwise disjoint.  Since the budget can reach at most $J$
neighborhoods, at least $J$ perturbation sets remain unobserved under every
assignment.  Restricted growth ensures that each set carries estimand weight
at least a constant multiple of $1/J$.  The joint-missingness bound therefore
gives risk of order
$ J(1/J)^2=1/J\asymp(d+1)/m.$
Appendix~\ref{subsec:app-hard-budget-tte} gives the explicit grouping construction and handles the case $|C|<2J$; it also gives the estimator and an extension to
budget-expanding graphs.

\begin{remark}[Why pairwise observability is not enough]
The budget creates a higher-order constraint that is not recorded by the
conflict graph.  If $m\geq2(d+1)$, any two treatment-arm potential outcomes
can be observed together by treating the union of their closed neighborhoods.
The induced conflict graph on the treatment arm therefore has no edges. A collection of
$r$ disjoint closed neighborhoods, however, requires $r(d+1)$ treated units and
cannot be observed jointly when $r(d+1)>m$.  The sampling hypergraph
records this restriction and gives the lower bound of order $(d+1)/m$.
\end{remark}

Appendix~\ref{subsec:app-budgeted-dte} gives an analogous result for
the minimax risk of the DTE under a treatment budget.  The budget and the
independence number jointly limit how many treatment-arm outcomes can be
observed. The treated-arm rank $\min\{m,\alpha(G)\}$ gives the lower bound.%

\subsection{A local complex-contagion contrast}
\label{subsec:complex-contagion-gap}

Consider the threshold exposure mapping
$f_i(z)=(z_i,\1\{C_i(z)\geq\theta\})$, where $C_i(z)=\sum_{j\in N_i}z_j$ is
the number of treated neighbors of unit $i$.  It records whether unit $i$ is
treated and whether at least $\theta$ of its neighbors are treated.  We study the minimax risk of a  contrast that compares treated units with at least $\theta$ treated neighbors to treated units with fewer that $\theta$ treated neighobors.  Writing
$b_{\rm hi}=(1,1)$ and $b_{\rm lo}=(1,0)$, the estimand is
\begin{equation}
 \tau_\theta(Y)=\frac1n\sum_{i=1}^n
       \{Y_i(b_{\rm hi})-Y_i(b_{\rm lo})\}.
 \label{eq:complex-contagion-target}
\end{equation}
All assignments are allowed.  In the graph families considered below, both
target exposures are attainable, and the relevant potential outcomes are
bounded by one.  At $\theta=1$, the contrast compares treated units with at
least one treated neighbor to treated units with none.  For $\theta\geq2$, it
measures the effect of crossing a threshold that requires several treated
neighbors, as in models of complex contagion
\citep{Granovetter1978,Watts2002,CentolaMacy2007}.  In the
terminology of \citet{ugander2023}, $b_{\rm hi}$ is absolute
$\theta$-neighborhood exposure to treatment.

\paragraph*{Observation constraint and regular-graph bounds}
Let $S$ be the set of units that reveal their low-exposure potential outcomes $Y_i(b_{\rm lo})$.
Every unit in $S$ is treated and has fewer than $\theta$ treated neighbors,
so it has at most $\theta-1$ neighbors within $S$.  Conversely, if every unit
in a set $S$ has at most $\theta-1$ neighbors within $S$, then treating all
units in $S$ and leaving all units in $S^c$ untreated reveals the low-exposure
potential outcome of every unit in $S$.  The
largest number of low-exposure potential outcomes that one assignment can
reveal is therefore the $(\theta-1)$-dissociation number
\[
 \operatorname{diss}_{\theta-1}(G)
 =\max\{|S|:S\subseteq V,\ \Delta(G[S])\leq\theta-1\}.
\]
This is the sampling rank of the low-exposure arm; at $\theta=1$,
it reduces to $\alpha(G)$.

\begin{theorem}[Threshold-contrast bounds on regular graphs]
\label{thm:contagion-main-bracket}
Fix $\theta\geq1$.  Let $G$ be $d$-regular on $n$ units, with
$d\geq16\theta$ and $n\geq4\sqrt d$.  Then
\[
 c_\theta\max\left\{\frac1{\operatorname{diss}_{\theta-1}(G)},
                        \frac{\sqrt d}{n}\right\}
 \leq \Rallstar(\cP(\cZ);f,a_\theta)
 \leq \Runbstar(\cP(\cZ);f,a_\theta)
 \leq C_\theta\frac dn .
\]
The constants depend only on $\theta$.  The Bernoulli design with
$p=\theta/d$ and the exposure Horvitz--Thompson estimator attain the upper
bound.
\end{theorem}

For the lower bounds, the reciprocal dissociation term follows from the sampling-rank lower bound.  The $\sqrt d/n$ term comes from balancing two packing-based lower bounds.  Let $\pi_v$ be the probability that the low-exposure outcome
of unit $v$ is observed.  For a threshold $\pi$, if $\pi_v\leq\pi$ for many
units, a packing of singleton low-exposure potential outcomes gives a lower bound
of order $1/(\pi n)$.  If $\pi_v>\pi$ for many units, observing the low
exposure at $v$ leaves all but at most $\theta-1$ of its $d$ neighbors
untreated.  Packing blocks of these neighbors gives a lower bound of order
$d\pi/n$.  The two bounds balance at $\pi\asymp d^{-1/2}$, giving
$\sqrt d/n$.

For the upper bound, under a Bernoulli design with $p=\theta/d$, both target
exposures have probabilities of order $1/d$.  Exposure indicators are
independent for units at distance greater than two; bounding the remaining
covariance terms gives the $d/n$ upper bound.
Appendix~\ref{subsec:app-complex-contagion} gives the lower-bound proofs and
the upper-bound variance calculation.

On the random regular graphs considered in Theorem~\ref{thm:dte-rates}(iii),
$\alpha(G)\lesssim n\log d/d$ with high probability.  Since
$\operatorname{diss}_{\theta-1}(G)\leq\theta\alpha(G)$, the lower bound
strengthens to order $d/(n\log d)$, while the Bernoulli upper bound remains
of order $d/n$ leaving a logarithmic gap.

The $d/n$ upper bound is sharp on some regular graph families for every fixed
$1\leq\theta\leq d-1$.  For example, consider disjoint cliques of size $d+1$.
Independently in each clique, treat either $\theta$ or $\theta+1$ units.  The
first choice gives the selected units the low exposure, and the second gives
them the high exposure.  Both minimax risks are then of order $d/n$.  The
same rate holds on disjoint copies of the complete bipartite graph $K_{d,d}$
(Appendix~\ref{subsec:app-complex-contagion}).  Whether a regular graph family
has minimax risk of order $\sqrt d/n$ remains open.

\paragraph*{What pairwise observability misses}

For $\theta\geq2$, any two low-exposure potential outcomes are jointly observable,
but a unit and $\theta$ of its neighbors cannot simultaneously receive the low exposure. The sampling hypergraph retains this higher-order restriction, whereas its pairwise projection discards it.

\begin{corollary}[Lower-bound certificate gap on graphs of girth at least five]
\label{cor:contagion-conflict-gap}
Fix $\theta\geq2$ and let $G$ be $d$-regular with $d\geq16\theta$ and
girth at least five.  The %
pairwise conflict graph gives a lower bound of order $1/n$, whereas the
sampling hypergraph yields
\[
 \Rallstar(\cP(\cZ);f,a_\theta)\geq c_\theta\frac{\sqrt d}{n}.
\]

\end{corollary}
Lemma~\ref{lem:contagion-conflict-projection} shows that the only pairwise
conflicts are between the two exposures of the same unit.  The hypergraph lower bound is therefore stronger by a factor of order
$\sqrt d$, reflecting higher-order observation constraints that pairwise
observability does not record. At $\theta=1$, the low-exposure conflict
graph is instead exactly $G$ (Remark~\ref{rem:simple-contagion-conflict-exact}).

\section{Bipartite interference}
\label{sec:bipartite-interference}
\label{subsec:bipartite-marketplace}

In a bipartite experiment, treatment is assigned to one population and
outcomes are measured on another.  Examples include interventions at
pollution sources with outcomes measured in nearby populations
\citep{ZiglerPapadogeorgou2021}, and marketplace experiments in which sellers
are randomized and buyer outcomes are measured
\citep{HarshawEtAlBipartite2022}.  

Let \(U\) be the treatment population, with
\(|U|=M\), and let \(V\) be the outcome population, with \(|V|=n\).  The
interference graph $G$ is bipartite, with edges only between $U$ and $V$.
Each outcome unit \(v\) has a nonempty neighborhood \(A_v\subseteq U\) of
treatment units that can affect its outcome.  Write \(d_v=|A_v|\) for its
\emph{outcome-unit degree}, the number of treatment units connected to outcome unit $v$. For an assignment \(z\in\{0,1\}^U\), the exposure of outcome unit \(v\) is
the number of treated units in its neighborhood,
$f_v(z)=\sum_{u\in A_v}z_u$.  For two doses
\(0\leq k'<k\leq\min_v d_v\), we study the corresponding
dose-response contrast:
\begin{equation}
 \tau_{k,k'}(Y)=\frac1n\sum_{v\in V}\{Y_v(k)-Y_v(k')\}.
 \label{eq:bipartite-dose-family}
\end{equation}
Write \(a_{k,k'}\) for the estimand coefficient vector of \(\tau_{k,k'}\).
The potential outcomes are
bounded by one, and the design class contains all distributions on
\(\{0,1\}^U\).  We call the doses \(0\) and \(d_v\) the \emph{extreme} doses of
outcome unit \(v\). %
Doses strictly between \(0\) and \(d_v\) are \emph{interior} doses. %

\label{subsec:bipartite-shared-main}
 We study the minimax risk of a contrast between two interior doses when many outcome units randomly draw their neighborhoods from a common pool of $M$ treatment units.
For example, each buyer may depend on several sellers randomly drawn from a common seller population.  %

\begin{definition}[Uniform \(d\)-out bipartite interference graph]
\label{def:d-out-incidence}
Fix a pool \(U\) of treatment units with \(|U|=M\), and let \(d\leq M\)
be the number of treatment neighbors of each outcome unit.  Each of the \(n\)
outcome units draws its neighborhood independently of the others, uniformly at
random among the \(\binom Md\) subsets of \(U\) of size \(d\).
Thus $d_v=d$ for every outcome unit $v$, and its realized dose ranges from
$0$ to $d$.
\end{definition}

For two interior doses $k'<k$, define $s=\sqrt{\min\{k,d-k\}+1}$.
Note that $s^2-1 = \min\{k,d-k\}$ is the distance from $k$ to the closest extreme dose.
We consider a treatment pool of size $M\asymp n/s$.

\begin{theorem}[Nearby interior doses on a shared pool]
\label{thm:bipartite-central-main}
Assume a uniform $d$-out bipartite interference graph
(Definition~\ref{def:d-out-incidence}), with $n\geq3$,
$M\geq16$, and $3\leq d\leq\sqrt M$.  Let $0<k'<k<d$ with $k-k'\leq c_0s$ for
a fixed $c_0$, and suppose $M\asymp n/s$.
Then, with probability at least $1-2/n^2$ over the graph draw,
\[
 c\frac{s}{n}
 \leq\Rallstar(\cP(\{0,1\}^U);f,a_{k,k'})
 \leq\Runbstar(\cP(\{0,1\}^U);f,a_{k,k'})
 \leq C\frac{s\log n}{n}.
\]
The constants depend only on $c_0$ and the constants implicit in $M\asymp n/s$.
\end{theorem}

\paragraph*{Lower bound}
Let $J_k$ be the maximum number of outcome units whose potential outcomes $Y_v(k)$ are jointly observed under one assignment.  %
With high probability over the graph draw, $J_k\leq C(n/s+M+\log n)$,
giving $J_k\lesssim n/s$ when $M\asymp n/s$. Thus, every assignment leaves
most dose-$k$ potential outcomes unobserved when $s$ is large.  A paired-arm
perturbation packing uses this limit on joint observation to give the
lower bound of order $s/n$.  The construction of the packing groups outcome units and pairs their sets of dose-$k$ and dose-$k'$ potential outcomes, so that every assignment leaves at least one set completely missing in sufficiently many pairs. The lower bound of Theorem~\ref{thm:arm-paired-certificate} converts
the resulting packing into the risk bound.  The details are in Appendix~\ref{subsec:app-bipartite-fixed-count}.

\paragraph*{Upper bound}
The upper bound is attained by a completely randomized design on the treatment
population paired with the exposure Horvitz--Thompson estimator.  Let $t$
be the nearest integer to $M(k+k')/(2d)$, and treat a uniformly chosen
$t$-subset of $U$.  Conditional on the realized interference graph, every outcome unit
has exposure probabilities
\[
 \pi_j=\frac{\binom tj\binom{M-t}{d-j}}{\binom Md},\qquad j\in\{k,k'\},
\]
The condition $k-k'\leq c_0s$ ensures that this design gives both target
doses exposure probabilities of order $1/s$, so
$\pi_k,\pi_{k'}\asymp1/s$.  The estimator
\[
 \widehat\tau_{\rm fix}=\frac1n\sum_{v\in V}
 \left\{\frac{\1\{f_v(Z)=k\}}{\pi_k}
       -\frac{\1\{f_v(Z)=k'\}}{\pi_{k'}}\right\}Y_v^{\rm obs}
\]
is unbiased for every realized interference graph.  The diagonal contribution to its
variance is of order $s/n$.  For the off diagonal contribution, note that under a Bernoulli assignment, the number of treated units in the pool would fluctuate and shift the doses of all outcome units together.  Fixing the number of treated units removes this common fluctuation, and Appendix~\ref{subsec:app-bipartite-fixed-count} controls the remaining covariance from overlapping neighborhoods, giving an upper bound rate of $s\log n/n$. The upper and lower bounds differ by a factor of $\log n$. %
The proofs are in Appendix~\ref{subsec:app-bipartite-fixed-count}.

The adjacent-dose contrasts $\tau_{k,k-1}$ are the increments of the
dose--response curve $k\mapsto n^{-1}\sum_vY_v(k)$.  By
Theorem~\ref{thm:bipartite-central-main}, the difficulty of estimating these
increments depends on the dose $k$.  For adjacent doses near either extreme,
with $\min\{k,d-k\}$ bounded, the lower and upper bounds are of order $1/n$ and
$\log n/n$, respectively.  For adjacent doses near the middle dose $d/2$, they are of order $\sqrt d/n$ and $\sqrt d\log n/n$.  In each case, the comparison assumes a
treatment pool of size $M\asymp n/s$.

For interior-dose contrasts, the conflict graph can lose information needed
for the lower bound. An interior dose specifies how many treatment
units in a neighborhood are treated, so the corresponding potential
outcomes may be pairwise observable even though the entire collection is
not jointly observable.
For the extreme doses $0$ and $d_v$, the induced conflict graph on
the two extreme-dose arms determines which collections are jointly
observable: each extreme exposure fixes all treatment assignments in its
neighborhood, so a pairwise-observable collection of extreme-dose potential
outcomes is jointly observable, by the same argument as in the proof of
Theorem~\ref{thm:conflict-hypergraph-equivalence}.

\begin{corollary}[A two-arm lower-bound certificate gap]
\label{cor:bipartite-full-projection-gap}
\label{rem:bipartite-dose-projection}
Under the uniform $d$-out bipartite interference graph model, take adjacent interior doses
\(k'=k-1\), with \(0<k'<k<d\), \(3\leq d\leq\sqrt M\), and
\(M\asymp n/s\).  If
\(\binom Md/n^2\to\infty\), then with probability tending to one the
pairwise conflict graph gives a lower bound of order $1/n$, whereas the
minimax risks satisfy
\[
 c\frac{s}{n}\leq R^*_{\rm all}\leq R^*_{\rm unb}
                    \leq C\frac{s\log n}{n}.
\]
\end{corollary}

The hypergraph lower bound is stronger by a factor of order $s$, which
grows when both $k$ and $d-k$ grow and is of order $\sqrt d$ for central
doses.  For adjacent doses, target exposures of distinct outcome units
conflict only when their neighborhoods are identical.  The probability of
any duplicate neighborhood is at most $n(n-1)/\{2\binom Md\}$, which
vanishes under the stated condition. Thus the induced conflict graph on
the two target-dose arms consists only
of the $n$ same-unit exposure pairs and misses the higher-order restrictions
on joint observation.  The proof is in
Appendix~\ref{subsec:app-bipartite-fixed-count}.

\section{Discussion}
\label{sec:discussion}

In this paper, we study design based minimax risk of estimating causal effects under interference for three problems: When a design is fixed, under a restricted class of designs, and when the design is unrestricted. In all three cases, we consider the minimax risk over the class of all and unbiased estimators.  These
problems correspond to practical situations faced by an experimenter:
analyzing an experiment whose design is fixed, choosing a design subject to
constraints, and assessing what is achievable without those constraints.
In all three problems, the minimax risk depends on which estimand-relevant
potential outcomes can be observed together and how often they are observed.
The minimax risks are determined jointly by the interference structure, exposure
mapping, and estimand.  Hence, an optimal design must be tailored to the
estimand and to the constraints of the experiment.  For example, the direct
and total treatment effects impose different joint observation requirements
on the same graph and lead to different optimal design constructions.  This
raises a broader question: when can one estimation strategy be simultaneously
minimax-rate optimal for several estimands?

Our framework also shows which part of the estimation strategy limits the
rate:  the estimator, the design, or the requirement of unbiasedness.
For example, under randomized saturation, changing the estimator under the same design improves the rate from $1/\xi$ to $1/n$
(Proposition~\ref{prop:fixed-ht-saturation-conditioning}).  By contrast,
constraints on the number of treated units, budgets, and lack of knowledge of cluster memberships at design stage can make the design the bottleneck. For instance, estimating the total treatment effect when cluster
memberships are unknown leads to an exponential dependence on cluster size. %
(Theorem~\ref{thm:invariant-tte}).  Unbiasedness usually has no
rate cost in most of our examples, and no first-order cost for the no-interference
average treatment effect. But it is the bottleneck in the core--periphery
family.  There, a biased strategy can omit a small, difficult-to-observe part
of the population, whereas an unbiased strategy cannot, and the ratio
of the two minimax risks can be of order up to $\sqrt n$
(Theorem~\ref{cor:dte-price-is-a-rate}). 
Several questions remain.   Characterizing when unbiasedness has no minimax-rate cost remains open. We have identified settings where pairwise
observability implies joint observability, but a broader characterization
is open.  The lower and upper bounds also leave gaps for complex contagion
and bipartite dose contrasts.  Another question is how estimation
changes when the exposure mapping is misspecified or the interference
structure is only partly known.

\begin{acks}[Acknowledgments]
This paper was written by the authors with the assistance of large language models, which included suggesting arguments, contributing to drafting and revision, and writing code for computational checks. We checked every result ourselves, and any errors are ours.
\end{acks}
\bibliographystyle{spr-ims-nameyear}
\bibliography{references}
\newpage
\appendix
This appendix contains proofs and supporting results for the main text.
For clarity, we repeat selected statements immediately before their proofs,
retaining their original numbering.

\addtocontents{toc}{\protect\setcounter{tocdepth}{2}}
\begingroup
\begin{center}
\bfseries Appendix contents
\end{center}
\makeatletter
\@starttoc{toc}
\makeatother
\endgroup
\clearpage

\section{Supporting results for Section~\ref{sec:setup-strategy}}
\label{app:setup-proofs}
\subsection{Proof of Rao--Blackwellization of the estimator}
In this section, we prove Lemma~\ref{lem:rao-blackwell}, which removes auxiliary  randomization $U$ in an estimator $T(U,D)$ by using Rao--Blackwellization and conditioning on $D$. Here, $D = (Z, Y_{O_f(Z)})$ is the observed data.

\begin{resultrestatement}{Lemma~\ref{lem:rao-blackwell}}{Rao--Blackwellized estimator}
Suppose $U$ takes values in a finite or countable set $\mathcal U$, and the joint law of $U$ and assignment $Z$
is known and does not depend on $Y$.  Let $D=(Z,Y_{O_f(Z)})$ and let
$T(U,D)$ have finite second moment for every admissible table.  Then
$\widetilde T(D)=\E[T(U,D)\mid D]$ depends only on the observed data $D$ and satisfies
\[
 \E[(\widetilde T-\tau_a(Y))^2\mid Y]
 \leq \E[(T-\tau_a(Y))^2\mid Y]
 \qquad\text{for every }Y\in\cY(f).
\]
If $T$ is unbiased, so is $\widetilde T$, and its variance
is no larger.
\end{resultrestatement}

\begin{proof}[Proof of Lemma~\ref{lem:rao-blackwell}]
Recall that $U$ is the auxiliary randomization and
$D=(Z,Y_{O_f(Z)})$ consists of the treatment assignment and the observed outcomes.
An estimator $T(U,D)$ uses both the auxiliary randomization and the observed
data.  We apply Rao--Blackwellization to $T(U,D)$ by taking its conditional
expectation given $D$, obtaining $\E[T(U,D)\mid D]$.  We show that
this conditional expectation can be computed from the observed data
without knowing the unobserved potential outcomes.

Fix an admissible table $y\in\cY(f)$ and write $D=(Z,y_{O_f(Z)})$.
All probabilities and expectations below are over the known joint
randomization law of $(U,Z)$, with the potential outcome table $y$ held fixed.
For this fixed table, $D$ is a deterministic function of $Z$ and contains
$Z$.  Hence $\sigma(D)=\sigma(Z)$, where $\sigma(X)$ denotes the
sigma-algebra generated by a random vector $X$ (since $D$ is a function
of $Z$, $\sigma(D)\subseteq\sigma(Z)$, and since $Z$ is the first
component of $D$, $\sigma(Z)\subseteq\sigma(D)$).
Thus conditioning on $D$ gives the same information as conditioning
on $Z$, and the
conditional law of $U$ given $D$ is its known conditional law given $Z$.
Hence
\[
 \widetilde T(D)=\E[T(U,D)\mid D]
 =\sum_{u\in\mathcal U}T(u,D)\Pr(U=u\mid Z).
\]
The formula for $\widetilde T$ uses only the observed data and the known
randomization probabilities, so it defines the same observed-data rule for every
admissible table.  Assignments of probability zero may be assigned any
value of $\widetilde T$. Conditional Jensen's inequality and iterated expectation now give
\begin{align*}
 \E[(\widetilde T(D)-\tau_a(y))^2]
 &=\E\left[\bigl(\E[T(U,D)-\tau_a(y)\mid D]\bigr)^2\right]\\
 &\leq\E[(T(U,D)-\tau_a(y))^2],
\end{align*}
which proves equation~\eqref{eq:rb-domination}. Moreover, if $T$ is design-unbiased, then iterated expectation gives
\[
 \E[\widetilde T(D)]
 =\E\{\E[T(U,D)\mid D]\}
 =\E[T(U,D)]
 =\tau_a(y).
\]
Thus design-unbiasedness is preserved.  For unbiased estimators,
mean squared error equals design variance, giving the variance inequality.
\end{proof}

\section{Supporting results for Section~\ref{sec:lower-bounds}}
\label{app:lower-bound-proofs}

In this section, we provide proofs for the lower-bound results in
Section~\ref{sec:lower-bounds} and give additional lower-bound analyses.
We first establish preliminary lemmas and prove the two-table packing lower bound
in Lemma~\ref{lem:two-table-general}, the lower bounds from packing disjoint
perturbations in Theorem~\ref{thm:fixed-joint-missingness}, and the singleton packing
lower bound in Corollary~\ref{cor:singleton-lower}.
We next prove Proposition~\ref{thm:three-problem-characterization}, which
represents the fixed-design, restricted-design, and all-design minimax
risks in terms of the weighted sampling hypergraphs. The example in
Section~\ref{subsec:cluster-with-noise} then shows that two designs with the
same sampling hypergraph and first-order inclusion probabilities can have
fixed-design unbiased minimax risks differing by a factor of order $n$.
Corollary~\ref{thm:restricted-joint-missingness} provides lower
bounds over a design class. For all design minimax problem, we derive the hypergraph based lower
bounds in Theorem~\ref{thm:global-hypergraph-lower} and use them to prove
the sampling-rank lower bound in Theorem~\ref{cor:rank-lower} and the
paired-set lower bound in Theorem~\ref{thm:paired-sets-main}.
The final subsection compares the sampling hypergraph with the
conflict graph of \citet{KandirosHarshawSavje2026}.

\subsection{Preliminary lemmas}
\label{subsec:lower-bound-preliminaries}
In this section, we collect some standard technical lemmas used throughout the proofs. 
Lemma~\ref{lem:posterior-risk} shows that the fixed-design minimax risk is
bounded below by the expected posterior variance of the estimand under any
finite-support prior $\nu$ on the potential outcomes. We then establish a second-moment
bound in Lemma~\ref{lem:orthogonal-second-moment}, conditional moments under a
hypercube prior in Lemma~\ref{lem:hypercube-moments}, and a weighted probability
bound in Lemma~\ref{lem:weighted-reciprocal}.

\begin{lemma}[Conditional bias--variance decomposition]
\label{lem:posterior-risk}
Fix a design $P$ and a finite-support prior $\nu$ on $\cY(f)$.
Draw $\mathsf Y\sim\nu$ and $Z\sim P$ independently, and let
$D=(Z,\mathsf Y_{O_f(Z)})$.
Write $\mathbb E_{\nu,P}$ for expectation under this joint law and
$m(D)=\mathbb E_{\nu,P}[\tau_a(\mathsf Y)\mid D]$.
For any estimator $T_P(D)$ with finite second moment under this law,
\begin{align*}
 \mathbb E_{\nu,P}[(T_P(D)-\tau_a(\mathsf Y))^2]
 &=\mathbb E_{\nu,P}[(T_P(D)-m(D))^2]\\
 &\quad+\mathbb E_{\nu,P}
 [\operatorname{Var}_{\nu,P}(\tau_a(\mathsf Y)\mid D)].
\end{align*}
Consequently, the fixed-design minimax risk satisfies
\begin{equation}
 \Rall(P;f,a)
 =\inf_{T_P\in\cT(P,f)}\sup_{y\in\cY(f)}
 \E_P[(T_P-\tau_a(y))^2\mid y]
 \geq\mathbb E_{\nu,P}
 [\operatorname{Var}_{\nu,P}(\tau_a(\mathsf Y)\mid D)].
 \label{eq:prior-posterior-risk}
\end{equation}
\end{lemma}

\begin{proof}
Conditional on $D$, the bias--variance decomposition gives
\[
 \mathbb E_{\nu,P}[(T_P(D)-\tau_a(\mathsf Y))^2\mid D]
 =(T_P(D)-m(D))^2
 +\operatorname{Var}_{\nu,P}(\tau_a(\mathsf Y)\mid D).
\]
Taking expectations proves the identity. The inequality follows because
worst-case risk is at least prior-averaged risk and the squared term is
nonnegative.
Taking the infimum over estimators gives the minimax lower bound, since
the right-hand side does not depend on the estimator.
\end{proof}

\begin{lemma}[A second-moment bound]
\label{lem:orthogonal-second-moment}
Let $T,V_1,\ldots,V_r$ be square-integrable random variables on the same
probability space, with $\E(V_bV_d)=0$ for $b\neq d$ and
$v_b=\E V_b^2>0$.  Then
\[
 \E T^2\geq\sum_{b=1}^r\frac{\{\E(TV_b)\}^2}{v_b}.
\]
\end{lemma}

\begin{proof}
Write $\alpha_b=\E(TV_b)/v_b$.
Expanding a nonnegative square and using the zero cross moments gives
\begin{align*}
 0&\leq\E\left(T-\sum_{b=1}^r\alpha_bV_b\right)^2\\
 &=\E T^2-2\sum_{b=1}^r\alpha_b\E(TV_b)
   +\sum_{b=1}^r\alpha_b^2v_b\\
 &=\E T^2-\sum_{b=1}^r\frac{\{\E(TV_b)\}^2}{v_b}.
\end{align*}
Rearranging proves the inequality.
\end{proof}

\begin{lemma}[Conditional moments under a hypercube prior]
\label{lem:hypercube-moments}
Fix a design $P$, pairwise disjoint nonempty perturbation sets
$B_1,\ldots,B_m\subseteq K_a$, and fixed signs
$s_{b,k}\in\{-1,1\}$ for $k\in B_b$.
Let $\theta=(\theta_1,\ldots,\theta_m)\in\{-1,1\}^m$, with independent
components, each taking values $-1$ and $1$ with probability $1/2$,
independently of $Z\sim P$.  Define the random potential outcome table
\[
 Y_k^\theta=
 \begin{cases}
 s_{b,k}\theta_b,&k\in B_b,\\
 0,&k\notin\bigcup_{b=1}^m B_b.
 \end{cases}
\]
This construction gives a hypercube prior on $\cY(f)$: each of the $2^m$
sign vectors has probability $2^{-m}$.
Let $D=(Z,Y^\theta_{O_f(Z)})$ and
$J_b=\1\{O_f(Z)\cap B_b\neq\varnothing\}$.
Under the joint law of this prior and design,
\begin{equation}
 \begin{gathered}
 \E(\theta_b\mid D)=
 \begin{cases}
 \theta_b,&J_b=1,\\
 0,&J_b=0,
 \end{cases}
 \qquad \Var(\theta_b\mid D)=1-J_b,\\
 \Cov(\theta_b,\theta_d\mid D)=0\quad(b\neq d).
 \end{gathered}
 \label{eq:packing-conditional-moments}
\end{equation}
In the first case, observing the potential outcome at an index $k\in B_b$
gives $\theta_b=s_{b,k}Y_k^\theta$, so $\theta_b$ is determined by $D$.
\end{lemma}

\begin{proof}
Since the perturbation sets $B_b$ are pairwise disjoint and all entries
lie in $[-1,1]$, each table $Y^\theta$ belongs to $\cY(f)$.
Fix data $D=(z,y_{O_f(z)})$ having positive probability.
If $J_b=1$, an observed potential outcome in $B_b$ determines
$\theta_b$, so $\E(\theta_b\mid D)=\theta_b$.
If $J_b=0$, changing only $\theta_b$ to $-\theta_b$ leaves
$(z,y_{O_f(z)})$ unchanged. Flipping $\theta_b$ pairs each sign vector
consistent with $(z,y_{O_f(z)})$ with
another consistent vector of the same prior probability $2^{-m}$ and
opposite value of $\theta_b$.  The assignment probability $P(z)$ is also
unchanged.  The paired terms cancel in the conditional expectation, giving
$\E(\theta_b\mid D)=0$.  Thus
\[
 \E(\theta_b\mid D)=
 \begin{cases}
 \theta_b,&J_b=1,\\
 0,&J_b=0.
 \end{cases}
\]
Since $\theta_b^2=1$, the two cases give
\[
 \Var(\theta_b\mid D)
 =\E(\theta_b^2\mid D)-\{\E(\theta_b\mid D)\}^2
 =1-J_b.
\]
For $b\neq d$, if either $J_b=0$ or $J_d=0$, flipping an unobserved
sign pairs equal-probability terms with opposite values of
$\theta_b\theta_d$.  Hence $\E(\theta_b\theta_d\mid D)=0$.
If both indicators equal one, both signs are determined by the data,
so $\E(\theta_b\theta_d\mid D)=\theta_b\theta_d$.
Combining the two cases gives
\[
 \Cov(\theta_b,\theta_d\mid D)
 =\E(\theta_b\theta_d\mid D)
  -\E(\theta_b\mid D)\E(\theta_d\mid D)
 =J_bJ_d\theta_b\theta_d-J_bJ_d\theta_b\theta_d=0.
\]
\end{proof}

\begin{remark}[Role of the random signs]
The perturbation sets $B_b$ specify which potential outcomes
vary together.  The independent uniform signs $\theta_b$ in
Lemma~\ref{lem:hypercube-moments}, also called Rademacher signs, turn this
packing into a finite prior on $\cY(f)$.  All potential outcomes in $B_b$ share the same random sign, so observing any one of them reveals
$\theta_b$.  If none is observed, the conditional variance of $\theta_b$
remains one. The lower-bound proofs use the conditional moments in
equation~\eqref{eq:packing-conditional-moments} to compute the
uncertainty in the estimand left by an assignment.
\end{remark}

\begin{lemma}[A weighted reciprocal-probability bound]
\label{lem:weighted-reciprocal}
Let $p_1,\ldots,p_r>0$ and $w_1,\ldots,w_r\in\mathbb R$.  Then
\begin{equation}
 \sum_{k=1}^r\frac{w_k^2}{p_k}
 \geq\frac{\left(\sum_{k=1}^r|w_k|\right)^2}
                {\sum_{k=1}^r p_k}.
 \label{eq:weighted-reciprocal}
\end{equation}
In particular, if $\sum_k p_k\leq C$ for some $C>0$, the left-hand side
is at least $(\sum_k|w_k|)^2/C$.
\end{lemma}

\begin{proof}
Cauchy--Schwarz gives
\[
 \left(\sum_{k=1}^r|w_k|\right)^2
 =\left(\sum_{k=1}^r\frac{|w_k|}{\sqrt{p_k}}\sqrt{p_k}\right)^2
 \leq\left(\sum_{k=1}^r\frac{w_k^2}{p_k}\right)
       \left(\sum_{k=1}^r p_k\right).
\]
Divide by $\sum_k p_k>0$.  The final assertion follows by replacing this
sum in the denominator by its upper bound $C$.
\end{proof}

\subsection{Fixed-design lower bounds and refinements}
\label{app:fixed-design-refinements}
\label{subsec:detection}

We first prove the two-table lower bound for unbiased and all estimator classes stated in
Lemma~\ref{lem:two-table-general} using the two-point prior, followed by the lower bounds from disjoint perturbations
in Theorem~\ref{thm:fixed-joint-missingness} which uses the hypercube prior.  We then derive its
single-set consequence in Lemma~\ref{lem:detection-principle} and prove
the singleton lower bound in Corollary~\ref{cor:singleton-lower}.
We finish with the mutually exclusive-set refinement in
Proposition~\ref{prop:crossblocklb}. Throughout this subsection the design is fixed.
Section~\ref{subsec:packing-design-class-bounds} treats optimization over a design class.

\subsubsection{The two-table lower bound}
We restate the lemma.
\begin{resultrestatement}{Lemma~\ref{lem:two-table-general}}{Two-table lower bound}
Let $Y,Y'\in\cY(f)$ differ exactly on a nonempty set $B\subseteq K_a$, and
let $\Delta=|\tau_a(Y)-\tau_a(Y')|>0$.  For every design $P$,
\[
 \Rall(P;f,a)\geq\frac{\Delta^2}{4}\rho_P(B), \mathrm{ and} \qquad
 \Runb(P;f,a)\geq\frac{\Delta^2}{4} \frac{\rho_P(B)}{1-\rho_P(B)}.
\]

The second bound is interpreted as $+\infty$ when $\rho_P(B)=1$;
in that case no unbiased estimator exists.
\end{resultrestatement}

\begin{proof}[Proof of Lemma~\ref{lem:two-table-general}]
We first prove the lower bound over all estimators.  This argument has
two parts: in part (i), we bound the worst-case design mean squared error
of any estimator from below by the expected posterior variance of the
estimand under a two-point prior; in part (ii), we compute this posterior
variance.

\emph{Part (i): Bounding risk by the expected posterior variance.}
Write $Y^{(0)}=Y$, $Y^{(1)}=Y'$, $t_j=\tau_a(Y^{(j)})$ for
$j\in\{0,1\}$, and $\bar t=(t_0+t_1)/2$.  Let $\nu$ be a prior
on $\cY(f)$ that assigns probability $1/2$ to each of $Y^{(0)}$ and
$Y^{(1)}$, and zero probability to all other tables.  For a random table
$\mathsf Y\sim\nu$, this means
\[
 \Pr_\nu\{\mathsf Y=Y^{(0)}\}
 =\Pr_\nu\{\mathsf Y=Y^{(1)}\}=\frac12.
\]
Draw $Z\sim P$ independently of $\mathsf Y$.  The prior $\nu$ and the
design $P$ therefore induce the product joint law $\nu\otimes P$ on
$(\mathsf Y,Z)$.  Write $\mathbb E_{\nu,P}$ for expectation under this
law.  The observed data are $D=(Z,\mathsf Y_{O_f(Z)})$.
Conditional on $\mathsf Y=y$, only the assignment is random, with law $P$.

Let
\[
 M=\{O_f(Z)\cap B=\varnothing\}
\]
be the event that all potential outcomes in a perturbation set $B$ are missing.  Write
$\rho=P(M)=\rho_P(B)$ and $q=P(M^c)=1-\rho$.

Let $T_P\in\cT(P,f)$ be an estimator with finite worst-case risk;
otherwise the lower bound is immediate.
Applying Lemma~\ref{lem:posterior-risk} to this two-point prior gives
\begin{align*}
 \sup_{y\in\cY(f)}\E_P[(T_P-\tau_a(y))^2\mid y]
 &\geq\frac12\sum_{j=0}^1\E_P[(T_P-t_j)^2\mid Y^{(j)}]\\
 &=\mathbb E_{\nu,P}[(T_P(D)-\tau_a(\mathsf Y))^2]\\
 &\geq\mathbb E_{\nu,P}
       [\operatorname{Var}_{\nu,P}(\tau_a(\mathsf Y)\mid D)].
\end{align*}

\emph{Part (ii): Computing the posterior variance.}
We now compute the posterior variance of $\tau_a(\mathsf Y)$ under the
two-point prior.  On $M$, no potential outcome in $B$ is
observed.  Since the two tables agree outside $B$, they give identical
observed data at each assignment $z$ for which $M$ occurs.  The prior
assigns equal probability to the tables, and the assignment probability
$P(z)$ is the same under both.  Hence, on $M$,
\[
 \Pr_{\nu,P}\{\mathsf Y=Y^{(j)}\mid D\}=\frac12,
 \qquad j=0,1.
\]
On $M^c$, an observed potential outcome in $B$ identifies the
table, since the tables differ in the potential outcome at every index in
$B$.  The posterior variance is therefore zero on $M^c$, whereas on $M$
the two estimand values retain equal posterior probabilities.
Consequently, almost surely,
\begin{align*}
 \operatorname{Var}_{\nu,P}(\tau_a(\mathsf Y)\mid D)
 &=\left\{\frac12(t_0-\bar t)^2+
           \frac12(t_1-\bar t)^2\right\}\mathbf 1_M\\
 &=\frac{\Delta^2}{4}\mathbf 1_M.
\end{align*}
Part (i) shows that the worst-case design risk is at least the expected
posterior variance.  Substituting the posterior variance computed in
part (ii) gives
\begin{align*}
 \sup_{y\in\cY(f)}\E_P[(T_P-\tau_a(y))^2\mid y]
 &\geq\mathbb E_{\nu,P}\!\left[
       \operatorname{Var}_{\nu,P}(\tau_a(\mathsf Y)\mid D)\right]\\
 &=\frac{\Delta^2}{4}\Pr_{\nu,P}(M)
  =\frac{\Delta^2}{4}\rho_P(B).
\end{align*}
Taking the infimum over all estimators proves the first bound.

We next prove the lower bound over unbiased estimators.
Suppose that $T_P$ is unbiased at every potential outcome table in $\cY(f)$.
First suppose $0<\rho_P(B)<1$.  Let $\mathbb E_j$ and
$\operatorname{Var}_j$ denote expectation and variance over $Z\sim P$
with the table fixed at $Y^{(j)}$.  By construction, the conditional
distribution of the observed data given $M$ is the same under both
tables: every potential outcome in $B$ is missing, and the
tables agree outside $B$.  Therefore define
\[
 u=\mathbb E_0[T_P\mid M]=\mathbb E_1[T_P\mid M],
 \qquad v_j=\mathbb E_j[T_P\mid M^c].
\]
Unbiasedness and conditioning on $M$ and $M^c$ give, for $j=0,1$,
\begin{align*}
 t_j=\mathbb E_j[T_P]
 &=\rho_P(B)\mathbb E_j[T_P\mid M]
   +(1-\rho_P(B))\mathbb E_j[T_P\mid M^c]\\
 &=\rho_P(B)u+(1-\rho_P(B))v_j.
\end{align*}
Subtracting $u$ from both sides and dividing by $1-\rho_P(B)>0$ yields
\[
 t_j-u=(1-\rho_P(B))(v_j-u),\qquad
 v_j-u=\frac{t_j-u}{1-\rho_P(B)}.
\]

For an unbiased estimator, the design mean squared error equals the
design variance.  The law of total variance, conditional on the indicator
$\mathbf 1_M$, gives
\begin{align*}
 \operatorname{Var}_j(T_P)
 &=\mathbb E_j\!\left[\operatorname{Var}_j(T_P\mid\mathbf 1_M)\right]
   +\operatorname{Var}_j\!\left(\mathbb E_j[T_P\mid\mathbf 1_M]\right)\\
 &=\rho_P(B)\operatorname{Var}_j(T_P\mid M)
   +(1-\rho_P(B))\operatorname{Var}_j(T_P\mid M^c)\\
 &\quad+\rho_P(B)(u-t_j)^2+(1-\rho_P(B))(v_j-t_j)^2\\
 &\geq\rho_P(B)(u-t_j)^2+(1-\rho_P(B))(v_j-t_j)^2.
\end{align*}
The inequality drops the two nonnegative conditional variances.
The unbiasedness identity also gives
$u-t_j=(1-\rho_P(B))(u-v_j)$ and
$v_j-t_j=\rho_P(B)(v_j-u)$.  Substituting $u-t_j=(1-\rho_P(B))(u-v_j)$ and
$v_j-t_j=\rho_P(B)(v_j-u)$, then using
$v_j-u=(t_j-u)/(1-\rho_P(B))$, gives
\begin{align*}
 \E_P[(T_P-t_j)^2\mid Y^{(j)}]
 &=\operatorname{Var}_j(T_P)\\
 &\geq\rho_P(B)(1-\rho_P(B))^2(u-v_j)^2
       +(1-\rho_P(B))\rho_P(B)^2(v_j-u)^2\\
 &=\rho_P(B)(1-\rho_P(B))(u-v_j)^2\\
 &=\frac{\rho_P(B)}{1-\rho_P(B)}(t_j-u)^2,
 \qquad j=0,1.
\end{align*}
We have obtained one risk inequality at $Y^{(0)}$ and one at $Y^{(1)}$.
The worst-case risk is at least the average of these two risks, so
averaging the inequalities for $j=0$ and $j=1$ yields
\begin{align*}
 \sup_{y\in\cY(f)}\E_P[(T_P-\tau_a(y))^2\mid y]
 &\geq\frac12\sum_{j=0}^1 \E_P[(T_P-t_j)^2\mid Y^{(j)}]\\
 &\geq\frac{\rho_P(B)}{1-\rho_P(B)}\,
        \frac{(t_0-u)^2+(t_1-u)^2}{2}\\
 &=\frac{\rho_P(B)}{1-\rho_P(B)}
        \left\{(u-\bar t)^2+\frac{\Delta^2}{4}\right\}\\
 &\geq\frac{\Delta^2}{4}\frac{\rho_P(B)}{1-\rho_P(B)}.
\end{align*}
Taking the infimum over unbiased estimators proves the second
bound for $0<\rho_P(B)<1$.  When $\rho_P(B)=0$, the bound is zero and holds by
nonnegativity of squared error.  When $\rho_P(B)=1$, the observed-data laws
under the two tables coincide, so every integrable estimator has the
same expectation under both.  Since $t_0\neq t_1$, no estimator can be
unbiased at both tables, and $\Runb(P;f,a)=+\infty$.
\end{proof}

\begin{remark}[Detection probability and total variation]
Recall that $\rho_P(B)$ is the probability that every potential outcome in $B$ is missing, whereas $q_P(B)=1-\rho_P(B)$ is the
probability of observing at least one such potential outcome. The latter event
distinguishes the two tables in Lemma~\ref{lem:two-table-general}.

For the tables $Y,Y'$ in Lemma~\ref{lem:two-table-general}, let
$\mathcal L_Y$ and $\mathcal L_{Y'}$ denote the laws of
$(Z,Y_{O_f(Z)})$ and $(Z,Y'_{O_f(Z)})$, respectively, under the same design
$P$. The detection probability equals the total-variation distance:
\begin{equation}
 \|\mathcal L_Y-\mathcal L_{Y'}\|_{\mathrm{TV}}
 =1-\rho_P(B)=q_P(B).
 \label{eq:tv-identity}
\end{equation}
Indeed, the laws agree on assignments observing no potential outcome in $B$.
On the complementary event, their observed outcomes differ.  Because the
assignment itself is observed and has the same law under both tables,
the total-variation distance equals the probability of observing at least
one potential outcome in $B$.
Thus the all-estimator bound in Lemma~\ref{lem:two-table-general} can be
written as
\[
 \Rall(P;f,a)\geq\frac{\Delta^2}{4}
 \bigl(1-\|\mathcal L_Y-\mathcal L_{Y'}\|_{\mathrm{TV}}\bigr).
\]
This is the same bound expressed through the overlap of the two data
laws, as in classical two-point lower-bound arguments. In this experiment,
the assignment is observed and outcomes are deterministic given the
assignment and table, so the overlap is exactly the joint missingness
probability. The total-variation formulation therefore does not give a
stronger bound here.
\end{remark}

\subsubsection{Lower bounds from disjoint perturbation sets}
We restate and prove the lower bound from disjoint perturbations. Recall that a perturbation packing $\cB=\{B_1,\ldots,B_m\}$ consists
of pairwise disjoint nonempty subsets of $K_a$.  For fixed signs
$s_{b,k}\in\{-1,1\}$, its estimand relevant coefficients are
$c_b=\sum_{k\in B_b}a_ks_{b,k}$.

\begin{resultrestatement}{Theorem~\ref{thm:fixed-joint-missingness}}{Lower bounds from disjoint perturbations}
Fix a design $P$, a perturbation packing $\cB=\{B_1,\ldots,B_m\}$, and signs
defining $c_b$.
\begin{enumerate}
\item[\textup{(i)}] Over all estimators,
\[
  \Rall(P;f,a)\geq\sum_{b=1}^m c_b^2\rho_P(B_b).
\]
\item[\textup{(ii)}] Over unbiased estimators,
\[
  \Runb(P;f,a)
  \geq\sum_{b:c_b\neq0}c_b^2
       \frac{\rho_P(B_b)}{1-\rho_P(B_b)}.
\]
If $c_b\neq0$ and $\rho_P(B_b)=1$ for some set, no unbiased
estimator exists and the bound is $+\infty$.
\end{enumerate}
\end{resultrestatement}

\begin{proof}[Proof of Theorem~\ref{thm:fixed-joint-missingness}]
Fix the design $P$ and the potential-outcome packing
$\cB=\{B_1,\ldots,B_m\}$, with the signs $s_{b,k}$ specified in the theorem.
We first construct a hypercube prior on $\cY(f)$.
For part (i) of the theorem, to prove the all-estimator lower bound, we
compute the expected posterior variance of the estimand under this prior.
For part (ii), to prove the unbiased-estimator lower bound, we establish
a second-moment bound for an unbiased estimator, then
subtract the prior second moment of the estimand to obtain a lower bound on the design variance.

\emph{Construction of the hypercube prior.}
Use the prior of Lemma~\ref{lem:hypercube-moments} with the perturbation
sets $B_1,\ldots,B_m$ and fixed signs $s_{b,k}$ from the theorem.
Thus $\theta=(\theta_1,\ldots,\theta_m)$ has one independent uniform
sign for each perturbation set $B_b$, and $Y_k^\theta=s_{b,k}\theta_b$
for $k\in B_b$, with zero entries outside $\bigcup_{b=1}^m B_b$.
The random sign $\theta_b$ determines whether the fixed sign pattern
on $B_b$ is used or reversed.
Substituting the table into the estimand gives
\begin{equation}
 \tau_a(Y^\theta)
 =\sum_{b=1}^m\sum_{k\in B_b}a_ks_{b,k}\theta_b
 =\sum_{b=1}^m c_b\theta_b.
 \label{eq:packing-prior-estimand}
\end{equation}
All expectations below are
over this finite prior and the design unless a subscript specifies otherwise.

\emph{Part (i): The lower bound over all estimators.}
Let
\[
 D=(Z,Y^\theta_{O_f(Z)}),\qquad
 J_b=\1\{O_f(Z)\cap B_b\neq\varnothing\},\qquad q_b=\E J_b=q_P(B_b).
\]
Here $J_b$ indicates whether the assignment reveals any potential outcome in $B_b$, and $q_b$ is its detection probability.
Lemma~\ref{lem:hypercube-moments} gives the conditional variances
$\Var(\theta_b\mid D)=1-J_b$ and zero conditional covariances for
distinct components of $\theta$.
Using equation~\eqref{eq:packing-prior-estimand} and expanding the conditional
variance of the estimand gives
\begin{align*}
 \Var\{\tau_a(Y^\theta)\mid D\}
 &=\sum_b c_b^2\Var(\theta_b\mid D)
   +2\sum_{b<d}c_bc_d\Cov(\theta_b,\theta_d\mid D)\\
 &=\sum_b c_b^2(1-J_b).
\end{align*}
Lemma~\ref{lem:posterior-risk} applies to this finite prior:
for any estimator $T_P$, its prior-average risk is at least the expected
posterior variance of the estimand.  Hence
\[
 \sup_{Y\in\cY(f)}\E_P[(T_P-\tau_a(Y))^2\mid Y]
 \geq\E\{\Var(\tau_a(Y^\theta)\mid D)\}
 =\sum_b c_b^2(1-q_b)
 =\sum_b c_b^2\rho_P(B_b).
\]
Taking the infimum over estimators proves part (i).

\emph{Part (ii): The lower bound for unbiased estimators.}
We now use unbiasedness to obtain the stronger bound in
part (ii).  Let $I=\{b\in\{1,\ldots,m\}:c_b\neq0\}$.
When $c_b=0$, the sign $\theta_b$ does not affect the estimand in
equation~\eqref{eq:packing-prior-estimand}, and the corresponding perturbation
set contributes no term to the claimed bound.  We therefore restrict the
sums in this part to $b\in I$, retaining the same prior on the full table.
If $I$ is empty, the bound is zero and follows from nonnegativity of
variance.  Otherwise, first suppose that $q_b=q_P(B_b)>0$ for every
$b\in I$.  We treat the case $q_b=0$ at the end of this part.

Let $T_P=T_P(D)$ be an unbiased estimator with finite
worst-case second moment.  We first show that unbiasedness implies
\begin{equation}
 \E T_P^2\geq\sum_{b\in I}\frac{c_b^2}{q_b}.
 \label{eq:packing-unbiased-second-moment}
\end{equation}
We will then subtract the second moment of the estimand to obtain the
average design variance under the prior.

For each $b\in I$, define $U_b=J_b\theta_b$.
When $J_b=1$, the data reveal $\theta_b$; when $J_b=0$, $U_b=0$.
Thus $U_b$ is a function of the observed data.
Since the assignment is independent of $\theta$, and the components of
$\theta$ are independent with mean zero and square one,
\begin{align*}
 \E U_b&=\E J_b\,\E\theta_b=0,\\
 \E U_b^2&=\E(J_b\theta_b^2)=\E J_b=q_b,\\
 \E(U_bU_d)&=\E(J_bJ_d)\E(\theta_b\theta_d)=0
 \qquad(b\neq d).
\end{align*}
The last identity does not require independence of $J_b$ and $J_d$.

We next compute $\E(T_PU_b)$.  Unbiasedness at each table
$Y^\theta$ gives $\E_P(T_P\mid\theta)=\sum_{d\in I}c_d\theta_d$.
Conditioning on $\theta$ and then averaging over the prior gives
\begin{align*}
 \E(T_P\theta_b)
 &=\E_\theta\{\theta_b\E_P(T_P\mid\theta)\}\\
 &=\sum_{d\in I}c_d\E_\theta(\theta_b\theta_d)
 =c_b.
\end{align*}
By the conditional-mean identity in
Lemma~\ref{lem:hypercube-moments},
$\E(\theta_b\mid D)=U_b$.  Since $T_P$ is a function of $D$,
\[
 \E(T_PU_b)
 =\E\{T_P\E(\theta_b\mid D)\}
 =\E\{\E(T_P\theta_b\mid D)\}
 =\E(T_P\theta_b)=c_b.
\]
The variables $U_b$, $b\in I$, have zero pairwise cross moments,
second moments $q_b>0$, and $\E(T_PU_b)=c_b$.
Applying Lemma~\ref{lem:orthogonal-second-moment} with $V_b=U_b$ gives
\[
 \E T_P^2\geq\sum_{b\in I}
 \frac{\{\E(T_PU_b)\}^2}{\E U_b^2}
 =\sum_{b\in I}\frac{c_b^2}{q_b},
\]
which proves equation~\eqref{eq:packing-unbiased-second-moment}.

We now convert this second-moment bound into a design variance bound.
The independent mean-zero signs give
\begin{align*}
 \E_\theta\{\tau_a(Y^\theta)^2\}
 &=\sum_{b\in I}c_b^2\E_\theta\theta_b^2
   +2\sum_{\substack{b,d\in I\\b<d}}
       c_bc_d\E_\theta(\theta_b\theta_d)\\
 &=\sum_{b\in I}c_b^2.
\end{align*}
For an unbiased estimator, design MSE equals design variance.
Using unbiasedness at each table and then
equation~\eqref{eq:packing-unbiased-second-moment}, we obtain
\begin{align*}
 \E_\theta\Var_P(T_P\mid Y^\theta)
 &=\E_\theta\left[\E_P(T_P^2\mid Y^\theta)
             -\{\E_P(T_P\mid Y^\theta)\}^2\right]\\
 &=\E T_P^2-\E_\theta\{\tau_a(Y^\theta)^2\}\\
 &\geq\sum_{b\in I}\frac{c_b^2}{q_b}-\sum_{b\in I}c_b^2
 =\sum_{b\in I}c_b^2\left(\frac1{q_b}-1\right).
\end{align*}
Every table in the prior belongs to $\cY(f)$, so the worst-case design
variance is at least this prior-average design variance:
\[
 \sup_{Y\in\cY(f)}\Var_P(T_P\mid Y)
 \geq\E_\theta\Var_P(T_P\mid Y^\theta)
 \geq\sum_{b\in I}c_b^2\left\{\frac1{q_P(B_b)}-1\right\}.
\]
This proves the claimed inequality for every unbiased
estimator with finite worst-case second moment.  If its worst-case second
moment is infinite, its worst-case variance is also infinite, since
unbiasedness fixes its mean at the bounded estimand.  The inequality then
holds as well.

Finally, suppose $q_P(B_b)=0$ for some $b\in I$.
Fix all the other components of $\theta$, and let $Y^+$ and $Y^-$ be
the tables obtained by setting $\theta_b=1$ and $\theta_b=-1$,
respectively. These tables agree outside $B_b$, and no potential outcome
in $B_b$ is observed with positive design probability.
Their observed-data laws therefore coincide.  By
equation~\eqref{eq:packing-prior-estimand}, however,
\[
 \tau_a(Y^+)-\tau_a(Y^-)=2c_b\neq0.
\]
An integrable estimator has the same expectation under the two identical
observed-data laws, so it cannot be unbiased at both tables.  Thus no
unbiased estimator exists in this case, as claimed.

\end{proof}

\subsubsection{Lower bounds from a single perturbation set and singleton packings}

The next lemma records the single-perturbation-set case of
the lower bound from disjoint perturbations (Theorem~\ref{thm:fixed-joint-missingness}).

\begin{lemma}[Lower bounds from joint missingness]
\label{lem:detection-principle}
\label{lem:detection}
For every fixed design $P$, every nonempty $B\subseteq K_a$, and every
choice of signs defining $c_B$,
\[
  \Rall(P;f,a)\geq c_B^2\rho_P(B).
\]
If $\rho_P(B)<1$, every unbiased estimator satisfies
\[
  \sup_{Y\in\cY(f)}\Var_P(T_P\mid Y)
  \geq c_B^2\frac{\rho_P(B)}{1-\rho_P(B)}.
\]
If $c_B\neq0$ and $\rho_P(B)=1$, no unbiased estimator exists.
\end{lemma}

\begin{proof}[Proof of Lemma~\ref{lem:detection-principle}]
Take the packing $\cB=\{B\}$, consisting of one perturbation set,
with the signs defining $c_B$.  Its coefficient is $c_B$, and its detection
probability is $q_P(B)=1-\rho_P(B)$.
Part (i) of the lower bound from disjoint perturbations (Theorem~\ref{thm:fixed-joint-missingness}) gives
\[
 \Rall(P;f,a)\geq c_B^2\{1-q_P(B)\}=c_B^2\rho_P(B).
\]
If $\rho_P(B)<1$, then $q_P(B)>0$.  Part (ii) of the same theorem gives,
for every unbiased estimator,
\[
 \sup_{Y\in\cY(f)}\Var_P(T_P\mid Y)
 \geq c_B^2\left\{\frac1{q_P(B)}-1\right\}
 =c_B^2\frac{\rho_P(B)}{1-\rho_P(B)}.
\]
When $c_B=0$, this inequality holds by nonnegativity of variance.
Finally, if $c_B\neq0$ and $\rho_P(B)=1$, then $q_P(B)=0$, so the
impossibility statement in part (ii) of the theorem applies.
\end{proof}

\begin{proof}[Proof of Corollary~\ref{cor:singleton-lower}]
Apply the lower bound from disjoint perturbations (Theorem~\ref{thm:fixed-joint-missingness}) to the singleton packing
\(B_k=\{k\}\), \(k\in K_a\), with \(s_k=\operatorname{sign}(a_k)\).  Then
\(\rho_P(B_k)=1-\pi_{P,k}\) and \(c_k=|a_k|\), so part (i) gives
equation~\eqref{eq:singleton-all}, part (ii) gives equation~\eqref{eq:singleton-unb}, and part
(ii) also shows that no unbiased estimator exists when
\(\pi_{P,k}=0\) for some \(k\in K_a\).
\end{proof}

The missingness lower bounds are attained when one unknown scalar is
either observed exactly or missing, as the following proposition shows.

\begin{proposition}[Exact price of unbiasedness for one missing component]
\label{prop:one-block-price}
Suppose we observe \(\theta\in[-C,C]\) exactly with probability
\(q\in(0,1)\). Conditional on \(\theta\) being missing, the law of the data
does not depend on \(\theta\).  Then
\begin{align}
  \inf_T\sup_{|\theta|\leq C}\E_\theta(T-\theta)^2
  &=C^2(1-q),
  \label{eq:one-block-all}\\
  \inf_{T:\,\E_\theta T=\theta}
  \sup_{|\theta|\leq C}\Var_\theta(T)
  &=C^2\frac{1-q}{q}.
  \label{eq:one-block-unb}
\end{align}
\end{proposition}

\begin{proof}
Let $J$ indicate that the scalar is observed, so $\Pr(J=1)=q$.
For unrestricted estimation, use $T_0=J\theta$, interpreted as zero
when $J=0$.  Its risk is
\[
 \E_\theta[(T_0-\theta)^2]=(1-q)\theta^2
 \leq(1-q)C^2.
\]
The two parameter values $\theta=C$ and $\theta=-C$ have separation
$2C$ and give identical data when $J=0$.  The two-table lower bound
therefore gives $(2C)^2(1-q)/4=C^2(1-q)$, matching this upper bound.

For unbiasedness, the same two values give the lower bound
$C^2(1-q)/q$ from the unbiased part of the two-table lemma.
Use $T_1=J\theta/q$.  Since $\E_\theta T_1=\theta$ and
$\E_\theta T_1^2=\theta^2/q$, its variance is
\[
 \Var_\theta(T_1)=\theta^2/q-\theta^2
 =\theta^2(1-q)/q\leq C^2(1-q)/q.
\]
This attains the second lower bound.  If $q=1$, both estimators recover
$\theta$ exactly and both risks are zero.
\end{proof}

\subsubsection{Lower bounds from mutually exclusive perturbation sets}
\label{subsec:contrast-refinements}
The lower bound from disjoint perturbations in
Theorem~\ref{thm:fixed-joint-missingness}(ii) can be
strengthened when two perturbation sets cannot be observed together.
Proposition~\ref{prop:crossblocklb} gives the refinement for one pair of perturbations;
Equation~\eqref{eq:crossblocklb-multipair} gives the extension to disjoint pairs of perturbations.
Write
\begin{equation}
 \lambda_P(B)=\frac{\rho_P(B)}{1-\rho_P(B)}
             =\frac{1-q_P(B)}{q_P(B)},
 \label{eq:joint-missingness-odds}
\end{equation}
with value $+\infty$ when $\rho_P(B)=1$.

The refinement adds a positive cross term for opposite-signed coefficient
sums. It is used in the no-interference benchmarks in
Appendix~\ref{app:no-interference-proofs}. For the no-interference ATE
under independent balanced assignment, it improves the unbiased lower
bound from $2/n$ to the exact value $4/n$
(Theorem~\ref{thm:ate-joint-benchmark}).

\begin{proposition}[Mutually exclusive, opposite-signed sets]
\label{prop:crossblocklb}
Let $B_1,B_2\subseteq K_a$ be disjoint packed sets.  For $r=1,2$, let
\[
  c_r=\sum_{k\in B_r}s_{r,k}a_k,\qquad
  J_r(Z)=\1\{O_f(Z)\cap B_r\neq\varnothing\},\qquad
  q_{P,r}=P\{J_r(Z)=1\},
\]
and suppose $J_1(Z)J_2(Z)=0$ almost surely.  If $q_{P,r}>0$ whenever
$c_r\neq0$, every unbiased estimator satisfies
\begin{equation}
  \sup_{Y\in\cY(f)}\Var_P(T_P\mid Y)
  \geq
  c_1^2\lambda_P(B_1)+c_2^2\lambda_P(B_2)-2c_1c_2.
  \label{eq:crossblocklb}
\end{equation}
If $c_r\neq0$ but $q_{P,r}=0$, the unbiased class is empty.
Products $c_r^2\lambda_P(B_r)$ with $c_r=0$ are interpreted as zero,
including when $q_{P,r}=0$.
\end{proposition}

The last term strengthens the lower bound from disjoint perturbations (Theorem~\ref{thm:fixed-joint-missingness}(ii)) when $c_1c_2<0$, as it
does for opposite-signed components of a causal contrast.

More generally, suppose all \(2J\) sets in the pairs
\(\{(B_{j1},B_{j2}):j=1,\ldots,J\}\) are pairwise disjoint and the two
observation indicators within each pair are mutually exclusive.  Then
\begin{equation}
  \sup_{Y\in\cY(f)}\Var_P(T_P\mid Y)
  \geq
  \sum_{j=1}^J\left\{
    c_{j1}^2\lambda_P(B_{j1})
    +c_{j2}^2\lambda_P(B_{j2})
    -2c_{j1}c_{j2}
  \right\}.
  \label{eq:crossblocklb-multipair}
\end{equation}
Here $q_{P,jr}=P\{O_f(Z)\cap B_{jr}\ne\varnothing\}$ and
$c_{jr}$ is the signed coefficient sum on $B_{jr}$.  No independence among
observation events for different pairs is required.

\begin{proof}[Proof of Proposition~\ref{prop:crossblocklb} and
equation~\eqref{eq:crossblocklb-multipair}]
We prove the multiple-pair statement; a single pair is its case $J=1$.
On block $B_{jr}$ put $Y_k=s_{jr,k}x_{jr}$, with
$x_{jr}\in[-1,1]$, and set all other entries to zero.  The target on this
submodel is $\sum_{j,r}c_{jr}x_{jr}$.  Write $t_z(x)$ for the estimator
at assignment $z$, after averaging any estimator randomization.
If $J_{jr}(z)=0$, the observed data, and hence $t_z$, do not depend on
$x_{jr}$.

For each pair draw a fair sign $\theta_j$, independently over pairs and
independently of $Z$.  Fix $j,r$.  For $\varepsilon\in\{-1,1\}$, let
$x^{(jr,\varepsilon)}(\theta_{-j})$ denote the table parameters with
$x_{jr}=\varepsilon$, the other parameter in pair $j$ set to zero, and
$x_{h1}=x_{h2}=\theta_h$ for $h\neq j$.  Define
\[
 A_{jr}(z)=\frac12\E_{\theta_{-j}}\left[
 t_z\{x^{(jr,1)}(\theta_{-j})\}
 -t_z\{x^{(jr,-1)}(\theta_{-j})\}\right].
\]
If $J_{jr}(z)=0$, these two tables give the same data, so $A_{jr}(z)=0$.
If the estimator has infinite worst-case variance, the required bound
holds immediately.  Otherwise the following moments are finite.
Unbiasedness at each of the two tables gives
\begin{align*}
 \E_P A_{jr}(Z)
 &=\frac12\E_{\theta_{-j}}\left[
 \tau\{Y(x^{(jr,1)})\}-\tau\{Y(x^{(jr,-1)})\}\right]\\
 &=\frac12(2c_{jr})=c_{jr}.
\end{align*}
Using the support restriction on $A_{jr}$, Cauchy--Schwarz gives
\[
 c_{jr}^2=(\E_P[A_{jr}J_{jr}])^2
 \leq\E_P A_{jr}^2\,\E_P J_{jr}^2
 =q_{P,jr}\E_P A_{jr}^2.
\]
Thus, for $q_{P,jr}>0$,
\[
 E_P A_{jr}(Z)=c_{jr},\qquad
 E_P A_{jr}(Z)^2\geq c_{jr}^2/q_{P,jr}.
\]
If $q_{P,jr}=0$ and $c_{jr}\ne0$, the first equality is impossible.
If both are zero, $A_{jr}=0$ almost surely and its contribution is zero.

Now evaluate the estimator on the tables $x_{j1}=x_{j2}=\theta_j$,
and write $t_z(\theta)$ for its value on these tables at assignment $z$.
At a fixed assignment, at most one potential outcome of each pair is observed.
The coefficient of $\theta_j$ in the sign expansion of $t_z(\theta)$ is
therefore $A_{j1}(z)+A_{j2}(z)$: when one potential outcome is observed, changing
the other from zero to $\theta_j$ does not change the data.
More explicitly,
\[
 \E_\theta[t_z(\theta)\theta_j]=A_{j1}(z)+A_{j2}(z).
\]
For this fixed assignment $z$, the signs satisfy
$\E_\theta\theta_j^2=1$ and $\E_\theta(\theta_j\theta_h)=0$ for $j\neq h$.
Applying Lemma~\ref{lem:orthogonal-second-moment} under the sign prior,
with $T=t_z(\theta)$ and $V_j=\theta_j$, gives
\[
 \E_\theta t_z(\theta)^2
 \geq\sum_j\{\E_\theta[t_z(\theta)\theta_j]\}^2
 =\sum_j\{A_{j1}(z)+A_{j2}(z)\}^2
 =\sum_{j,r}A_{jr}(z)^2.
\]
The last equality uses mutual exclusivity within each pair.  On these tables
the target is $\sum_j(c_{j1}+c_{j2})\theta_j$, whose prior second moment
is $\sum_j(c_{j1}+c_{j2})^2$.  Since the estimator is unbiased at every
one of these tables,
\begin{align*}
 \E_\theta\Var_P(T\mid Y(\theta))
 &=\E_P\E_\theta t_Z(\theta)^2
   -\sum_j(c_{j1}+c_{j2})^2\\
 &\geq\sum_{j,r}\frac{c_{jr}^2}{q_{P,jr}}
       -\sum_j(c_{j1}+c_{j2})^2\\
 &=\sum_j\left\{
 c_{j1}^2\left(\frac1{q_{P,j1}}-1\right)
 +c_{j2}^2\left(\frac1{q_{P,j2}}-1\right)-2c_{j1}c_{j2}\right\}.
\end{align*}
Terms with $c_{jr}=q_{P,jr}=0$ are interpreted as zero.
This proves equation~\eqref{eq:crossblocklb-multipair}. The worst-case
design variance is at least the prior-average design variance
$\E_\theta\Var_P(T\mid Y(\theta))$.  The argument used no independence between pairs' observation
events.
\end{proof}

\subsection{Sampling-design reduction of the three minimax problems}
\label{app:observation-law-proofs}

In this subsection, we prove Proposition~\ref{thm:three-problem-characterization}
using two lemmas. Lemma~\ref{lem:revelation-sufficiency} shows that the
fixed-design minimax risks are unchanged when estimators
$T(Z,Y_{O_f(Z)})$ are restricted to functions of
$(O_f(Z),Y_{O_f(Z)})$. The proof constructs such an estimator by
averaging over assignments that reveal the same indices; this averaging
cannot increase risk and preserves unbiasedness. This reduces the fixed-design problem
to the induced sampling design, or equivalently the weighted sampling hypergraph, and the restricted-design problem to the collection of
attainable induced sampling designs. Lemma~\ref{lem:achievable-laws} shows that,
when all designs on $\cZ_0$ are allowed, every probability distribution
on the feasible hyperedges is attainable. Optimizing over these
probabilities then leaves only the unweighted sampling hypergraph.

\begin{lemma}[Sufficiency of the induced sampling design]
\label{lem:revelation-sufficiency}
Fix a design $P$. For squared-error loss, both the minimax risk over all
estimators, $\Rall(P;f,a)$, and the minimax risk over unbiased
estimators, $\Runb(P;f,a)$, are unchanged if estimators are restricted to
measurable functions of $O_f(Z)$ and $Y_{O_f(Z)}$.
\end{lemma}

\begin{proof}
Fix $P$ and an admissible table $Y$. Let $T(D)$ be an estimator based
on the full observed data $D=(Z,Y_{O_f(Z)})$. Put $A=O_f(Z)$ and let
$D_0=(A,Y_A)$ be the data retained after discarding the assignment.
For an observed data pattern $A$ with $Q_{P,f}(A)>0$, the conditional probability of
assignment $z$ given this pattern is
\[
 P(Z=z\mid O_f(Z)=A)
 =\frac{P(z)\1\{O_f(z)=A\}}{Q_{P,f}(A)}.
\]
For the fixed table $Y$, the event $\{D_0=(A,Y_A)\}$ is precisely the
event $\{O_f(Z)=A\}$.  Thus the same conditional probabilities apply
when we condition on $D_0$.  Define
\[
 \overline T(A,y_A)
 =\sum_{z:O_f(z)=A}\frac{P(z)}{Q_{P,f}(A)}T(z,y_A).
\]
The probabilities in this sum are known and independent of the table, so
the formula defines a single estimator using only $(A,y_A)$.
For each fixed $Y$, $\overline T(D_0)=\E_P[T(D)\mid D_0,Y]$.
This is Rao--Blackwellization of the design-based estimator with respect
to the retained data.  Conditional Jensen and iterated expectation give
\begin{align*}
 \E_P[(\overline T(D_0)-\tau_a(Y))^2\mid Y]
 &\leq\E_P\{\E_P[(T(D)-\tau_a(Y))^2\mid D_0,Y]\mid Y\}\\
 &=\E_P[(T(D)-\tau_a(Y))^2\mid Y],\\
 \E_P[\overline T(D_0)\mid Y]&=\E_P[T(D)\mid Y].
\end{align*}
Thus averaging cannot increase risk and preserves unbiasedness.
Conversely, every function of $D_0$ is also a function of the full data
$D$, so it is allowed in the original estimator class.  The infima over full-data estimators and retained-data estimators
therefore agree, both for all estimators and for unbiased estimators.  Patterns of probability
zero may be assigned any value.
\end{proof}

\begin{lemma}[Every hypergraph law is implementable]
\label{lem:achievable-laws}
Every probability distribution on the hyperedges of
\(\cH_{\mathrm{obs}}(f;\cZ_0)\) is induced by some design in
\(\cP(\cZ_0)\).
\end{lemma}

\begin{proof}
Write $\cH_{\mathrm{obs}}(f;\cZ_0)=(K,\mathcal E)$.
For each hyperedge $A\in\mathcal E$, choose one witness
assignment \(z_A\in\cZ_0\) such that \(O_f(z_A)=A\).  Given a probability law
\(Q\) on the hyperedges, define $P(z_A)=Q(A)$ and assign probability
zero to all other assignments.  Witnesses for distinct hyperedges are
distinct, since one assignment has only one observation pattern.
Consequently $P$ is a probability law on $\cZ_0$ and, for every hyperedge $A$,
\[
 Q_{P,f}(A)=\sum_{z:O_f(z)=A}P(z)=P(z_A)=Q(A).
\]
\end{proof}

\begin{resultrestatement}{Proposition~\ref{thm:three-problem-characterization}}{Sampling-design reduction of the three minimax problems}
Fix $f$, $a$, and $\cY(f)$.
\begin{enumerate}
\item[\textup{(i)}] For a fixed design $P$, both $\Rall(P)$ and
$\Runb(P)$ depend on $P$ only through $Q_{P,f}$.
\item[\textup{(ii)}]
For a restricted design class $\cC$, the joint minimax risk $\Rallstar(\cC)$ depends on the class only through $\mathcal Q_f(\cC)$. The corresponding risk $\Runbstar(\cC)$ depends only on $\mathcal Q_f(\cC_+)$, where $\cC_+=\cC\cap\cP_+(f,a)$.
\item[\textup{(iii)}]
When all designs on $\cZ_0$ are allowed, both $\Rallstar(\cP(\cZ_0))$ and $\Runbstar(\cP(\cZ_0))$ depend on the assignment and exposure model only through $\cH_{\mathrm{obs}}(f;\cZ_0)$.
\end{enumerate}
\end{resultrestatement}

\begin{proof}[Proof of Proposition~\ref{thm:three-problem-characterization}]
By Lemma~\ref{lem:revelation-sufficiency}, the fixed-design risks can be
computed using estimators of $(A,Y_A)$, where $A\sim Q_{P,f}$.
For every fixed table, the distribution of $(A,Y_A)$ is
specified by $Q_{P,f}$.
Recall from Definition~\ref{def:weighted-co-observation-hypergraph} that,
writing $\cH_{\mathrm{obs}}(f;\cZ_0)=(K,\mathcal E)$, the weighted
sampling hypergraph is $(\cH_{\mathrm{obs}}(f;\cZ_0),w_P)$,
where $w_P(A)=Q_{P,f}(A)$ for $A\in\mathcal E$.
In particular,
\[
 \operatorname{supp}(Q_{P,f})
 =\{A\in\mathcal E:w_P(A)>0\}.
\]
Conversely, the weighted hypergraph determines the induced sampling design by
$Q_{P,f}(A)=w_P(A)$ for $A\in\mathcal E$ and $Q_{P,f}(A)=0$ otherwise.
Thus two designs inducing the same weighted sampling hypergraph
induce the same distribution of $(A,Y_A)$ for every fixed table and
hence have the same minimax risks, proving part (i).

For part (ii), recall that
\[
 \Rallstar(\cC)=\inf_{P\in\cC}\Rall(P),\qquad
 \Runbstar(\cC)=\inf_{P\in\cC_+}\Runb(P).
\]
By part (i), designs with the same induced sampling design contribute the same
value of $\Rall(P)$ to the first infimum and the same value of $\Runb(P)$
to the second. Thus $\Rallstar(\cC)$ is determined by the set
$\mathcal Q_f(\cC)=\{Q_{P,f}:P\in\cC\}$, and $\Runbstar(\cC)$ by
$\mathcal Q_f(\cC_+)=\{Q_{P,f}:P\in\cC_+\}$. These are precisely the
collections of weighted sampling hypergraphs attainable in the two
classes. This proves part (ii).

For part (iii), write $\cH_{\mathrm{obs}}(f;\cZ_0)=(K,\mathcal E)$
and let $\Delta(\mathcal E)$ denote all probability laws on $\mathcal E$.
Every design $P\in\cP(\cZ_0)$ induces a sampling design
$Q_{P,f}\in\Delta(\mathcal E)$, or equivalently a weighted
sampling hypergraph. Conversely, Lemma~\ref{lem:achievable-laws}
constructs a design inducing any chosen $Q\in\Delta(\mathcal E)$.
Consequently,
\[
 \mathcal Q_f(\cP(\cZ_0))=\Delta(\mathcal E).
\]
Thus optimization over all designs is equivalent to optimization over
all probability laws on the feasible hyperedges. In the unbiased problem,
positivity requires $\sum_{A\in\mathcal E:k\in A}Q(A)>0$ for each
$k\in K_a$, which is a condition on $Q$ alone. Hence the attainable
induced sampling designs for this problem are precisely the laws in
$\Delta(\mathcal E)$ satisfying these inequalities.
In both cases, the hyperedge probabilities are optimized over, and
both risks therefore depend only on the unweighted sampling hypergraph, with $f$, $a$, and $\cY(f)$ held fixed.
This proves part (iii).
\end{proof}

\subsection{An example with same hypergraph, but different fixed design minimax rates}
\label{subsec:app-cluster-with-noise}

\label{subsec:cluster-with-noise}

Proposition~\ref{thm:three-problem-characterization}(i) shows that the the fixed-design minimax risks depend on
the full induced sampling design $Q_{P,f}$, or equivalently, the weighted sampling hypergraph. In this subsection, we construct an example with two different designs and same unweighted sampling hypergraph and first-order inclusion probabilities, but different unbiased minimax rates. These examples show that the fixed design minimax risks depend on the weighted sampling hypergraph and not just the unweighted hypergraph, or the marginal probabilities.
We use the no-interference setup of Section~\ref{subsec:ate-benchmark}, with
\(\cZ=\{0,1\}^n\), \(f_i(z)=z_i\), potential outcomes in \([-1,1]\),
and ATE coefficients \(a_{i,1}=1/n\) and \(a_{i,0}=-1/n\).

The first design is the independent Bernoulli law $P_{1/2}$ with $P(Z_i=1) = 1/2$, for which
\begin{equation}
  \Runb(P_{1/2};f,a)=\frac4n
  \label{eq:cluster-bernoulli}
\end{equation}
by Theorem~\ref{thm:ate-joint-benchmark}.  The second is a clustered
design built from a common shock and independent noise. Let
$C\sim{\rm Bernoulli}(1/2)$ and let
$\varepsilon_1,\ldots,\varepsilon_n$ be independent
${\rm Bernoulli}(\eta)$ variables, independent of $C$, where
$\eta\in(0,1)$. Set
\begin{equation}
  Z_i=C\mathbin{\oplus}\varepsilon_i,\qquad i=1,\ldots,n.
  \label{eq:cluster-with-noise}
\end{equation}
where $\oplus$ denotes addition modulo two. Write $P_{\rm cw}$ for this design.
Both designs are laws on $\cZ=\{0,1\}^n$ that give positive probability to
every assignment, so they have the same sampling hypergraph
$\cH_{\mathrm{obs}}(f)$, the family of all one-index-per-unit
transversals.  Moreover $\Pr_{P_{\rm cw}}(Z_i=1)=1/2$ for every $i$, so every
first-order inclusion probability matches $P_{1/2}$.  The designs differ only in
their joint dependence.  For $\eta=c/n$, the probabilities that all units are untreated and that
all units are treated are equal:
\[
 P_{\rm cw}(Z=0^n)=P_{\rm cw}(Z=1^n)
 =\frac12\{(1-c/n)^n+(c/n)^n\}\longrightarrow\frac{e^{-c}}2.
\]

\begin{theorem}[Cluster randomization with vanishing noise]
\label{thm:cluster-with-noise}
Let $\eta=c/n$ for a fixed $c>0$ and $n>c$, and let $P_{\rm cw}$ denote the
design defined by equation~\eqref{eq:cluster-with-noise}.  Then $P_{\rm cw}$ and $P_{1/2}$ have
the same sampling hypergraph and the same first-order inclusion probabilities, and as $n\to\infty$,
\begin{align}
  \Runb(P_{\rm cw};f,a)
  &\geq\frac{2e^{-c}}{2-e^{-c}}\,\{1-o(1)\},
  \label{eq:cluster-lower}\\
  \Runb(P_{\rm cw};f,a)
  &\leq
  4\left\{\frac1n+\left(1-\frac1n\right)(1-2\eta)^2\right\}
  \leq 4+o(1).
  \label{eq:cluster-upper}
\end{align}
Hence $\Runb(P_{\rm cw};f,a)=\Theta(1)$, whereas
$\Runb(P_{1/2};f,a)=4/n$.
\end{theorem}

\begin{proof}[Proof of Theorem~\ref{thm:cluster-with-noise}]
For the lower bound, we use the lower bound from disjoint perturbations
with the two sets
\[
 B_1=\{(i,1):i=1,\ldots,n\},\qquad
 B_0=\{(i,0):i=1,\ldots,n\}.
\]
These sets partition the potential-outcome table into its treatment and
control arms. To see the associated perturbations explicitly, let
$\theta_1,\theta_0$ be independent uniform signs and set
$Y_i(1)=\theta_1$ and $Y_i(0)=\theta_0$ for every unit $i$.
This gives four admissible tables, and their ATE is
$\tau_a(Y^\theta)=\theta_1-\theta_0$. Thus the perturbation coefficients
are $c_1=1$ and $c_0=-1$.

The treatment-arm set $B_1$ is entirely missing exactly when every unit
receives control. Under the common-shock construction, this occurs either
when $C=0$ and all $\varepsilon_i=0$, or when $C=1$ and all
$\varepsilon_i=1$. Its probability is therefore
\[
 \rho_n=\frac12\{(1-\eta)^n+\eta^n\}
       =\frac{e^{-c}}2\{1-o(1)\}
\]
along the sequence $\eta=c/n$. By symmetry, the control-arm set $B_0$
is entirely missing with the same probability, on the event $Z=1^n$.
Although every unit has treatment probability $1/2$, each entire arm
therefore remains unobserved with probability bounded away from zero.
Applying Theorem~\ref{thm:fixed-joint-missingness}(ii) to the packing
$\{B_1,B_0\}$ gives
\[
 \Runb(P_{\rm cw};f,a)
 \geq c_1^2\frac{\rho_n}{1-\rho_n}
       +c_0^2\frac{\rho_n}{1-\rho_n}
 =2\frac{\rho_n}{1-\rho_n}
 =\frac{2e^{-c}}{2-e^{-c}}\,\{1-o(1)\},
\]
which proves equation~\eqref{eq:cluster-lower}.

For the upper bound, use the HT estimator with the marginal treatment
probabilities $1/2$. Write $S_i=2Z_i-1$ and
$b_i=Y_i(1)+Y_i(0)$. Direct expansion gives
\[
 \widehat\tau_{\rm HT}
 =\frac2n\sum_i\{Z_iY_i(1)-(1-Z_i)Y_i(0)\}
 =\tau_a(Y)+\frac1n\sum_i b_iS_i.
\]
Since $\E S_i=0$, this estimator is unbiased for every admissible table.
Moreover, $S_i=(2C-1)(1-2\varepsilon_i)$, so
$\Var(S_i)=1$ and, for $i\ne j$,
$\Cov(S_i,S_j)=\E(S_iS_j)=(1-2\eta)^2$.
Consequently,
\[
 \Var_{P_{\rm cw}}(\widehat\tau_{\rm HT})
 =\frac1{n^2}\left\{\sum_i b_i^2
       +(1-2\eta)^2\sum_{i\ne j}b_ib_j\right\}
 \leq4\left\{\frac1n+
       \left(1-\frac1n\right)(1-2\eta)^2\right\},
\]
where the inequality uses $|b_i|\leq2$.
This bound holds uniformly over the potential-outcome tables and proves
equation~\eqref{eq:cluster-upper}.

\end{proof}

For fixed $\eta\in(0,1)$, each arm block is detected with probability
tending to one, so the resulting missingness lower bound tends to zero.
The scaling $\eta=c/n$ keeps the joint missingness probabilities bounded
away from zero and gives the constant-order lower bound above. The example shows that the unweighted sampling hypergraph and first-order inclusion probabilities do not
suffice to determine the unbiased fixed-design minimax rate. %

\subsection{From fixed-design lower bounds to lower bounds on restricted designs}
\label{subsec:packing-design-class-bounds}
In this subsection, we derive lower bounds for the joint minimax risks
over a design class $\cC$. The lower bounds from disjoint perturbations
(Theorem~\ref{thm:fixed-joint-missingness}) bound the risk under a specified
design and perturbation packing. To obtain a lower bound that
optimizes over the design, we first take the supremum of the fixed design lower bounds
over perturbation packings $\cB$ and the fixed perturbation signs
$s_{b,k}\in\{-1,1\}$,
and then take the infimum over designs in $\cC$.
The packing and signs giving the strongest lower bound may therefore
differ across designs. We express the resulting bounds in terms of the
attainable induced sampling designs $\mathcal Q_f(\cC)$, using
Proposition~\ref{thm:three-problem-characterization}.
Corollary~\ref{thm:restricted-joint-missingness} states these bounds for
both estimator classes.

For an induced sampling design $Q$, define
\[
 q_Q(B)=\sum_{A:A\cap B\neq\varnothing}Q(A)
\]
to be the probability that at least one potential outcome in $B$ is observed.
For a perturbation packing $\cB=\{B_1,\ldots,B_m\}$, define the
detection vector $q_Q=(q_Q(B_1),\ldots,q_Q(B_m))$.
Its $b$th entry is the probability that at least one potential outcome in $B_b$ is observed under $Q$. Fix perturbation signs
$s=(s_{b,k})$, giving coefficients $c_b=\sum_{k\in B_b}a_ks_{b,k}$.
The all-estimator and unbiased-estimator lower bounds in
Theorem~\ref{thm:fixed-joint-missingness} depend on $Q$ only through
$q_Q$. Write their respective lower-bound expressions as
\begin{equation}
  \Phi_{\rm all}(q;\cB,s)
  =\sum_b c_b^2(1-q_b),
  \qquad
  \Phi_{\rm unb}(q;\cB,s)
  =\sum_{b:c_b\neq0}c_b^2\left(\frac1{q_b}-1\right).
  \label{eq:phi-all}
\end{equation}
A term in $\Phi_{\rm unb}$ is $+\infty$ if $q_b=0$ and $c_b\neq0$.

\begin{corollary}[Lower bounds under design restrictions]
\label{thm:restricted-joint-missingness}
\label{cor:three-problems}
\label{prop:price-restricted}
For a design class $\cC$,
\begin{equation}
  \Rallstar(\cC)
  \geq
  \inf_{Q\in\mathcal Q_f(\cC)}
  \sup_{(\cB,s)}\Phi_{\rm all}(q_Q;\cB,s),
  \label{eq:joint-cert-all}
\end{equation}
and for unbiased estimation,
\begin{equation}
  \Runbstar(\cC)
  \geq
  \inf_{Q\in\mathcal Q_f(\cC_+)}
  \sup_{(\cB,s)}\Phi_{\rm unb}(q_Q;\cB,s),
  \label{eq:joint-cert-unb}
\end{equation}
where in both displays the supremum is over perturbation packings and
choices of signs.
\end{corollary}

\begin{proof}[Proof of Corollary~\ref{thm:restricted-joint-missingness}]
For a fixed design $P$, let $Q=Q_{P,f}$ be its induced sampling design.
The lower bounds from disjoint perturbations
(Theorem~\ref{thm:fixed-joint-missingness}) hold for every perturbation
packing $\cB$ and choice of perturbation signs $s$.
Taking the supremum over these choices gives
\begin{align*}
 \Rall(P)&\geq\sup_{(\cB,s)}\Phi_{\rm all}(q_Q;\cB,s),\\
 \Runb(P)&\geq\sup_{(\cB,s)}\Phi_{\rm unb}(q_Q;\cB,s).
\end{align*}
Taking the infimum over designs $P\in\cC$ in the first inequality gives
\[
 \inf_{P\in\cC}\Rall(P)
 \geq\inf_{P\in\cC}\sup_{(\cB,s)}
       \Phi_{\rm all}(q_{Q_{P,f}};\cB,s)
 =\inf_{Q\in\mathcal Q_f(\cC)}\sup_{(\cB,s)}
       \Phi_{\rm all}(q_Q;\cB,s).
\]
The last equality follows from the definition of $\mathcal Q_f(\cC)$:
it is exactly the set of laws $Q_{P,f}$ induced by designs in $\cC$.
The same calculation with $\Phi_{\rm unb}$ and $\cC_+$ proves the
unbiased statement.  Designs outside $\cC_+$ have an unobservable
relevant potential outcome and admit no unbiased estimator.
\end{proof}

When all designs on $\cZ_0$ are allowed, take $\cC=\cP(\cZ_0)$ in the
same bounds.  The infimum then ranges over all probability laws on the
hyperedges of $\cH_{\mathrm{obs}}(f;\cZ_0)$.  For the unbiased bound,
restrict to laws that give every relevant potential outcome positive probability
of observation.

For a restricted class, the feasible observed sets are not sufficient to specify
the attainable probabilities.  Section~\ref{subsec:cluster-with-noise}
gives two designs with the same support and first-order inclusion probabilities whose minimax risks differ by a factor of $n$.

\begin{remark}[Order of optimization]
\label{rem:order-of-optimization}
Corollary~\ref{thm:restricted-joint-missingness} places the supremum
over perturbation packings and signs \emph{inside} the infimum over
induced sampling designs. The inequality
$\inf_Q\sup_{(\cB,s)}\Phi\geq\sup_{(\cB,s)}\inf_Q\Phi$
holds for either lower-bound expression $\Phi$, but the two orders
have different meanings. In the stated order, the perturbation packing
and signs may be chosen separately for each allowed induced sampling design $Q$,
to give the strongest lower bound for that law. Thus the worst-case
potential-outcome perturbations may differ across designs in the class.
In the swapped order, a single packing and choice of signs must give
a lower bound uniformly over all allowed induced sampling designs; the resulting
lower bound can be zero even when the stated order gives a positive
lower bound. Any fixed packing $\cB$ and signs $s$ nevertheless give
the valid, and in general weaker, lower bound
$\inf_{Q\in\mathcal Q_f(\cC)}\Phi(q_Q;\cB,s)$.
For the unbiased problem, the same comparison uses
$Q\in\mathcal Q_f(\cC_+)$.
\end{remark}

\subsection{Lower bounds for the all-design minimax risk from the sampling hypergraph}
\label{app:hypergraph-certificates}

The preceding subsection optimized the fixed-design packing lower bounds
over attainable induced sampling designs. When all designs on $\cZ_0$ are
allowed, every probability law on the feasible hyperedges is attainable
(Lemma~\ref{lem:achievable-laws}). In this subsection, we derive explicit
lower-bound functionals of the unweighted sampling hypergraph for
the joint minimax problems, optimizing over all designs and either all
estimators or unbiased estimators. We obtain these lower bounds by
controlling which perturbation sets of potential outcomes a feasible assignment can reveal and
which must remain unobserved. The functionals use the hypergraph
restricted to the estimand-relevant indices and the estimand coefficients.

We begin by proving Theorem~\ref{thm:global-hypergraph-lower}, which is the general result underlying all the lower bound results in this subsection. It applies to arbitrary linear estimands and gives lower bounds for both estimator classes.
We then develop three specializations of this theorem, used throughout
the applications to exploit different constraints on which relevant
potential outcomes can be observed together. First, in Section~\ref{subsec:one-arm-rank-bounds}, Theorem~\ref{cor:rank-lower} gives a reciprocal
\emph{sampling-rank} lower bound. This lower bound uses only
one arm of an estimand and bounds minimax risks of both estimator classes in terms of the largest number of that arm's potential outcomes observable under a single assignment. Second, in Section~\ref{subsec:observed-arm-count-bounds}, Corollary~\ref{cor:two-arm-coverage} gives lower bounds based on exposure counts by using
two arms together, and Corollary~\ref{cor:multi-arm-coverage} extends the result to finitely many arms. These bounds apply to unbiased estimators and depend on the expected exposure counts jointly attainable across the arms. Third, the specialization in Section~\ref{subsec:block-certificates} groups potential-outcomes within two arms and pairs the groups across arms. Theorem~\ref{thm:arm-paired-certificate} allows unequal weights across
pairs and gives a lower bound valid for all estimators. When all pairs
have equal weight, it yields Theorem~\ref{thm:paired-sets-main}.

Before developing these bounds, we define the restrictions of the
sampling hypergraph to the estimand-relevant potential outcomes and to
individual arms.

\begin{definition}[Estimand and arm traces]
\label{def:trace-hypergraphs}
For an estimand $a$, let $K_a=\{k\in K:a_k\neq0\}$ be its relevant
set of potential outcomes.  The estimand trace of the sampling hypergraph is
\begin{equation}
  \cH_{\mathrm{obs}}^a(f;\cZ_0)
  =
  \bigl(K_a,\{O_f(z)\cap K_a:z\in\cZ_0\}\bigr).
  \label{eq:target-relevant-trace}
\end{equation}
If $A\subseteq K_a$ is a set of potential outcomes, the trace on $A$ is
\[
  \cH_{\mathrm{obs}}^A(f;\cZ_0)
  =
  \bigl(A,\{O_f(z)\cap A:z\in\cZ_0\}\bigr).
\]

\end{definition}

When the sampling hypergraph is abbreviated by $\cH$, we write these traces as $\cH^a$ and $\cH^A$. Thus $\cH^a$ is the sampling hypergraph restricted
to the estimand-relevant potential outcomes.
The trace $\cH^A$ retains only the vertices in $A$ and intersects each original
hyperedge with $A$. When $A$ consists of the potential outcomes in one exposure arm of a two-arm contrast, we call it an \emph{arm trace}.
 The estimand trace retains constraints on observing different arms
together, even when each arm can be observed in full under a separate
assignment. To derive a lower bound using only a subset $A\subseteq K_a$, we restrict attention to admissible tables whose potential outcomes
outside $A$ are fixed at zero. A lower bound over this smaller class is
also a lower bound for the original problem. We use this restriction
below to obtain one-arm lower bounds.

\subsubsection{General lower bounds for linear estimands from the sampling hypergraph}
In this subsection, we derive general lower bounds for the joint minimax
risks of an arbitrary linear estimand $\tau_a$, optimizing over all
designs on $\cZ_0$ and either all estimators or unbiased estimators.
The lower bounds are expressed as functionals of the estimand trace
$\cH_{\mathrm{obs}}^a$ and the coefficients of the estimand.

Fix a perturbation packing $\cB=\{B_1,\ldots,B_m\}$ and perturbation
signs $s_{b,k}\in\{-1,1\}$, defining
$c_b=\sum_{k\in B_b}a_ks_{b,k}$. Each $B_b$ indexes a set
of potential outcomes that vary with one unknown sign $\theta_b$
under the hypercube prior. Observing any outcome in $B_b$ reveals
$\theta_b$; observing no outcome in $B_b$ leaves $\theta_b$ unknown. Thus, an assignment $z$ detects the perturbation on $B_b$ exactly when the corresponding hyperedge $O_f(z)$ intersects $B_b$. For the all-estimator lower bound, each undetected set $B_b$ contributes $c_b^2$ to the posterior variance of the estimand. The expected posterior variance lower bounds the minimax risk. We therefore bound the posterior variance of the estimand remaining under every feasible assignment.

Write $\cH_{\mathrm{obs}}^a=(K_a,\mathcal E_a)$, where
$\mathcal E_a=\{O_f(z)\cap K_a:z\in\cZ_0\}$ is its set of hyperedges.
To quantify how much a feasible assignment detects or misses a perturbation, assign a
nonnegative weight $u_b$ to each perturbation set $B_b$, and write $u=(u_b)$.
The \emph{maximum observed weight} is the largest total weight of perturbation
sets detected by a feasible hyperedge,
\begin{equation}
  \nu_u(\cH_{\rm obs}^a,\cB)
  =\max_{E\in\mathcal E_a}
    \sum_{b:E\cap B_b\neq\varnothing}u_b ,
  \label{eq:max-weighted-coverage}
\end{equation}
and the \emph{minimum missing weight} is the smallest total weight of
perturbation sets missed by a feasible hyperedge,
\begin{equation}
  \mu_u(\cH_{\rm obs}^a,\cB)
  =\min_{E\in\mathcal E_a}\sum_{b:E\cap B_b=\varnothing}u_b
  =\sum_bu_b-\nu_u(\cH_{\rm obs}^a,\cB).
  \label{eq:packing-deficit-functional}
\end{equation}
For fixed packing and weights, both quantities depend only on the feasible
hyperedges. For the all-estimator bound, we take $u_b=c_b^2$. The minimum
missing weight, written $\mu_{c^2}$, then bounds the posterior variance remaining
under every assignment. For the unbiased bound, the weights can be chosen
freely. In this case, the maximum observed weight bounds $\sum_b u_bq_b$, where $q_b$ is the probability of detecting $B_b$. The proof uses this bound with
Cauchy--Schwarz to control the reciprocal detection probabilities.

\begin{theorem}[General lower bounds from the sampling hypergraph]
\label{thm:global-hypergraph-lower}
\label{prop:joint-price}
\label{thm:design-free-hitting}
\label{prop:design-free-unbiased}
For every perturbation packing and every choice of signs,
\begin{equation}
  \Rallstar\bigl(\cP(\cZ_0);f,a\bigr)
  \geq
  \mu_{c^2}(\cH_{\rm obs}^a,\cB).
  \label{eq:design-free-hitting}
\end{equation}
For every such packing and signs, and every $u\geq0$ with
$\nu_u(\cH_{\rm obs}^a,\cB)>0$,
\begin{equation}
  \Runbstar\bigl(\cP(\cZ_0);f,a\bigr)
  \geq
  \frac{\bigl(\sum_b|c_b|\sqrt{u_b}\bigr)^2}
       {\nu_u(\cH_{\rm obs}^a,\cB)}
  -\sum_bc_b^2.
  \label{eq:design-free-unbiased}
\end{equation}
Both right-hand sides may be maximized over packings and signs.
\end{theorem}

When every perturbation set has weight $u_0$, the minimum missing weight
is $u_0$ times the smallest number of sets missed by a feasible hyperedge:
\begin{equation}
  \mu_u(\cH_{\rm obs}^a,\cB)
  =u_0\Bigl\{
    |\cB|-\max_{E\in\mathcal E_a}\#\{b:E\cap B_b\neq\varnothing\}
  \Bigr\},
  \label{eq:equal-block-capacity-lower}
\end{equation}
\begin{proof}[Proof of Theorem~\ref{thm:global-hypergraph-lower}]
Fix a perturbation packing, its signs, and a design $P$. We first apply
the fixed-design lower bounds from disjoint perturbations. For each
estimator class, we then obtain a lower bound independent of $P$, which
also bounds the infimum over all designs.
Write $Q=Q_{P,f}$, $q_b=q_Q(B_b)=q_P(B_b)$, and $E=O_f(Z)$ for the
observed hyperedge. Since each $B_b\subseteq K_a$, replacing $E$ by its
trace $E\cap K_a$ does not change whether it intersects $B_b$.

For all estimators, the fixed-design lower bound from disjoint
perturbations (Theorem~\ref{thm:fixed-joint-missingness}(i)) gives
\[
  \Rall(P)\geq\sum_b c_b^2(1-q_b).
\]
For each feasible observed hyperedge $E$, the definition of the maximum
observed weight gives
\[
  \sum_{b:E\cap B_b\neq\varnothing}c_b^2
  \leq\nu_{c^2}(\cH_{\rm obs}^a,\cB).
\]
Therefore,
\[
  \sum_b c_b^2q_b
  =\E_Q\left[\sum_{b:E\cap B_b\neq\varnothing}c_b^2\right]
  \leq\nu_{c^2}(\cH_{\rm obs}^a,\cB).
\]
Substituting into the fixed-design lower bound yields
\[
  \Rall(P)
  \geq\sum_b c_b^2-\nu_{c^2}(\cH_{\rm obs}^a,\cB)
  =\mu_{c^2}(\cH_{\rm obs}^a,\cB).
\]
The right-hand side is the same for every $P$. Taking the infimum over
all designs proves equation~\eqref{eq:design-free-hitting}.

For unbiased estimators, only perturbation sets with nonzero set effect
$c_b$ contribute to the fixed-design lower bound. Define their index set by
\[
  J=\{b\in\{1,\ldots,m\}:c_b\neq0\}.
\]
For $b\notin J$, changing the sign on $B_b$ does not change the estimand,
so these sets contribute zero to the lower bound. If $q_b=0$ for some
$b\in J$, the design never detects a perturbation that changes the
estimand; no unbiased estimator exists and the claimed lower bound holds.
Otherwise, the fixed-design lower bound from
disjoint perturbations (Theorem~\ref{thm:fixed-joint-missingness}(ii)) gives
\begin{equation}
  \Runb(P)\geq\sum_{b\in J}\frac{c_b^2}{q_b}
                     -\sum_{b\in J}c_b^2.
  \label{eq:hypergraph-unb-proof-fixed}
\end{equation}
Fix nonnegative weights $u=(u_b)_{b=1}^m$ with
$\nu_u(\cH_{\rm obs}^a,\cB)>0$. Cauchy--Schwarz gives
\begin{equation}
  \left(\sum_{b\in J}|c_b|\sqrt{u_b}\right)^2
  =\left(\sum_{b\in J}\frac{|c_b|}{\sqrt{q_b}}\sqrt{u_bq_b}\right)^2
  \leq\left(\sum_{b\in J}\frac{c_b^2}{q_b}\right)
       \left(\sum_{b\in J}u_bq_b\right).
  \label{eq:hypergraph-unb-proof-cs}
\end{equation}
For every feasible hyperedge $E$, nonnegativity of the weights and the
definition of $\nu_u$ imply
\begin{equation}
  \sum_{\substack{b\in J:\ E\cap B_b\neq\varnothing}}u_b
  \leq\sum_{\substack{1\leq b\leq m:\ E\cap B_b\neq\varnothing}}u_b
  \leq\nu_u(\cH_{\rm obs}^a,\cB).
  \label{eq:hypergraph-unb-proof-hyperedge}
\end{equation}
Taking expectations in equation~\eqref{eq:hypergraph-unb-proof-hyperedge},
the expected detected weight of the sets indexed by $J$ satisfies
\begin{equation}
  \sum_{b\in J}u_bq_b
  =\E_Q\left[\sum_{\substack{b\in J:\ E\cap B_b\neq\varnothing}}u_b\right]
  \leq\nu_u(\cH_{\rm obs}^a,\cB).
  \label{eq:hypergraph-unb-proof-expected}
\end{equation}
Combining equation~\eqref{eq:hypergraph-unb-proof-expected} with the
Cauchy--Schwarz inequality in equation~\eqref{eq:hypergraph-unb-proof-cs}
and dividing by the positive quantity
$\nu_u(\cH_{\rm obs}^a,\cB)$ gives
\begin{equation}
  \sum_{b\in J}\frac{c_b^2}{q_b}
  \geq\frac{\bigl(\sum_{b\in J}|c_b|\sqrt{u_b}\bigr)^2}
                {\nu_u(\cH_{\rm obs}^a,\cB)}.
  \label{eq:hypergraph-unb-proof-reciprocal}
\end{equation}
Substituting equation~\eqref{eq:hypergraph-unb-proof-reciprocal} into the
fixed-design unbiased lower bound in equation~\eqref{eq:hypergraph-unb-proof-fixed}
yields
\begin{equation}
  \begin{aligned}
  \Runb(P)
  &\geq\frac{\bigl(\sum_{b\in J}|c_b|\sqrt{u_b}\bigr)^2}
                {\nu_u(\cH_{\rm obs}^a,\cB)}
       -\sum_{b\in J}c_b^2\\
  &=\frac{\bigl(\sum_{b=1}^m|c_b|\sqrt{u_b}\bigr)^2}
                {\nu_u(\cH_{\rm obs}^a,\cB)}
       -\sum_{b=1}^m c_b^2,
  \end{aligned}
  \label{eq:hypergraph-unb-proof-conclusion}
\end{equation}
where the equality uses $c_b=0$ for $b\notin J$. The right-hand side is
the same for every $P$. Taking the infimum over all designs proves
equation~\eqref{eq:design-free-unbiased}.
\end{proof}

\subsubsection{Lower bounds based on one arm of an estimand}
\label{subsec:one-arm-rank-bounds}
In this subsection, we derive lower bounds for the joint minimax risks
using the \emph{sampling rank} of one equal-weight arm of an estimand. The sampling rank $r(V)$
is the largest number of potential outcomes in $V$ that any
feasible assignment can observe. For unbiased estimators, we use singleton perturbation sets to get the rank based lower bounds. For lower bounds on all estimators, we group potential outcomes into
perturbation sets so that every assignment leaves sufficiently many sets of potential outcomes
unobserved. We obtain both lower bounds by applying the general lower bounds from the sampling hypergraph (Theorem~\ref{thm:global-hypergraph-lower}).
We restate the one-arm theorem before giving the two constructions.

For the chosen arm $V$, define the observed exposure count
$n_V(z)=|O_f(z)\cap V|$. $n_V(z)$ is the number of units revealing a potential outcome in $V$ under assignment $z$. For a full exposure-$d$
arm, $n_V(z)=n_d(z)$ from Section~\ref{subsec:observed-data}.
Since $r(V)=\max_{z\in\cZ_0}n_V(z)$, every design $P$ on $\cZ_0$
satisfies
\[
 \E_P[n_V(Z)]=\sum_{k\in V}\pi_{P,k}\leq r(V).
\]
Thus no design can make the expected exposure count within $V$ $\E_P[n_V(Z)]$ exceed the rank $r(V)$. For unbiased estimation, the singleton lower
bound and Cauchy--Schwarz give
$\Runb(P)\geq A^2\{1/\E_P[n_V(Z)]-1/N_V\}$ whenever all indices
in $V$ have positive observation probability. Substituting
$\E_P[n_V(Z)]\leq r(V)$ gives the reciprocal-rank lower bound below;
if a relevant index has zero observation probability, no unbiased
estimator exists. The all-estimator argument instead uses the grouped
perturbations described above.

\begin{resultrestatement}{Theorem~\ref{cor:rank-lower}}{Reciprocal sampling-rank lower bound}
Let $V\subseteq K_a$ contain $N_V=|V|$ relevant potential outcomes from one
estimand arm, each with estimand coefficient satisfying $|a_k|=A/N_V$.
Write $\cH_{\mathrm{obs}}(f;\cZ_0)=(K,\mathcal E)$ and define the
sampling rank
\[
  r(V)=\max_{E\in\mathcal E}|E\cap V|.
\]
If $r(V)\geq1$, then
\begin{equation*}
  \Runbstar\bigl(\cP(\cZ_0)\bigr)
  \geq A^2\left\{\frac1{r(V)}-\frac1{N_V}\right\}.
  \tag{\ref*{eq:rank-unb}}
\end{equation*}
For all estimators, if $1\leq r(V)\leq(1-\epsilon)N_V$ for some
$0<\epsilon<1$, then
\begin{equation*}
  \Rallstar\bigl(\cP(\cZ_0)\bigr)
  \geq \frac{A^2c_\epsilon}{r(V)}
  \tag{\ref*{eq:rank-all}}
\end{equation*}
for a constant $c_\epsilon>0$ depending only on $\epsilon$.
\end{resultrestatement}

\begin{proof}[Proof of Theorem~\ref{cor:rank-lower}]
Here $A=\sum_{k\in V}|a_k|$ is the total absolute coefficient weight
on $V$. Throughout the proof, use the restricted outcome class
\[
 \cY_V(f)=\{Y\in\cY(f):Y_k=0\text{ for every }k\notin V\}.
\]
Since $\cY_V(f)\subseteq\cY(f)$, the worst-case risk of any fixed
design and estimator over $\cY(f)$ is at least its worst-case risk over
$\cY_V(f)$. Thus a lower bound on the minimax risk over $\cY_V(f)$ is
also a lower bound for the original problem. The same conclusion holds
for unbiased estimation: every estimator unbiased on $\cY(f)$ is also
unbiased on $\cY_V(f)$.
If $r(V)=0$, every assignment misses $V$. Apply the lower bounds from
disjoint perturbations in Theorem~\ref{thm:fixed-joint-missingness} with
the one-set packing $\cB=\{V\}$ and signs $s_{1,k}=\operatorname{sign}(a_k)$
for $k\in V$. Then $c_1=\sum_{k\in V}|a_k|=A$ and $\rho_P(V)=1$
for every design $P$. Part (i) of
Theorem~\ref{thm:fixed-joint-missingness} gives the all-estimator lower
bound $A^2$, and part (ii) of the same theorem shows that no unbiased
estimator exists when
$A>0$. In the remainder assume $r(V)\geq1$.
To prove equation~\eqref{eq:rank-unb}, use the singleton packing of \(V\) and take
\(u_b=1\).  Since every assignment reveals at most \(r(V)\) potential outcomes in \(V\), \(\nu_u=r(V)\), while \(|c_b|=A/N_V\).  Theorem
~\ref{thm:global-hypergraph-lower} gives
\[
  \Runbstar\bigl(\cP(\cZ_0)\bigr)
  \geq\frac{(N_V A/N_V)^2}{r(V)}-N_V(A/N_V)^2
  =
  A^2\left\{\frac1{r(V)}-\frac1{N_V}\right\}.
\]
To prove equation~\eqref{eq:rank-all}, write $r=r(V)\geq1$ and choose
$m=\min\{2r,N_V\}$ perturbation sets of approximately equal size.
Since any assignment observes at most $r$ indices in $V$, it can
intersect at most $r$ sets. This leaves at least $m-r$ perturbations
unresolved, while keeping each perturbation's estimand contribution
sufficiently large. Specifically, partition $V$ into $m$ disjoint sets
of sizes $\lfloor N_V/m\rfloor$ and $\lceil N_V/m\rceil$, and align
the perturbation signs with the estimand coefficients.  Each perturbation set then has coefficient magnitude at least
$(A/N_V)\lfloor N_V/m\rfloor\geq A/(2m)$, since $N_V/m\geq1$
and $\lfloor x\rfloor\geq x/2$ for $x\geq1$. Every assignment reveals
at most $r$ potential outcomes in $V$, so it intersects at most $r$ perturbation sets.
At least $m-r$ perturbation sets are therefore unobserved under every assignment.
The general hypergraph lower bound (Theorem~\ref{thm:global-hypergraph-lower}) gives
\[
  \Rallstar\bigl(\cP(\cZ_0);f,a\bigr)
  \geq\frac{A^2(m-r)}{4m^2}.
\]
If $2r\leq N_V$, then $m=2r$ and
\[
  \Rallstar\bigl(\cP(\cZ_0);f,a\bigr)\geq\frac{A^2}{16r}.
\]
If $2r>N_V$, use the singleton partition $m=N_V$, with exact coefficient
magnitude $A/N_V$ in each set, to obtain
\[
  \Rallstar\bigl(\cP(\cZ_0);f,a\bigr)
  \geq\frac{A^2(N_V-r)}{N_V^2}
  \geq\frac{A^2\epsilon}{N_V}
  \geq\frac{A^2\epsilon}{2r}.
\]
Thus $c_\epsilon=\min\{1/16,\epsilon/2\}$ is valid for
$0<\epsilon<1$.
\end{proof}

In some cases, the lower bounds based on each individual arm are not sufficient. The following two subsections obtains stronger lower bounds by
considering both arms together.  For example, for the total treatment effect under neighborhood interference, each
arm can be observed in full under a separate assignment: $z=1^n$ reveals
$Y_i(1_{N_i^+})$ for all units, while $z=0^n$ reveals
$Y_i(0_{N_i^+})$ for all units. Thus both arms have rank $n$, and the
one-arm rank lower bounds in Theorem~\ref{cor:rank-lower} are
uninformative. We first use singleton perturbations
and unbiasedness to obtain bounds based on expected exposure counts (Section \ref{subsec:observed-arm-count-bounds}).
We then allow grouped perturbations to obtain bounds valid for all
estimators. (Section \ref{subsec:block-certificates})

\subsubsection{Lower bounds for unbiased estimators based on expected exposure counts}
\label{subsec:observed-arm-count-bounds}

In this subsection, we derive lower bounds for the joint minimax risk
over all designs and unbiased estimators in terms of the expected
exposure counts achievable by a single design. The one-arm rank bounds from the previous section
maximize each arm's observed exposure count $n_V(z)$ separately over
assignments, giving $r(V)=\max_{z\in\cZ_0}n_V(z)$.
Here we instead consider the expected exposure counts of both arms jointly attainable under a single design. This retains the tradeoff between counts in the two arms. We first consider two arms and then extend the argument to multiple arms.

Consider an equal-weight contrast between two arms. Let
$V_+,V_-\subseteq K_a$ be the positive and negative arms, which
partition $K_a$, and write $N_+=|V_+|$ and $N_-=|V_-|$.
For an assignment $z$, let $n_\pm(z)=|O_f(z)\cap V_\pm|$ denote
the number of units whose potential outcome in the positive or negative
arm, respectively, is revealed. %
Since $E=O_f(z)\cap K_a$, we have
$n_\pm(z)=|E\cap V_\pm|$. Thus define
$r_\pm(E)=(n_+(z),n_-(z))=(|E\cap V_+|,|E\cap V_-|)$,
the pair of observed exposure counts under assignment $z$.
We call $r_\pm(E)$ the \emph{exposure-count vector}. Writing
$\cH_{\rm obs}^a=(K_a,\mathcal E_a)$, define the region of attainable expected exposure counts by
\[
  \mathfrak C_\pm(\cH_{\rm obs}^a)
  =\operatorname{conv}\{r_\pm(E):E\in\mathcal E_a\}.
\]
Under a design $P$, the expected exposure counts in the two arms are
$x_\pm(P)=\E_P[n_\pm(Z)]$. The pair $(x_+(P),x_-(P))$ is an average
of the feasible exposure-count vectors, with probabilities supplied by $P$,
and therefore belongs to $\mathfrak C_\pm$. Conversely, choose a point
$x=\sum_{j=1}^J\lambda_j r_\pm(E_j)\in\mathfrak C_\pm$, where
$E_j\in\mathcal E_a$, $\lambda_j\geq0$, and $\sum_{j=1}^J\lambda_j=1$.
For each $E_j$, choose an assignment $z_j$ satisfying
$O_f(z_j)\cap K_a=E_j$. A design that selects $z_j$ with probability
$\lambda_j$ has expected exposure-count vector
$\sum_{j=1}^J\lambda_j r_\pm(E_j)=x$.
Hence $\mathfrak C_\pm$ is exactly the set of expected exposure-count
pairs achievable by a design. Define
\begin{equation}
  \Gamma_\pm(\cH_{\rm obs}^a;V_+,V_-)
  =
  \min_{(x_+,x_-)\in\mathfrak C_\pm(\cH_{\rm obs}^a)}
  \left\{
    \frac1{x_+}+\frac1{x_-}
    -\frac1{N_+}-\frac1{N_-}
  \right\},
  \label{eq:two-arm-coverage-functional}
\end{equation}
The objective inside the minimum is assigned value $+\infty$ when
$x_+=0$ or $x_-=0$.

The lower bound involves the reciprocals of the expected exposure counts
in both arms, with the finite-arm corrections $1/N_+$ and $1/N_-$.
Optimizing over the attainable region allows the design to balance these
counts, subject to the constraints imposed by feasible assignments.

\begin{corollary}[Two-arm coverage]
\label{cor:two-arm-coverage}
\label{thm:equal-weight-covering}
Let
\[
  \tau
  =
  \frac1{N_+}\sum_{k\in V_+}Y_k
  -
  \frac1{N_-}\sum_{k\in V_-}Y_k
\]
be an equal-weight two-arm contrast.  Then
\begin{equation}
  \Runbstar\bigl(\cP(\cZ_0)\bigr)
  \geq
  \Gamma_\pm(\cH_{\rm obs}^a;V_+,V_-).
  \label{eq:two-arm-coverage-lower}
\end{equation}
\end{corollary}

\begin{proof}
Fix a design $P$ and let $q_k=P\{k\in O_f(Z)\}$. If $q_k=0$ for any
$k\in K_a$, no unbiased estimator exists and the claimed lower bound holds.
Otherwise, apply the fixed-design lower bound from disjoint perturbations
(Theorem~\ref{thm:fixed-joint-missingness}(ii)) to the singleton packing
of $V_+\cup V_-$, using aligned signs. This gives
\begin{equation}
  \Runb(P)\geq
  \frac1{N_+^2}\sum_{k\in V_+}\left(\frac1{q_k}-1\right)
  +\frac1{N_-^2}\sum_{k\in V_-}\left(\frac1{q_k}-1\right).
  \label{eq:two-arm-proof-fixed}
\end{equation}
The expected exposure counts in the two arms are
\begin{equation}
  x_\pm(P)=\sum_{k\in V_\pm}q_k
       =\E_P[n_\pm(Z)].
  \label{eq:two-arm-proof-counts}
\end{equation}
Apply Lemma~\ref{lem:weighted-reciprocal} separately to each arm, with
$w_k=1$ and $p_k=q_k$. Since there are $N_\pm$ indices in the respective
arm and their detection probabilities sum to $x_\pm(P)$, the lemma gives
\begin{equation}
  \sum_{k\in V_\pm}\frac1{q_k}\geq\frac{N_\pm^2}{x_\pm(P)}.
  \label{eq:two-arm-proof-reciprocal}
\end{equation}
Substituting equation~\eqref{eq:two-arm-proof-reciprocal} into
 equation~\eqref{eq:two-arm-proof-fixed} yields
\begin{equation}
  \Runb(P)\geq\frac1{x_+(P)}+\frac1{x_-(P)}-\frac1{N_+}-\frac1{N_-}.
  \label{eq:two-arm-proof-count-bound}
\end{equation}
The vector $(x_+(P),x_-(P))$ belongs to $\mathfrak C_\pm$, so the right-hand
side of equation~\eqref{eq:two-arm-proof-count-bound} is at least
$\Gamma_\pm$. Taking the infimum over all designs proves
 equation~\eqref{eq:two-arm-coverage-lower}.
\end{proof}

\paragraph{Estimands with more than two arms}
The two-arm lower bound above extends to an estimand with any finite
number of arms, with equal absolute coefficients within each arm.
Corollary~\ref{cor:multi-arm-coverage} gives this extension in terms of
the expected exposure counts.

\begin{corollary}[Multi-arm coverage]
\label{cor:multi-arm-coverage}
Let the estimand be a linear contrast whose nonempty arms
\(V_1,\dots,V_L\) form a disjoint partition of its index set and each carry
equal within-arm weight \(\lvert a_k\rvert=A_\ell/N_\ell\) for
\(k\in V_\ell\), where \(N_\ell=|V_\ell|\).
Write $\cH_{\rm obs}^a=(K_a,\mathcal E_a)$. For an assignment $z$,
define $n_\ell(z)=|O_f(z)\cap V_\ell|$, the observed exposure count in arm
$\ell$. For $E=O_f(z)\cap K_a$, define the exposure-count vector
and the region of attainable expected exposure counts by
\[
 \begin{aligned}
  r(E)&=\bigl(n_1(z),\ldots,n_L(z)\bigr)
       =\bigl(|E\cap V_1|,\ldots,|E\cap V_L|\bigr),\\
  \mathfrak C_{\rm arm}(\cH_{\rm obs}^a)
  &=\operatorname{conv}\{r(E):E\in\mathcal E_a\}.
 \end{aligned}
\]
Then
\begin{equation}
  \Runbstar\bigl(\cP(\cZ_0);f,a\bigr)
  \geq
  \min_{x\in\mathfrak C_{\rm arm}(\cH_{\rm obs}^a)}
  \sum_{\ell=1}^L\frac{A_\ell^2}{x_\ell}
  -\sum_{\ell=1}^L\frac{A_\ell^2}{N_\ell},
  \label{eq:multi-arm-coverage-lower}
\end{equation}
Each summand in the objective is assigned value $+\infty$ if
$A_\ell>0$ and $x_\ell=0$, and value zero if $A_\ell=0$.
\end{corollary}

\begin{proof}
Fix a design $P$, and let
\[
  I=\{\ell\in\{1,\ldots,L\}:A_\ell>0\}
\]
index the arms with nonzero estimand weight. Define
$q_k=P\{k\in O_f(Z)\}$ and
$x_\ell(P)=\sum_{k\in V_\ell}q_k=\E_P[n_\ell(Z)]$.
If $q_k=0$ for any $k\in\bigcup_{\ell\in I}V_\ell$, no unbiased
estimator exists and the claimed lower bound holds. Otherwise, the
fixed-design lower bound from disjoint perturbations
(Theorem~\ref{thm:fixed-joint-missingness}(ii)), applied to the singleton
packing of $\bigcup_{\ell\in I}V_\ell$, gives
\begin{equation}
  \Runb(P)\geq\sum_{\ell\in I}\frac{A_\ell^2}{N_\ell^2}
              \left(\sum_{k\in V_\ell}\frac1{q_k}-N_\ell\right).
  \label{eq:multi-arm-proof-fixed}
\end{equation}
For each $\ell\in I$, apply Lemma~\ref{lem:weighted-reciprocal} with
$w_k=1$ and $p_k=q_k$ on $V_\ell$. This gives
\begin{equation}
  \sum_{k\in V_\ell}\frac1{q_k}\geq\frac{N_\ell^2}{x_\ell(P)}.
  \label{eq:multi-arm-proof-reciprocal}
\end{equation}
Substituting equation~\eqref{eq:multi-arm-proof-reciprocal} into
 equation~\eqref{eq:multi-arm-proof-fixed} yields
\begin{equation}
  \Runb(P)\geq\sum_{\ell\in I}\frac{A_\ell^2}{x_\ell(P)}
                    -\sum_{\ell\in I}\frac{A_\ell^2}{N_\ell}.
  \label{eq:multi-arm-proof-count-bound}
\end{equation}
Arms outside $I$ have $A_\ell=0$, so under the stated convention the
right-hand side equals the objective in
 equation~\eqref{eq:multi-arm-coverage-lower}. The vector $x(P)=(x_1(P),\ldots,x_L(P))$ is an average
of feasible exposure-count vectors and therefore belongs to
$\mathfrak C_{\rm arm}$. Consequently, the right-hand side of
 equation~\eqref{eq:multi-arm-proof-count-bound} is at least the minimum
in equation~\eqref{eq:multi-arm-coverage-lower}. Taking the infimum over
all designs proves the claim.
\end{proof}

\subsubsection{Lower bounds for all estimators based on paired perturbation sets}
\label{subsec:block-certificates}

The exposure-count lower bounds of
Section~\ref{subsec:observed-arm-count-bounds} need not capture the full
difficulty of estimation. Consider again the total treatment effect,
now on a complete graph. Each arm has sampling rank $n$, so the
one-arm rank lower bounds in Theorem~\ref{cor:rank-lower} give a trivial lower bound of zero. On the other hand, Corollary~\ref{cor:two-arm-coverage} gives the lower bound $2/n$ for
unbiased estimators: every assignment reveals either an entire arm or
neither arm, so the expected exposure counts satisfy $x_++x_-\leq n$.
The reciprocal-count expression is minimized at $x_+=x_-=n/2$,
achieved by choosing the all-treated or all-control assignment with
equal probability, and its value is $2/n$.

However, no assignment reveals outcomes from both TTE arms. Put a prior
on potential-outcome tables by setting every treatment-arm outcome equal
to one uniform random sign and every control-arm outcome equal to another
independent uniform random sign, with all other potential outcomes fixed
at zero. Every assignment leaves at least one of these signs unobserved.
The resulting posterior variance of the TTE is therefore at least $1$,
giving a lower bound of $1$ over all designs and estimators.

Theorem~\ref{thm:arm-paired-certificate} below generalizes this
construction to paired groups from the two arms and yields the stronger
lower bound in this example. Specifically, we group
potential outcomes within each arm and pair the groups across
arms, assigning an independent unknown sign to each group. Unlike the
singleton construction used for the exposure-count lower bounds,
each group may contain several potential outcomes. An assignment
reveals both signs in a pair only if it observes at least one outcome
in each group. Otherwise, at least one sign remains unknown and
contributes to the posterior variance of the estimand.

Bounding how many pairs a single assignment can observe on both sides
gives a lower bound for the joint minimax risk over all designs and
estimators. We obtain this bound by applying the general hypergraph
lower bound (Theorem~\ref{thm:global-hypergraph-lower}) to the disjoint
sets forming the pairs. Theorem~\ref{thm:arm-paired-certificate} allows
unequal weights across pairs; its equal-weight specialization gives
Theorem~\ref{thm:paired-sets-main}, used in the neighborhood and
bipartite examples. The same bounds hold for unbiased estimators because
they form a smaller estimator class.

\begin{definition}[Arm-paired packing and maximum jointly observed pair count]
\label{def:arm-paired-packing}
An \emph{arm-paired packing} is a perturbation packing
$\cB=\{B_1,\ldots,B_m\}$ in which every set carries a splitting
$B_b=B_b^+\cup B_b^-$ into two disjoint nonempty parts, one drawn from each
arm of the contrast, together with perturbation signs $s_{b,k}$ for
which the two parts carry equal weight,
\begin{equation}
  w_b
  =\Bigl|\sum_{k\in B_b^+}s_{b,k}a_k\Bigr|
  =\Bigl|\sum_{k\in B_b^-}s_{b,k}a_k\Bigr| .
  \label{eq:arm-paired-weight}
\end{equation}
The two parts carry separate unknown signs under the hypercube prior.
A hyperedge $E\in\cH_{\rm obs}^a$ reveals both signs, and hence
\emph{resolves} the pair $B_b$, when it meets both parts, $E\cap B_b^+\neq\varnothing$ and $E\cap B_b^-\neq\varnothing$.
The \emph{maximum jointly observed pair count} of the packing is the largest number of pairs
any one feasible hyperedge resolves,
\begin{equation}
  \mathcal M(\cB)
  =
  \max_{E\in\cH_{\rm obs}^a}
  \#\{b\leq m:\ E\text{ resolves }B_b\}.
  \label{eq:resolution-capacity}
\end{equation}
\end{definition}

\begin{theorem}[Arm-paired set certificate]
\label{thm:arm-paired-certificate}
Let $\cB=\{B_1,\ldots,B_m\}$ be an arm-paired packing with per-arm
weights $w_b$ and maximum jointly observed pair count $\mathcal M=\mathcal M(\cB)$.
Then the joint minimax risks satisfy
\begin{equation}
  \Rallstar\bigl(\cP(\cZ_0)\bigr)
  \ \geq\
  \sum_{b=1}^m w_b^2
  \ -\
  \max_{S\subseteq[m],\,|S|\leq\mathcal M}\ \sum_{b\in S}w_b^2 ,
  \label{eq:arm-paired-lower}
\end{equation}
that is, the total squared per-arm weight less its $\mathcal M$ largest terms.
The same bound holds for $\Runbstar\bigl(\cP(\cZ_0)\bigr)$.  In particular:
\begin{enumerate}
\item[\textup{(i)}] if no feasible hyperedge resolves any set, then
$\Rallstar\bigl(\cP(\cZ_0)\bigr)\geq\sum_bw_b^2$;
\item[\textup{(ii)}] if the per-arm weights are equal, $w_b=A/m$ for every $b$, then
\begin{equation}
  \Rallstar\bigl(\cP(\cZ_0)\bigr)
  \ \geq\ A^2\,\frac{m-\mathcal M}{m^2}.
  \label{eq:arm-paired-equal}
\end{equation}
Suppose in addition that arm-paired packings with equal per-arm weights exist for every
$m\leq N$ and that their maximum jointly observed pair counts obey a common bound
$\mathcal M(\cB_m)\leq\mathcal M$.  Then optimizing $m$ gives $A^2/(4\mathcal M)$ when
$1\leq\mathcal M\leq N/2$, and $A^2$ when $\mathcal M=0$.
\end{enumerate}
\end{theorem}

\begin{proof}[Proof of Theorem~\ref{thm:arm-paired-certificate}]
Fix an arm-paired packing $\cB$. We apply the general hypergraph lower
bound to its $2m$ disjoint parts, then use the maximum jointly observed pair count to
bound the squared effect missed by every feasible hyperedge.

The collection
\[
  \widetilde{\cB}=\{B_b^+,B_b^-:b=1,\ldots,m\}
\]
is a perturbation packing. By equation~\eqref{eq:arm-paired-weight},
each part of pair $b$ has squared set effect $w_b^2$. For a feasible
hyperedge $E\in\cH_{\rm obs}^a$, let
\[
  R(E)=\{b\in\{1,\ldots,m\}:E\text{ resolves }B_b\}.
\]
If $b\notin R(E)$, the hyperedge misses at least one of $B_b^+$ and
$B_b^-$. The total squared effect of the missed parts therefore satisfies
\begin{equation}
 \begin{aligned}
  &\sum_{b=1}^m w_b^2\left(
    \1\{E\cap B_b^+=\varnothing\}
   +\1\{E\cap B_b^-=\varnothing\}\right)\\
  &\qquad\geq\sum_{b\notin R(E)}w_b^2
   =\sum_{b=1}^m w_b^2-\sum_{b\in R(E)}w_b^2\\
  &\qquad\geq\sum_{b=1}^m w_b^2
   -\max_{S\subseteq[m],\,|S|\leq\mathcal M}\sum_{b\in S}w_b^2.
 \end{aligned}
 \label{eq:arm-paired-proof-missing}
\end{equation}
The last inequality uses $|R(E)|\leq\mathcal M$, by the definition of
maximum jointly observed pair count in equation~\eqref{eq:resolution-capacity}.

Equation~\eqref{eq:arm-paired-proof-missing} holds for every feasible
hyperedge. Taking the minimum over hyperedges bounds the minimum missing
squared effect for $\widetilde{\cB}$ by the same expression. Applying
Theorem~\ref{thm:global-hypergraph-lower} to this packing gives
\begin{equation}
  \Rallstar\bigl(\cP(\cZ_0)\bigr)
  \geq\sum_{b=1}^m w_b^2
   -\max_{S\subseteq[m],\,|S|\leq\mathcal M}\sum_{b\in S}w_b^2.
  \label{eq:arm-paired-proof-risk}
\end{equation}
This proves equation~\eqref{eq:arm-paired-lower}. The same bound holds
for $\Runbstar$, since restricting to unbiased estimators cannot reduce
the minimax risk.

If $\mathcal M=0$, only the empty set is allowed in the maximum in
equation~\eqref{eq:arm-paired-lower}.  Its weight is zero, giving part (i).
For part (ii), suppose $w_b=A/m$ for every $b$.
Since $0\leq\mathcal M\leq m$, the maximum is
$\mathcal M(A/m)^2$, and equation~\eqref{eq:arm-paired-lower} gives
\begin{equation}
 \Rallstar\bigl(\cP(\cZ_0)\bigr)
 \geq m(A/m)^2-\mathcal M(A/m)^2
 =A^2\frac{m-\mathcal M}{m^2}.
  \label{eq:arm-paired-proof-equal}
\end{equation}
Now suppose equal-weight packings exist for every $m\leq N$, with
$\mathcal M(\cB_m)\leq\mathcal M$.
For $1\leq\mathcal M\leq N/2$, the choice $m=2\mathcal M$ is available.
Using the common capacity bound in
equation~\eqref{eq:arm-paired-proof-equal} gives
\begin{equation}
 \Rallstar\bigl(\cP(\cZ_0)\bigr)
 \geq A^2\frac{2\mathcal M-\mathcal M}{(2\mathcal M)^2}
 =\frac{A^2}{4\mathcal M}.
  \label{eq:arm-paired-proof-optimized}
\end{equation}
If $\mathcal M=0$, the choice $m=1$ gives $A^2$.
\end{proof}

For the next restatement, let $B_b^+$ and $B_b^-$, $b=1,\ldots,m$,
be nonempty, pairwise disjoint sets of potential outcomes from
the two respective arms, and let $\mathcal M$ be the largest number of
pairs for which one feasible observed set meets both parts, as in
equation~\eqref{eq:resolution-capacity}.

\begin{resultrestatement}{Theorem~\ref{thm:paired-sets-main}}{Lower bound from paired perturbation sets}
For the arm-paired packing $\{B_b^+,B_b^-:b\leq m\}$, suppose
that each of the $2m$ parts carries the same
estimand weight,
\[
 \sum_{k\in B_b^+}|a_k|=\sum_{k\in B_b^-}|a_k|=\frac{A}{m},
 \qquad b=1,\ldots,m.
\]
Then
\[
 \Rallstar\bigl(\cP(\cZ_0)\bigr)\geq A^2\frac{m-\mathcal M}{m^2}.
\]
The same bound holds for $\Runbstar\bigl(\cP(\cZ_0)\bigr)$.
\end{resultrestatement}

\begin{proof}[Proof of Theorem~\ref{thm:paired-sets-main}]
The pairs in the theorem form an arm-paired packing.  Choose the fixed
signs to align with the estimand coefficients on each part.  Then
\[
 w_b=\sum_{k\in B_b^+}|a_k|
    =\sum_{k\in B_b^-}|a_k|=\frac{A}{m}.
\]
Its maximum jointly observed pair count is the $\mathcal M$ in the theorem.
Part (ii) of Theorem~\ref{thm:arm-paired-certificate} therefore gives
\[
 \Rallstar\bigl(\cP(\cZ_0)\bigr)
 \geq A^2\frac{m-\mathcal M}{m^2},
\]
and the same bound for $\Runbstar\bigl(\cP(\cZ_0)\bigr)$.
\end{proof}

Theorem~\ref{thm:arm-paired-certificate} is applied in
Appendices~\ref{app:neighborhood-proofs} and~\ref{app:bipartite-proofs}.
When no assignment resolves any pair, part (i) gives the sum of squared
per-arm weights as a lower bound, as in the packed-neighborhood bound
for the total effect (Lemma~\ref{lem:packed-neighborhood-lower}).
When at most $\mathcal M$ pairs can be resolved, equal-weight packings
and part (ii) give a reciprocal-capacity lower bound under the stated
packing-size condition. This is used in the block-resolution bound for
the total effect (Lemma~\ref{lem:block-resolution}) and in the bipartite
dose example.

\subsection{Comparison with conflict-graph lower bounds}
\label{sec:pairwise}

In this subsection, we compare the lower bounds based on the sampling hypergraph with
lower bounds based on the conflict graph of \citet{KandirosHarshawSavje2026} for the all-design minimax problem. We first show that the independence number of the conflict graph is an upper bound of the sampling rank, so replacing the rank by the independence number can weaken the resulting lower bound. We then show that under the setting of neighborhood interference with the identity exposure
mapping and unrestricted treatment assignments, the conflict graph
and sampling hypergraph determine one another. We also prove the corresponding statement for the restriction of the hypergraph and the conflict graph to estimand-relevant potential outcomes. As a consequence, in this specific setting, the independence number of the conflict graph is equal to the sampling rank, and the conflict graph contains all the information needed for the all-design minimax problem. Finally, we consider cases when assignment restrictions or coarsened exposure mappings break this equivalence, and the sampling hypergraph strictly contains more information than the conflict graph.

\subsubsection{Pairwise observability and sampling rank}
\label{subsec:pairwise-relaxation}
Fix an exposure mapping $f$ and a nonempty set of feasible treatment assignments $\cZ_0$. Recall that the sampling hypergraph is
$\cH_{\rm ops}=(K,\mathcal E)$, where $K$ is the set of potential
outcomes and $\mathcal E=\{O_f(z):z\in\cZ_0\}$ is the collection of hyperedges recording the maximal sets of potential outcomes that can be jointly observed. For a given assignment $z\in\cZ_0$, write $E=O_f(z)$
for the corresponding hyperedge. Fix a nonempty set $T\subseteq K$ of
potential outcomes, each individually attainable under $\cZ_0$.
In the sampling hypergraph,
a collection of nodes is contained in a hyperedge if and only if the
corresponding potential outcomes can be jointly observed under a single
treatment assignment. The sampling rank of $T$ is defined as
\begin{equation}
  r(T)=\max_{E\in\mathcal E}|E\cap T|.
  \label{eq:pairwise-rank}
\end{equation}
$r(T)$ is the largest number of potential outcomes in $T$ that one feasible
assignment reveals.

\begin{definition}[Pairwise conflicts and conflict graph]
\label{def:conflict-graph}
\label{def:pairwise-conflict}
 For an exposure mapping $f$ and nonempty feasible assignment set $\cZ_0$, two distinct potential outcomes $Y_k,Y_\ell$ have a \emph{pairwise conflict} if $I_k(z)I_\ell(z)=0$ for every $z\in\cZ_0$: at least one outcome in the pair must remain missing under every feasible assignment.
The conflict graph $H_{\rm conf}$ has vertex set $K$. Distinct vertices $k,\ell\in K$ are adjacent in the conflict graph if and only if their potential outcomes have a pairwise conflict.
For $T\subseteq K$, write $H_{\rm conf}[T]$ for the induced subgraph on
$T$, retaining exactly the conflict edges with both endpoints in $T$.
\end{definition}

Recall from Section~\ref{subsec:packings-and-the-three-problems} that a collection $S\subseteq T$ is \emph{pairwise observable} if each pair of distinct potential outcomes in $S$ can be observed under some feasible assignment; the assignment may differ across pairs. The collection is \emph{jointly observable} if one feasible assignment observes every potential outcome in $S$.
To express this distinction in terms of assignments, recall the revelation sets
\begin{equation}
  \cZ_k=\{z\in\cZ_0:I_k(z)=1\}
  =\{z\in\cZ_0:k\in O_f(z)\},\qquad k\in K.
  \label{eq:revelation-sets}
\end{equation}
Thus $S$ is pairwise observable if and only if $\cZ_k\cap\cZ_\ell\ne\varnothing$ for every distinct $k,\ell\in S$, equivalently $S$ is an independent set in $H_{\rm conf}[T]$. In contrast, $S$ is jointly observable if and only if $\bigcap_{k\in S}\cZ_k\ne\varnothing$, equivalently $S$ is contained in a hyperedge of the sampling hypergraph. Pairwise intersections need not imply a common intersection, so pairwise observability of $S$ need not imply joint observability.

However, every jointly observable collection $S$ is pairwise observable. In particular, the set of potential outcomes $E\cap T$ is an independent set in $H_{\rm conf}[T]$ for every hyperedge $E\in\mathcal E$, giving
\begin{equation}
  r(T)\leq\alpha\{H_{\rm conf}[T]\}.
  \label{eq:pairwise-inequality}
\end{equation}
Here $r(T)$ is the largest size of a jointly observable collection in $T$, whereas $\alpha\{H_{\rm conf}[T]\}$ is the largest size of a pairwise observable collection in $T$. When $T$ lies in one estimand arm, the rank-based lower bounds of Theorem~\ref{cor:rank-lower} have $r(T)$ in the denominator. Replacing $r(T)$ by the potentially larger independence number can therefore weaken these lower bounds.

\begingroup\sloppy
\subsubsection{Equivalence of the conflict graph and sampling hypergraph under identity neighborhood exposures}
\label{subsec:pairwise-exact}
\label{subsec:conflict-comparison}
\par\endgroup

We now show that under identity neighborhood exposures and unrestricted treatment assignments, pairwise observability implies joint observability. Consequently, the rank inequality in \eqref{eq:pairwise-inequality} is an equality. In particular, the next theorem shows that under neighborhood interference with identity exposure mapping, and unrestricted treatment assignments, the hyperedges of the sampling hypergraph with the maximal independent sets in $H_{\rm conf}$, showing that the two representations determine one another.

\begin{theorem}[Equivalence of the conflict graph and sampling hypergraph]
\label{thm:conflict-hypergraph-equivalence}
\label{prop:subcube-collapse}
Suppose $\cZ_0=\{0,1\}^n$ and, for every $k=(i,e)\in K$,
observing $Y_i(e)$ requires precisely the neighborhood configuration
$z_{N_i^+}=e$, with treatment assignments outside $N_i^+$ unrestricted.
Thus
\[
 \cZ_{(i,e)}=\{z\in\{0,1\}^n:z_{N_i^+}=e\}.
\]
Then the hyperedge set $\mathcal E$ is exactly the family of maximal
independent sets of the conflict graph $H_{\rm conf}$:
\begin{equation}
 \mathcal E
 =\{S\subseteq K:S\text{ is maximal independent in }H_{\rm conf}\}.
 \label{eq:conflict-hypergraph-characterization}
\end{equation}
Conversely, two distinct vertices are adjacent in $H_{\rm conf}$ if and
only if no hyperedge contains both. Thus the conflict graph and
sampling hypergraph determine one another.
\end{theorem}

\begin{proof}
First, note that the sampling hypergraph determines the conflict graph: two distinct vertices are adjacent in the conflict graph whenever no hyperedge contains both vertices.

Next, we prove that every independent set in $H_{\rm conf}$ is jointly observable. In terms of the revelation sets, the key implication is that, for every $S\subseteq K$,
\[
 \bigl[\cZ_k\cap\cZ_\ell\ne\varnothing
 \text{ for every distinct }k,\ell\in S\bigr]
 \quad\Longrightarrow\quad
 \bigcap_{k\in S}\cZ_k\ne\varnothing.
\]

To prove this implication, suppose $\cZ_k\cap\cZ_\ell\ne\varnothing$ for every distinct $k,\ell\in S$. For $k=(i,e)$ and $\ell=(j,e')$, an assignment in $\cZ_k\cap\cZ_\ell$ satisfies both neighborhood configurations. Therefore,
\[
 e_v=e'_v\qquad\text{for every }v\in N_i^+\cap N_j^+.
\]
Define the assignment $z^S=(z^S_1,\ldots,z^S_n)\in\{0,1\}^n$ by
\[
 z^S_v=\begin{cases}
 e_v,&\text{if }v\in N_i^+\text{ for some }(i,e)\in S,\\
 0,&\text{otherwise}.
 \end{cases}
\]
The agreement above ensures that this assignment is well-defined. By construction, $z^S_{N_i^+}=e$ for every $(i,e)\in S$. Since all binary treatment assignments are feasible,
\[
 z^S\in\bigcap_{k\in S}\cZ_k.
\]
Hence the common intersection is nonempty. Thus every independent set in $H_{\rm conf}$ is jointly observable. Conversely, every jointly observable collection is an independent set in $H_{\rm conf}$, because an assignment observing the entire collection also observes each pair.

We now show that the hyperedges are exactly the maximal independent sets in $H_{\rm conf}$. Let $S$ be a maximal independent set in $H_{\rm conf}$. By the preceding argument, there exists an assignment $z$ such that $S\subseteq O_f(z)$. The set $O_f(z)$ is an independent set in $H_{\rm conf}$, because all its potential outcomes are observed together under $z$. Maximality of $S$ therefore gives $S=O_f(z)$, so $S$ is a hyperedge.

Conversely, let $E=O_f(z)$ be a hyperedge. It is an independent set in $H_{\rm conf}$ because its potential outcomes are jointly observable. To show maximality, take any vertex $(i,e')\notin E$. Since $E$ contains exactly one potential outcome for each unit, it contains $(i,e)$ for some $e\ne e'$. The two distinct potential outcomes $Y_i(e)$ and $Y_i(e')$ cannot be observed together, so their vertices are adjacent in $H_{\rm conf}$. Consequently, no vertex outside $E$ can be added while preserving the property of being an independent set in $H_{\rm conf}$. Hence $E$ is a maximal independent set in $H_{\rm conf}$, proving \eqref{eq:conflict-hypergraph-characterization}.

\end{proof}

 As a consequence of the above theorem, the rank inequality in \eqref{eq:pairwise-inequality} is an equality. Indeed, for any
nonempty $T\subseteq K$, taking $S$ to be a maximum independent set
of $H_{\rm conf}[T]$ gives
$r(T)\geq|S|=\alpha\{H_{\rm conf}[T]\}$. Combining this inequality with
equation~\eqref{eq:pairwise-inequality} gives
$r(T)=\alpha\{H_{\rm conf}[T]\}$.

\paragraph*{Restriction to estimand-relevant potential outcomes}
Theorem~\ref{thm:conflict-hypergraph-equivalence} establishes equivalence of the conflict graph and the sampling hypergraph defined on $K$ under neighborhood interference with the identity exposure mapping and unrestricted treatment assignments. We next use this result to show that the induced conflict graph on the estimand-relevant potential outcomes  $K_a=\{k\in K:a_k\neq0\}$ retains enough information to determine the all-design minimax risks over all estimators and over unbiased estimators. The relevant graph is the induced subgraph $H_{\rm conf}[K_a]$, and the corresponding hypergraph is the estimand trace $\cH_{\rm obs}^a=(K_a,\mathcal E_a)$, where $\mathcal E_a=\{O_f(z)\cap K_a:z\in\cZ_0\}$. Unlike hyperedges of $\cH_{\rm obs}$, trace hyperedges can have different sizes and can strictly contain one another.

Under the assumptions of Theorem~\ref{thm:conflict-hypergraph-equivalence}, the maximal trace hyperedges of $\cH_{\rm obs}^a$ are exactly the maximal independent sets of $H_{\rm conf}[K_a]$. The next theorem shows, for any exposure mapping, that maximal trace hyperedges suffice for all-design minimax for both all-estimator and unbiased estimator classes. Corollary~\ref{cor:induced-conflict-minimax} then combines this result with Theorem~\ref{thm:conflict-hypergraph-equivalence} to establish the induced-conflict-graph consequence.

\begin{theorem}[Maximal estimand trace hyperedges suffice for all-design minimax risk]
\label{thm:trace-conflict-minimax}
Fix an exposure mapping $f$, a nonempty feasible assignment set $\cZ_0$, estimand coefficients $a$, and the bounded-table class $\cY(f)=[-1,1]^K$. Let $\mathcal E_a^{\max}$ be the inclusion-maximal members of $\mathcal E_a$, and define
\[
 \mathcal C_{\max}
 =\left\{P\in\cP(\cZ_0):
 P\bigl(O_f(Z)\cap K_a\in\mathcal E_a^{\max}\bigr)=1\right\}.
\]
Then, for both estimator classes,
\begin{equation}
 R^*_{\eta}\bigl(\cP(\cZ_0);f,a\bigr)
 =R^*_{\eta}(\mathcal C_{\max};f,a),
 \qquad \eta\in\{\mathrm{all},\mathrm{unb}\}.
 \label{eq:trace-maximal-minimax}
\end{equation}
\end{theorem}

\begin{proof}
We show that any design--estimator pair $(P,\hat\tau)$ can be replaced by a pair whose design belongs to $\mathcal C_{\max}$, without increasing worst-case risk or losing unbiasedness. Write
\[
 r_Y(P,\hat\tau)=\E_P\bigl[(\hat\tau(Z,Y_{O_f(Z)})-\tau_a(Y))^2\bigr]
\]
for its squared-error risk at a fixed table $Y$.

First, we remove the dependence on outcomes outside the estimand revelant potential outcomes $K_a$. Given $Y$, define $Y^0$ by $Y^0_k=Y_k$ for $k\in K_a$ and $Y^0_k=0$ otherwise. Then $Y^0\in[-1,1]^K$ and $\tau_a(Y^0)=\tau_a(Y)$. Define
\[
 \hat\tau_0(z,y_{O_f(z)})=\hat\tau(z,y^0_{O_f(z)}),
\]
where $y^0_{O_f(z)}$ is obtained by replacing the observed coordinates outside $K_a$ by zero. Consequently,
\[
 r_Y(P,\hat\tau_0)=r_{Y^0}(P,\hat\tau),\qquad
 \sup_Y r_Y(P,\hat\tau_0)\leq\sup_Y r_Y(P,\hat\tau).
\]
If $\hat\tau$ is unbiased at every admissible table, then
$\E_P[\hat\tau_0]=\tau_a(Y^0)=\tau_a(Y)$, so $\hat\tau_0$ is also unbiased.

Next, for each $E\in\mathcal E_a$, let $q(E)=P\{O_f(Z)\cap K_a=E\}$ be the probability that the observed estimand-relevant indices are exactly $E$. For $q(E)>0$, define an estimator depending only on $(E,Y_E)$ by
\[
 \hat\tau_1(E,y_E)
 =\sum_{z:O_f(z)\cap K_a=E}\frac{P(z)}{q(E)}
   \hat\tau_0(z,\widetilde y^{\,z}_{O_f(z)}),
\]
where $\widetilde y^{\,z}_k=y_k$ for $k\in E$ and $\widetilde y^{\,z}_k=0$ for $k\in O_f(z)\setminus K_a$. For every $E\in\mathcal E_a$ with $q(E)=0$, define $\hat\tau_1(E,y_E)=0$ for all $y_E$. Such traces are feasible but have probability zero under $P$, so this choice affects neither risk nor unbiasedness. For $q(E)>0$, the weights in the displayed sum are the conditional assignment probabilities given $O_f(Z)\cap K_a=E$ and do not depend on the outcome table. Thus, for every fixed $Y$, conditional averaging gives
\begin{align*}
 \E_q[\hat\tau_1(E,Y_E)]&=\E_P[\hat\tau_0(Z,Y_{O_f(Z)})],\\
 \E_q[(\hat\tau_1(E,Y_E)-\tau_a(Y))^2]
 &\leq r_Y(P,\hat\tau_0).
\end{align*}
This preserves unbiasedness and cannot increase worst-case risk.

Finally, for each $E\in\mathcal E_a$, choose a maximal trace hyperedge $M(E)\supseteq E$. Put
\[
 q'(M)=\sum_{E:M(E)=M}q(E),\qquad M\in\mathcal E_a^{\max}.
\]
Choose one assignment with trace $M$ for each $M$ and assign it probability $q'(M)$. This defines a design $P'\in\mathcal C_{\max}$. After observing $Y_M$, use
\[
 \hat\tau_2(M,Y_M)
 =\sum_{E:M(E)=M}\frac{q(E)}{q'(M)}\hat\tau_1(E,Y_E),
 \qquad q'(M)>0.
\]
Every term is available because $E\subseteq M$. Summing over $M$ and applying Jensen's inequality give, for every $Y$,
\begin{align*}
 \E_{q'}[\hat\tau_2(M,Y_M)]&=\E_q[\hat\tau_1(E,Y_E)],\\
 \E_{q'}[(\hat\tau_2(M,Y_M)-\tau_a(Y))^2]
 &\leq\E_q[(\hat\tau_1(E,Y_E)-\tau_a(Y))^2].
\end{align*}
Hence the pair $(P',\hat\tau_2)$ has no larger worst-case risk than $(P,\hat\tau)$ and is unbiased whenever $\hat\tau$ is unbiased. Estimators on events of probability zero may be defined arbitrarily. Taking infima proves that restricting designs to $\mathcal C_{\max}$ cannot increase either minimax risk. The reverse inequality follows from $\mathcal C_{\max}\subseteq\cP(\cZ_0)$, proving \eqref{eq:trace-maximal-minimax}.
\end{proof}

\begin{proof}
We show that any design--estimator pair $(P,\hat\tau)$ can be replaced by a pair whose design belongs to $\mathcal C_{\max}$, without increasing worst-case risk or losing unbiasedness.

Fix a design $P$ and an estimator $\hat\tau(Z,Y_{O_f(Z)})$. For a table $Y$, let $Y^0$ agree with $Y$ on $K_a$ and be zero elsewhere. Then $Y^0\in\cY(f)$ and $\tau_a(Y^0)=\tau_a(Y)$. Replacing the observed outcomes outside $K_a$ by zero therefore gives the original estimator's risk at $Y^0$, so it cannot increase worst-case risk and preserves unbiasedness.

For $E\in\mathcal E_a$, put $q(E)=P\{O_f(Z)\cap K_a=E\}$, the probability under $P$ of observing exactly the estimand-relevant indices in $E$. For each fixed table $Y$ and each $E$ with $q(E)>0$, define
\[
 \hat\tau_1(E,Y_E)
 :=\E_P\!\left[
 \hat\tau(Z,Y^0_{O_f(Z)})
 \,\middle|\,O_f(Z)\cap K_a=E
 \right].
\]
This averages over assignments producing the same trace $E$. Since the outcomes outside $K_a$ have been replaced by zero, the conditional expectation depends on $Y$ only through $Y_E$. Equivalently, for arbitrary input $y_E$, let $(y_E,0)_{O_f(z)}$ retain $y_E$ and fill the coordinates in $O_f(z)\setminus K_a$ with zeros. Then
\[
 \hat\tau_1(E,y_E)=
 \begin{cases}
 \displaystyle\sum_{z:O_f(z)\cap K_a=E}\frac{P(z)}{q(E)}
 \hat\tau\bigl(z,(y_E,0)_{O_f(z)}\bigr),&q(E)>0,\\[6pt]
 0,&q(E)=0.
 \end{cases}
\]
When $q(E)=0$, we define $\hat\tau_1(E,y_E)=0$ for every $y_E$. Such traces never occur under $P$, so this choice affects neither expectations nor risks. The law of iterated expectation and conditional Jensen's inequality give, for every $Y$,
\begin{align*}
 \E_q[\hat\tau_1(E,Y_E)]
 &=\E_P[\hat\tau(Z,Y^0_{O_f(Z)})],\\
 \E_q[(\hat\tau_1(E,Y_E)-\tau_a(Y))^2]
 &\leq\E_P[(\hat\tau(Z,Y^0_{O_f(Z)})-\tau_a(Y^0))^2].
\end{align*}
Taking the supremum over $Y$ proves the worst-case risk comparison with the original estimator; the expectation identity proves preservation of unbiasedness.

For each $E\in\mathcal E_a$, choose $M(E)\in\mathcal E_a^{\max}$ with $E\subseteq M(E)$, and set
\[
 q'(M)=\sum_{E:M(E)=M}q(E),
\]
where $q'(M)$ is the total probability of the original traces assigned to the maximal trace $M$. This law is implemented by assigning probability $q'(M)$ to one feasible assignment with trace $M$, giving a design $P'\in\mathcal C_{\max}$. For each fixed table $Y$ and each $M$ with $q'(M)>0$, define
\[
 \hat\tau_2(M,Y_M)
 :=\E_q\!\left[
 \hat\tau_1(E,Y_E)\,\middle|\,M(E)=M
 \right].
\]
Here the expectation is over the original trace $E$ with law $q$, conditional on its chosen maximal extension being $M$. Since every such $E$ is contained in $M$, this conditional expectation depends on $Y$ only through $Y_M$. Equivalently, for arbitrary input $y_M$, define
\[
 \hat\tau_2(M,y_M)=
 \begin{cases}
 \displaystyle\sum_{E:M(E)=M}\frac{q(E)}{q'(M)}\hat\tau_1(E,y_E),&q'(M)>0,\\[6pt]
 0,&q'(M)=0.
 \end{cases}
\]
For $q'(M)>0$, the weight $q(E)/q'(M)$ is the conditional probability of the original trace $E$ given $M(E)=M$. Here $y_E$ is available as a restriction of $y_M$, since $E\subseteq M$. When $q'(M)=0$, we define $\hat\tau_2(M,y_M)=0$ for every $y_M$. Such maximal traces never occur under $P'$, so this choice affects neither expectations nor risks. The law of iterated expectation and conditional Jensen's inequality give
\begin{align*}
 \E_{q'}[\hat\tau_2(M,Y_M)]&=\E_q[\hat\tau_1(E,Y_E)],\\
 \E_{q'}[(\hat\tau_2(M,Y_M)-\tau_a(Y))^2]
 &\leq\E_q[(\hat\tau_1(E,Y_E)-\tau_a(Y))^2],
\end{align*}
where the inequality follows from Jensen's inequality. Thus $(P',\hat\tau_2)$ preserves unbiasedness and has no larger worst-case risk than $(P,\hat\tau)$. Restricting designs to $\mathcal C_{\max}$ therefore cannot increase either minimax risk. The reverse inequality follows from $\mathcal C_{\max}\subseteq\cP(\cZ_0)$, proving \eqref{eq:trace-maximal-minimax}.
\end{proof}

\begin{corollary}[Induced conflict graph determines all-design minimax risks]
\label{cor:induced-conflict-minimax}
Under the assumptions of Theorem~\ref{thm:conflict-hypergraph-equivalence}, the maximal trace hyperedges $\mathcal E_a^{\max}$ are exactly the maximal independent sets in $H_{\rm conf}[K_a]$. Consequently, with the estimand coefficients $a$ and bounded-table outcome class $\cY(f)=[-1,1]^K$ fixed, the induced conflict graph determines both all-design minimax risks.
\end{corollary}

\begin{proof}
Let $S\subseteq K_a$ be an independent set in $H_{\rm conf}[K_a]$. It is also an independent set in $H_{\rm conf}$, so Theorem~\ref{thm:conflict-hypergraph-equivalence} gives an assignment $z$ with $S\subseteq O_f(z)$. Hence $S\subseteq O_f(z)\cap K_a\in\mathcal E_a$. Conversely, every trace hyperedge is an independent set in $H_{\rm conf}[K_a]$, since its potential outcomes are observed together. Thus the independent sets in $H_{\rm conf}[K_a]$ are exactly the subsets of trace hyperedges. Taking inclusion-maximal members identifies the maximal independent sets with $\mathcal E_a^{\max}$.

By Theorem~\ref{thm:trace-conflict-minimax}, designs supported on these maximal trace hyperedges suffice for both all-design minimax problems. Every probability distribution on these hyperedges is attainable by choosing a corresponding assignment for each hyperedge. The sampling-design reduction of Proposition~\ref{thm:three-problem-characterization} then shows that, with the estimand and outcome class fixed, the induced conflict graph determines both risks.
\end{proof}

The graph of \citet{KandirosHarshawSavje2026} is defined on the estimand-relevant potential outcomes and corresponds to $H_{\rm conf}[K_a]$. Their maximal-independent-set realization result \citep[Proposition~A.6]{KandirosHarshawSavje2026} identifies the same maximal observation sets. Theorem~\ref{thm:trace-conflict-minimax} shows explicitly why these sets suffice for the all-design minimax problem, both over all estimators and over unbiased estimators. For one equal-weight arm, our sampling-rank lower bound takes the independence-number form corresponding to the uniformly bounded case $q=\infty$ of their Theorem~3.2.

\subsubsection{When pairwise observability does not imply joint observability}
\label{subsec:coarsened-conflict-gap}
Assignment restrictions can make a pairwise-observable collection
impossible to observe jointly. For example, with three units, no
interference, and a budget of at most two treated units, any pair of
$Y_1(1),Y_2(1),Y_3(1)$ can be observed together, but all three cannot.
For the treatment arm $T$, the induced conflict graph $H_{\rm conf}[T]$
is edgeless, so
$\alpha\{H_{\rm conf}[T]\}=3$, whereas $r(T)=2$.
The set $T$ contains only one potential outcome per unit, so the
same-unit conflicts of $H_{\rm conf}$ have no edges in this induced
subgraph.

Coarsened exposures can impose joint constraints that are not determined
by pairwise conflicts, even when every binary treatment assignment is
feasible. Consider six units: three focal units $1,2,3$ and three other
units $a,b,c$, with focal-unit neighbor sets
\[
 N_1=\{a,b\},\qquad N_2=\{b,c\},\qquad N_3=\{a,c\}.
\]
For each focal unit, let the exposure record its own treatment and the
number of treated neighbors:
$f_i(z)=(z_i,\sum_{j\in N_i}z_j)$.
For $d\in\{1,2\}$, let
$T_d=\{(i,(1,d)):i=1,2,3\}$ index the arm in which the focal unit is
treated and has exactly $d$ treated neighbors.

Any two potential outcomes in $T_1$ can be observed together.
Observing all three, however, would require
\[
 z_a+z_b=1,\qquad z_b+z_c=1,\qquad z_a+z_c=1.
\]
Adding these equations gives $2(z_a+z_b+z_c)=3$, which is impossible for
binary assignments. Thus the induced graph $H_{\rm conf}[T_1]$ is edgeless and has
independence number $3$, whereas $r(T_1)=2$.
In contrast, treating $a,b,c$ and all three focal units reveals every
potential outcome in $T_2$. Hence the induced graph $H_{\rm conf}[T_2]$ is also
edgeless, but $r(T_2)=3$. The induced conflict graphs on these two arms have the
same pairwise structure, while their sampling ranks differ.
The sampling hypergraph records the additional fact that the three
outcomes in the first arm cannot be observed jointly.

Our applications exhibit the same distinction. In the complex-contagion
example, pairs of treated-low
potential outcomes can be jointly observed even though some larger
collections cannot; Lemma~\ref{lem:contagion-conflict-projection} describes the pairwise
projection, and Corollary~\ref{cor:contagion-conflict-gap} gives the resulting
lower-bound gap. The bipartite interior-dose example also has
pairwise-observable collections that are not jointly observable;
Corollary~\ref{cor:bipartite-full-projection-gap} gives the
corresponding lower-bound gap.

\section{Supporting results for Section~\ref{sec:fixed-design-ht}}
\label{supp:ht-upper-bounds}

This appendix first derives the exact variance of exposure HT for a fixed
potential-outcome table, then obtains two upper bounds for its worst-case
variance. Combining these upper bounds with the lower bounds from
Appendix~\ref{app:lower-bound-proofs} proves the fixed-design minimax
comparisons in Theorems~\ref{thm:fixed-bounded-positive-dependence}
and~\ref{thm:fixed-ht-comparison}. We then verify the main-text examples.
Finally, we prove Proposition~\ref{prop:fixed-ht-saturation-conditioning},
which uses an example to show that the H-T estimator is not always minimax optimal. Throughout this appendix, the variance calculations fix a design $P$,
an exposure mapping $f$, and
estimand coefficients $a$, so that the estimand is
$\tau_a(Y)=\sum_{(i,e)\in K_a}a_{i,e}Y_i(e)$, and every potential outcome
table $Y\in\cY(f)$ satisfies $|Y_i(e)|\leq1$.  Recall from
Section~\ref{sec:fixed-design-ht} the observation indicators
$I_{i,e}=\1\{f_i(Z)=e\}$ and the inclusion probabilities
$\pi_{P,i}(e)=P(I_{i,e}=1)$.  We assume that $P$ is positive on $K_a$, that
is, $\pi_{P,i}(e)>0$ for every $(i,e)\in K_a$.  The variance calculations
also use the pairwise inclusion probabilities
\[
 \pi_{P,ij}(e,e')=P\{I_{i,e}=1,I_{j,e'}=1\}.
\]

\subsection{Exact variance identity}
\label{subsec:app-ht-variance-identity}

We derive the exact variance of the exposure HT estimator in
equation~\eqref{supp:eq:general-ht-variance} and split it into four terms in
equation~\eqref{supp:eq:four-part-variance} below.  The exposure-wise bound in Proposition~\ref{supp:prop:exposure-wise-ht} and the cross-exposure bound in Proposition~\ref{supp:prop:cross-exposure-ht} control these terms. Recall that the exposure HT estimator is

\[
 \widehat\tau_{\mathrm{HT},P}
 =\sum_{(i,e)\in K_a}
 a_{i,e}\frac{I_{i,e}Y_i(e)}{\pi_{P,i}(e)}.
\]
For a fixed potential outcome table $Y$, the estimator is a linear
combination of the variables $I_{i,e}/\pi_{P,i}(e)$ with coefficients
$a_{i,e}Y_i(e)$.  For any two relevant potential outcome indices
$(i,e),(j,e')\in K_a$, the covariance of these variables is
\[
 \Cov_P\left\{
   \frac{I_{i,e}}{\pi_{P,i}(e)},
   \frac{I_{j,e'}}{\pi_{P,j}(e')}
 \right\}
 =\frac{\pi_{P,ij}(e,e')}
 {\pi_{P,i}(e)\pi_{P,j}(e')}-1
 =:\phi_{P,ij}(e,e').
\]
This includes the case $(i,e)=(j,e')$, where
$\pi_{P,ii}(e,e)=\pi_{P,i}(e)$ and the covariance is the variance
$\pi_{P,i}(e)^{-1}-1$.  Therefore

\begin{equation}
 \Var_P(\widehat\tau_{\mathrm{HT},P}\mid Y)
 =\sum_{\substack{(i,e)\in K_a\\(j,e')\in K_a}}
 a_{i,e}a_{j,e'}Y_i(e)Y_j(e')\phi_{P,ij}(e,e').
 \label{supp:eq:general-ht-variance}
\end{equation}
All sums over pairs of indices in this section use ordered pairs.  We split the sum in
equation~\eqref{supp:eq:general-ht-variance} according to whether the two
potential outcomes belong to the same unit and whether they have
the same exposure:

\begin{equation}
 \Var_P(\widehat\tau_{\mathrm{HT},P}\mid Y)
 =D_a(Y;P)+\mathrm{OD}_a(Y;P)
 +\mathrm{CD}_{\mathrm{self},a}(Y;P)
 +\mathrm{CD}_{\mathrm{pair},a}(Y;P),
 \label{supp:eq:four-part-variance}
\end{equation}
where
\begin{align}
 D_a(Y;P)
 &=\sum_{(i,e)\in K_a}a_{i,e}^2Y_i(e)^2
   \left\{\frac1{\pi_{P,i}(e)}-1\right\},
 \notag\\ %
 \mathrm{OD}_a(Y;P)
 &=\sum_{\substack{i\ne j\\(i,e),(j,e)\in K_a}}
   a_{i,e}a_{j,e}Y_i(e)Y_j(e)\phi_{P,ij}(e,e),
 \notag\\ %
 \mathrm{CD}_{\mathrm{self},a}(Y;P)
 &=\sum_i\sum_{\substack{e\ne e'\\(i,e),(i,e')\in K_a}}
   a_{i,e}a_{i,e'}Y_i(e)Y_i(e')\phi_{P,ii}(e,e'),
 \notag\\ %
 \mathrm{CD}_{\mathrm{pair},a}(Y;P)
 &=\sum_{\substack{i\ne j,\ e\ne e'\\
       (i,e),(j,e')\in K_a}}
   a_{i,e}a_{j,e'}Y_i(e)Y_j(e')\phi_{P,ij}(e,e').
 \notag %
\end{align}
The diagonal term $D_a$ collects the pairs with $(i,e)=(j,e')$; the
off-diagonal term $\mathrm{OD}_a$ collects pairs of distinct units with the
same exposure; and the cross-exposure terms $\mathrm{CD}_{\mathrm{self},a}$
and $\mathrm{CD}_{\mathrm{pair},a}$ collect pairs with different exposures,
for the same unit and for two distinct units, respectively.  A unit cannot
have two exposures at once, so $\pi_{P,ii}(e,e')=0$ and
$\phi_{P,ii}(e,e')=-1$ for $e\ne e'$.  Since $|Y_i(e)|\leq1$, the diagonal term satisfies
$D_a(Y;P)\leq D_a(P)$ for every $Y\in\cY(f)$, where $D_a(P)$ is defined in equation~\eqref{eq:fixed-diagonal}.
The quantities with an argument $Y$ are contributions to the exact
variance for that table. The corresponding quantities without $Y$ bound
these contributions uniformly over $\cY(f)$.
For the off-diagonal term, taking absolute values gives
$|\mathrm{OD}_a(Y;P)|\leq\mathrm{OD}_a^+(P)+\mathrm{OD}_a^-(P)$,
where $\mathrm{OD}_a^-(P)$ is defined in
Section~\ref{subsec:app-exposure-wise}. Lemma~\ref{supp:lem:same-exposure-sign}
bounds the contribution from negative observation-indicator covariances
by $D_a(P)+\mathrm{OD}_a^+(P)$. Note that the signs in the superscripts refer to covariances between sampling indicators; the signs of the variance contributions
also depend on the coefficients and potential outcomes.

\begin{remark}[Independent observation indicators]
\label{rem:independent-indicators}
Suppose that each unit has at most one relevant potential outcome, and that the observation indicators $I_{i,e}$, $(i,e)\in K_a$, are
pairwise independent.  Then $\phi_{P,ij}(e,e')=0$ for every pair of distinct
potential outcome indices, and equation~\eqref{supp:eq:general-ht-variance}
reduces to
\[
 \Var_P(\widehat\tau_{\mathrm{HT},P}\mid Y)
 =\sum_{k\in K_a}a_k^2Y_k^2(\pi_{P,k}^{-1}-1).
\]
Its maximum over $\cY(f)$, attained when $|Y_k|=1$ for every $k\in K_a$, equals
the unbiased singleton lower bound in equation~\eqref{eq:singleton-unb}.  Thus
exposure HT is minimax among unbiased estimators.  The first
assumption is needed because two exposures of the same unit are mutually
exclusive, so their indicators are dependent whenever both have positive
probability.

Sampling from a finite population of $N$ units is the special case with one
potential outcome per unit, observed when the unit is sampled
with inclusion probability $\pi_i$.  For the population mean of outcomes in
$[-1,1]$, equation~\eqref{eq:singleton-unb} gives the lower bound
$N^{-2}\sum_i(\pi_i^{-1}-1)$.  \citet{AronowLopatto2026} prove the
corresponding lower bound for the population total, with outcomes in known
intervals, for all design-unbiased estimators. Their population-total
lower bound is attained exactly when the inclusion indicators are pairwise
independent. The calculation above shows that pairwise independence
suffices for HT to attain the population-mean lower bound in our setting.
Theorem~\ref{thm:fixed-bounded-positive-dependence} extends this sufficient condition, at the level of rates, to designs with bounded positive dependence, including negatively dependent indicators, and to causal contrasts under interference, where the exposures of a unit are mutually exclusive and pairwise independence cannot hold.
\end{remark}

\subsection{Upper bounds for the variance of exposure HT}
\label{subsec:app-exposure-wise}

We use two approaches to bound the worst-case variance of exposure HT.
Proposition~\ref{supp:prop:exposure-wise-ht} bounds the variance within
each exposure level and combines the exposure-specific errors by
Cauchy--Schwarz, introducing a factor equal to the number of exposure
levels. This bound is used in the fixed-design minimax comparisons.
Proposition~\ref{supp:prop:cross-exposure-ht} instead bounds the
cross-exposure covariance terms directly, avoiding the extra factor on number of exposure levels; we use
this bound for the bounded-local-interference example.

For the first approach, Lemma~\ref{supp:lem:same-exposure-sign} controls
the negative same-exposure covariance terms. Let $\cD_a=\{e:(i,e)\in K_a\text{ for some }i\}$ and
$L_a=|\cD_a|$. Define
\begin{align*}
 D_{a,e}(P)
 &=\sum_{i:(i,e)\in K_a}a_{i,e}^2
   \left\{\frac1{\pi_{P,i}(e)}-1\right\},\\
 \mathrm{OD}_{a,e}^+(P)
 &=\sum_{\substack{i\ne j\\(i,e),(j,e)\in K_a}}
   |a_{i,e}a_{j,e}|\{\phi_{P,ij}(e,e)\}^{+},\\
 \mathrm{OD}_{a,e}^-(P)
 &=\sum_{\substack{i\ne j\\(i,e),(j,e)\in K_a}}
   |a_{i,e}a_{j,e}|\{\phi_{P,ij}(e,e)\}^{-}.
\end{align*}
Write $D_a(P)=\sum_eD_{a,e}(P)$ and
$\mathrm{OD}_a^\pm(P)=\sum_e\mathrm{OD}_{a,e}^\pm(P)$.

\begin{lemma}[Sign bound for same-exposure association]
\label{supp:lem:same-exposure-sign}
For every $e\in\cD_a$,
\[
 \mathrm{OD}_{a,e}^-(P)\leq D_{a,e}(P)+\mathrm{OD}_{a,e}^+(P).
\]
Consequently, $\mathrm{OD}_a^-(P)\leq D_a(P)+\mathrm{OD}_a^+(P)$.
\end{lemma}

\begin{proof}
For fixed $e$, consider
\[
 W_e=\sum_{i:(i,e)\in K_a}|a_{i,e}|
 \frac{I_{i,e}}{\pi_{P,i}(e)}.
\]
Expanding $\Var_P(W_e)$, we get
\[
 \begin{aligned}
 0\leq\Var_P(W_e)
 &=D_{a,e}(P)+
   \sum_{\substack{i\ne j\\(i,e),(j,e)\in K_a}}
   |a_{i,e}a_{j,e}|\phi_{P,ij}(e,e)\\
 &=D_{a,e}(P)+\mathrm{OD}_{a,e}^+(P)
   -\mathrm{OD}_{a,e}^-(P).
 \end{aligned}
\]
The second line splits $\phi_{P,ij}(e,e)=\{\phi_{P,ij}(e,e)\}^+-\{\phi_{P,ij}(e,e)\}^-$. Rearranging gives
$\mathrm{OD}_{a,e}^-(P)\leq
D_{a,e}(P)+\mathrm{OD}_{a,e}^+(P)$. Summing this inequality over
$e\in\cD_a$ proves the second claim.
\end{proof}

\begin{proposition}[Exposure-wise HT upper bound]
\label{supp:prop:exposure-wise-ht}
For every positive design $P$,
\[
 \sup_{Y\in\cY(f)}\Var_P(\widehat\tau_{\mathrm{HT},P}\mid Y)
 \leq2\left\{\sum_{e\in\cD_a}
 \sqrt{D_{a,e}(P)+\mathrm{OD}_{a,e}^+(P)}\right\}^{2}.
\]
In particular,
\[
 \sup_{Y\in\cY(f)}\Var_P(\widehat\tau_{\mathrm{HT},P}\mid Y)
 \leq2L_a\{D_a(P)+\mathrm{OD}_a^+(P)\}.
\]
\end{proposition}

\begin{proof}
Write $\widehat\tau_{\mathrm{HT},P}-\tau_a(Y)=\sum_{e\in\cD_a}U_e$, where
\[
 U_e=\sum_{i:(i,e)\in K_a}a_{i,e}Y_i(e)
 \left\{\frac{I_{i,e}}{\pi_{P,i}(e)}-1\right\}.
\]
Here $U_e$ is the HT error for the estimand with coefficients $a_{i,e}$ at exposure $e$ and zero elsewhere, so equation~\eqref{supp:eq:general-ht-variance} gives
\[
 \Var_P(U_e\mid Y)=\sum_{\substack{i,j\\(i,e),(j,e)\in K_a}}a_{i,e}a_{j,e}Y_i(e)Y_j(e)\phi_{P,ij}(e,e).
\]
Bounding $|Y_i(e)|\leq1$ and $|\phi_{P,ij}(e,e)|=\{\phi_{P,ij}(e,e)\}^++\{\phi_{P,ij}(e,e)\}^-$ gives
\begin{align*}
 \Var_P(U_e\mid Y)
 &\leq D_{a,e}(P)
   +\sum_{\substack{i\ne j\\(i,e),(j,e)\in K_a}}
    |a_{i,e}a_{j,e}|\,|\phi_{P,ij}(e,e)|\\
 &=D_{a,e}(P)+\mathrm{OD}_{a,e}^+(P)
   +\mathrm{OD}_{a,e}^-(P)\\
 &\leq2\{D_{a,e}(P)+\mathrm{OD}_{a,e}^+(P)\},
\end{align*}
where the last inequality follows from
Lemma~\ref{supp:lem:same-exposure-sign}.

We now combine the exposure-specific contributions. The Cauchy--Schwarz inequality for
covariances gives
\[
 \begin{aligned}
 \Var_P\left(\sum_{e\in\cD_a}U_e\middle|Y\right)
 &=\sum_{e,e'\in\cD_a}\Cov_P(U_e,U_{e'}\mid Y)\\
 &\leq\sum_{e,e'\in\cD_a}
   \sqrt{\Var_P(U_e\mid Y)\Var_P(U_{e'}\mid Y)}\\
 &=\left\{\sum_{e\in\cD_a}
   \sqrt{\Var_P(U_e\mid Y)}\right\}^{2}\\
 &\leq2\left\{\sum_{e\in\cD_a}
   \sqrt{D_{a,e}(P)+\mathrm{OD}_{a,e}^+(P)}\right\}^{2}.
 \end{aligned}
\]
Taking the supremum over $Y\in\cY(f)$ proves the first result. Finally, the Cauchy--Schwarz inequality $(\sum_ex_e)^2\leq L_a\sum_ex_e^2$, applied with $x_e=\{D_{a,e}(P)+\mathrm{OD}_{a,e}^+(P)\}^{1/2}$, gives
\[
 \left\{\sum_{e\in\cD_a}
 \sqrt{D_{a,e}(P)+\mathrm{OD}_{a,e}^+(P)}\right\}^{2}
 \leq L_a\{D_a(P)+\mathrm{OD}_a^+(P)\},
\]
which proves the second result.
\end{proof}

The factor $L_a$ in Proposition~\ref{supp:prop:exposure-wise-ht} comes from combining the exposure levels by the Cauchy--Schwarz inequality.  We also use an HT upper bound that keeps the cross-exposure covariance terms and has no factor $L_a$.  Since $|\phi_{P,ii}(e,e')|=1$ for $e\ne e'$, the same-unit term below carries no covariance factor.

To bound cross-exposure terms directly, define
\begin{align*}
 \mathrm{CD}_{\mathrm{self},a}(P)
 &=\sum_i\sum_{\substack{e\ne e'\\(i,e),(i,e')\in K_a}}
   |a_{i,e}a_{i,e'}|,\\
 \mathrm{CD}_{\mathrm{pair},a}(P)
 &=\sum_{\substack{i\ne j,\ e\ne e'\\
       (i,e),(j,e')\in K_a}}
   |a_{i,e}a_{j,e'}|\,|\phi_{P,ij}(e,e')|.
\end{align*}

\begin{lemma}[Same-unit cross-exposure bound]
\label{supp:lem:self-cross-domination}
$\mathrm{CD}_{\mathrm{self},a}(P)\leq D_a(P)$.
\end{lemma}

\begin{proof}
Fix a unit $i$ with at least one relevant exposure, and sum over $e$ such that $(i,e)\in K_a$. Write $x_e=|a_{i,e}|$. Since the exposure events are disjoint,
$\sum_e\pi_{P,i}(e)\leq1$.  Applying
Lemma~\ref{lem:weighted-reciprocal} with $w_e=x_e$ and
$p_e=\pi_{P,i}(e)$ gives
\[
 \sum_e\frac{x_e^2}{\pi_{P,i}(e)}
 \geq\frac{(\sum_ex_e)^2}{\sum_e\pi_{P,i}(e)}
 \geq(\sum_ex_e)^2.
\]
Subtracting $\sum_ex_e^2$ gives
\[
 \sum_{e\ne e'}|a_{i,e}a_{i,e'}|
 \leq\sum_e a_{i,e}^2\{\pi_{P,i}(e)^{-1}-1\}.
\]
Summing over the units proves the bound.
\end{proof}

\begin{proposition}[HT upper bound with cross-exposure component]
\label{supp:prop:cross-exposure-ht}
For every positive design $P$,
\begin{equation}
 \sup_{Y\in\cY(f)}\Var_P(\widehat\tau_{\mathrm{HT},P}\mid Y)
 \leq3D_a(P)+2\mathrm{OD}_a^+(P)
 +\mathrm{CD}_{\mathrm{pair},a}(P).
 \label{supp:eq:cross-exposure-upper}
\end{equation}
Consequently, if
$\mathrm{CD}_{\mathrm{pair},a}(P)\leq
C_{\mathrm{cross}}\{D_a(P)+\mathrm{OD}_a^+(P)\}$, then
\[
 \sup_{Y\in\cY(f)}\Var_P(\widehat\tau_{\mathrm{HT},P}\mid Y)
 \leq(3+C_{\mathrm{cross}})\{D_a(P)+\mathrm{OD}_a^+(P)\}.
\]
\end{proposition}

\begin{proof}
Taking absolute values in equation~\eqref{supp:eq:four-part-variance} yields
\[
 \sup_Y\Var_P(\widehat\tau_{\mathrm{HT},P}\mid Y)
 \leq D_a(P)+\mathrm{OD}_a^+(P)+\mathrm{OD}_a^-(P)
 +\mathrm{CD}_{\mathrm{self},a}(P)
 +\mathrm{CD}_{\mathrm{pair},a}(P).
\]
Substituting the bounds
$\mathrm{OD}_a^-(P)\leq D_a(P)+\mathrm{OD}_a^+(P)$ and
$\mathrm{CD}_{\mathrm{self},a}(P)\leq D_a(P)$ from Lemmas~\ref{supp:lem:same-exposure-sign} and~\ref{supp:lem:self-cross-domination} gives equation~\eqref{supp:eq:cross-exposure-upper}. Under the stated bound on $\mathrm{CD}_{\mathrm{pair},a}(P)$, the right-hand side is at most
$(3+C_{\mathrm{cross}})D_a(P)+(2+C_{\mathrm{cross}})\mathrm{OD}_a^+(P)$, which proves the consequence.
\end{proof}

\subsection{Proofs of the fixed-design comparisons}
\label{subsec:app-fixed-ht-comparison}
Each proof combines the exposure-wise HT upper bound with fixed-design
lower bounds from Appendix~\ref{app:lower-bound-proofs}. The first
comparison uses singleton lower bounds. The second also uses a packing
lower bound to account for stronger positive same-exposure dependence.

\begin{resultrestatement}{Theorem~\ref{thm:fixed-bounded-positive-dependence}}{Bounded positive dependence}
Suppose that \(P\) has bounded positive same-exposure dependence with
\(\mathrm{OD}_a^+(P)\leq C_0D_a(P)\), and that \(L_a\leq L_0\).  Then
\begin{equation*}
 D_a(P)\leq\Runb(P;f,a)
 \leq2L_0(1+C_0)D_a(P).
 \tag{\ref{eq:fixed-bpd-unbiased}}
\end{equation*}
If, in addition,
\(\pi_{P,k}\geq c>0\) for every \(k\in K_a\), with \(c\) independent of the
population size, then
\begin{equation*}
 cD_a(P)\leq\Rall(P;f,a)\leq\Runb(P;f,a)
 \leq2L_0(1+C_0)D_a(P).
 \tag{\ref{eq:fixed-bpd-all}}
\end{equation*}
\end{resultrestatement}

\begin{proof}[Proof of Theorem~\ref{thm:fixed-bounded-positive-dependence}]
Corollary~\ref{cor:singleton-lower} gives
\(\Runb(P;f,a)\geq D_a(P)\).  Exposure HT is unbiased, so $\Runb(P;f,a)\leq V_{\mathrm{HT}}(P)$.  Proposition~\ref{supp:prop:exposure-wise-ht} and the bounded-positive-dependence condition
give
\[
 \Runb(P;f,a)
 \leq V_{\mathrm{HT}}(P)
 \leq2L_a\{D_a(P)+\mathrm{OD}_a^+(P)\}
 \leq2L_0(1+C_0)D_a(P),
\]
where the last step uses $L_a\leq L_0$ and $\mathrm{OD}_a^+(P)\leq C_0D_a(P)$.
If \(\pi_{P,k}\geq c\) for every \(k\in K_a\), the all-estimator singleton
bound also gives
\[
 \Rall(P;f,a)
 \geq\sum_{k\in K_a}a_k^2(1-\pi_{P,k})
 \geq cD_a(P),
\]
because
\(1-\pi_{P,k}\geq c(1-\pi_{P,k})/\pi_{P,k}\).  Finally,
\(\Rall(P;f,a)\leq\Runb(P;f,a)\), which completes the proof.
\end{proof}

\begin{resultrestatement}{Theorem~\ref{thm:fixed-ht-comparison}}{HT optimality under packing-matched dependence}
Suppose that \(L_a\leq L_0\), \(\pi_{P,k}\geq c\) for every \(k\in K_a\),
and
\begin{equation*}
 \mathrm{OD}_a^+(P)\leq C_0L_{\cB}(P),
 \tag{\ref{eq:packing-matches-positive-dependence}}
\end{equation*}
where \(L_0\), \(c>0\), and \(C_0\) are independent of the population size.
Then
\begin{equation*}
 \begin{aligned}
 \max\{cD_a(P),L_{\cB}(P)\}
 &\leq\Rall(P;f,a)\leq\Runb(P;f,a)\\
 &\leq V_{\mathrm{HT}}(P)
 \leq2L_0\{c^{-1}+C_0\}\Rall(P;f,a).
 \end{aligned}
 \tag{\ref{eq:fixed-ht-comparison}}
\end{equation*}
Thus exposure HT remains minimax-rate optimal when the positive dependence is
matched by joint missingness.
\end{resultrestatement}

\begin{proof}[Proof of Theorem~\ref{thm:fixed-ht-comparison}]
Choose $s_k=\operatorname{sign}(a_k)$ in each perturbation set, so that $c_b=\sum_{k\in B_b}|a_k|$ for each $B_b\in\cB$. Theorem
\ref{thm:fixed-joint-missingness}(i) gives $\Rall(P;f,a)\geq\sum_bc_b^2\rho_P(B_b)=L_{\cB}(P)$.
As in the proof of Theorem~\ref{thm:fixed-bounded-positive-dependence}, the condition $\pi_{P,k}\geq c$ gives $\Rall(P;f,a)\geq cD_a(P)$.  The
unbiased estimators form a subclass, so
$\Rall(P;f,a)\leq\Runb(P;f,a)$, and exposure HT is unbiased, so
$\Runb(P;f,a)\leq V_{\mathrm{HT}}(P)$.  This proves all but the final
inequality in equation~\eqref{eq:fixed-ht-comparison}.

For the final inequality, Proposition~\ref{supp:prop:exposure-wise-ht} and $L_a\leq L_0$ give the first line below, the packing-matched condition gives the second, and the lower bounds $D_a(P)\leq c^{-1}\Rall(P;f,a)$ and $L_{\cB}(P)\leq\Rall(P;f,a)$ established above give the third:
\begin{align*}
 V_{\mathrm{HT}}(P)
 &\leq2L_0\{D_a(P)+\mathrm{OD}_a^+(P)\}\\
 &\leq2L_0\{D_a(P)+C_0L_{\cB}(P)\}\\
 &\leq2L_0(c^{-1}+C_0)\Rall(P;f,a).
\end{align*}
\end{proof}

\subsection{Calculations for minimax risk of the fixed-design examples}
\label{subsec:app-fixed-ht-examples}
We give the calculations for the examples in Section~\ref{sec:fixed-design-ht}, in the order in which they appear in the main text.

\subsubsection{No-interference ATE}

Under no interference, $f_i(Z)=Z_i$, so the observation indicators are $I_{i,1}=Z_i$ and $I_{i,0}=1-Z_i$. For the ATE, $a_{i,1}=1/n$ and $a_{i,0}=-1/n$ for every unit $i$.

Under independent Bernoulli assignment with treatment probability $p\in(0,1)$, $Z_i$ and $Z_j$ are independent for $i\ne j$. Hence $I_{i,e}$ and $I_{j,e'}$ are independent for every $e,e'\in\{0,1\}$, giving
\[
 P\{I_{i,e}=1,I_{j,e'}=1\}
 =\pi_{P,i}(e)\pi_{P,j}(e'),\qquad i\ne j.
\]
In particular, $\phi_{P,ij}(e,e)=0$ for both $e=0$ and $e=1$, so $\mathrm{OD}_a^+(P)=0$. Since $\pi_{P,i}(1)=p$ and $\pi_{P,i}(0)=1-p$,
\[
 D_a(P)=\sum_{i=1}^n\frac1{n^2}
 \bigl\{(p^{-1}-1)+((1-p)^{-1}-1)\bigr\}
 =\frac1n\{p^{-1}+(1-p)^{-1}-2\}.
\]

Under complete randomization with exactly $n_1$ treated units, put $n_0=n-n_1$ and $p=n_1/n$, where $1\leq n_1\leq n-1$. For distinct units $i\ne j$ and $e\in\{0,1\}$,
\[
 \pi_{P,i}(e)=\frac{n_e}{n},\qquad
 P\{I_{i,e}=1,I_{j,e}=1\}
 =\frac{n_e(n_e-1)}{n(n-1)}
 \leq\left(\frac{n_e}{n}\right)^2.
\]
Thus $\phi_{P,ij}(e,e)\leq0$ for both exposure levels, and again $\mathrm{OD}_a^+(P)=0$. The same diagonal calculation gives
\[
 D_a(P)=\frac1n\{p^{-1}+(1-p)^{-1}-2\},\qquad p=\frac{n_1}{n}.
\]
In either design, if $p$ is bounded away from zero and one, then the inclusion probabilities are bounded below, $L_a=2$, and $D_a(P)\asymp1/n$. Theorem~\ref{thm:fixed-bounded-positive-dependence} therefore gives
\[
 \Rall(P;f,a)\asymp\Runb(P;f,a)\asymp\frac1n.
\]

\subsubsection{Bounded local interference}
\label{app:fixed-design-regimes}
\label{sec:neighborhood-regimes}

We now control the covariance terms under independent assignment and bounded local interference.

Suppose treatment is assigned independently and interference is restricted to closed neighborhoods $N_i^+$. Thus every relevant observation indicator $I_{i,e}$ is a function of $Z_{N_i^+}$. Assume
\begin{equation}
  \max_i|N_i^+|\leq r,
  \qquad
  \max_i|\{j:N_i^+\cap N_j^+\ne\varnothing\}|\leq M,
  \label{eq:main-local-overlap}
\end{equation}
where $r$ and $M$ are bounded uniformly in $n$. Suppose also that each unit
has at most $L_{\max}$ relevant exposures, with $L_{\max}$ bounded uniformly
in $n$, and that,
uniformly on $K_a$,
\begin{equation}
  0<\underline\pi\leq\pi_{P,i}(e)
  \leq\overline\pi<1.
  \label{eq:main-local-positivity}
\end{equation}
This covers, for example, identity neighborhood exposures and count or fraction exposures computed from $Z_{N_i^+}$.

\begin{proposition}[Bounded local interference]
\label{prop:main-bounded-local}
Under independent assignment, equations~\eqref{eq:main-local-overlap}--\eqref{eq:main-local-positivity}, and the bound of $L_{\max}$ relevant exposures per unit,
\[
  \mathrm{OD}_a^+(P)+\mathrm{CD}_{\mathrm{pair},a}(P)\leq C D_a(P),
\]
where $C$ depends only on
$\underline\pi,\overline\pi,M$, and $L_{\max}$. Consequently,
\[
  \Runb(P;f,a)\asymp D_a(P),
\]
and HT is fixed-design minimax-rate optimal among unbiased estimators.
\end{proposition}

\begin{proof}
\label{app:bounded-local-proof}  The unbiased
singleton lower bound in equation~\eqref{eq:singleton-unb} gives
\[
  \Runb(P;f,a)\geq
  D_a(P)
  =
  \sum_{(i,e)\in K_a}a_{i,e}^2
  \left\{\frac1{\pi_{P,i}(e)}-1\right\}.
\]
It remains to control the off-diagonal terms in the HT upper bound.

If $N_i^+\cap N_j^+=\varnothing$, independent assignment gives $Z_{N_i^+}\perp Z_{N_j^+}$. Since $I_{i,e}$ and $I_{j,e'}$ are functions of these respective vectors,
\[
 P\{I_{i,e}=1,I_{j,e'}=1\}
 =\pi_{P,i}(e)\pi_{P,j}(e'),\qquad
 \phi_{P,ij}(e,e')=0.
\]
Thus only pairs with overlapping coordinate sets can contribute. For any relevant pair, positivity gives the uniform bound
\[
 |\phi_{P,ij}(e,e')|
 =\left|\frac{P\{I_{i,e}=1,I_{j,e'}=1\}}
 {\pi_{P,i}(e)\pi_{P,j}(e')}-1\right|
 \leq B,\qquad B=\underline\pi^{-2}+1.
\]
For $k=(i,e)\in K_a$, define
\[
 \mathcal N(k)=\{\ell=(j,e')\in K_a:j\ne i,\ N_i^+\cap N_j^+\ne\varnothing\}.
\]
The overlap and exposure-count assumptions imply $|\mathcal N(k)|\leq ML_{\max}$. Each term in $\mathrm{OD}_a^+(P)+\mathrm{CD}_{\mathrm{pair},a}(P)$ corresponds to an ordered pair $(k,\ell)$ with $\ell\in\mathcal N(k)$. Using $2|a_ka_\ell|\leq a_k^2+a_\ell^2$ and the symmetry of this relation,
\begin{align*}
 \mathrm{OD}_a^+(P)+\mathrm{CD}_{\mathrm{pair},a}(P)
 &\leq B\sum_{k\in K_a}\sum_{\ell\in\mathcal N(k)}|a_ka_\ell|\\
 &\leq\frac B2\sum_{k\in K_a}\sum_{\ell\in\mathcal N(k)}(a_k^2+a_\ell^2)\\
 &=B\sum_{k\in K_a}|\mathcal N(k)|a_k^2
 \leq BML_{\max}\sum_{k\in K_a}a_k^2.
\end{align*}
Set $C_1=BML_{\max}$.
The upper bound \(\pi_{P,k}\leq\overline\pi<1\) gives
\[
  D_a(P)\geq
  \left(\frac1{\overline\pi}-1\right)\sum_{k\in K_a}a_k^2 ,
\]
so \(\mathrm{OD}_a^+(P)+\mathrm{CD}_{\mathrm{pair},a}(P)\leq C D_a(P)\) with $C=C_1/(\overline\pi^{-1}-1)$.  Equation~\eqref{supp:eq:cross-exposure-upper} then gives $V_{\mathrm{HT}}(P)\leq3D_a(P)+2\{\mathrm{OD}_a^+(P)+\mathrm{CD}_{\mathrm{pair},a}(P)\}\leq(3+2C)D_a(P)$, which matches the singleton lower bound in equation~\eqref{eq:singleton-unb} up to the factor $3+2C$.

\end{proof}

We now verify Example~2 of Section~\ref{sec:fixed-design-ht}.  Let $P$ be
Bernoulli assignment with treatment probability $p$ bounded away from zero and
one, on a graph whose maximum degree is at most $\Delta$, with $\Delta$
bounded uniformly in $n$.  For the DTE or the full-neighborhood TTE, write
$\pi_i^+$ and $\pi_i^-$ for the exposure probabilities of the two target
exposures of unit $i$.  Each target observation indicator is a function of
$Z_{N_i^+}$, so equation~\eqref{eq:main-local-overlap} holds with
$r=\Delta+1$, and $M=\Delta^2+1$, the largest number of units
within distance two of a given unit, including the unit itself. Indeed,
$N_i^+\cap N_j^+\ne\varnothing$ if and only if $\operatorname{dist}_G(i,j)\leq2$, and there are at most $1+\Delta+\Delta(\Delta-1)=1+\Delta^2$ such units.  %
Each unit has $L_{\max}=2$ relevant exposures. Writing $d_i=|N_i|$, their probabilities are
\[
 \begin{array}{c|cc}
 &\pi_i^+&\pi_i^-\\ \hline
 \mathrm{DTE}&p(1-p)^{d_i}&(1-p)^{d_i+1}\\
 \mathrm{TTE}&p^{d_i+1}&(1-p)^{d_i+1}.
 \end{array}
\]
Since $d_i\leq\Delta$, every displayed probability satisfies
\[
 \min(p,1-p)^{\Delta+1}\leq\pi_i^\pm\leq\max(p,1-p).
\]
Thus equation~\eqref{eq:main-local-positivity} holds with $\underline\pi=\min(p,1-p)^{\Delta+1}$ and $\overline\pi=\max(p,1-p)$, uniformly in $n$.
Proposition~\ref{prop:main-bounded-local} gives
$V_{\mathrm{HT}}(P)\lesssim D_a(P)$, and since every relevant potential
outcome has coefficient $\pm1/n$,
\[
 D_a(P)=\frac1{n^2}\sum_{i=1}^n\bigl\{(\pi_i^+)^{-1}+(\pi_i^-)^{-1}-2\bigr\}
 \asymp\frac1n.
\]
The two target exposures of each unit are mutually exclusive, so
$\pi_i^++\pi_i^-\leq1$, and the all-estimator singleton lower bound in
equation~\eqref{eq:singleton-all} gives
\[
 \Rall(P;f,a)
 \geq\frac1{n^2}\sum_{i=1}^n(2-\pi_i^+-\pi_i^-)
 \geq\frac1n.
\]
Hence $\Rall(P;f,a)\asymp\Runb(P;f,a)\asymp V_{\mathrm{HT}}(P)\asymp1/n$.

\subsubsection{A rare exposure on a regular graph}

We verify Example~3 of Section~\ref{sec:fixed-design-ht}.  Let $G$ be a
$d$-regular graph, let $P$ be Bernoulli assignment with $p=1/(d+2)$, and
consider the DTE, whose relevant potential outcomes have coefficients
$\pm1/n$.  The treated exposure $(1,0_{N_i})$ has probability $p(1-p)^d$ and
the control exposure $(0,0_{N_i})$ has probability $(1-p)^{d+1}$.  With
$x=(1-p)^{-1}=(d+2)/(d+1)$,
\[
 D_a(P)=\frac1n\bigl\{(d+2)x^d+x^{d+1}-2\bigr\}
 =\frac{(d+2)x^{d+1}-2}{n}.
\]
Since $1\leq x^{d+1}\leq e$, this gives $(d+1)/n\leq D_a(P)\leq e(d+2)^2/\{(d+1)n\}$, so
$D_a(P)\asymp(d+1)/n$.

We next bound $\mathrm{OD}_a^+(P)$. Write $e_i^+=(1,0_{N_i})$ and $e_i^-=(0,0_{N_i})$ for unit $i$'s treated and control target exposures. For every pair of distinct units, define
\[
 s_{ij}=|N_i\cap N_j|,\qquad t_{ij}=|N_i^+\cap N_j^+|.
\]
The treated-arm potential outcomes of adjacent units have a pairwise conflict, so $\phi_{P,ij}(e_i^+,e_j^+)=-1$. For nonadjacent $i,j$, observing both treated-arm outcomes requires $Z_i=Z_j=1$ and $Z_v=0$ for all $v\in N_i\cup N_j$, giving $\phi_{P,ij}(e_i^+,e_j^+)=x^{s_{ij}}-1$. For the control arm, observing both outcomes requires $Z_v=0$ for all $v\in N_i^+\cup N_j^+$, giving $\phi_{P,ij}(e_i^-,e_j^-)=x^{t_{ij}}-1$. For $0\leq u\leq d+1$, convexity of the
exponential gives
\[
 x^u-1\leq e^{u/(d+1)}-1\leq(e-1)\frac{u}{d+1}.
\]
For a fixed unit $i$, counting the two-step walks from $i$ gives
$\sum_{j\ne i}s_{ij}=d(d-1)$, and counting the units $k\in N_i^+$ together with
the $d$ other closed neighborhoods containing $k$ gives
$\sum_{j\ne i}t_{ij}=d(d+1)$.  Each pair carries weight $1/n^2$, so
\[
 \mathrm{OD}_a^+(P)
 \leq\frac{e-1}{n^2(d+1)}\sum_{i=1}^n\{d(d-1)+d(d+1)\}
 \leq\frac{2(e-1)(d+1)}{n}
 \leq2(e-1)D_a(P).
\]
Thus the design has bounded positive dependence with $C_0=2(e-1)$ and
$L_a=2$, and Theorem~\ref{thm:fixed-bounded-positive-dependence} gives
\[
 D_a(P)\leq\Runb(P;f,a)\leq V_{\mathrm{HT}}(P)\leq4(2e-1)D_a(P).
\]
The treated exposure probability is of order $1/(d+1)$, so the all-estimator
part of Theorem~\ref{thm:fixed-bounded-positive-dependence} does not apply.
Proposition~\ref{prop:dte-regular-all-estimator} gives the matching
all-estimator lower bound directly.

\begin{proposition}[A fixed design lower bound for the regular-graph DTE]
\label{prop:dte-regular-all-estimator}
Let $G$ be a $d$-regular graph on $n$ units, let $P$ be Bernoulli assignment
with $p=1/(d+2)$, and let $a$ be the DTE coefficients.  If
$n\geq6(d+2)^3$, then
\[
 \Rall(P;f,a)\geq\frac{d+2}{8n}.
\]

\end{proposition}

\begin{proof}
Let $\pi=p(1-p)^d$, where $p=1/(d+2)$, be the common observation probability of the treated-arm DTE outcomes $Y_i(1,0_{N_i})$. We group approximately $1/\pi$ such outcomes with independent observation indicators. Each group then has a joint-missingness probability bounded
away from zero, allowing us to apply the disjoint-perturbation lower bound.
The treated-arm potential outcome $Y_i(1,0_{N_i})$ is observed exactly when
$Z_i=1$ and $Z_j=0$ for every $j\in N_i$. Since
$1/\pi=(d+2)\{1+1/(d+1)\}^d$ and $1\leq\{1+1/(d+1)\}^d\leq e$,
\[
 d+2\leq1/\pi\leq e(d+2).
\]
Let $k=\lfloor1/\pi\rfloor$, so that $d+2\leq k\leq e(d+2)$ and $k\pi\leq1$.

The square graph $G^2$, which joins units at distance at most two, has
maximum degree at most $d+d(d-1)=d^2$, so a greedy coloring partitions the
units into at most $d^2+1$ classes.  Two units in the same class are at
distance at least three, so their closed neighborhoods are disjoint.  Split
each color class into disjoint groups of exactly $k$ units, and denote all the resulting groups by $A_1,\ldots,A_m$. Fewer than $k$ units in each color class are not included in these groups. Since the groups contain $mk$ units, the number of units not included in the groups satisfies
\[
 n-mk\leq(d^2+1)(k-1)\leq e(d+2)^3\leq\frac n2.
\]
Thus $mk\geq n/2$.

For each group $A_b$, define the corresponding set of treated-arm potential-outcome indices by
\[
 B_b=\{(i,e_i^+):i\in A_b\},\qquad e_i^+=(1,0_{N_i}).
\]
Within $A_b$, the closed neighborhoods $N_i^+$ are pairwise disjoint. Hence the observation indicators $\{I_{i,e_i^+}:i\in A_b\}$ are mutually independent under Bernoulli assignment. Independence between different groups is not required. Thus
\[
 \rho_P(B_b)=(1-\pi)^k\geq(1-\pi)^{1/\pi}\geq\frac14.
\]
The last inequality holds because $\pi\leq1/(d+2)\leq1/2$ and
$x\mapsto(1-x)^{1/x}$ is decreasing on $(0,1)$ with value $1/4$ at $x=1/2$.
The sets $B_1,\ldots,B_m$ are nonempty and pairwise disjoint, so they form a perturbation packing. Each treated-arm potential outcome has DTE coefficient $1/n$, so with all
perturbation signs equal to $+1$ each set has $c_b=k/n$.
The lower bound from disjoint perturbations (Theorem~\ref{thm:fixed-joint-missingness}(i)) gives
\[
 \Rall(P;f,a)\geq\sum_{b=1}^m\Bigl(\frac kn\Bigr)^2\frac14
 =\frac{(mk)k}{4n^2}\geq\frac{k}{8n}\geq\frac{d+2}{8n}.
\]
\end{proof}
Combining Proposition~\ref{prop:dte-regular-all-estimator} with the upper
bound above gives
\[
 \Rall(P;f,a)\asymp\Runb(P;f,a)\asymp V_{\mathrm{HT}}(P)
 \asymp(d+1)/n.
\]

\subsubsection{Cluster randomization}

We verify Example~4 of Section~\ref{sec:fixed-design-ht} using
Theorem~\ref{thm:fixed-ht-comparison}.  The $n=\xi m$ units form $\xi$
clusters of size $m$; each cluster is treated with probability $p$,
independently across clusters, and $f_i(Z)=Z_i$.  Every relevant potential
outcome of the ATE has coefficient $\pm1/n$, and the exposure probabilities
are $p$ and $1-p$, so $\pi_{P,k}\geq c=\min(p,1-p)$ and $L_a=2$.

Within a cluster, all units have the same treatment, so the treated-arm
indicators of two units $i\ne j$ in the same cluster coincide and
$\phi_{P,ij}(1,1)=p/p^2-1=p^{-1}-1$; similarly
$\phi_{P,ij}(0,0)=(1-p)^{-1}-1$.  Indicators in different clusters are
independent.  Hence
\[
 \mathrm{OD}_a^+(P)=\xi\,\frac{m(m-1)}{n^2}\Bigl\{\frac{1-p}{p}+\frac{p}{1-p}\Bigr\}
 \leq\frac1\xi\Bigl\{\frac{1-p}{p}+\frac{p}{1-p}\Bigr\},
\]
while $D_a(P)=n^{-1}\{(1-p)/p+p/(1-p)\}$, so $\mathrm{OD}_a^+(P)/D_a(P)$ grows
like $m$.

For the packing, take one perturbation set for each cluster $g$ and arm $t$,
$B_{g,t}=\{(i,t):i\in g\}$.  These $2\xi$ sets are disjoint, each has
$c_b=\sum_{k\in B_{g,t}}|a_k|=m/n=1/\xi$, and $B_{g,1}$ and $B_{g,0}$ are
missing exactly when the cluster is in control and treated, respectively.
Therefore
\[
 L_{\cB}(P)=\xi\Bigl(\frac1\xi\Bigr)^2\{(1-p)+p\}=\frac1\xi,
\]
and equation~\eqref{eq:packing-matches-positive-dependence} holds with
$C_0=(1-p)/p+p/(1-p)$.  When $p$ is bounded away from zero and one,
Theorem~\ref{thm:fixed-ht-comparison} gives
$\Rall(P;f,a)\asymp\Runb(P;f,a)\asymp V_{\mathrm{HT}}(P)\asymp1/\xi$.

\subsubsection{Randomized saturation for a policy contrast}

We verify Example~5 of Section~\ref{sec:fixed-design-ht} through
Theorem~\ref{thm:fixed-ht-comparison}. Each of the $\xi$ clusters of size $m$
independently draws $Q_g\in\{q,q'\}$ and then treats a uniformly chosen set of $mQ_g$ units. Assume that $q\ne q'$ are bounded away from zero and one uniformly in $n$, and that $mq$ and $mq'$ are integers between $1$ and $m-1$.
For $u\in\{q,q'\}$, write $r_u=P(Q_g=u)$, with both probabilities bounded away from zero, and define $s_1(u)=u$ and $s_0(u)=1-u$.
The exposure is $f_i(Z)=(Z_i,Q_{g(i)})$. The estimand is the policy effect
$\mu_{\Gamma_q}(Y)-\mu_{\Gamma_{q'}}(Y)$ of equation~\eqref{eq:policy-mean}, where $\Gamma_q$ treats exactly $mq$ units in every cluster, so that
$\mu_{\Gamma_q}(Y)=n^{-1}\sum_i\{qY_i(1,q)+(1-q)Y_i(0,q)\}$.
The coefficients are
\[
 a_{i,(t,q)}=\frac{s_t(q)}n,\qquad
 a_{i,(t,q')}=-\frac{s_t(q')}n,\qquad t\in\{0,1\}.
\]
Thus $L_a=4$ and every coefficient satisfies $|a_k|\leq1/n$. For $u\in\{q,q'\}$, the exposure $(t,u)$ has probability
\[
 \pi_{P,i}(t,u)=P\{I_{i,(t,u)}=1\}=r_us_t(u),
\]
which is bounded away from zero.

For distinct units $i,j$ in the same cluster, conditional on $Q_g=u$, exactly $ms_t(u)$ units receive exposure $(t,u)$. Uniform selection therefore gives
\begin{align*}
 P\{I_{i,(t,u)}=1,I_{j,(t,u)}=1\}
 &=r_u\frac{ms_t(u)}m\frac{ms_t(u)-1}{m-1}\\
 &=r_us_t(u)\frac{ms_t(u)-1}{m-1}.
\end{align*}
Consequently,
\[
 \phi_{P,ij}\bigl((t,u),(t,u)\bigr)
 =\frac{ms_t(u)-1}{(m-1)r_us_t(u)}-1
 \leq\frac1{r_u}-1,
\]
because $ms_t(u)-1\leq(m-1)s_t(u)$. Observation indicators for units in different clusters are independent. Summing over the four exposures and the $\xi m(m-1)$ ordered pairs within clusters gives
\[
 \mathrm{OD}_a^+(P)
 \leq4\xi\frac{m^2}{n^2}
 \max_{u\in\{q,q'\}}\left(\frac1{r_u}-1\right)
 =\frac4\xi\max_{u\in\{q,q'\}}\left(\frac1{r_u}-1\right).
\]
Also, $D_a(P)=O(1/n)$ because every exposure probability is bounded away from zero.

We now prove the lower bound using a perturbation packing. Take one perturbation set for each cluster $g$, saturation $u\in\{q,q'\}$, and arm $t\in\{0,1\}$:
\[
 B_{g,t,u}=\{(i,(t,u)):i\in g\}.
\]
These sets are pairwise disjoint. Within each set the coefficients have a common sign; choosing perturbation signs to agree with the coefficients gives block weight $c_b=ms_t(u)/n$. Since $ms_t(u)\geq1$, at least one unit in $g$ has exposure $(t,u)$ whenever $Q_g=u$. Hence the entire block is missing exactly when $Q_g\ne u$, giving $\rho_P(B_{g,t,u})=1-r_u$. Using $s_1(u)^2+s_0(u)^2\geq1/2$ and $\sum_{u\in\{q,q'\}}(1-r_u)=1$,
\[
 L_{\cB}(P)
 =\xi\frac{m^2}{n^2}
 \sum_{u\in\{q,q'\}}(1-r_u)\sum_{t\in\{0,1\}}s_t(u)^2
 \geq\frac1{2\xi}.
\]
Thus equation~\eqref{eq:packing-matches-positive-dependence} holds with
$C_0=8\max_{u\in\{q,q'\}}(r_u^{-1}-1)$, and Theorem~\ref{thm:fixed-ht-comparison} gives
$\Rall(P;f,a)\asymp\Runb(P;f,a)\asymp V_{\mathrm{HT}}(P)\asymp1/\xi$.
For comparison, fixing the saturation at $q$ in every cluster gives rate
$1/n$ for the within-saturation contrast between $(1,q)$ and $(0,q)$,
by the calculation for Example~1. This is a different estimand from the
policy contrast between saturation levels considered here.

\subsection{Improving exposure HT by conditioning on randomized saturation}
\label{subsec:app-fixed-ht-saturation-conditioning}

In this section we present an example design with a different exposure mapping and estimand where the exposure HT is not minimax optimal. Instead, a conditional HT estimator that conditions on the realized saturation improves estimation.

Recall the design in Section~\ref{subsec:fixed-ht-saturation-conditioning}: the $n=\xi m$ units form $\xi$ administrative clusters, each containing $m/2$ two-person households. Interference is confined to households, and the estimand contrasts exposures $(1,1)$ and $(0,0)$. Each cluster independently draws a saturation level $Q_g\in\{q_L,q_H\}$ with $P(Q_g=q_H)=r$, then treats exactly $mQ_g$ units uniformly.  Closely related multilevel designs, with households nested in administrative clusters and saturation randomized at the cluster level, are used in field experiments on voter mobilization \citep{SinclairMcConnellGreen2012} and cash transfers \citep{EggerEtAl2022}; see also \citet{HudgensHalloran2008}, \citet{BairdEtAl2018}, and \citet{VazquezBare2023}.

Write $I_i(e)=\1\{f_i(Z)=e\}$ and
$\pi_i(e\mid q)=P\{f_i(Z)=e\mid Q_{g(i)}=q\}$. The conditional estimator is
\[
 \widehat\tau_{\mathrm{cond}}
 =\frac1n\sum_i\left\{
 \frac{I_i(1,1)Y_i(1,1)}{\pi_i(1,1\mid Q_{g(i)})}
 -\frac{I_i(0,0)Y_i(0,0)}{\pi_i(0,0\mid Q_{g(i)})}\right\}.
\]
Exposure HT uses $\pi_i(e)=(1-r)\pi_i(e\mid q_L)+r\pi_i(e\mid q_H)$
in the denominators. We compute the conditional probabilities in the
proof below.

\begin{lemma}[Conditional unbiasedness and variance]
\label{lem:app-ht-conditional-variance}
Fix a potential outcome table $Y$, and let $Q$ be the vector $(Q_1,\ldots,Q_\xi)$ of realized cluster saturations, which is a function of the assignment $Z$ and is therefore observed. If $T$ has finite second moment and
$\E_P(T\mid Q,Y)=\tau_a(Y)$ almost surely, then $T$ is design-unbiased and
\[
\Var_P(T\mid Y)=\E_P\{\Var_P(T\mid Q,Y)\mid Y\}.
\]
For any estimator $H$ with finite second moment,
\[
\Var_P(H\mid Y)\geq
\Var_P\{\E_P(H\mid Q,Y)\mid Y\}.
\]
\end{lemma}
\begin{proof}

Iterated expectation gives
$\E_P(T\mid Y)=\E_P\{\E_P(T\mid Q,Y)\mid Y\}=\tau_a(Y)$.  The law of total
variance gives, for any $H$ with finite second moment,
\[
 \Var_P(H\mid Y)=\E_P\{\Var_P(H\mid Q,Y)\mid Y\}+\Var_P\{\E_P(H\mid Q,Y)\mid Y\}.
\]
For $H=T$ the second term vanishes because $\E_P(T\mid Q,Y)$ is the constant
$\tau_a(Y)$.  For general $H$, dropping the nonnegative first term gives the
inequality.

\end{proof}

\begin{resultrestatement}{Proposition~\ref{prop:fixed-ht-saturation-conditioning}}{Conditioning on randomized saturation}
Under the design above, assume that $r,q_L,q_H$ are bounded away from zero and one, that $q_L<q_H$ with $q_H-q_L$ bounded away from zero, and that $mq$ is an integer satisfying $2\leq mq\leq m-2$ for each $q\in\{q_L,q_H\}$, with all bounds uniform in the population size. Then
\[
 \Rall(P;f,a)\asymp\Runb(P;f,a)\asymp\frac1n,
 \qquad
 V_{\mathrm{HT}}(P)\asymp\frac1\xi.
\]
\end{resultrestatement}

\begin{proof}[Proof of Proposition~\ref{prop:fixed-ht-saturation-conditioning}]
All constants below may depend on the fixed lower bounds on $r$, $1-r$,
$q_L$, $1-q_L$, $q_H$, $1-q_H$, and $q_H-q_L$, but not on $\xi$ or $m$.
We first show that the conditional estimator is unbiased and has
worst-case variance $O(1/n)$. Conditional on $Q_g=q$, the target exposure probabilities are
\[
 \pi_i(1,1\mid Q_g=q)=\frac{(mq)_2}{(m)_2},
 \qquad
 \pi_i(0,0\mid Q_g=q)=\frac{(m-mq)_2}{(m)_2},
\]
where $(x)_k=x(x-1)\cdots(x-k+1)$ denotes the falling factorial.
Under the count conditions in the proposition, $mq\geq2$ and $m-mq\geq2$ for both
saturations. Since $(t)_2\geq t^2/2$ for integers $t\geq2$,
the conditional exposure probabilities are at least $q^2/2$ and
$(1-q)^2/2$, respectively, and hence are uniformly bounded away from zero.

Therefore, for every fixed
potential outcome table,
\[
 \E_P(\widehat\tau_{\mathrm{cond}}\mid Q,Y)
 =\frac1n\sum_i\{Y_i(1,1)-Y_i(0,0)\}=\tau(Y).
\]
Iterated expectation therefore gives $\E_P(\widehat\tau_{\mathrm{cond}}\mid Y)=\tau(Y)$.

We next bound the conditional variance.  Write
$\widehat\tau_{\mathrm{cond}}=n^{-1}\sum_iX_i$.  The variables $X_i$ are
uniformly bounded because the outcomes and inverse conditional exposure probabilities are bounded. Two units in the same household contribute a bounded
covariance, and there are $O(n)$ such ordered pairs.

Conditional on $Q$, the design within cluster $g$ is a uniform
$mq$-subset of the $m$ units. Put $t=mq$ and $u=m-t$. For two distinct
households, the joint probabilities of their target exposure events are
\[
 \frac{(t)_4}{(m)_4},\qquad
 \frac{(u)_4}{(m)_4},\qquad
 \frac{(t)_2(u)_2}{(m)_4},
\]
for two treated households, two control households, and one of each,
respectively. Their marginal products are
$\{(t)_2/(m)_2\}^2$, $\{(u)_2/(m)_2\}^2$, and
$(t)_2(u)_2/(m)_2^2$. Since $q$ is bounded away from zero and one,
\begin{align*}
 \frac{(t)_4}{(m)_4}-\left\{\frac{(t)_2}{(m)_2}\right\}^2
   &=O(m^{-1}),\\
 \frac{(u)_4}{(m)_4}-\left\{\frac{(u)_2}{(m)_2}\right\}^2
   &=O(m^{-1}),\\
 \frac{(t)_2(u)_2}{(m)_4}
   -\frac{(t)_2(u)_2}{(m)_2^2}
   &=O(m^{-1}).
\end{align*}
These are the covariances of the two treated-household indicators, the two
control-household indicators, and a treated- and control-household indicator.
Because the potential outcomes and inverse exposure probabilities are
uniformly bounded, expanding $\Cov_P(X_i,X_j\mid Q,Y)$ in these three terms
gives
\[
 \bigl|\Cov_P(X_i,X_j\mid Q,Y)\bigr|\lesssim m^{-1}
\]
for units in distinct households in the same administrative cluster.
There are $O(\xi m^2)$ ordered pairs of units in distinct households
within the same administrative cluster.  Units in different
administrative clusters are conditionally independent.  It follows that
\[
 \sup_{Y\in\cY(f)}
 \Var_P(\widehat\tau_{\mathrm{cond}}\mid Q,Y)
 \lesssim \frac{n+\xi m^2/m}{n^2}
 \asymp\frac1n.
\]
Since $\E_P(\widehat\tau_{\mathrm{cond}}\mid Q,Y)=\tau(Y)$,
Lemma~\ref{lem:app-ht-conditional-variance} gives
\[
\Var_P(\widehat\tau_{\mathrm{cond}}\mid Y)
=\E_P\{\Var_P(\widehat\tau_{\mathrm{cond}}\mid Q,Y)\mid Y\}
\lesssim1/n.
\]

We next obtain a matching lower bound for all estimators from singleton
perturbations. The two target exposures of each unit are mutually exclusive, so $\pi_i(1,1)+\pi_i(0,0)\leq1$.  Every
relevant potential outcome has coefficient $\pm1/n$, so the all-estimator
singleton lower bound in equation~\eqref{eq:singleton-all} gives
\[
 \Rall(P;f,a)
 \geq\frac1{n^2}\sum_{i=1}^n
 \{2-\pi_i(1,1)-\pi_i(0,0)\}
 \geq\frac1n.
\]

Together with the conditional estimator, this proves
$\Rall(P;f,a)\asymp\Runb(P;f,a)\asymp1/n$.

It remains to show that marginal exposure HT has worst-case variance
of order $1/\xi$. For the upper bound, write
$\widehat\tau_{\mathrm{HT},P}=\sum_{g=1}^\xi W_g$, where $W_g$ collects the
terms of the units in cluster $g$.  The marginal exposure probabilities
$\pi_i(1,1)$ and $\pi_i(0,0)$ are averages of the conditional ones and hence
bounded away from zero, so each of the $m$ terms in $W_g$ is at most a
constant times $1/n$ in absolute value, and $|W_g|\leq C/\xi$.  The clusters
are assigned independently, so the $W_g$ are independent and
\[
 \Var_P(\widehat\tau_{\mathrm{HT},P}\mid Y)
 =\sum_{g=1}^\xi\Var_P(W_g\mid Y)
 \leq\xi\frac{C^2}{\xi^2}=\frac{C^2}\xi.
\]
For the matching lower bound on the worst-case HT variance, choose
a particular table: set $Y_i(1,1)=1$ and $Y_i(0,0)=0$ for every unit, and set the other potential
outcomes to zero.  Let
\[
 h_m(q)=\frac{\pi_i(1,1\mid Q_g=q)}{\pi_i(1,1)}.
\]
Conditional on the realized saturations,
\[
 \E_P(\widehat\tau_{\mathrm{HT},P}\mid Q,Y)
 =\frac1\xi\sum_{g=1}^\xi h_m(Q_g).
\]

As $m\to\infty$, $(mq)_2/(m)_2\to q^2$, so
\[
 h_m(q_H)-h_m(q_L)\to\frac{q_H^2-q_L^2}{(1-r)q_L^2+rq_H^2},
\]
which is bounded away from zero.  Hence there is an $m_0$ such that
$|h_m(q_H)-h_m(q_L)|$ is bounded away from zero for $m\geq m_0$.  The
saturations $Q_g$ are independent, and
$\Var\{h_m(Q_g)\}=r(1-r)\{h_m(q_H)-h_m(q_L)\}^2$, so the second part of
Lemma~\ref{lem:app-ht-conditional-variance} gives, for $m\geq m_0$,
\[
 \Var_P(\widehat\tau_{\mathrm{HT},P}\mid Y)
 \geq\Var_P\Bigl\{\frac1\xi\sum_{g=1}^\xi h_m(Q_g)\Bigm|Y\Bigr\}
 =\frac{\Var\{h_m(Q_g)\}}\xi\gtrsim\frac1\xi.
\]
For $m<m_0$, the bound $V_{\mathrm{HT}}(P)\geq\Rall(P;f,a)\geq1/n\geq1/(m_0\xi)$
suffices.

Thus $V_{\mathrm{HT}}(P)\asymp1/\xi$, completing the proof.
\end{proof}

\subsubsection{Conditioning on randomized saturation without interference}
\label{subsec:app-fixed-ht-no-interference-saturation}
Use the randomized-saturation design from Proposition~\ref{prop:fixed-ht-saturation-conditioning}, but now let $f_i(Z)=Z_i$ and consider the finite-population average treatment
effect
\[
 \tau(Y)=\frac1n\sum_i\{Y_i(1)-Y_i(0)\}.
\]
Write $\bar q=\E(Q_g)$.  Exposure HT uses the unconditional probabilities
$\bar q$ and $1-\bar q$.  The estimator that conditions on the observed
saturations is
\[
 \widehat\tau_{\mathrm{cond}}
 =\frac1n\sum_i\left\{
 \frac{Z_iY_i(1)}{Q_{g(i)}}
 -\frac{(1-Z_i)Y_i(0)}{1-Q_{g(i)}}\right\}.
\]

Conditional on $Q$, this estimator is unbiased.  Its variance is bounded
as in the proof of Proposition~\ref{prop:fixed-ht-saturation-conditioning},
with households replaced by single units.  Conditional on $Q$, the treated
units of cluster $g$ form a uniform $mQ_g$-subset, so for distinct units in
the same cluster $\Cov_P(Z_i,Z_j\mid Q)=-Q_g(1-Q_g)/(m-1)$, of order $1/m$,
while units in different clusters are conditionally independent.  The $n$
individual variance terms and the $O(\xi m^2)$ pairs within clusters
therefore give

\[
 \sup_{Y\in\cY(f)}
 \Var_P(\widehat\tau_{\mathrm{cond}}\mid Q,Y)
 \lesssim \frac{n+\xi m^2/m}{n^2}
 \asymp\frac1n.
\]
The conditional mean is $\tau(Y)$ for every $Q$, so
Lemma~\ref{lem:app-ht-conditional-variance} gives the same variance bound

without conditioning on $Q$.  The treatment and control exposures of each
unit are mutually exclusive, so the all-estimator singleton lower bound gives
$\Rall(P;f,a)\geq1/n$, as in the proof of
Proposition~\ref{prop:fixed-ht-saturation-conditioning}.  Hence both
fixed-design minimax risks are of order $1/n$.

Exposure HT has worst-case variance of order $1/\xi$.  For the
table $Y_i(1)=1$ and $Y_i(0)=0$,
\[
 \E_P(\widehat\tau_{\mathrm{HT},P}\mid Q,Y)
 =\frac1\xi\sum_{g=1}^\xi\frac{Q_g}{\bar q}.
\]

The second part of Lemma~\ref{lem:app-ht-conditional-variance} gives
\[
 \Var_P(\widehat\tau_{\mathrm{HT},P}\mid Y)
 \geq\frac{\Var(Q_g)}{\xi\bar q^2}
 =\frac{r(1-r)(q_H-q_L)^2}{\xi\bar q^2}\gtrsim\frac1\xi,
\]
and the cluster decomposition used in the proof of
Proposition~\ref{prop:fixed-ht-saturation-conditioning} gives the matching
upper bound.
  Thus, when $m=n/\xi$ grows, the worst-case variance of exposure HT
exceeds the fixed-design minimax risk by a factor of order $m$.

\section{Supporting results for Section~\ref{sec:no-interference}}
\label{app:no-interference-proofs}
We prove the all-desgin minimax risk for the average treatment effect under no interference, in Theorem~\ref{thm:ate-joint-benchmark} using the weighted sampling hypergraph. We then establish the exact
finite-sample optimality of independent Bernoulli assignment with treatment
probability $1/2$ for the all-estimator problem and prove the constrained-design
benchmark in Theorem~\ref{thm:ate-fixed-count}.

\subsection{No-interference ATE}
\label{subsec:app-ate-benchmark}

The sharp lower bound for all estimators follows from a Bayes-risk
calculation for a shifted binomial count. We first give that calculation
and then show how the same count arises on every sampling hyperedge.

\begin{lemma}[Bayes risk for a shifted binomial count]
\label{lem:ate-shifted-binomial-bayes}
Fix $M>0$ and $L>M+1$. Let $\Theta$ be uniform on $[L-M,L+M]$
and put $m=\lfloor\sqrt n\,\Theta\rfloor$. For sufficiently large $n$,
suppose we observe $C_n=m+B_n$, where, conditional on $\Theta$,
$B_n\sim\operatorname{Binomial}(n-m,1/2)$. The conditional law of
$B_n$ depends on $\Theta$ only through $m$.
For every $0<A<M$,
\[
 \liminf_{n\to\infty}\inf_\delta
 \mathbb E\bigl[(\delta(C_n)-m/\sqrt n)^2\bigr]
 \geq (1-A/M)\int_{-A}^{A} z^2\phi(z)\,dz,
\]
where $\phi$ is the standard normal density.
\end{lemma}

\begin{proof}
We first show that adding uniform noise to the rescaled count preserves
the minimum Bayes risk for estimating $m/\sqrt n$. Define
\[
 V_n=\frac{2C_n-n}{\sqrt n}+\eta_n,\qquad
 \eta_n\sim\operatorname{Uniform}[-1/\sqrt n,1/\sqrt n),
 \qquad \eta_n\perp(\Theta,C_n).
\]
Rounding $V_n$ to the nearest point of the rescaled lattice recovers
the count. In terms of the original count,
\[
 C_n=\left\lfloor\frac{\sqrt n\,V_n+n}{2}+\frac12\right\rfloor
 \quad\text{almost surely}.
\]
Thus any estimator based on $C_n$ can be evaluated from $V_n$.
Conversely, for any estimator $\delta(V_n)$ of finite risk, define
$\bar\delta(C_n)=\mathbb E[\delta(V_n)\mid C_n]$. Independence of the
added noise gives $\mathbb E[\delta(V_n)\mid C_n,\Theta]=\bar\delta(C_n)$.
Conditional Jensen's inequality then gives
\[
 \mathbb E[(\bar\delta(C_n)-m/\sqrt n)^2]
 \leq \mathbb E[(\delta(V_n)-m/\sqrt n)^2].
\]
Taking infima in both directions shows that the two observations give
the same minimum Bayes risk.

Next we show that the conditional law of $V_n$ given $\Theta$ converges
in total variation to $N(\Theta,1)$, uniformly over the prior support
$\mathcal I=[L-M,L+M]$. For $c=0,\ldots,n$, define
\[
 s_{n,c}=\frac{2c-n}{\sqrt n},\qquad
 J_{n,c}=[s_{n,c}-1/\sqrt n,s_{n,c}+1/\sqrt n).
\]
For $\theta\in\mathcal I$, put $m=\lfloor\sqrt n\,\theta\rfloor$,
$b=n-m$, and $j=c-m$. On $J_{n,c}$ the conditional density of $V_n$ is
$f_{n,\theta}(v)=(\sqrt n/2)2^{-b}\binom bj$, with zero density when
$j\notin\{0,\ldots,b\}$. On intervals meeting any fixed bounded interval,
$j-b/2=O(\sqrt b)$ uniformly in $\theta$, since
$2j-b=\sqrt n(s_{n,c}-m/\sqrt n)$ and $b/n\to1$ uniformly.
In this range Stirling's formula gives
\[
 2^{-b}\binom bj
 =\frac{2}{\sqrt b}\phi\!\left(\frac{2j-b}{\sqrt b}\right)
 \{1+O(b^{-1})\}.
\]
For $v\in J_{n,c}$, we have $|v-s_{n,c}|\leq1/\sqrt n$ and
$|m/\sqrt n-\theta|\leq1/\sqrt n$. Hence, for each fixed $H>0$,
\[
 \sup_{\theta\in\mathcal I,\ |v|\leq H}
 |f_{n,\theta}(v)-\phi(v-\theta)|\longrightarrow0.
\]
The conditional mean and variance of $V_n$ are $m/\sqrt n$ and
$1-m/n+1/(3n)$, respectively. Thus for
$X\mid\Theta=\theta\sim N(\theta,1)$,
\[
 \mathbb E[V_n^2\mid\Theta=\theta]
 =\frac{m^2}{n}+1-\frac mn+\frac{1}{3n}\leq D,
 \qquad
 \mathbb E[X^2\mid\Theta=\theta]=1+\theta^2\leq D,
 \qquad D=(L+M)^2+2.
\]
Integrating on $[-H,H]$ and bounding both tail probabilities by their
second moments gives
\[
 \sup_{\theta\in\mathcal I}\int_{\mathbb R}
 |f_{n,\theta}(v)-\phi(v-\theta)|\,dv
 \leq 2H\sup_{\theta\in\mathcal I,\ |v|\leq H}
 |f_{n,\theta}(v)-\phi(v-\theta)|+\frac{2D}{H^2}.
\]
Taking $n\to\infty$ for fixed $H$ and then $H\to\infty$ proves the
uniform total variation convergence.

We now compare the Bayes risks. We clip estimates to
$[L-M-1,L+M+1]$, which contains both targets. Write $W=2M+2$
for the length of this interval and $\theta_n=m/\sqrt n$. For any
clipped estimator $\delta$, the conditional risks satisfy
\begin{align*}
 &\left|\mathbb E[(\delta(V_n)-\theta_n)^2\mid\Theta=\theta]
       -\mathbb E[(\delta(X)-\theta)^2\mid\Theta=\theta]\right|\\
 &\qquad\leq W^2\int_{\mathbb R}|f_{n,\theta}(v)-\phi(v-\theta)|\,dv
            +2W|\theta_n-\theta|.
\end{align*}
The first term bounds the change in the observation law. The second
comes from factoring the difference of the squared losses at the two
targets. The bound is $o(1)$ uniformly in $\theta$ and $\delta$.
Integrating over the prior and taking infima shows that the minimum
Bayes risks differ by $o(1)$.

Finally, we lower-bound the Bayes risk in the Gaussian experiment.
For $0<A<M$, put $\mathcal I_A=[L-M+A,L+M-A]$. For $x\in\mathcal I_A$,
the interval $[x-A,x+A]$ lies inside $\mathcal I$. Restricting the
prior-averaged loss to these observations and parameters gives
\begin{align*}
 \mathbb E[(\delta(X)-\Theta)^2]
 &=\frac{1}{2M}\int_{\mathbb R}\int_{\mathcal I}
   (\delta(x)-\theta)^2\phi(x-\theta)\,d\theta\,dx\\
 &\geq \frac{1}{2M}\int_{\mathcal I_A}\int_{-A}^{A}
   (\delta(x)-x-z)^2\phi(z)\,dz\,dx\\
 &\geq (1-A/M)\int_{-A}^{A} z^2\phi(z)\,dz.
\end{align*}
The second line sets $z=\theta-x$. For the last inequality,
$\int_{-A}^{A}z\phi(z)\,dz=0$, so the cross term vanishes, and
$|\mathcal I_A|=2(M-A)$. Taking the infimum over $\delta$ and using the
Bayes-risk comparison proves the lemma.
\end{proof}

\begin{proof}[Proof of Theorem~\ref{thm:ate-joint-benchmark}]
In the no-interference experiment, the expanded finite population is the set
$K=[n]\times\{0,1\}$, and the ATE coefficients are $a_{i,1}=1/n$ and
$a_{i,0}=-1/n$. Each sampling hyperedge has the form
\[
 E=\{(i,d_i):i\in[n]\},\qquad d_i\in\{0,1\},
\]
where $d_i$ is the treatment level assigned to unit $i$.
For a design $P$, write $Q=Q_{P,f}$ for its induced sampling design.
By Proposition~\ref{thm:three-problem-characterization}, we can work with
the observed data $(\mathsf E,\mathsf Y_{\mathsf E})$, where $\mathsf E\sim Q$.

\textit{Sharp all-estimator lower bound.}
We construct a prior $\nu$ on the space $\cY(f)$ of admissible potential-outcome tables and compute its Bayes risk. Fix $M>0$ and $L>M+1$. Draw
$\Theta\sim\operatorname{Uniform}[L-M,L+M]$ and set
$m=\lfloor\sqrt n\,\Theta\rfloor$. Thus, for sufficiently large $n$, throughout the prior interval, we have $1\le m\le n$. Conditional on $m$, draw $S$ uniformly from the $m$-subsets of $[n]$. Given $S$, draw independent signs
$\varepsilon_i\in\{-1,1\}$ with probability $1/2$ for $i\in[n]\setminus S$ and set
\[
 (\mathsf Y_i(0),\mathsf Y_i(1))=
 \begin{cases}
 (-1,1), & i\in S,\\
 (\varepsilon_i,\varepsilon_i), & i\in[n]\setminus S.
 \end{cases}
\]
Let $\nu$ be the law of the resulting table $\mathsf Y$. It has finite
support because each coordinate of $\mathsf Y$ is in $\{-1,1\}$.
Draw the assignment independently of all the variables used to
construct the table.
Each unit in $S$ contributes $2/n$ to the ATE, and each unit outside
$S$ contributes zero. Thus
\[
 \tau_a(\mathsf Y)=\frac1n\sum_{i=1}^n
 \{\mathsf Y_i(1)-\mathsf Y_i(0)\}=\frac{2|S|}{n}=\frac{2m}{n}.
\] The conditional construction also allows
$m=0$. In that case, $S=\varnothing$ and each unit has
$(\mathsf Y_i(0),\mathsf Y_i(1))=(\varepsilon_i,\varepsilon_i)$,
with signs drawn independently across units.

Fix an edge $E\in\mathcal E$ and consider the outcomes it reveals. Let
\[
 S_+(E,\mathsf Y_E)=\{i\in[n]:a_{i,d_i}\mathsf Y_i(d_i)>0\},
 \qquad C_n(E)=|S_+(E,\mathsf Y_E)|.
\]
Every $i\in S$ belongs to $S_+(E,\mathsf Y_E)$ for either value of $d_i$.
For $i\in[n]\setminus S$, write
$I_i=\mathbf 1\{(2d_i-1)\varepsilon_i=1\}$ for the indicator that unit
$i$ contributes positively on $E$. For this fixed edge, conditional on $S$, these
indicators are independent Bernoulli$(1/2)$ variables. Since
$C_n(E)=m+\sum_{i\in[n]\setminus S}I_i$, for every fixed $E$,
\[
 C_n(E)\mid m\sim m+\operatorname{Binomial}(n-m,1/2).
\]
For fixed $E$ and $m$, the set $S_+(E,\mathsf Y_E)$ is exchangeable across unit
labels because $S$ is uniform and the membership indicators outside
$S$ are independent and fair. Conditional also on $C_n(E)=c$, this set
is uniform among the $\binom nc$ subsets of size $c$. In particular,
for each $T\subseteq[n]$ with $|T|=c$,
\[
 \Pr_\nu\{S_+(E,\mathsf Y_E)=T\mid C_n(E)=c,m\}
 =\binom nc^{-1}
\]
whenever the conditioning event has positive probability. For a fixed $E$, the set determines $\mathsf Y_E$, since membership
specifies the sign of each observed contribution and every observed
outcome is in $\{-1,1\}$.

For each fixed $E$, the posterior law of the target given $\mathsf Y_E$
depends only on $C_n(E)$, and the joint law of $(m,C_n(E))$ is the same
for every edge. Write $C_n=C_n(\mathsf E)$ for the observed count. By Lemma~\ref{lem:posterior-risk},
\[
 \Rall(P;f,a)
 \geq \sum_{E\in\mathcal E}Q(E)\,
   \mathbb E_\nu\!\left[\operatorname{Var}_\nu\{\tau_a(\mathsf Y)
     \mid \mathsf Y_E\}\right]
 =\mathbb E\!\left[\operatorname{Var}(2m/n\mid C_n)\right].
\]
The equality uses the fact that every hyperedge contributes the same expected
posterior variance and the weights sum to one. Write $r_n(M,L)$ for the
Bayes risk in Lemma~\ref{lem:ate-shifted-binomial-bayes}. The conditional
mean minimizes squared-error Bayes risk, by
Lemma~\ref{lem:posterior-risk}, so
\[
 \mathbb E[\operatorname{Var}(2m/n\mid C_n)]
 =\frac4n\,\mathbb E[\operatorname{Var}(m/\sqrt n\mid C_n)]
 =\frac4n\,r_n(M,L).
\]
The first equality uses $2m/n=(2/\sqrt n)(m/\sqrt n)$.

The prior is the same for every design. Taking the infimum over $P$
and applying Lemma~\ref{lem:ate-shifted-binomial-bayes} gives, for
every $0<A<M$ and $L>M+1$,
\[
 \liminf_{n\to\infty} n\inf_P\Rall(P;f,a)
 \geq 4\liminf_{n\to\infty}r_n(M,L)
 \geq 4(1-A/M)\int_{-A}^{A}z^2\phi(z)\,dz.
\]
The left-hand side is independent of $M$, $L$, and $A$. For each fixed
$A>0$, choose $L=M+2$ and let $M\to\infty$ in this inequality. This gives
\[
 \liminf_{n\to\infty} n\inf_P\Rall(P;f,a)
 \geq 4\int_{-A}^{A}z^2\phi(z)\,dz
 \qquad\text{for every } A>0.
\]
Letting $A\to\infty$ then yields
\[
 \liminf_{n\to\infty} n\inf_P\Rall(P;f,a)
 \geq 4\int_{\mathbb R}z^2\phi(z)\,dz=4.
\]

\textit{Unbiased lower bound.} Pair the two mutually exclusive singleton
potential outcomes of each unit. Their coefficients are $1/n$ and $-1/n$.
Let $q_i=P(Z_i=1)$. If some $q_i\in\{0,1\}$, no estimator is unbiased for
the ATE over the full parameter space, so $\Runb(P;f,a)=+\infty$.
Otherwise Proposition~\ref{prop:crossblocklb} gives
\[
 \Runb(P;f,a)
 \geq \frac1{n^2}\sum_{i=1}^n
 \left\{\frac{1-q_i}{q_i}+\frac{q_i}{1-q_i}+2\right\}
 =\frac1{n^2}\sum_{i=1}^n\frac{1}{q_i(1-q_i)}
 \geq\frac4n.
\]
The last inequality uses $q_i(1-q_i)\leq 1/4$ for each unit.
Taking the infimum over $P$ gives $\Runbstar(\cP(\cZ))\geq 4/n$.

\emph{Upper bounds.} Under independent Bernoulli assignment $P_{1/2}$ with treatment probability $1/2$, consider the estimator
\[
  \widehat\tau_\alpha
  =\frac\alpha n\sum_{i=1}^n
  \{Z_iY_i(1)-(1-Z_i)Y_i(0)\}.
\]
The unbiased HT estimator corresponds to $\alpha=2$. Under $P_{1/2}$,
each summand can be written as
\[
 Z_iY_i(1)-(1-Z_i)Y_i(0)
 =\tfrac12\{Y_i(1)-Y_i(0)\}+
 (Z_i-\tfrac12)\{Y_i(1)+Y_i(0)\}.
\]
The centered assignment indicators are independent and each has
variance $1/4$. Therefore,
\[
 E(\widehat\tau_\alpha\mid Y)=\frac\alpha2\tau(Y),
 \qquad
 \Var(\widehat\tau_\alpha\mid Y)
 =\frac{\alpha^2}{4n^2}\sum_i\{Y_i(1)+Y_i(0)\}^2.
\]
The squared bias plus the variance therefore gives
\begin{align*}
  E(\widehat\tau_\alpha-\tau)^2
  ={}&\left(1-\frac\alpha2\right)^2
  \left\{\frac1n\sum_i(Y_i(1)-Y_i(0))\right\}^2\\
  &+\frac{\alpha^2}{4n^2}\sum_i\{Y_i(1)+Y_i(0)\}^2 .
\end{align*}
Put $u_i=Y_i(1)-Y_i(0)$ and $v_i=Y_i(1)+Y_i(0)$.  For each unit,
\[
  |u_i|+|v_i|\leq2.
\]
The risk is convex in each pair $(u_i,v_i)$. Its maximum over $|u_i|+|v_i|\leq2$
is attained at one of $(\pm2,0)$ and $(0,\pm2)$.  Suppose $k$ units use a
vertex of the first type.  Choosing their signs to agree gives the largest
squared bias, while the remaining $n-k$ units each contribute $4$ to the
variance term.  The maximum risk for a given $k$ is
\[
 4\left(1-\frac\alpha2\right)^2\frac{k^2}{n^2}
 +\frac{\alpha^2(n-k)}{n^2},
 \qquad 0\leq k\leq n.
\]
This expression is convex in $k$, so its maximum occurs at $k=0$ or $k=n$.
Consequently,
\[
  \sup_Y E(\widehat\tau_\alpha-\tau)^2
  =
  \max\left\{
    4\left(1-\frac\alpha2\right)^2,\,
    \frac{\alpha^2}{n}
  \right\}.
\]
Equating the two terms gives
\(\alpha=2\sqrt n/(\sqrt n+1)\) and worst-case risk
\(4/(\sqrt n+1)^2\).  Hence the shrinkage upper bound gives
\[
  \limsup_{n\to\infty} n\inf_P\Rall(P;f,a)\leq4.
\]

This matches the all-estimator lower bound to first order. Setting $\alpha=2$ gives the unbiased HT estimator with worst-case risk $4/n$, which attains the unbiased lower bound exactly.

\end{proof}
\subsection{Exact optimality of independent Bernoulli assignment}
Theorem~\ref{thm:ate-joint-benchmark} gives the first-order joint minimax
risk over all estimators. The following design reduction identifies an
optimal design for this problem at each finite $n$.

\begin{lemma}[Design reduction for the no-interference ATE]
\label{lem:ate-fair-design-reduction}
For every design--estimator pair $(P,T)$, there is an estimator $T^*$ under independent Bernoulli assignment $P_{1/2}$ with treatment probability $1/2$ whose worst-case squared-error risk is no larger than that of $(P,T)$. In particular, $P_{1/2}$ minimizes the fixed-design minimax risk over all designs.
\end{lemma}

\begin{proof}
We use within-unit swaps to turn any design $P$ into $P_{1/2}$.
For a swapped unit, we also reverse both outcome signs so that
its contribution to the ATE is preserved.

\emph{Assignment law.}
Use the population index set $K=[n]\times\{0,1\}$ and the sampling
edges $O_f(z)=\{(i,z_i):i\in[n]\}$. For $g\in\{0,1\}^n$, write
$\oplus$ for addition modulo two, applied coordinatewise to assignments.
Define the transformed table $Y^g$ on the same index set by
\[
 Y_i^g(d)=(-1)^{g_i}Y_i(d\oplus g_i).
\]
For $g_i=0$ the two entries are unchanged. For $g_i=1$ they are
swapped and negated, so $Y_i^g(1)-Y_i^g(0)=Y_i(1)-Y_i(0)$ in
either case. Thus $Y^g\in\cY(f)$ and $\tau_a(Y^g)=\tau_a(Y)$.

Draw $Z\sim P$ and independently draw $G$ uniformly on $\{0,1\}^n$.
Set $W=Z\oplus G$. For each $w\in\{0,1\}^n$,
\[
 \Pr(W=w)=2^{-n}\sum_{z\in\{0,1\}^n}P(Z=z)=2^{-n}.
\]
Hence $W$ has law $P_{1/2}$.

\emph{Estimator construction.}
Fix an estimator $T$ under $P$. Given $W=w$, we observe $Y_{O_f(w)}$.
For each $z$ in the support of $P$, set $g=w\oplus z$. The outcomes
that $T$ would use on $O_f(z)$ at table $Y^g$ can be computed from
$Y_{O_f(w)}$, since
\[
 Y_i^g(z_i)=(-1)^{g_i}Y_i(w_i),\qquad i\in[n].
\]
Define the estimator under $P_{1/2}$ by
\[
 T^*(w,Y_{O_f(w)})
 =\sum_{z:\,P(Z=z)>0}P(Z=z)\,
   T\bigl(z,Y^{w\oplus z}_{O_f(z)}\bigr).
\]
Since $\Pr(Z=z,W=w)=2^{-n}P(Z=z)$, the weights in this sum are
the conditional probabilities of $Z$ given $W=w$.

\emph{Risk comparison.}
Write $r(P,T;Y)$ for the squared-error risk of $T$ at table $Y$.
For fixed $Y$, the estimator $T^*(W,Y_{O_f(W)})$ is the conditional
mean of $T(Z,Y^G_{O_f(Z)})$ given $W$. Conditional on $G=g$, we
have $Z\sim P$ and $\tau_a(Y^g)=\tau_a(Y)$, so the conditional risk
of that term is $r(P,T;Y^g)$. Conditional Jensen's inequality and
averaging over $G$ give
\[
 \begin{aligned}
 r(P_{1/2},T^*;Y)
 &\leq \mathbb E[(T(Z,Y^G_{O_f(Z)})-\tau_a(Y))^2]\\
 &=2^{-n}\sum_{g\in\{0,1\}^n}r(P,T;Y^g)
 \leq \sup_{Y'\in\cY(f)}r(P,T;Y').
 \end{aligned}
\]
Taking the supremum over $Y$ proves the comparison for each pair
$(P,T)$. Since $P_{1/2}$ is itself admissible, taking infima yields
\[
 \inf_P\Rall(P;f,a)=\Rall(P_{1/2};f,a).
\]
\end{proof}

\subsection{No-interference ATE under a constraint on the number of treated units}
\label{subsec:app-fixed-count-ate}

\begin{resultrestatement}{Theorem~\ref{thm:ate-fixed-count}}{Minimax risk of ATE under constraints on the number of treated units}
For the ATE in equation~\eqref{eq:ate-target} and every \(m\in\{1,\ldots,n-1\}\), the
following statements hold.
\begin{enumerate}[(i)]
\item The joint minimax risk over all designs with exactly $m$ treated units is
\begin{equation*}
  \Rallstar(\cP(\cZ_m))
  \asymp
  \Runbstar(\cP(\cZ_m))
  \asymp
  \frac1m+\frac1{n-m}.
\end{equation*}
\item Let $k_m=\min\{m,\lfloor n/2\rfloor\}$. The joint minimax risk over all designs with at most $m$ treated units is 
\begin{equation*}
  \Rallstar(\cP(\cZ_{\le m}))
  \asymp
  \Runbstar(\cP(\cZ_{\le m}))
  \asymp
  \frac1{k_m}+\frac1{n-k_m}
  \asymp
  \frac1m.
\end{equation*}
\end{enumerate}
The implied constants are universal and do not depend on \(m\) or \(n\).
\end{resultrestatement}

\begin{proof}[Proof of Theorem~\ref{thm:ate-fixed-count}]
Write $V_d=\{(i,d):i\in[n]\}$ for the arm indexed by $d\in\{0,1\}$.
We first give a fixed-design unbiased lower bound that applies to
both assignment spaces.

\emph{Unbiased lower bound for a fixed design.} For any design $P$, let $q_i=P(Z_i=1)$. If any $q_i$ is zero or one, no estimator is unbiased for the ATE over $\mathcal Y(f)$, so $\Runb(P;f,a)=+\infty$. Otherwise, apply Proposition~\ref{prop:crossblocklb} to the $n$ pairs of same-unit potential outcomes, with coefficients $1/n$ and $-1/n$, to obtain
\begin{align}
 \Runb(P;f,a)
 &\geq\frac1{n^2}\sum_{i=1}^n
 \left\{\frac{1-q_i}{q_i}+\frac{q_i}{1-q_i}+2\right\}\notag\\
 &=\frac1{n^2}\sum_{i=1}^n
 \left\{\frac1{q_i}+\frac1{1-q_i}\right\}.
 \label{eq:app-ate-paired-unbiased}
\end{align}

\emph{Exactly $m$ treated units.} On $\cZ_m$, the two arm
ranks are $r(V_1)=m$ and $r(V_0)=n-m$. Each arm contains $n$ outcomes
with coefficient magnitude $1/n$.
We apply the sampling-rank lower bound (Theorem~\ref{cor:rank-lower}) to the arm with
rank $\min\{m,n-m\}\leq n/2$ to obtain
\[
  \Rallstar(\cP(\cZ_m);f,a)
  \geq c_0\max\{m^{-1},(n-m)^{-1}\}
  \geq c\,\frac{n}{m(n-m)}
\]
for a universal constant \(c>0\).

For the all-design lower bound on the class of unbiased estimators, note that every design on $\cZ_m$ satisfies $\sum_iq_i=m$. Cauchy--Schwarz therefore gives
\[
 \sum_i\frac1{q_i}\geq\frac{n^2}{m},
 \qquad
 \sum_i\frac1{1-q_i}\geq\frac{n^2}{n-m}.
\]
Substituting these inequalities into equation~\eqref{eq:app-ate-paired-unbiased} gives
\[
 \Runbstar(\cP(\cZ_m);f,a)
 \geq\frac1m+\frac1{n-m}
 =\frac{n}{m(n-m)}.
\]

\emph{Difference-in-means upper bound.} Let \(S\) be a uniform \(m\)-subset of
\(\{1,\ldots,n\}\), and use
\[
  \widehat\tau_{\rm DIM}
  =
  \frac1m\sum_{i\in S}Y_i(1)
  -
  \frac1{n-m}\sum_{i\notin S}Y_i(0).
\]
This estimator is design-unbiased.  Put
\[
  C_i=\frac{Y_i(1)}{m}+\frac{Y_i(0)}{n-m}.
\]
Then
\[
 \widehat\tau_{\rm DIM}
 =\sum_{i\in S}C_i-\frac1{n-m}\sum_{i=1}^nY_i(0),
\]
and the second term is nonrandom.  The simple-random-sampling identity gives
\[
  \Var\left(\sum_{i\in S}C_i\right)
  =
  \frac{m(n-m)}{n(n-1)}
  \sum_{i=1}^n(C_i-\bar C)^2
  \leq
  \frac{m(n-m)}{n(n-1)}
  \sum_{i=1}^n C_i^2 .
\]
Since \(|C_i|\leq 1/m+1/(n-m)=n/\{m(n-m)\}\),
\[
  \sum_{i=1}^n C_i^2
  \leq \frac{n^3}{m^2(n-m)^2},
  \qquad
  \sup_Y
  \Var(\widehat\tau_{\rm DIM}\mid Y)
  \leq
  \frac{n^2}{(n-1)m(n-m)}.
\]
Since the difference-in-means estimator is unbiased, this upper
bound applies to both minimax risks. Together with the two lower
bounds, it proves the claims for $\cZ_{m}$.

\emph{At most $m$ treated units.} Now consider the assignment space $\cZ_{\le m}$.
We first bound the unbiased minimax risk from above and below.
For the upper bound, the uniform design on $\cZ_k$ is feasible for each
$1\leq k\leq m$. The preceding variance bound for this design is
$n^2/\{(n-1)k(n-k)\}$. To minimize this bound, we maximize $k(n-k)$
over the feasible counts. Since $k(n-k)=n^2/4-(k-n/2)^2$, we can take
$k_m=\min\{m,\lfloor n/2\rfloor\}$. Using the uniform
design on \(\cZ_{k_m}\) with the same difference-in-means estimator gives
\[
  \Runbstar(\cP(\cZ_{\le m});f,a)
  \leq
  \frac{n^2}{(n-1)k_m(n-k_m)}.
\]

For the unbiased lower bound, put
$s_m=\min\{m,n/2\}$. It suffices to consider designs with $0<q_i<1$
for all $i$. Equation~\eqref{eq:app-ate-paired-unbiased} and convexity of
$x\mapsto1/\{x(1-x)\}$ give
\[
 \Runb(P;f,a)\geq\frac1{n^2}\sum_i\frac1{q_i(1-q_i)}
 \geq\frac1{n\bar q(1-\bar q)}
 \geq\frac{n}{s_m(n-s_m)},
\]
where $\bar q=n^{-1}\sum_iq_i\leq m/n$ and $x(1-x)$ is maximized at
$1/2$.  Combining this lower bound with the difference-in-means upper bound gives the finite-sample bracket for $\cP(\cZ_{\le m})$ in Remark~\ref{rem:fixed-count-fpc}.

The same upper bound applies to $\Rallstar(\cP(\cZ_{\le m});f,a)$,
since the all-estimator class includes the difference-in-means
estimator. For the all-estimator lower bound, the treated-arm rank on
$\cZ_{\le m}$ is $r(V_1)=m$. If $m\leq n/2$,
the sampling-rank lower bound (Theorem~\ref{cor:rank-lower}) applied to the treated arm gives a constant multiple of
\(1/m\), which is a constant multiple of \(n/\{k_m(n-k_m)\}\) because
\(k_m=m\).  If \(m>n/2\), the singleton lower bound in equation~\eqref{eq:singleton-all} gives
\(1/n\), which is a constant multiple of
\(n/\{k_m(n-k_m)\}\).  This proves the minimax rates for $\cP(\cZ_{\le m})$.

Finally, the rates for $\cP(\cZ_{\le m})$ simplify to $1/m$. If $m\leq n/2$, then $k_m=m$ and
\[
 \frac1m\leq\frac{n}{m(n-m)}\leq\frac2m.
\]
If $m>n/2$, then $k_m=\lfloor n/2\rfloor$, so
$n/\{k_m(n-k_m)\}\asymp1/n\asymp1/m$.
The comparison constants are independent of $m$ and $n$, as required in equation~\eqref{eq:ate-budget-rate}.
\end{proof}

\begin{remark}[Finite-population correction and sharp constants]
\label{rem:fixed-count-fpc}
For $\cP(\cZ_m)$, the proof gives
\[
  \frac{n}{m(n-m)}
  \leq\Runbstar(\cP(\cZ_m))
  \leq\frac{n^2}{(n-1)m(n-m)},
\]
whose two ends differ by the factor \(n/(n-1)\).  The factor $n/(n-1)$ comes from the variance calculation for sampling without replacement. It tends to one and does not affect the minimax rate.  For $\cP(\cZ_{\le m})$, put $s_m=\min\{m,n/2\}$.  The corresponding bracket is
\[
 \frac{n}{s_m(n-s_m)}\leq\Runbstar(\cP(\cZ_{\leq m}))
 \leq\frac{n^2}{(n-1)k_m(n-k_m)}.
\]
Here $s_m=k_m$ when $m\leq n/2$ or $n$ is even.  For odd $n$ and
$m>n/2$, the lower bound is $4/n$ and the upper bound uses a nearly balanced
integer number of treated units.
\end{remark}

\section{Supporting results for Section~\ref{sec:partial-interference}}
\label{app:invariant-proofs}

This appendix proves Lemma~\ref{lem:invariant-representation} and Proposition~\ref{prop:unknown-partition}, then establishes the cluster-count probability and HT variance bounds used to prove Theorem~\ref{thm:invariant-tte}.

Remark~\ref{rem:peer-not-total} explains why the direct effects are indexed by the number of treated peers and not by the number of treated units in the cluster.

\begin{remark}[Number of treated peers versus number of treated units in the cluster]
\label{rem:peer-not-total}

For an assignment $z$, let $k_g(z)=\sum_{j:g(j)=g}z_j$ be the number of treated units in cluster $g$. Under a random assignment $Z$, write $K_g=k_g(Z)$. Since $k_g(z)=C_i(z)+z_i$, the pairs $(z_i,C_i(z))$ and $(z_i,k_g(z))$ contain the same information. They give two ways to index the potential outcomes. However, holding the second coordinate fixed gives different contrasts under the two indexings.

To see this, write $\widetilde Y_i(t,k)=Y_i(t,k-t)$ for the table indexed by own treatment and the number of treated units in the cluster. A direct effect with $a$ treated peers is $Y_i(1,a)-Y_i(0,a)$. With a fixed number $1\leq k\leq m-1$ of treated units in the cluster, the contrast is instead
\[
 \widetilde Y_i(1,k)-\widetilde Y_i(0,k)
 =Y_i(1,k-1)-Y_i(0,k).
\]
This changes both own treatment and the number of treated peers. It can be nonzero even if every direct effect is zero. For example, if $m\geq2$ and $Y_i(t,a)=a/(m-1)$, every direct effect is zero, but the contrast at a fixed number of treated units in the cluster is $-1/(m-1)$. This contrast includes a spillover effect.

A contrast that averages $Y_i(1,a)$ and $Y_i(0,a)$ over different distributions of $a$ can likewise be nonzero when all direct effects vanish. To average direct effects, use the same peer-count weights for both treatment arms.
\end{remark}

\subsection{Permutation-invariant designs}
\label{subsec:app-invariant-designs}

Throughout, the population is
partitioned into $\xi$ clusters of common size $m$ with $n=\xi m$ and $\xi\geq2$, and
$\cC_{\mathrm{inv}}$ is the class of Definition~\ref{def:label-invariant}.

\begin{proof}[Proof of Lemma~\ref{lem:invariant-representation}]
Let $M=\sum_{i=1}^n Z_i$ be the number of treated units. Any two
assignments with the same treated count are related by a permutation of
the unit labels, so a permutation-invariant design gives them the same
probability. For each $k$ with $P(M=k)>0$, it follows that
\[
 P(Z=z\mid M=k)=\binom nk^{-1}
 \qquad\text{for every }z\text{ with }\sum_{i=1}^n z_i=k.
\]
Thus we can generate the design by drawing $M$ from its marginal law
and then choosing the treated set uniformly among sets of size $M$.
Conversely, every permutation preserves the count and maps the uniform
law on treated sets of that size to itself.
\end{proof}

\begin{proof}[Proof of Proposition~\ref{prop:unknown-partition}]
The design $P$ must be chosen before the partition is known.
The estimator can use the partition once it is revealed.
Write $R(P,\pi)=\Rall(P;f^\pi,a^\pi)$.

\emph{Relabelling the experiment.}
For a permutation $\sigma\in\mathfrak S_n$, define
$(\sigma z)_i=z_{\sigma^{-1}(i)}$ and let $\sigma\pi$ be the
partition obtained by relabelling each unit $i$ as $\sigma(i)$.
Write $P_\sigma(z)=P(\sigma z)$ for the law of $\sigma^{-1}Z$
when $Z\sim P$. Relabelling gives a bijection between admissible
tables for $\pi$ and $\sigma\pi$, preserving the estimand and
relabelling the observed data. Transporting estimators in both
directions gives
\[
 R(P,\sigma\pi)=R(P_\sigma,\pi).
\]
Every partition in $\Pi$ is a relabelling of $\pi_0$. Hence a
permutation-invariant $P$ has $\sup_{\pi\in\Pi}R(P,\pi)=R(P,\pi_0)$.
Taking the infimum over such designs gives the upper bound on
the left-hand side of equation~\eqref{eq:unknown-partition}.

\emph{Averaging the relabelled designs.}
For an arbitrary $P$, put $N=|\mathfrak S_n|$ and
$\bar P=N^{-1}\sum_\sigma P_\sigma$. This design is
permutation invariant because relabelling permutes the mixture
components. A uniform $\sigma$ maps $\pi_0$ to a uniform
partition in $\Pi$, since each partition has the same number
of preimages. Thus
\[
 \sup_{\pi\in\Pi}R(P,\pi)
 \geq \frac1{|\Pi|}\sum_{\pi\in\Pi}R(P,\pi)
 =\frac1N\sum_{\sigma\in\mathfrak S_n}R(P_\sigma,\pi_0).
\]
We show that the last average bounds $R(\bar P,\pi_0)$ from above.

\emph{Averaging out the component label.}
Let $U$ be uniform on $\mathfrak S_n$ and draw $Z\mid U=\sigma$
from $P_\sigma$. The joint law of $(U,Z)$ is fixed by $P$ and does not depend on the outcome table. At the reference partition $\pi_0$, write
$D=(Z,Y_{O_{f^{\pi_0}}(Z)})$ and $\tau=\tau_{a^{\pi_0}}(Y)$.
If the average component risk is infinite, the desired bound is
immediate. Otherwise, for $\epsilon>0$, choose $T_\sigma$ with
worst-case risk at most $R(P_\sigma,\pi_0)+\epsilon$. Averaging over $U$ gives
\[
 \begin{aligned}
 \sup_{Y\in\cY(f^{\pi_0})}\E[(T_U(D)-\tau)^2\mid Y]
 &=\sup_{Y\in\cY(f^{\pi_0})}\frac1N\sum_{\sigma\in\mathfrak S_n}
   \E_{P_\sigma}[(T_\sigma(D)-\tau)^2\mid Y]\\
 &
 \leq \frac1N\sum_{\sigma\in\mathfrak S_n}R(P_\sigma,\pi_0)+\epsilon.
 \end{aligned}
\]
By Lemma~\ref{lem:rao-blackwell}, averaging over $U$ conditional
on the observed data gives an estimator under $\bar P$ with no
greater risk at any table. The conditional weights depend only
on $Z$. For $z$ in the support of $\bar P$, the resulting
estimator is
\[
 \bar T(z,y_{O_{f^{\pi_0}}(z)})
 =\sum_{\sigma:\,P_\sigma(z)>0}
   \frac{P_\sigma(z)}{\sum_{\rho\in\mathfrak S_n}P_\rho(z)}
   T_\sigma(z,y_{O_{f^{\pi_0}}(z)}).
\]
Letting $\epsilon\downarrow0$ gives
$R(\bar P,\pi_0)\leq N^{-1}\sum_\sigma R(P_\sigma,\pi_0)$.
Since $\bar P$ is invariant, taking the infimum over $P$ gives
the reverse inequality in equation~\eqref{eq:unknown-partition}.

For the unbiased risk, relabelling preserves unbiasedness, so
the same risk identity holds. When the average component risk
is finite, choose each $T_\sigma$ to be unbiased under
$P_\sigma$ and within $\epsilon$ of its unbiased minimax risk.
Then $\E[T_U(D)\mid U,Y]=\tau$, so mixing preserves unbiasedness.
Averaging out $U$ preserves the mean and weakly decreases the
squared-error risk at each table. The same inequalities and
infima therefore give the unbiased identity.
\end{proof}

\subsection{Bounds on cluster-count probabilities}
\label{subsec:app-cluster-count-bounds}

Conditional on the number of treated units in the population, permutation invariance gives a hypergeometric cluster count. We need bounds on the probability of an interior count and on the two extreme counts.

\begin{lemma}[Bounds on cluster-count probabilities]
\label{lem:count-bottleneck}
Fix a cluster $g$ and let $n\geq2m$.
\begin{enumerate}
\item[\textup{(i)}] For every $\delta\in(0,1/2)$ there is a constant $C_\delta$
such that every $P\in\cC_{\mathrm{inv}}$ satisfies
\[
  \Pr(K_g=\ell)\leq\frac{C_\delta}{\sqrt m},
  \qquad \delta m\leq \ell \leq(1-\delta)m .
\]
\item[\textup{(ii)}] For every $P\in\cC_{\mathrm{inv}}$ and every value of the
number $M$ of treated units in the population,
\[
  \min\bigl\{\Pr(K_g=m\mid M),\ \Pr(K_g=0\mid M)\bigr\}\leq2^{-m}.
\]
\end{enumerate}
\end{lemma}
\begin{proof}

By Lemma~\ref{lem:invariant-representation} it is enough to treat a design that
treats a uniformly chosen set of $M$ units, and then to average over $M$.
Under such a design
\begin{equation}
  h(M)=\Pr(K_g=\ell\mid M)
  =\frac{\binom M\ell\binom{n-M}{m-\ell}}{\binom nm}
  =\frac{\binom m\ell\binom{n-m}{M-\ell}}{\binom nM},
  \label{eq:hypergeometric-count}
\end{equation}
which is positive exactly for $\ell \leq M\leq n-m+\ell$.

\emph{(i)}  On that range,
\[
  \frac{h(M+1)}{h(M)}
  =\frac{(M+1)(n-M-m+\ell)}{(M+1-\ell)(n-M)},
\]
and this ratio is at least one if and only if $mM\leq \ell(n+1)-m$.  Hence $h$ is
unimodal in $M$ and attains its maximum at
$M_*=\lfloor \ell(n+1)/m\rfloor$.  Write $p=\ell/m$ and $M_*=np+e$ with $|e|<1$, so
that
\[
  \theta=\frac{M_*}{n}=p+\frac en,
  \qquad
  \lambda=\frac{M_*-\ell}{n-m}=p+\frac e{n-m}.
\]
Since $n-m\geq m$, for $m$ larger than a threshold depending only on $\delta$
all three of $p$, $\theta$ and $\lambda$ lie in $[\delta/2,1-\delta/2]$, so Stirling's bounds apply uniformly at $M_*$.
Stirling's bounds give, uniformly over proportions in that range,
$\binom Nk\asymp_\delta N^{-1/2}\exp\{NH(k/N)\}$ with
$H(x)=-x\log x-(1-x)\log(1-x)$, so applying the upper bound to the two
binomial coefficients in the numerator of equation~\eqref{eq:hypergeometric-count} and
the lower bound to the denominator,
\[
  h(M_*)\leq C_\delta\sqrt{\frac n{m(n-m)}}\,
  \exp\bigl\{mH(p)+(n-m)H(\lambda)-nH(\theta)\bigr\}.
\]
Because $n\theta=mp+(n-m)\lambda$, concavity of $H$ makes the exponent at most
zero, and $n/(n-m)\leq2$ because $n\geq2m$.  Therefore
$h(M)\leq h(M_*)\leq C'_\delta m^{-1/2}$ for every $M$, and the same bound
holds after averaging over $M$.  Enlarging the constant covers the finitely
many $m$ below the threshold.

\emph{(ii)}  Conditionally on the number $M$ of treated units in the population, write
$u_M=\Pr(K_g=m\mid M)$ and $v_M=\Pr(K_g=0\mid M)$.  If either vanishes the
claim is immediate.  Otherwise
\[
  u_Mv_M
  =\prod_{i<m}\frac{(M-i)(n-M-i)}{(n-i)^2}
  \leq\prod_{i<m}\frac{\{(n-2i)/2\}^2}{(n-i)^2}
  \leq4^{-m},
\]
by $xy\leq\{(x+y)/2\}^2$, so $\min(u_M,v_M)\leq2^{-m}$.
\end{proof}

\subsection{Upper bound under independent cluster assignments}
The upper-bound designs assign treatments independently across clusters. The following lemma bounds the variance of exposure HT through the probabilities of the two cluster counts required by the estimand.

\begin{lemma}[Independent-cluster upper bound]
\label{lem:invariant-cluster-upper}
Assign clusters independently, drawing their numbers of treated units from a common law
$q$, where $q_j=P(K_g=j)$ and choosing a uniform subset conditional on each count. For either
contrast in equation~\eqref{eq:invariant-estimands}, let $j_+$ and $j_-$ be the
counts required by its positive and negative target exposures. If both
$q_{j_+}$ and $q_{j_-}$ are positive, the exposure HT contrast satisfies
\[
 \sup_Y\operatorname{Var}(\widehat\tau_{\mathrm{HT}}\mid Y)
 \leq\frac1\xi\left(\frac1{q_{j_+}}+\frac1{q_{j_-}}\right).
\]
\end{lemma}
\begin{proof}
We write $\widehat\tau_{\mathrm{HT}}=\xi^{-1}\sum_{g=1}^{\xi} T_g$,
where $T_g$ is the difference between the two cluster-level HT means.
For a target exposure $e=(t,a)$, put $j=t+a$. On $K_g=j$,
the number of units with exposure $e$ is $s_e=j$ if $t=1$
and $s_e=m-j$ if $t=0$. Since the treated subset is uniform
conditional on $K_g=j$, each unit $i$ in cluster $g$ has
\[
 \pi_{P,i}(e)=P\{f_i(Z)=e\}=q_j\,s_e/m.
\]
Write $T_{g,e}$ for the cluster-level HT mean for exposure $e$.
The factor $1/m$ in that mean cancels the $m$ in the inverse
exposure probability. Thus
\[
 T_{g,e}=
 \begin{cases}
 \displaystyle \frac1{q_j s_e}
 \sum_{i:g(i)=g,\, f_i(Z)=e}Y_i(e),&K_g=j,\\[4pt]
 0,&K_g\ne j.
 \end{cases}
\]
The $s_e$ outcomes in the sum lie in $[-1,1]$. Hence
$\E[T_{g,e}^2\mid Y]\leq q_j/q_j^2=1/q_j$.

Write $e_+$ and $e_-$ for the positive and negative target exposures.
Their required counts are $(m,0)$ for the TTE and $(c+1,c)$ for the
DTE with $c$ treated peers. Since the counts differ, at most one
of $T_{g,e_+}$ and $T_{g,e_-}$ is nonzero. Consequently, for
$T_g=T_{g,e_+}-T_{g,e_-}$,
\[
 \E[T_g^2\mid Y]
 =\E[T_{g,e_+}^2\mid Y]+\E[T_{g,e_-}^2\mid Y]
 \leq \frac1{q_{j_+}}+\frac1{q_{j_-}}.
\]
The $T_g$ are independent conditional on $Y$ because each
depends only on its own cluster's assignment. Therefore
\[
 \operatorname{Var}(\widehat\tau_{\mathrm{HT}}\mid Y)
 =\frac1{\xi^2}\sum_{g=1}^{\xi}\operatorname{Var}(T_g\mid Y)
 \leq \frac1\xi\left(\frac1{q_{j_+}}+\frac1{q_{j_-}}\right).
\]
\end{proof}

\subsection{Proof of Theorem~\ref{thm:invariant-tte}}

\begin{resultrestatement}{Theorem~\ref{thm:invariant-tte}}{The price of graph-unaware assignment under partial interference}

Let $\xi\geq2$.
\begin{enumerate}
\item[\textup{(i)}] \emph{Graph-aware designs.}  The total treatment effect
and every direct treatment effect $\tau_{\mathrm{DTE}}(c)$,
$0\leq c\leq m-1$, satisfy
\begin{equation*}
  \frac1\xi\leq\Rallstar(\cP(\cZ))\leq\frac4\xi,
  \qquad
  \frac2\xi\leq\Runbstar(\cP(\cZ))\leq\frac4\xi .
\end{equation*}
\item[\textup{(ii)}] \emph{Graph-unaware designs: total treatment effect.}
Over permutation-invariant designs, the minimax risks satisfy
\begin{equation*}
  \frac18\min\Bigl\{1,\frac{2^m}\xi\Bigr\}
  \leq\Rallstar(\cC_{\mathrm{inv}})
  \leq\Runbstar(\cC_{\mathrm{inv}})
  \leq\min\Bigl\{4,\frac{2^{m+1}}\xi\Bigr\}.
\end{equation*}
\item[\textup{(iii)}] \emph{Graph-unaware designs: direct treatment effect
away from the extremes.}  Fix $\delta\in(0,1/2)$.  If
$\delta m\leq c\leq(1-\delta)m-1$ and $m\geq m_0(\delta)$, then
the minimax risks satisfy
\begin{equation*}
  \Runbstar(\cC_{\mathrm{inv}})\asymp_\delta\frac{\sqrt m}{\xi},
  \qquad
  \Rallstar(\cC_{\mathrm{inv}})\asymp_\delta
  \min\Bigl\{1,\frac{\sqrt m}{\xi}\Bigr\}.
\end{equation*}
\item[\textup{(iv)}] \emph{Graph-unaware designs: direct treatment effect at
the two extremes.}  For $c=0$ and $c=m-1$, the minimax risks satisfy
\begin{equation*}
  \Rallstar(\cC_{\mathrm{inv}})\asymp
  \Runbstar(\cC_{\mathrm{inv}})\asymp\frac1\xi,
  \qquad c\in\{0,m-1\}.
\end{equation*}
\end{enumerate}
\end{resultrestatement}

\begin{proof}
\emph{(i) Graph-aware designs.}
The two contrasts have the same block calculation because their target
exposures require distinct cluster counts. Write $e_+,e_-$ for the
positive and negative exposures and $j_+,j_-$ for their required
counts. For the TTE these counts are $(m,0)$, and for the DTE with
$c$ treated peers they are $(c+1,c)$.

For any design $P$, let $q_{g,j}=P(K_g=j)$ and take the two blocks
$B_{g,\pm}=\{(i,e_\pm):g(i)=g\}$ in each cluster $g$.
The sampled edge meets $B_{g,\pm}$ exactly when $K_g=j_\pm$.
The blocks are disjoint and each has total coefficient magnitude
$m/n=1/\xi$. Since $q_{g,j_+}+q_{g,j_-}\leq1$, the lower bound
from disjoint perturbations (Theorem~\ref{thm:fixed-joint-missingness}(i))
gives
\[
 \Rall(P)\geq\frac1{\xi^2}\sum_{g=1}^{\xi}
 (2-q_{g,j_+}-q_{g,j_-})\geq\frac1\xi.
\]
Part (ii) of the same theorem gives
\[
 \Runb(P)\geq\frac1{\xi^2}\sum_{g=1}^{\xi}
 \left(\frac1{q_{g,j_+}}+\frac1{q_{g,j_-}}-2\right)\geq\frac2\xi.
\]
For positive probabilities, the reciprocal sum is at least
$4/(q_{g,j_+}+q_{g,j_-})\geq4$. A zero probability gives infinite
unbiased risk. Both lower bounds hold after taking the infimum over $P$.

For the upper bound, independently in each cluster draw the treated
count from $\{j_+,j_-\}$ with probability $1/2$ each and choose a
uniform treated subset of that size. Lemma~\ref{lem:invariant-cluster-upper}
gives worst-case variance at most $(2+2)/\xi$ for the unbiased HT
estimator. The same upper bound applies over all estimators.

\textit{(ii) Graph-unaware designs: total treatment effect.}  Fix $P\in\cC_{\mathrm{inv}}$ and an integer $1\leq L\leq\xi$.
Form $\lfloor\xi/L\rfloor$ disjoint groups of $L$ clusters.
For each group $H$ of cluster labels and each target exposure
$e\in\{(1,m-1),(0,0)\}$, take the block
$B_{H,e}=\{(i,e):g(i)\in H\}$. Each block has total coefficient
magnitude $L/\xi$. Let $\rho_{H,1}$ and $\rho_{H,0}$ be the probabilities
that the sampled edge misses the treated and control blocks, respectively.
Conditional on the population treated count $M$, the union bound gives
\[
 \Pr\{O_f(Z)\cap B_{H,(1,m-1)}\ne\varnothing\mid M\}\leq Lu_M,
 \qquad
 \Pr\{O_f(Z)\cap B_{H,(0,0)}\ne\varnothing\mid M\}\leq Lv_M,
\]
where $u_M$ and $v_M$ are the conditional probabilities of the
cluster counts $m$ and $0$ from Lemma~\ref{lem:count-bottleneck}(ii).
At least one is at most $2^{-m}$ for each $M$. Thus, for every group $H$,
\[
 \rho_{H,1}+\rho_{H,0}
 \geq \E[(1-Lu_M)_++(1-Lv_M)_+]\geq1-L2^{-m}.
\]
The lower bound from disjoint perturbations (Theorem~\ref{thm:fixed-joint-missingness}(i)) applied to the
$2\lfloor \xi/L\rfloor$ sets gives
\[
  \Rall(P)\geq\Bigl\lfloor\frac{\xi}{L}\Bigr\rfloor\frac{L^2}{\xi^2}
  \bigl(1-L2^{-m}\bigr).
\]
Take $L=\min\{\xi,2^{m-1}\}$. If $L=\xi$ the bound is at least $1/2$; otherwise
$\lfloor \xi/L\rfloor\geq \xi/(2L)$ and the bound is at least $L/(4\xi)=2^{m-3}/\xi$.
Both cases give $\frac18\min\{1,2^m/\xi\}$.

For the upper bound, independent Bernoulli assignment with treatment probability $1/2$ is permutation invariant and
makes the clusters independent with $q_0=q_m=2^{-m}$, so
Lemma~\ref{lem:invariant-cluster-upper} gives an unbiased estimator
with worst-case variance at most $(1/q_0+1/q_m)/\xi=2^{m+1}/\xi$.  Alternatively,
treat the entire population with probability $1/2$ and treat no units otherwise.  This design is permutation invariant, and the estimator returning
$2n^{-1}\sum_iY_i^{\rm obs}$ when everyone is treated and
$-2n^{-1}\sum_iY_i^{\rm obs}$ when nobody is, is unbiased for
$\tau_{\mathrm{TTE}}$ with variance $(\mu_{1,m-1}+\mu_{0,0})^2\leq4$, where
$\mu_{t,a}=n^{-1}\sum_iY_i(t,a)$.  Taking
whichever of the two designs has the smaller bound gives
$\min\{4,2^{m+1}/\xi\}$. 

\textit{(iii) Interior direct treatment effect: unbiased estimators.}  Write
$\beta_m=C_\delta m^{-1/2}$ for the bound of
Lemma~\ref{lem:count-bottleneck}(i), which applies to both $c$ and $c+1$ under
the range $\delta m\leq c\leq(1-\delta)m-1$ and holds uniformly over $\cC_{\mathrm{inv}}$ and over
clusters.  Using the two perturbation sets per cluster from part~\textup{(i)} for the direct treatment effect, the lower bound from disjoint perturbations (Theorem~\ref{thm:fixed-joint-missingness}(ii)) gives
\[
  \Runb(P)\geq\frac1{\xi^2}\sum_{g=1}^\xi
  \Bigl(\frac1{q_{g,c}}+\frac1{q_{g,c+1}}-2\Bigr)
  \geq\frac1\xi\Bigl(\frac2{\beta_m}-2\Bigr),
\]
which is bounded below by a positive constant times $\sqrt m/\xi$ once $m$ exceeds a threshold depending on $\delta$. The calculation uses the marginal cluster-count probabilities and permits dependence across clusters.  For the upper bound,
independent Bernoulli assignment with treatment probability $p=(c+\frac12)/m$ is permutation invariant and makes the
clusters independent with count law $\operatorname{Binomial}(m,p)$.
For $j=c,c+1$, we have $|j-mp|=1/2$, and $p$ is bounded away from
zero and one in terms of $\delta$. Stirling's formula therefore gives
$q_j=\binom mj p^j(1-p)^{m-j}\asymp_\delta m^{-1/2}$, so
Lemma~\ref{lem:invariant-cluster-upper} gives
$(1/q_c+1/q_{c+1})/\xi$, again of order $\sqrt m/\xi$.

\textit{(iii) Interior direct treatment effect: all estimators.}  The grouping is taken
using the uniform probability bound $\beta_m$, so the same packing applies to every permutation-invariant design.  Assume $m$ large
enough that $\beta_m<1/2$, put $L=\min\{\xi,\lfloor1/(2\beta_m)\rfloor\}$ and form
$\lfloor \xi/L\rfloor$ disjoint groups of $L$ clusters.  In each group take one
perturbation block $B_{H,(1,c)}=\{(i,(1,c)):g(i)\in H\}$, with
total coefficient magnitude $L/\xi$. The sampled edge meets this block
only if some cluster in the group has $c+1$
treated units, which by the union bound has probability at most
$L\beta_m\leq1/2$, irrespective of dependence among clusters. The lower bound from disjoint perturbations (Theorem~\ref{thm:fixed-joint-missingness}(i)) gives
\[
  \Rall(P)\geq\Bigl\lfloor\frac{\xi}{L}\Bigr\rfloor\frac{L^2}{\xi^2}\cdot\frac12 .
\]
If $L=\xi$, the bound is at least $1/2$.  Otherwise $\lfloor \xi/L\rfloor\geq \xi/(2L)$ and
$\lfloor1/(2\beta_m)\rfloor\geq1/(4\beta_m)$, so the bound is at least
$L/(4\xi)\geq1/(16\xi \beta_m)$. Since $\beta_m=C_\delta/\sqrt m$, the two cases give a lower bound of order $\min\{1,\sqrt m/\xi\}$. For the upper bound, the Bernoulli design with rate $p=(c+1/2)/m$ and exposure HT has
worst-case risk of order $\sqrt m/\xi$, and the estimator that always returns
zero has worst-case risk at most four; use whichever bound is smaller.

\textit{(iv) Direct treatment effect at the extremes.}  The graph-aware lower bounds in part~\textup{(i)} apply to the smaller design class $\cC_{\mathrm{inv}}$ at both endpoints.  For $c=0$, use independent
Bernoulli assignment with rate $p=1/(m+1)$.  The two required cluster counts
have probabilities
\[
 q_0=(1-p)^m=\left(\frac{m}{m+1}\right)^m,
 \qquad
 q_1=mp(1-p)^{m-1}=q_0,
\]
Both are at least $e^{-1}$ because $m\log(1+1/m)\leq 1$.  Lemma
\ref{lem:invariant-cluster-upper} therefore gives $\Runbstar(\cC_{\mathrm{inv}})\leq 2e/\xi$.  The case $c=m-1$ follows by exchanging treatment and control and
using rate $p=m/(m+1)$.
\end{proof}

\section{Supporting results for Section~\ref{sec:neighborhood-interference}}
\label{app:neighborhood-proofs}\label{app:additional-examples}
We first prove the ordering-design bounds used for the DTE and TTE constructions.
The DTE proofs and core--periphery example follow, then the TTE lower and
upper bounds and their comparison with the conflict-graph bounds. We next
prove the treatment-budget results for both contrasts. The appendix ends
with the threshold-contrast bounds and the clique and complete-bipartite
community examples.

\subsection{General lemmas for ordering designs}

We define a design that selects units whose requested assignments are compatible, then derive its selection probabilities and an estimator variance bound. The DTE and distance-two TTE constructions specialize this design.

Throughout, \(S_i\subseteq[n]\) denotes the block of assignment coordinates that
unit \(i\)'s exposure depends on, so that \(f_i(z)\) depends only on \(z_{S_i}\); for the DTE and TTE constructions, \(S_i=N^+_i\).  We write
the estimand as a two-arm contrast
\begin{equation}
  \tau(Y)=\frac1n\sum_{i=1}^n\{Y_i(+)-Y_i(-)\},
  \label{eq:designs-generic-target}
\end{equation}
with \(+\) and \(-\) naming the two exposure levels being compared.

\subsubsection{Ordering designs}
\label{subsec:ordering-designs}

For each unit $i$ and arm $a\in\{+,-\}$, choose an assignment on $S_i$ that puts unit $i$ at exposure $a$. We call this assignment a \emph{request}. More generally, one may specify a distribution over such assignments for each unit and arm. For the treated DTE arm, the request treats $i$ and leaves its neighbors untreated; the control request leaves all of $N_i^+$ untreated. For the TTE, the two requests treat all of $B_1(i)$ or leave all of it untreated.

Two requests can be implemented together if they agree on the intersection of their assignment blocks. We say that units $i$ and $j$ \emph{conflict} if some choice of their arms and supported requests assigns different treatments to a unit in $S_i\cap S_j$. Define
\begin{equation}
 \mathcal N_i=\{j\ne i:i\text{ and }j\text{ conflict}\},
 \qquad \mathcal N_i^+=\mathcal N_i\cup\{i\}.
 \label{eq:designs-clash-set}
\end{equation}
The sets $\mathcal N_i$ are fixed by the blocks and request distributions. The design selects units with no conflicts between them, then draws their arms.

\begin{definition}[Ordering design]
\label{def:ordering-design}
Fix assignment blocks \(S_i\) and their requests, proposal probabilities
\(p_i\in(0,1]\), and a total order \(\prec\) on the units, either fixed in
advance or drawn at random independently of everything else.  Assign treatment
as follows.
\begin{enumerate}
\item \emph{Propose.}  Each unit \(i\) proposes with probability \(p_i\),
independently across units, to join a set \(\cI\).  Let \(B_i\) denote the
proposal indicator.
\item \emph{Select.}  A proposing unit is added to \(\cI\) when no conflicting unit precedes it in \(\prec\) and proposed:
\begin{equation}
  R_i=\1\{i\in\cI\}
  =B_i\prod_{j\in\mathcal N_i:\,j\prec i}(1-B_j).
  \label{eq:ordering-membership}
\end{equation}
\item \emph{Assign.}  Each unit of \(\cI\) draws an arm uniformly from
\(\{+,-\}\), independently across units and independently of the proposals and order, and the design implements its request on \(S_i\). If requests are random, draw them from the specified distributions independently of the proposals and order. Units outside the selected blocks are left untreated.
\end{enumerate}
\end{definition}

The assignment is well defined on the full cube. If two conflicting units both propose, the later unit is excluded. Thus selected units have requests that agree on every overlapping block, whatever their arms, and each selected unit attains its assigned exposure. Under a restricted assignment space, one must also check that the completed assignment belongs to that space.

Write
\begin{equation}
 q_i=\Pr(i\in\cI),
 \qquad
 C_{ij}=\Cov\left(\frac{R_i}{q_i},\frac{R_j}{q_j}\right)
 =\frac{\E[R_iR_j]}{q_iq_j}-1,
 \label{eq:ordering-q-cov}
\end{equation}
where $C_{ij}$ is defined when $q_i,q_j>0$. The event that unit $i$ is selected and assigned a specified arm has probability $q_i/2$. A unit outside $\cI$ may also attain the exposure, so its full exposure probability can exceed $q_i/2$.

\begin{lemma}[Selection probabilities]
\label{lem:ordering-propensity}
Under Definition~\ref{def:ordering-design},
\(q_i=p_i\,\E\prod_{j\in\mathcal N_i:\,j\prec i}(1-B_j)\), and in particular:
\begin{enumerate}
\item[\textup{(i)}] if \(\mathcal N_i=\varnothing\), then \(q_i=p_i\);
\item[\textup{(ii)}] for a fixed order,
\(q_i=p_i\prod_{j\in\mathcal N_i:\,j\prec i}(1-p_j)\);
\item[\textup{(iii)}] for a uniformly random order,
\begin{equation}
  q_i=p_i\int_0^1\prod_{j\in\mathcal N_i}(1-tp_j)\,dt,
  \label{eq:ordering-uniform-q}
\end{equation}
which for a common proposal probability \(p\) equals
\(\{1-(1-p)^{|\mathcal N^+_i|}\}/|\mathcal N^+_i|\), and at \(p=1\) equals
\(1/|\mathcal N^+_i|\) exactly.
\end{enumerate}
\end{lemma}

\begin{proof}
Equation~\eqref{eq:ordering-membership} gives the first expression, because $B_i$ is independent of the order and the other proposal indicators. If $\mathcal N_i$ is empty, the product equals one, proving part (i). For a fixed order, independence of the proposal indicators gives the product in part (ii). For part (iii), attach independent uniform marks \(U_j\) on \([0,1]\) and order by them.
Conditionally on \(U_i=t\), each \(j\in\mathcal N_i\) blocks \(i\) exactly when
\(B_j=1\) and \(U_j<t\), which has probability \(tp_j\), independently across
\(j\).  Integrating over \(t\) gives equation~\eqref{eq:ordering-uniform-q}.  With
\(p_j\equiv p\) the integral is
\(p\int_0^1(1-tp)^{m}dt=\{1-(1-p)^{m+1}\}/(m+1)\) with
\(m=|\mathcal N_i|\); at \(p=1\) the design keeps \(i\) exactly when \(i\)
precedes every conflicting unit, and the earliest element of
\(\mathcal N^+_i\) is uniform on that set.
\end{proof}

A uniformly random order gives $q_i>0$ for every unit, including when every proposal probability equals one. For a fixed order, $q_i>0$ precisely when every earlier conflicting unit has proposal probability less than one.

\paragraph*{The estimator}
Because selection together with either specified arm has probability
\(q_i/2\), define the Horvitz--Thompson estimator based on
\(\cI\),
\begin{equation}
  \widehat\tau_{\cI}
  =\frac2n\sum_{i\in\cI}\frac1{q_i}
   \bigl\{\1\{f_i(Z)=+\}-\1\{f_i(Z)=-\}\bigr\}Y_i^{\rm obs},
  \label{eq:ordering-thinned-estimator}
\end{equation}
which is unbiased whenever $q_i>0$ for every unit. The estimator uses the selection indicators $R_i$, which need not be determined by the observed assignment $Z$. To obtain an estimator in the observed-data class $\cT(P,f)$, define
\[
 \widehat\tau_{\rm RB}
 =\E[\widehat\tau_{\cI}\mid Z,Y_{O_f(Z)}].
\]
Conditional on $Z$, the law of the auxiliary randomization is known and does not depend on the potential outcome table. Lemma~\ref{lem:rao-blackwell} therefore shows that $\widehat\tau_{\rm RB}$ is a valid observed-data estimator, remains unbiased, and has no larger variance. We bound the variance of $\widehat\tau_{\cI}$ and use the bound for $\widehat\tau_{\rm RB}$.

\begin{proposition}[Variance upper bound for ordering designs]
\label{prop:ordering-envelope}
Suppose \(|Y_i(\pm)|\leq1\).  Under any ordering design with \(q_i>0\) for
every unit, the estimator in equation~\eqref{eq:ordering-thinned-estimator} and
\(\widehat\tau_{\rm RB}\) are unbiased and both obey
\begin{equation}
  \sup_{Y\in\cY(f)}\Var(\widehat\tau_{\rm RB}\mid Y)
  \leq
  \sup_{Y\in\cY(f)}\Var(\widehat\tau_{\cI}\mid Y)
  \leq
  \frac4{n^2}
  \left\{2\sum_{i=1}^n\frac1{q_i}
  +\sum_{i\neq j}|C_{ij}|\right\}.
  \label{eq:ordering-envelope}
\end{equation}
\end{proposition}

\begin{proof}
Put
\[
 d_i=Y_i(+)-Y_i(-),
 \qquad s_i=Y_i(+)+Y_i(-).
\]
Fix $Y$ and write $R=(R_1,\ldots,R_n)$. Let $A_i$ be $1$ for
the positive arm and $-1$ for the negative arm, extending the arm
draws to unselected units by independent signs. Conditional on $R$,
the $A_i$ are independent with mean zero and variance one. Then
\[
 \widehat\tau_{\cI}=\frac1n\sum_i\frac{R_i}{q_i}(d_i+A_i s_i).
\]
Since $R_i^2=R_i$, this gives
\[
 \E(\widehat\tau_{\cI}\mid R,Y)
 =\frac1n\sum_i\frac{R_i}{q_i}d_i,
 \qquad
 \Var(\widehat\tau_{\cI}\mid R,Y)
 =\frac1{n^2}\sum_i\frac{R_i}{q_i^2}s_i^2.
\]
The conditional mean identity and $\E R_i=q_i$ prove unbiasedness. Taking expectations in the conditional variance identity and using $|s_i|\leq2$ gives
\[
 \E\{\Var(\widehat\tau_{\cI}\mid R,Y)\}
 \leq\frac4{n^2}\sum_i\frac1{q_i}.
\]
For the conditional mean, the definition of $C_{ij}$ gives
\begin{align*}
 \Var\{\E(\widehat\tau_{\cI}\mid R,Y)\}
 &=\frac1{n^2}\sum_{i,j}d_id_jC_{ij}\\
 &\leq\frac4{n^2}
 \left\{\sum_i\left(\frac1{q_i}-1\right)
       +\sum_{i\ne j}|C_{ij}|\right\},
\end{align*}
because $C_{ii}=q_i^{-1}-1$ and $|d_i|\leq2$.  Combining the conditional variance and conditional mean bounds by the law of total variance gives
\[
 \Var(\widehat\tau_{\cI}\mid Y)
 \leq\frac4{n^2}\left\{2\sum_iq_i^{-1}-n+\sum_{i\ne j}|C_{ij}|\right\}.
\]
Dropping $-n$ and taking the supremum over $Y$ yields
equation~\eqref{eq:ordering-envelope}.  Finally,
$\widehat\tau_{\rm RB}=\E(\widehat\tau_{\cI}\mid Z,Y_{O_f(Z)})$ is unbiased and cannot have larger variance.
\end{proof}

The first term in equation~\eqref{eq:ordering-envelope} depends on the selection probabilities. The second depends on covariances between selection indicators, including pairs with shared conflicting units. The following lemma gives these covariances when all units propose and the order is uniform.

\begin{lemma}[Membership covariances under a uniformly random order]
\label{lem:ordering-covariance}
Take \(p_i\equiv1\) and a uniformly random order, so that
\(q_i=1/|\mathcal N^+_i|\).  For \(i\neq j\),
\begin{enumerate}
\item[\textup{(i)}] if $i$ and $j$ conflict, then \(C_{ij}=-1\);
\item[\textup{(ii)}] if \(\mathcal N^+_i\cap\mathcal N^+_j=\varnothing\), then
\(C_{ij}=0\);
\item[\textup{(iii)}] otherwise
\(C_{ij}=\frac{|\mathcal N_i^+\cap\mathcal N_j^+|}{|\mathcal N_i^+\cup\mathcal N_j^+|}\leq
|\mathcal N^+_i\cap\mathcal N^+_j|\,/\,|\mathcal N^+_j\setminus\mathcal N^+_i|\).
\end{enumerate}
Consequently, if the second-order sets satisfy
\(|\bigcup_{j\in\mathcal N^+_i}\mathcal N^+_j|\leq\kappa|\mathcal N^+_i|\) for
every \(i\), then
\begin{equation}
  \sum_{i\neq j}|C_{ij}|\ \leq\ \kappa\sum_{i=1}^n|\mathcal N^+_i| .
  \label{eq:ordering-covariance-sum}
\end{equation}
\end{lemma}

\begin{proof}
Here \(R_i=\1\{i\text{ is the earliest element of }\mathcal N^+_i\}\).  If
\(j\in\mathcal N_i\) then \(R_iR_j=0\), giving (i).  If \(\mathcal N^+_i\) and
\(\mathcal N^+_j\) are disjoint, the pair of earliest elements of the two sets is
uniform on \(\mathcal N^+_i\times\mathcal N^+_j\) by symmetry of the uniform
order, so the two events are independent and (ii) holds.  For units $i,j$ that do not conflict, put $a=|\mathcal N_i^+|$, $b=|\mathcal N_j^+|$,
and $c=|\mathcal N_i^+\cap\mathcal N_j^+|$.  Integrating the independent
uniform priority marks over $U_i=x<U_j=y$ gives
\[
 \int_0^1\int_x^1(1-x)^{a-c-1}(1-y)^{b-1}\,dy\,dx
 =\frac1{b(a+b-c)}.
\]
Here the $a-c-1$ marks belonging only to the first set must exceed $x$;
the $b-1$ marks other than $j$ in the second set must exceed $y$.
The reverse ordering contributes $1/\{a(a+b-c)\}$, so
\[
  \Pr(R_i=R_j=1)=\frac{1}{b(a+b-c)}+\frac{1}{a(a+b-c)}.
\]
Dividing by $q_iq_j=1/(ab)$ and subtracting one gives
\[
 C_{ij}=\frac{a+b}{a+b-c}-1
 =\frac{c}{a+b-c}
 \leq\frac{c}{|\mathcal N_j^+\setminus\mathcal N_i^+|},
\]
proving part (iii). The denominator on the right is positive because it contains $j$.

For the sum bound, the exact formula in part (iii) gives $C_{ij}\leq1$
because the intersection is contained in the union. Such a covariance vanishes unless the two conflict neighborhoods intersect. For fixed $i$, at most
\[
 \left|\bigcup_{j\in\mathcal N_i^+}\mathcal N_j^+\right|
 -|\mathcal N_i^+|
 \leq(\kappa-1)|\mathcal N_i^+|
\]
units outside $\mathcal N_i^+$ have a nonzero covariance with $R_i$. The $|\mathcal N_i^+|-1$ units that conflict with $i$ each contribute $|C_{ij}|=1$. Therefore
\[
 \sum_{i\ne j}|C_{ij}|
 \leq\sum_i\{(\kappa-1)|\mathcal N_i^+|+|\mathcal N_i^+|-1\}
 \leq\kappa\sum_i|\mathcal N_i^+|,
\]
which proves equation~\eqref{eq:ordering-covariance-sum}.
\end{proof}

\subsection{Proofs for the DTE rates and the ordered proposal design}
\label{subsec:app-dte-ordered-thinning}

The treated DTE potential outcomes can be jointly observed exactly on independent sets of $G$. We use this observation constraint for the lower bound, then apply the ordering-design variance bound to prove Theorem~\ref{thm:dte-rates}.

\begin{proposition}[General DTE lower bound]
\label{prop:dte-general-rank}
For the DTE on any graph $G$,
\begin{equation}
  \Runbstar\bigl(\cP(\cZ);f,a_{\rm DTE}\bigr)
  \geq
  \max\left\{\frac4n,\frac1{\alpha(G)}-\frac1n\right\}
  \geq\frac1{2\alpha(G)},
  \label{eq:dte-general-unb-lower}
\end{equation}
and, for a universal constant $c>0$,
\begin{equation}
  \Rallstar\bigl(\cP(\cZ);f,a_{\rm DTE}\bigr)
  \geq\frac{c}{\alpha(G)}.
  \label{eq:dte-general-all-lower}
\end{equation}
\end{proposition}

\begin{proof}
Write $V_+=\{(i,(1,0_{N_i})):i\in[n]\}$ for the treated arm.
Its sampling rank is $r(V_+)=\alpha(G)$, since the units observed
in this arm form an independent set and treating exactly an independent
set observes its treated-arm outcomes.

We first prove the unbiased lower bounds in equation~\eqref{eq:dte-general-unb-lower}.
If any target exposure has probability zero under $P$, then
$\Runb(P;f,a_{\rm DTE})=+\infty$. For the unbiased calculation below,
fix a design with positive target-exposure probabilities and a
design-unbiased estimator $T$. For the term
\(1/\alpha(G)-1/n\), write \(\pi_{i,+}\) for the probability of observing
unit \(i\)'s treated-arm potential outcome under $P$.  Since every
assignment observes at most \(\alpha(G)\) such potential outcomes,
\(\sum_i\pi_{i,+}\leq\alpha(G)\).  The unbiased singleton lower
bound in equation~\eqref{eq:singleton-unb} gives the first inequality,
and Lemma~\ref{lem:weighted-reciprocal}, with $w_i=1$ and
$p_i=\pi_{i,+}$, gives the second:
\[
 \sup_Y\Var(T\mid Y)\geq
 \frac1{n^2}\sum_i\left(\frac1{\pi_{i,+}}-1\right)
 \geq\frac1{\alpha(G)}-\frac1n
\]
for every design-unbiased estimator $T$.

The term \(4/n\) instead uses the \(n\) pairs
\(\{Y_i(1,0_{N_i}),Y_i(0,0_{N_i})\}\), one pair for each unit.
No assignment observes both members of any pair.  If their inclusion probabilities are \(\pi_{i,+},\pi_{i,-}\), then
\(\pi_{i,+}+\pi_{i,-}\leq1\).  The paired lower bound in equation~\eqref{eq:crossblocklb-multipair}, applied to coefficients $1/n$ and $-1/n$, gives
\[
 \sup_Y\Var(T\mid Y)\geq
 \frac1{n^2}\sum_i
 \left(\frac1{\pi_{i,+}}+\frac1{\pi_{i,-}}\right)
 \geq\frac4n.
\]
Taking the infimum over positive designs gives the maximum in equation~\eqref{eq:dte-general-unb-lower}. If $\alpha(G)\leq n/2$, then $1/\alpha(G)-1/n\geq1/\{2\alpha(G)\}$. If $\alpha(G)>n/2$, the bound $4/n$ is at least $1/\{2\alpha(G)\}$.

For the all-estimator risk, the sampling-rank lower bound (Theorem~\ref{cor:rank-lower}) applied to $V_+$ gives
\(c_{1/2}/\alpha(G)\) when \(\alpha(G)\leq n/2\).  If $\alpha(G)>n/2$, the singleton lower bound in equation~\eqref{eq:singleton-all} gives
\(n^{-2}\sum_i(2-\pi_{i,+}-\pi_{i,-})\geq1/n
\geq1/\{2\alpha(G)\}\).  Taking the infimum over designs proves equation~\eqref{eq:dte-general-all-lower}.
\end{proof}

\begin{proposition}[Variance of the ordered proposal design for the DTE]
\label{prop:dte-ordered-proposal-thinning}
Suppose \(|Y_i(e)|\leq1\) and consider the DTE on a graph \(G\).  Let
\(R_i\) indicate that unit \(i\) is included under the ordered proposal
design \(P_{\prec,p}\), write $q_i=P(R_i=1)$, and define the normalized selection covariances
\begin{equation}
  C_{ij}=\Cov\left(\frac{R_i}{q_i},\frac{R_j}{q_j}\right).
  \label{eq:ordered-thinning-q-cov}
\end{equation}
There is a universal constant \(c>0\) such that, whenever \(q_i>0\) for every
unit,
\begin{equation}
  \sup_{Y\in\cY(f)}\Var(\widehat\tau_{\rm RB}\mid Y)
  \leq
  V_{\rm DTE}(\prec,p):=\frac4{n^2}
  \left\{2\sum_{i=1}^n\frac1{q_i}
  +\sum_{i\neq j}|C_{ij}|\right\},
  \label{eq:ordered-thinning-upper}
\end{equation}
and consequently, for every choice of order and proposal probabilities with positive selection probabilities,
\begin{equation}
  \frac{c}{\alpha(G)}
  \leq
  \Rallstar(\cP(\cZ);f,a_{\rm DTE})
  \leq
  \Runbstar(\cP(\cZ);f,a_{\rm DTE})
  \leq
  V_{\rm DTE}(\prec,p).
  \label{eq:ordered-thinning-sandwich}
\end{equation}
\end{proposition}

\begin{proof}[Proof of Proposition~\ref{prop:dte-ordered-proposal-thinning}]
The two DTE requests on $N_i^+$ treat or control unit $i$ and leave its neighbors untreated. Requests at nonadjacent units agree on shared neighbors; requests at adjacent units can conflict. Thus the conflict neighborhood in Definition~\ref{def:ordering-design} is $\mathcal N_i=N_i$, and the selected set $\cI$ is independent in $G$. Each selected unit attains its assigned DTE arm, with selection-and-arm probability $q_i/2$.

Proposition~\ref{prop:ordering-envelope} gives equation~\eqref{eq:ordered-thinning-upper} for the observed-data estimator $\widehat\tau_{\rm RB}$. Since this estimator is unbiased,
\[
 \Rallstar(\cP(\cZ);f,a_{\rm DTE})
 \leq\Runbstar(\cP(\cZ);f,a_{\rm DTE})
 \leq\sup_Y\Var(\widehat\tau_{\rm RB}\mid Y).
\]
Combining these inequalities with Proposition~\ref{prop:dte-general-rank} proves equation~\eqref{eq:ordered-thinning-sandwich}.
\end{proof}

\emph{Rates for representative graph families.}
The table gives the common minimax rates and the inputs to the ordered
proposal design that attains them with $\widehat\tau_{\rm RB}$.
Here $\eta\in(0,1)$ is fixed and $m_j$ is the size of star component $j$.
\begin{center}
\begin{tabular}{llll}
\toprule
Family & Order & Proposal probabilities & Rate\\
\midrule
Bounded degree & any & \(p_i=\eta\) & \(1/n\)\\
Disjoint stars & center first & center: \(\eta/m_j\); leaves: \(\eta\) & \(1/n\)\\
Cliques of size \(m\) & any & \(p_i=\eta/m\) & \(m/n\)\\
\bottomrule
\end{tabular}
\end{center}
A large maximum degree alone need not make the DTE difficult. The star
construction keeps each component's contribution to the variance bound
of order $m_j$ before the $n^{-2}$ normalization. Since the component
sizes sum to $n$, the risk is $O(1/n)$ even for a single large star.

For bounded-degree graphs, take any fixed order and proposal probabilities
bounded away from zero and one.  Then all \(q_i\)'s are bounded below and each $R_i$ depends on proposals in $N_i^+$. Selection indicators with disjoint proposal sets are independent, and each $N_i^+$ intersects only a bounded number of other such sets. Thus $\sum_{i\ne j}|C_{ij}|=O(n)$ and $\sum_iq_i^{-1}=O(n)$. Proposition~\ref{prop:dte-ordered-proposal-thinning} gives an upper bound of order $1/n$, matching Proposition~\ref{prop:dte-general-rank} because $\alpha(G)\asymp n$.

For a disjoint union of stars, fix \(\eta\in(0,1)\).  In a component of
size \(m\), put the center \(h\) first, take \(p_h=\eta/m\), and set
\(p_i=\eta\) on its leaves.  Then
\[
 q_h=\eta/m,\qquad q_i=\eta(1-\eta/m),\qquad
 C_{hi}=-1,\qquad C_{ij}=\frac{\eta/m}{1-\eta/m}
\]
for distinct leaves \(i,j\).  Consequently \(\sum_iq_i^{-1}=O(m)\)
and \(\sum_{i\ne j}|C_{ij}|=O(m)\), with constants depending only on
\(\eta\); the formulas also cover an isolated center, with empty leaf sums.
Different components use independent proposal coins, so their covariance
terms vanish.  Summing over component sizes gives \(O(n)\) for both sums
and hence the upper rate \(1/n\).  Choosing all leaves in nontrivial stars
and all isolated vertices shows \(n/2\leq\alpha(G)\leq n\), so the lower
rate matches uniformly over component sizes.

For disjoint \(m\)-cliques, take any order inside each clique and proposal
probabilities \(p_i\asymp1/m\).  Then \(q_i\asymp1/m\), so $\sum_iq_i^{-1}=O(nm)$. At most one unit is selected from each clique, giving $C_{ij}=-1$ within a clique; different cliques have independent proposal coins. Hence $\sum_{i\ne j}|C_{ij}|=O(nm)$ and equation~\eqref{eq:ordered-thinning-upper} gives an upper bound of order $m/n$.  Since
\(\alpha(G)=n/m\), the lower bound matches.

\subsubsection{Proof of the regular-graph DTE bounds}
\label{subsec:app-dte-regular}

Use a uniformly random order and proposal probability one. A unit is selected exactly when it is first in its closed neighborhood, so $q_i=1/(d+1)$. Lemma~\ref{lem:ordering-covariance} gives $C_{ij}=-1$ for adjacent units and
\[
 C_{ij}=\frac{|B_1(i)\cap B_1(j)|}{|B_1(i)\cup B_1(j)|}
 \qquad\text{for nonadjacent units}.
\]
The covariance sum counts ordered pairs of distinct units, so the
adjacent pairs contribute $nd$. For nonadjacent units,
$|B_1(i)\cup B_1(j)|\geq d+1$. Hence their contribution is at most
\[
 \frac1{d+1}\sum_{i,j}|B_1(i)\cap B_1(j)|=n(d+1),
\]
because
\[
 \sum_{i,j}|B_1(i)\cap B_1(j)|
 =\sum_v|\{i:v\in B_1(i)\}|^2=n(d+1)^2.
\]
Substituting $\sum_iq_i^{-1}=n(d+1)$ and $\sum_{i\ne j}|C_{ij}|\leq nd+n(d+1)$ into equation~\eqref{eq:ordered-thinning-upper} gives
\begin{equation}
 \Runbstar(\cP(\cZ);f,a_{\rm DTE})\leq\frac{4(4d+3)}n.
 \label{eq:dte-regular-app-sandwich}
\end{equation}

\begin{resultrestatement}{Theorem~\ref{thm:dte-rates}}{Minimax rates for the neighborhood DTE}
There are universal constants $c,C>0$ such that the following hold.
\begin{enumerate}
\item[\textup{(i)}] For every graph $G$,
\[
  \frac{c}{\alpha(G)}
  \leq\Rallstar(\cP(\cZ);f,a_{\rm DTE})
  \leq\Runbstar(\cP(\cZ);f,a_{\rm DTE}).
\]
\item[\textup{(ii)}] If $G$ is $d$-regular, then
\begin{equation*}
  \frac{c}{\alpha(G)}
  \leq
  \Rallstar(\cP(\cZ);f,a_{\rm DTE})
  \leq
  \Runbstar(\cP(\cZ);f,a_{\rm DTE})
  \leq
  C\,\frac{d+1}{n}.
\end{equation*}
\item[\textup{(iii)}] If $G$ is a uniformly random simple \(d\)-regular
graph, there are universal
constants \(d_0,c_0>0\) such that, whenever
\(d_0\leq d<n-c_0n/\log n\) and \(nd\) is even, with probability
\(1-o(1)\) as \(n\to\infty\),
\begin{equation*}
  \frac{d}{n\log d}
  \lesssim
  \Rallstar(\cP(\cZ);f,a_{\rm DTE})
  \leq
  \Runbstar(\cP(\cZ);f,a_{\rm DTE})
  \lesssim
  \frac dn .
\end{equation*}
\end{enumerate}
\end{resultrestatement}

\begin{proof}[Proof of Theorem~\ref{thm:dte-rates}]
Part~\textup{(i)} is Proposition~\ref{prop:dte-general-rank}, together with inclusion of the unbiased estimator class in the all-estimator class. For part~\textup{(ii)}, use the uniformly random order with proposal probability one. Equation~\eqref{eq:dte-regular-app-sandwich} gives
\[
 \Runbstar(\cP(\cZ);f,a_{\rm DTE})
 \leq\frac{4(4d+3)}n\leq\frac{16(d+1)}n.
\]
The lower bound is again Proposition~\ref{prop:dte-general-rank}.

Under the random-regular assumptions in part~\textup{(iii)}, the independence-number bound gives
$\alpha(G)\lesssim n\log d/d$ with probability $1-o(1)$
\citep{FriezeLuczak1992,CooperEtAl2002,KrivelevichEtAl2001}.
On this event, Proposition~\ref{prop:dte-general-rank} gives
$\Rallstar\gtrsim d/(n\log d)$, and part~\textup{(ii)} gives
$\Runbstar\lesssim d/n$ because $d\geq d_0>0$. This proves part~\textup{(iii)}.
\end{proof}

\subsection{Proof for the core-periphery DTE}
\label{subsec:app-core-periphery-dte}

Let $C$ be the core, with $|C|=S$, and let $P=V\setminus C$ be the periphery. The core is a clique, every core unit is adjacent to every peripheral unit, and $\Delta(G[P])\leq r$.

\begin{proof}[Proof of Theorem~\ref{cor:dte-price-is-a-rate}]
Write $B_C=\{(i,(1,0_{N_i})):i\in C\}$ and
$B_P=\{(i,(1,0_{N_i})):i\in P\}$ for the treated-arm blocks
in the core and periphery. We first give the lower bounds, then
construct a strategy attaining each rate.

\textit{All-estimator lower bound.} Proposition~\ref{prop:dte-general-rank} gives $\Rallstar\gtrsim1/n$. For $S\geq2$, restrict to tables whose only nonzero coordinates are in
$B_C$. Since $C$ is a clique, $r(B_C)=1$. Their contribution is $S/n$ times their mean, so the sampling-rank lower bound (Theorem~\ref{cor:rank-lower}) gives $\Rallstar\gtrsim S^2/n^2$. For $S=1$, this term is absorbed by $1/n$. Taking the maximum gives the lower bound of order $1/n+S^2/n^2$.

\textit{Unbiased lower bound.} Fix a design $P_Z$. If a target exposure has zero probability, then
$\Runb(P_Z;f,a_{\rm DTE})=+\infty$. Otherwise write $\pi_{i,+}$ for the probability of observing unit $i$'s treated DTE potential outcome. Since this potential outcome can be observed for at most one core unit, $\sum_{i\in C}\pi_{i,+}\leq1$. For $S\geq2$, the singleton lower bound and Cauchy--Schwarz give
\[
 \Runb(P_Z;f,a_{\rm DTE})
 \geq\frac1{n^2}\sum_{i\in C}(\pi_{i,+}^{-1}-1)
 \geq\frac{S^2-S}{n^2}\geq\frac{S^2}{2n^2}.
\]
When $B_P$ is nonempty, we use $B_C$ and $B_P$ as disjoint perturbation
blocks. Define
\[
 u=P_Z\{O_f(Z)\cap B_C\ne\varnothing\},
 \qquad
 v=P_Z\{O_f(Z)\cap B_P\ne\varnothing\}.
\]
Both are positive under the exposure condition above. The sampled edge
cannot meet both blocks because every core unit is adjacent to every
peripheral unit, so $u+v\leq1$. The lower bound from disjoint perturbations (Theorem~\ref{thm:fixed-joint-missingness}(ii)) gives
\[
 \Runb(P_Z;f,a_{\rm DTE})
 \geq\frac1{n^2}\left\{S^2\frac{1-u}{u}+(n-S)^2\frac{1-v}{v}\right\}
 \geq\frac{2S(n-S)}{n^2}.
\]
For fixed $u$, the expression decreases as $v$ increases, so it is minimized at $v=1-u$. The last inequality then follows from the arithmetic--geometric mean inequality. Combining the core and core--periphery bounds gives a lower bound of order
\[
 \frac{S^2+S(n-S)}{n^2}=\frac Sn.
\]
If $B_P$ is empty, the core bound alone gives the required rate for $S\geq2$. For $S=1$, Proposition~\ref{prop:dte-general-rank} supplies the same order $1/n$. These bounds hold for every design, so they also hold after taking the infimum over designs.

\textit{All-estimator upper bound.} If the periphery is empty, the zero estimator has risk at most
$4=4S^2/n^2$. Otherwise, we keep the core untreated, run the bounded-degree ordered proposal design on $G[P]$, and use the HT estimator of Proposition~\ref{prop:dte-ordered-proposal-thinning} on $G[P]$ with weights $1/n$ in place of $1/|P|$. The bounded-degree calculation gives variance $O_r(|P|/n^2)\leq O_r(1/n)$. Omitting the core contribution incurs absolute bias at most $2S/n$, giving risk of order $1/n+S^2/n^2$.

\textit{Unbiased upper bound.} Use Proposition~\ref{prop:dte-ordered-proposal-thinning}, ordering the core before the periphery. Give each core unit proposal probability $\eta/(n+1)$ and each peripheral unit proposal probability $\eta/(r+2)$, where $0<\eta\leq1/4$ is fixed. For unit $i$, let $D_i$ be its earlier neighbors. Their total proposal probability is at most $2\eta\leq1/2$, so
\[
 q_i=p_i\prod_{h\in D_i}(1-p_h)
 \geq p_i\left(1-\sum_{h\in D_i}p_h\right)\geq p_i/2.
\]
Consequently,
\[
 \sum_iq_i^{-1}\lesssim_r S(n+1)+(n-S)\lesssim_r Sn.
\]
Adjacent units cannot both be selected, so $C_{ij}=-1$. There are $S(S-1)+2S(n-S)+O_r(n)$ ordered adjacent pairs, which is $O_r(Sn)$. For nonadjacent units, a common earlier neighbor's nonproposal probability
appears once in the joint selection probability and twice in the product
$q_iq_j$. All other factors cancel, giving
\[
 C_{ij}=\prod_{h\in D_i\cap D_j}(1-p_h)^{-1}-1
 \lesssim\sum_{h\in D_i\cap D_j}p_h.
\]
Indeed, $-\log(1-p_h)\leq2p_h$ and the sum in the exponent is at most one, so $e^x-1\leq(e-1)x$ applies. Each core unit has at most $n$ later neighbors, and each peripheral unit has at most $r$ later neighbors. Counting pairs for each common earlier neighbor therefore gives
\[
 \sum_{i\ne j:\,ij\notin E}|C_{ij}|
 \lesssim S\frac{\eta}{n+1}n^2+(n-S)\frac{\eta}{r+2}r^2
 \lesssim_r Sn.
\]
Substituting these bounds in equation~\eqref{eq:ordered-thinning-upper}
gives worst-case variance $O_r(S/n)$ for the unbiased observed-data
estimator $\widehat\tau_{\rm RB}$.

Combining the bounds gives the two rates. Their ratio satisfies
\[
 \frac{S/n}{1/n+S^2/n^2}=\frac{Sn}{n+S^2}
 \asymp\min\{S,n/S\},
\]
which is largest, of order $\sqrt n$, when $S\asymp\sqrt n$.
\end{proof}

\subsection{Proofs for the TTE bounds}
\label{subsec:app-tte-packed-clusters}

We prove Theorem~\ref{thm:restricted-growth-gcr-tte}. The first lower bound uses disjoint closed neighborhoods, within each of which treatment- and control-arm potential outcomes cannot be observed together. The second bounds how many pairs of blocks of potential outcomes an assignment can resolve. The restricted-growth condition then gives a matching graph-cluster upper bound.

Write $1$ and $0$ for the fully treated and fully untreated neighborhood exposures, respectively. Define
\[
  \mathsf{Pack}_1(G)
  =
  \max_{C\subseteq V:\operatorname{dist}_G(c,c')\geq3}
  \sum_{c\in C}|B_1(c)|^2,
\]
where the distance condition applies to distinct centers $c,c'\in C$. For a treated set
\(T\subseteq V\), put
\begin{equation}
  A_T=\{i:B_1(i)\subseteq T\},
  \qquad
  B_T=\{i:B_1(i)\cap T=\varnothing\},
  \qquad
  M(G)=\max_{T\subseteq V}\min\{|A_T|,|B_T|\}.
  \label{eq:tte-capacity}
\end{equation}
The sets $A_T$ and $B_T$ index the treatment- and control-arm potential outcomes observed under $T$. Thus $M(G)$ is the largest number of potential outcomes that can be observed from each arm in one assignment. Combine the two lower bounds by setting
\begin{equation}
  L_{\rm TTE}(G)
  =
  \max\left\{
    \frac{\mathsf{Pack}_1(G)}{n^2},\
    \max_{1\leq K\leq n}(K-M(G))_+\frac{\lfloor n/K\rfloor^2}{n^2}
  \right\}.
  \label{eq:tte-generic-lower-functional}
\end{equation}
Lemma~\ref{lem:block-resolution} shows that the second term equals \(1\) when \(M(G)=0\) and is at least
\(1/(16M(G))\) when \(M(G)\geq1\).

\begin{lemma}[Packed-neighborhood TTE lower bound]
\label{lem:packed-neighborhood-lower}
For the full-neighborhood TTE on every graph \(G\),
\begin{equation}
  \Rallstar(\cP(\cZ);f,a_{\rm TTE})
  \geq
  \frac{\mathsf{Pack}_1(G)}{n^2}.
  \label{eq:packed-neighborhood-lower}
\end{equation}
Consequently the same lower bound holds for \(\Runbstar\).
\end{lemma}

\begin{proof}[Proof of Lemma~\ref{lem:packed-neighborhood-lower}]
Choose a set \(C\) of units attaining \(\mathsf{Pack}_1(G)\).  These units
serve as neighborhood centers.  Take the arm-paired
packing indexed by \(C\),
\[
  B_c^+=\{(i,1):i\in B_1(c)\},
  \qquad
  B_c^-=\{(i,0):i\in B_1(c)\},
\]
with positive perturbation signs. Each pair has per-arm weight $w_c=|B_1(c)|/n$, so $\sum_cw_c^2=\mathsf{Pack}_1(G)/n^2$.  The centers are at pairwise distance at
least three, so the balls \(B_1(c)\) are disjoint and the blocks are a packing.
No feasible hyperedge meets both $B_c^+$ and $B_c^-$. For any $i,j\in B_1(c)$, both closed neighborhoods contain $c$. Observing the treatment-arm potential outcome of $i$ requires $Z_c=1$, while observing the control-arm potential outcome of $j$ requires $Z_c=0$. Thus no assignment observes both parts of the pair, and the resolution capacity is $\mathcal M(\cB)=0$.
Theorem~\ref{thm:arm-paired-certificate}(i) gives
equation~\eqref{eq:packed-neighborhood-lower}.
\end{proof}

\begin{lemma}[Lower bound from block resolution]
\label{lem:block-resolution}
For the full-neighborhood TTE, with \(M(G)\) defined in
equation~\eqref{eq:tte-capacity}, the joint minimax risk satisfies
\begin{equation}
  \Rallstar(\cP(\cZ);f,a_{\rm TTE})
  \ \geq\ \max_{1\leq K\leq n}(K-M(G))_+\frac{\lfloor n/K\rfloor^2}{n^2}
  \ \geq\ \begin{cases}
  1, & M(G)=0,\\[4pt]
  \dfrac{1}{16M(G)}, & M(G)\geq1 .
  \end{cases}
  \label{eq:block-resolution-lower}
\end{equation}
\end{lemma}

\begin{proof}[Proof of Lemma~\ref{lem:block-resolution}]
Partition \(V\) into \(K\) blocks \(W_1,\dots,W_K\) of size at least
\(\lfloor n/K\rfloor\) and take the arm-paired packing
\(W_k^+=\{(i,1):i\in W_k\}\), \(W_k^-=\{(i,0):i\in W_k\}\) with positive signs,
so the per-arm weight is $w_k=|W_k|/n$.  Under the assignment \(z\) with \(T=\{i:z_i=1\}\), the sampled edge
$O_f(z)$ resolves pair $k$ when it meets both $W_k^+$ and $W_k^-$.
This happens exactly when $W_k$ meets both $A_T$ and $B_T$ in
equation~\eqref{eq:tte-capacity}. Each resolved pair uses at least one unit of $A_T$ and one unit of $B_T$. Since the $W_k$ are disjoint, the number of resolved pairs is at most $\min\{|A_T|,|B_T|\}\leq M(G)$. Theorem~\ref{thm:arm-paired-certificate}(ii) therefore gives
\[
 \Rallstar(\cP(\cZ);f,a_{\rm TTE})
 \geq(K-M(G))_+\frac{\lfloor n/K\rfloor^2}{n^2}.
\]
Maximizing over $K$ proves the first inequality.
If \(M(G)=0\), take \(K=1\).  Otherwise \(1\leq M(G)\leq n/2\), so
\(K=2M(G)\) is admissible and
$\lfloor n/K\rfloor\geq n/(2K)$ gives
\[
 (K-M(G))\frac{\lfloor n/K\rfloor^2}{n^2}
 \geq\frac{M(G)}{4K^2}=\frac1{16M(G)}.
\]
\end{proof}

Combining Lemmas~\ref{lem:packed-neighborhood-lower} and~\ref{lem:block-resolution} gives
\begin{equation}
 \Rallstar(\cP(\cZ);f,a_{\rm TTE})\geq L_{\rm TTE}(G).
 \label{eq:tte-generic-bracket}
\end{equation}
We next bound the variance of exposure HT under independent cluster assignments.

\begin{proposition}[TTE upper bound for graph-cluster randomization]
\label{prop:generic-gcr-tte-upper}
For the full-neighborhood TTE with \(|Y_i(1)|,|Y_i(0)|\leq1\), fix a partition
$\mathcal C$ and assign clusters independently to treatment with probability
$\rho\in(0,1)$. Exposure HT is unbiased and satisfies the following bound. Here
$\mathcal C_i$ is the set of clusters meeting $B_1(i)$ and
$\pi_i^\star=\min\{\rho^{|\mathcal C_i|},(1-\rho)^{|\mathcal C_i|}\}$:
\begin{equation}
  \sup_{Y\in\cY(f)}
  \Var(\widehat\tau_{\rm GCR,HT}\mid Y)
  \leq
  \frac4{n^2}
  \sum_{i,j:\mathcal C_i\cap\mathcal C_j\neq\varnothing}
  \frac1{\pi_i^\star\pi_j^\star}.
  \label{eq:generic-gcr-tte-upper}
\end{equation}
Consequently,
\begin{equation}
  L_{\rm TTE}(G)
  \leq
  \Rallstar(\cP(\cZ);f,a_{\rm TTE})
  \leq
  \Runbstar(\cP(\cZ);f,a_{\rm TTE})
  \leq
  \frac4{n^2}
  \sum_{i,j:\mathcal C_i\cap\mathcal C_j\neq\varnothing}
  \frac1{\pi_i^\star\pi_j^\star}.
  \label{eq:generic-gcr-tte-sandwich}
\end{equation}
\end{proposition}

\begin{proof}[Proof of Proposition~\ref{prop:generic-gcr-tte-upper}]
The exposure HT estimator is
\[
  \widehat\tau_{\rm GCR,HT}
  =
  \frac1n\sum_i
  \left[
  \frac{Y_i(1)\1\{f_i(Z)=1\}}{\pi_i^1}
  -
  \frac{Y_i(0)\1\{f_i(Z)=0\}}{\pi_i^0}
  \right].
\]
Write the estimator as $n^{-1}\sum_iX_i$, where $X_i$ is the expression in brackets. The exposure probabilities are
$\pi_i^1=\rho^{|\mathcal C_i|}$ and $\pi_i^0=(1-\rho)^{|\mathcal C_i|}$, so both are positive and $\E(X_i\mid Y)=Y_i(1)-Y_i(0)$. This proves unbiasedness.

For fixed $Y$, only the two exposure indicators in $X_i$ are random,
and both depend only on the assignments of clusters in $\mathcal C_i$.
Since at most one indicator is nonzero, $|X_i|\leq1/\pi_i^\star$. If $\mathcal C_i\cap\mathcal C_j=\varnothing$, the variables $X_i$ and $X_j$ are independent. For the remaining pairs, Cauchy--Schwarz gives
\[
 |\Cov(X_i,X_j\mid Y)|
 \leq\sqrt{\Var(X_i\mid Y)\Var(X_j\mid Y)}
 \leq\frac4{\pi_i^\star\pi_j^\star}.
\]
Summing over ordered pairs proves equation~\eqref{eq:generic-gcr-tte-upper}. Equation~\eqref{eq:tte-generic-bracket} and the unbiasedness of exposure HT give equation~\eqref{eq:generic-gcr-tte-sandwich}.
\end{proof}

\begin{lemma}[Restricted-growth packing]
\label{lem:restricted-growth-packing}
If \(G\) has restricted growth coefficient \(\kappa\) from radius one, then
\begin{equation}
  \mathsf{Pack}_1(G)
  \geq
  \frac1\kappa\sum_{i=1}^n |B_1(i)|
  =
  \frac{n(\bar d+1)}{\kappa}.
  \label{eq:restricted-growth-packing}
\end{equation}
\end{lemma}

\begin{proof}[Proof of Lemma~\ref{lem:restricted-growth-packing}]
Let \(w_i=|B_1(i)|\).  Greedily choose a remaining unit \(c\) of maximal
\(w_i\), add it to \(C\), and exclude the units in \(B_2(c)\).  The selected units have
pairwise distance at least three.  When \(c\) is selected, every newly excluded
unit $v\in D_c=B_2(c)\cap\{\text{remaining units}\}$ has \(w_v\leq w_c\), so the sum of neighborhood sizes over these excluded units is at
most
\[
  \sum_{v\in D_c}w_v
  \leq |B_2(c)|w_c
  \leq \kappa |B_1(c)|w_c
  =\kappa w_c^2 .
\]
The sets $D_c$ partition the units, since each unit is excluded once. Summing over the selected centers therefore gives
\(\sum_i w_i\leq\kappa\sum_{c\in C}w_c^2\), which is
equation~\eqref{eq:restricted-growth-packing}.
\end{proof}

\begin{proposition}[Restricted-growth TTE bounds]
\label{prop:restricted-growth-detail}
For the full-neighborhood TTE, if \(G\) has restricted growth
coefficient \(\kappa\), then
\begin{equation}
  \Rallstar(\cP(\cZ);f,a_{\rm TTE})
  \geq
  \frac{\bar d+1}{\kappa n}.
  \label{eq:restricted-growth-gcr-lower}
\end{equation}
Moreover, choose a maximal set of centers at pairwise distance at least
three and assign each unit to a nearest center, breaking ties by a fixed
rule. Assign clusters independently to treatment with probability
$\rho\in(0,1)$. The HT estimator under this design obeys
\begin{equation}
  \sup_{Y\in\cY(f)}
  \Var(\widehat\tau_{\rm GCR,HT}\mid Y)
  \leq
  C(\kappa,\rho)\frac{\bar d+1}{n}.
  \label{eq:restricted-growth-gcr-upper}
\end{equation}
Consequently, for fixed \(\kappa\) and fixed \(\rho\),
\[
  \Rallstar(\cP(\cZ);f,a_{\rm TTE})
  \asymp
  \Runbstar(\cP(\cZ);f,a_{\rm TTE})
  \asymp
  \frac{\bar d+1}{n}.
\]
For $d$-regular graphs and fixed $\kappa$, the rate is $(d+1)/n$.
\end{proposition}

\begin{proof}[Proof of Proposition~\ref{prop:restricted-growth-detail}]
Lemmas~\ref{lem:packed-neighborhood-lower} and~\ref{lem:restricted-growth-packing} give
\[
 \Rallstar(\cP(\cZ);f,a_{\rm TTE})
 \geq\frac{\mathsf{Pack}_1(G)}{n^2}
 \geq\frac{\bar d+1}{\kappa n}.
\]

For the upper bound, choose a maximal set of units at pairwise distance at
least three to serve as cluster centers.  Assign each unit to a nearest center, breaking ties by a fixed rule. Maximality ensures that every unit is within distance two of a center, so each
cluster is contained in the two-hop ball around its center.
Assign clusters independently to treatment with probability \(\rho\).

Let \(\mathcal C_i\) be the set of clusters intersecting \(B_1(i)\).  If a
cluster with center \(c\) intersects \(B_1(i)\), then
\(\operatorname{dist}_G(c,i)\leq3\).  The closed neighborhoods of these centers are disjoint. For each contributing center $c$, $B_1(i)\subseteq B_4(c)$ and restricted growth gives
$|B_1(i)|\leq|B_4(c)|\leq\kappa^3|B_1(c)|$. All these center neighborhoods lie in $B_4(i)$.  Hence
\[
  |\mathcal C_i|\kappa^{-3}|B_1(i)|
  \leq |B_4(i)|
  \leq \kappa^3|B_1(i)|,
  \qquad\text{so}\qquad
  |\mathcal C_i|\leq\kappa^6 .
\]
Therefore
\[
  \pi_i^1=\rho^{|\mathcal C_i|}\geq \rho^{\kappa^6},
  \qquad
  \pi_i^0=(1-\rho)^{|\mathcal C_i|}
  \geq (1-\rho)^{\kappa^6}.
\]
Let $\pi_\star=\min\{\rho^{\kappa^6},(1-\rho)^{\kappa^6}\}$. If $\mathcal C_i\cap\mathcal C_j\ne\varnothing$, a shared cluster has its center within distance three of both $i$ and $j$, so $\operatorname{dist}_G(i,j)\leq6$. Proposition~\ref{prop:generic-gcr-tte-upper} therefore gives
\begin{align*}
 \sup_Y\Var(\widehat\tau_{\rm GCR,HT}\mid Y)
 &\leq\frac4{\pi_\star^2n^2}\sum_i|B_6(i)|\\
 &\leq\frac{4\kappa^5}{\pi_\star^2n^2}\sum_i|B_1(i)|
 =\frac{4\kappa^5}{\pi_\star^2}\frac{\bar d+1}{n}.
\end{align*}
This proves equation~\eqref{eq:restricted-growth-gcr-upper} and the matching rates.
\end{proof}

\begin{resultrestatement}{Theorem~\ref{thm:restricted-growth-gcr-tte}}{Minimax rates for the neighborhood TTE}

The following hold.
\begin{enumerate}
\item[\textup{(i)}] For every graph $G$ with $M(G)\geq1$,
\[
  \max\left\{\frac{\mathsf{Pack}_1(G)}{n^2},
             \frac{1}{16M(G)}\right\}
  \leq\Rallstar(\cP(\cZ);f,a_{\rm TTE})
  \leq\Runbstar(\cP(\cZ);f,a_{\rm TTE}).
\]
If $M(G)=0$, then $\Rallstar(\cP(\cZ);f,a_{\rm TTE})\geq1$.
\item[\textup{(ii)}] If $G$ has restricted-growth coefficient $\kappa$, then
\[
 \frac{\bar d+1}{\kappa n}
 \leq\Rallstar(\cP(\cZ);f,a_{\rm TTE})
 \leq\Runbstar(\cP(\cZ);f,a_{\rm TTE})
 \leq C(\kappa)\frac{\bar d+1}{n},
\]
where $C(\kappa)$ depends only on $\kappa$.  For fixed $\kappa$, both minimax
risks therefore have rate $(\bar d+1)/n$.
\end{enumerate}
\end{resultrestatement}

\begin{proof}[Proof of Theorem~\ref{thm:restricted-growth-gcr-tte}]
Part~\textup{(i)} follows from Lemmas~\ref{lem:packed-neighborhood-lower} and~\ref{lem:block-resolution}. For part~\textup{(ii)}, apply Proposition~\ref{prop:restricted-growth-detail} with $\rho=1/2$. The lower bound has constant $1/\kappa$, and the upper-bound constant depends only on $\kappa$.
\end{proof}

\subsubsection{Comparison with the conflict-graph bounds}
\label{subsec:tte-khs-comparison}

\begin{remark}[Neighborhood packing in the conflict graph]
\label{rem:tte-conflict-biclique}
For full-neighborhood TTE, treated copy $i$ and control copy $j$ conflict
exactly when $B_1(i)\cap B_1(j)\ne\varnothing$.  There are no conflicts
within one arm.  If $W=B_1(c)$, every treated copy from $W$ conflicts with
every control copy from $W$, since the corresponding neighborhoods contain
$c$.  The paired-set lower bound therefore assigns squared per-arm mass
$|W|^2/n^2$ to this complete bipartite subgraph.  Disjoint neighborhoods
provide disjoint sets of potential outcomes, so the paired-set lower bounds add.
Thus the neighborhood-packing lower bound uses the complete bipartite
subgraphs associated with disjoint neighborhoods.
\end{remark}

All comparisons here use the uniformly bounded endpoint \(q=\infty\) of
\citet{KandirosHarshawSavje2026}.  For \(n\geq5\), let \(H\) denote their conflict graph on
\(2n\) target exposures.  Their global and local lower theorems give,
respectively,
\[
 L_{\rm glob}(H)=\sup_{T\subseteq V(H),\,|T|\geq5}
 \frac{(|T|/n)^2}{\alpha(H[T])},
\]
\[
 L_{\rm loc}(H)=\frac1{n^2}
 \sup_{T\subseteq V(H),\,|T|\geq5}
 \min\left\{d_{\min}(H,T)^2,
 \frac{|T|}{\sqrt{\sum_{v\in T}d_H(v)^{-2}}},\ |T|^2\right\}.
\]
Here \(d_{\min}(H,T)\) is the minimum over \(T\) of the degrees in \(H\).  Their unbiased upper theorem
specializes to \(d_{\rm avg}(H)/n\).  Thus, up to universal constants,
\[
 \max\{L_{\rm glob}(H),L_{\rm loc}(H)\}
 \lesssim R^*_{\rm all}\leq R^*_{\rm unb}
 \lesssim\frac{d_{\rm avg}(H)}n.
\]
For full-neighborhood TTE, put $\Phi(G)=n^{-2}\sum_i|B_2(i)|$. Then
\(d_{\rm avg}(H)=d_{\rm avg}(G^2)+1=n\Phi(G)\), including the own-unit
opposite-arm conflict.  Their \(q=\infty\) upper therefore gives
\(R^*_{\rm unb}\lesssim\Phi(G)\) on every graph.  Proposition~\ref{prop:restricted-growth-detail} gives a direct variance proof
for the restricted-growth class.

Suppose now that $G$ is $d$-regular with restricted-growth coefficient
$\kappa$.  Each vertex of $H$ has degree $|B_2(i)|$, and hence
$d+1\leq d_H(i)\leq\kappa(d+1)$.  For every $T\subseteq V(H)$, bipartiteness gives
$\alpha(H[T])\geq |T|/2$. Consequently,
\[
 L_{\rm glob}(H)\leq\frac4n.
\]
Taking $T=V(H)$ gives $L_{\rm glob}(H)=4/n$: one arm is an
independent set of size $n$, and the $n$ own-unit conflicts form a perfect
matching, so $\alpha(H)=n$.
The degree upper bound also gives, for every nonempty $T$,
\[
 \frac{|T|}{\sqrt{\sum_{v\in T}d_H(v)^{-2}}}
 \leq\kappa(d+1)\sqrt{|T|}
 \leq\kappa(d+1)\sqrt{2n}.
\]
The degree lower bound gives
\[
 L_{\rm loc}(H)\geq\frac1{n^2}
 \min\{(d+1)^2,(d+1)\sqrt{2n},4n^2\}.
\]
When $d+1\leq\sqrt n$, the global bound supplies $1/n$; when
$d+1>\sqrt n$, the local bound supplies order $(d+1)/n^{3/2}$.
Thus, for $n\geq5$,
\[
 \max\{L_{\rm glob}(H),L_{\rm loc}(H)\}
 \asymp_\kappa
 \max\left\{\frac1n,\frac{d+1}{n^{3/2}}\right\},
\]
while their unbiased upper bound is of order $(d+1)/n$. The resulting
lower and upper bounds differ by a factor of order
$\min\{d+1,\sqrt n\}$. Lemma~\ref{lem:packed-neighborhood-lower}
and the restricted-growth packing bound raise the lower bound to order
$(d+1)/n$, matching the upper bound.

\subsubsection{General upper bound and spectral lower bound}
\label{app:tte-general-refinements}
\label{sec:tte-one-rate}

We give an unconditional upper bound using small proposal probabilities, then bound $M(G)$ to obtain the spectral-expander lower bound.

\begin{proposition}[A general TTE upper bound]
\label{prop:tte-small-proposals}
Let \(D=\max_i(|B_2(i)|-1)\).  For the full-neighborhood TTE and target
potential outcomes bounded by one,
\[
 \Rallstar\leq\Runbstar\leq
 \min\left\{4,\frac{28(D+1)}n\right\}.
\]
In particular, every \(d\)-regular graph has
\(\Runbstar\lesssim\min\{1,(d^2+1)/n\}\).
\end{proposition}

\begin{proof}[Proof of Proposition~\ref{prop:tte-small-proposals}]
Apply ordered proposal thinning to \(H=G^2\), whose maximum degree is \(D\),
with fixed order and \(p=1/\{2(D+1)\}\).  Write \(E_i\) for the earlier
neighbors of \(i\) in \(H\).  The included units have disjoint closed
\(G\)-neighborhoods, so the TTE potential outcome requested by each arm
coin is observed.  Moreover,
\[
 q_i=p(1-p)^{|E_i|}\geq p(1-Dp)\geq p/2,
 \qquad \sum_iq_i^{-1}\leq4n(D+1).
\]
Units adjacent in $H$ have normalized selection covariance \(-1\).  For distinct
units nonadjacent in $H$, independence of the proposal coins gives
\[
 C_{ij}=(1-p)^{-|E_i\cap E_j|}-1
 \leq4p|E_i\cap E_j|.
\]
Indeed \(-\log(1-p)\leq2p\), \(2pD<1\), and
\(e^x-1\leq2x\) for \(0\leq x\leq1\).  Each vertex belongs to at most \(D\) of the sets \(E_i\), so
\[
 \sum_{i,j}|E_i\cap E_j|
 =\sum_h|\{i:h\in E_i\}|^2\leq nD^2.
\]
There are at most \(nD\) ordered adjacent pairs in \(H\). Therefore
\[
 \sum_{i\ne j}|C_{ij}|\leq nD+4pnD^2\leq3nD.
\]
Write $X_i$ for the selected-unit contribution before division by $n$:
\[
 X_i=\frac{2R_i}{q_i}\{\1\{A_i=1\}Y_i(1)-\1\{A_i=0\}Y_i(0)\},
\]
where $R_i$ indicates inclusion and $A_i$ is the independent fair arm coin.
The estimator $n^{-1}\sum_iX_i$ is unbiased, since
$\E[X_i\mid Y]=Y_i(1)-Y_i(0)$. Conditional averaging given the observed
data gives the unbiased estimator $\widehat\tau_{\rm RB}$, with no larger
variance, as in Proposition~\ref{prop:ordering-envelope}. To obtain the
stated constant, we bound the diagonal terms directly by
$\Var(X_i\mid Y)\leq \E[X_i^2\mid Y]\leq4/q_i$.
For the off-diagonal terms, when $i\ne j$,
$\Cov(X_i,X_j\mid Y)=C_{ij}\{Y_i(1)-Y_i(0)\}\{Y_j(1)-Y_j(0)\}$.
Summing these variance and covariance terms gives
\[
 \sup_Y\Var(\widehat\tau_{\rm RB}\mid Y)
 \leq\frac4{n^2}\left\{\sum_iq_i^{-1}+\sum_{i\ne j}|C_{ij}|\right\}
 \leq\frac{28(D+1)}n.
\]
Alternatively, assign the entire population to treatment or control with
one fair coin.  The resulting HT estimator is unbiased and has variance at
most \(4\).  Choose the strategy with the smaller bound.  Finally,
\(D\leq d^2\) on a \(d\)-regular graph.
\end{proof}

\begin{lemma}[Identification of $M(G)$]
\label{lem:m-identification}
Write $N[S]=\bigcup_{i\in S}B_1(i)$ and
$N^2[S]=\bigcup_{i\in S}B_2(i)$. Define the vertex-neighborhood profile by
$\psi(s)=\min\{|N[S]|:S\subseteq V,\ |S|=s\}$ for $s=0,\ldots,n$,
so that $\psi(0)=0$. Then
\begin{equation}
  M(G)
  =\max_{A\subseteq V}\min\bigl(|A|,\,n-|N^2[A]|\bigr)
  \ \leq\ \max\bigl\{m\in\{0,\ldots,\lfloor n/2\rfloor\}:\,m+\psi(\psi(m))\leq n\bigr\}.
  \label{eq:m-isoperimetric}
\end{equation}
Equivalently $M(G)$ is the largest $m$ for which $G^2$ contains two disjoint
$m$-sets with no edge between them---the bipartite hole number of $G^2$.
\end{lemma}

\begin{proof}[Proof of Lemma~\ref{lem:m-identification}]
If $A_T$ is the set of units whose closed neighborhoods lie in the treated
set $T$, then every vertex of $B_T$ has distance greater than two from
$A_T$.  Otherwise the two closed neighborhoods would share a vertex
required to be both treated and controlled.  Hence
$B_T\subseteq V\setminus N^2[A_T]$, giving one direction of the equality.
Conversely, for any $A\subseteq V$, treat $T=N[A]$.  Then $A\subseteq A_T$
and $B_T=V\setminus N^2[A]$, proving the reverse direction.  Taking subsets
of the larger of the two sets gives the bipartite-hole characterization.

The profile $\psi$ is nondecreasing: every larger set contains a smaller
one whose closed neighborhood it contains.  Thus, for $|A|=m$,
$|N^2[A]|\geq\psi(|N[A]|)\geq\psi(\psi(m))$.  Two separated $m$-sets
therefore require $m+\psi(\psi(m))\leq n$, proving the inequality.
\end{proof}

\begin{corollary}[Spectral upper bounds on \(M\)]
\label{cor:m-spectral}
If \(G^2\) is \(d_2\)-regular with \(d_2>0\) and largest absolute
nonconstant adjacency eigenvalue \(\lambda_*\), then
\(M(G)\leq\lambda_*n/d_2\).
Alternatively, if \(G\) is \(d\)-regular with \(d>0\) and largest absolute
nonconstant adjacency eigenvalue \(\lambda\), then
\[
 M(G)\leq\frac{n\lambda^2}{d^2+\lambda^2}.
\]
In particular, \(\lambda\leq C\sqrt d\) implies
\(\Rallstar\gtrsim d/n\), without a girth assumption.
\end{corollary}
\begin{proof}
For two disjoint \(m\)-sets with no edge between them in \(G^2\), expander mixing gives
\(d_2m^2/n\leq\lambda_*m\), proving the first bound.
For the second, let \(A,B\) be two \(m\)-sets at distance greater than two
in \(G\), and let \(Q\) be the adjacency matrix of \(G\).  There are no
length-two walks between them, so \(\1_A^\top Q^2\1_B=0\).
Write \(x=\1_A-(m/n)\1\) and \(y=\1_B-(m/n)\1\).  Regularity gives
\[
 0=\frac{d^2m^2}{n}+x^\top Q^2y,
 \qquad
 |x^\top Q^2y|\leq\lambda^2\|x\|\|y\|
 =\lambda^2m(1-m/n).
\]
Rearranging yields \(m\leq n\lambda^2/(d^2+\lambda^2)\).
Lemma~\ref{lem:m-identification} now gives both upper bounds on $M(G)$.
If $\lambda\leq C\sqrt d$ and $M(G)\geq1$,
Lemma~\ref{lem:block-resolution} yields
\[
 \Rallstar\geq\frac1{16M(G)}
 \geq\frac{d^2+\lambda^2}{16n\lambda^2}
 \geq\frac{d}{16C^2n}.
\]
If $M(G)=0$, the same lemma gives $\Rallstar\geq1\geq d/n$.
\end{proof}

\begin{corollary}[TTE on spectral expanders]
\label{cor:tte-spectral-bracket}
Suppose \(G\) is \(d\)-regular, \(d\geq1\), and every nonconstant
adjacency eigenvalue has absolute value at most \(C_0\sqrt d\).  Then
\[
 c\frac dn\leq\Rallstar\leq\Runbstar
 \leq C\min\left\{1,\frac{d^2+1}{n}\right\},
\]
where constants depend only on \(C_0\).
\end{corollary}
\begin{proof}
The lower bound follows from Corollary~\ref{cor:m-spectral} with $C=C_0$.
For the upper bound, Proposition~\ref{prop:tte-small-proposals} and $D\leq d^2$
give
\[
 \Runbstar\leq\min\left\{4,\frac{28(d^2+1)}n\right\}
 \leq28\min\left\{1,\frac{d^2+1}n\right\}.
\]
\end{proof}

\subsection{Proofs for the hard-budget TTE example}
\label{subsec:app-hard-budget-tte}

The lit-ball design samples a fixed number of units and treats their closed
neighborhoods. We account for sampling without replacement when bounding
the covariance of the control-arm observation indicators.
We first prove the covariance bound in Lemma~\ref{lem:lit-ball-covariance-control}, then apply it to prove Theorem~\ref{thm:hard-budget-regular-tte}.  Proposition~\ref{prop:hard-budget-expanding-tte} gives the subsequent expansion-based refinement.

Recall the lit-ball design in Definition~\ref{def:hard-budget-lit-ball}, and
write $P_\ell$ for its law.  Since
$|T|\leq\ell(d+1)\leq m$, it respects the treatment budget.

For each unit $i$, define
\[
 A_i=\{j:B_1(j)\cap B_1(i)\neq\varnothing\}=B_2(i),
 \qquad
 q_i=\Pr(C_{\rm lit}\cap A_i=\varnothing).
\]
Each sampled unit attains the treated arm, while unit $i$ attains the
control arm exactly when $C_{\rm lit}\cap A_i=\varnothing$.  The corresponding
augmented-data estimator is
\begin{equation}
  \widehat\tau_{\rm ball}
  =\frac1n\sum_{i=1}^n
  \left\{\frac{\1\{i\in C_{\rm lit}\}}{\ell/n}\,Y_i^{\rm obs}
        -\frac{\1\{C_{\rm lit}\cap A_i=\varnothing\}}{q_i}\,Y_i^{\rm obs}\right\}.
  \label{eq:lit-ball-estimator}
\end{equation}
It is unbiased whenever every $q_i>0$.  Because $C_{\rm lit}$ cannot
generally be reconstructed from $Z$, the observed-data estimator is
$\widehat\tau_{{\rm ball},{\rm RB}}
=\E[\widehat\tau_{\rm ball}\mid Z,Y_{O_f(Z)}]$; Lemma~\ref{lem:rao-blackwell}
preserves unbiasedness and the variance bound.

\begin{remark}[The budgeted design and thinning]
\label{rem:lit-ball-no-thinning}
The ordering design on $G^2$ keeps selected neighborhoods disjoint
because each selected neighborhood receives an independent arm assignment.  If
two selected neighborhoods overlapped and received different arms, their
assignments would conflict on the shared units.

The lit-ball design uses selected neighborhoods only to observe treated-arm
potential outcomes.  Overlapping balls therefore make the same treatment request on
their shared units.  Control-arm potential outcomes are observed separately for units
whose neighborhoods avoid the entire treated union.  Thus the lit-ball design can use overlapping sampled neighborhoods.
\end{remark}

\begin{lemma}[Covariance control for lit balls]
\label{lem:lit-ball-covariance-control}
Let \(A_1,\ldots,A_n\subseteq\{1,\ldots,n\}\), let \(B=\max_i|A_i|\), and let
\(C_{\rm lit}\) be a uniform \(\ell\)-subset of \(\{1,\ldots,n\}\), where $1\leq\ell\leq n$.
Assume \(\sum_{i,j}|A_i\cap A_j|\leq nB^2\). Put
\(q_i=\Pr(C_{\rm lit}\cap A_i=\varnothing)\) and
\(W_i=\1\{C_{\rm lit}\cap A_i=\varnothing\}/q_i\).  If
\(\ell B\leq c_0n\) for a
sufficiently small universal constant \(c_0\), then
\begin{equation}
  \sum_{i,j=1}^n|\Cov(W_i,W_j)|\leq C\,\ell B^2,
  \qquad\text{so}\qquad
  \frac1{n^2}\sum_{i,j=1}^n|\Cov(W_i,W_j)|\leq\frac C\ell,
  \label{eq:lit-ball-covariance-control}
\end{equation}
with \(C\) universal.
\end{lemma}

\begin{proof}[Proof of Lemma~\ref{lem:lit-ball-covariance-control}]
For \(A,B'\subseteq\{1,\ldots,n\}\) of size at most $B$, write \(a=|A|\), \(b=|B'|\), and
\(t=|A\cap B'|\).  If \(q(A)=\Pr(C_{\rm lit}\cap A=\varnothing)\), then
\[
  q(A)=\prod_{r=0}^{\ell-1}
  \left(1-\frac{|A|}{n-r}\right).
\]
The case $B=0$ is trivial; assume $B,\ell\geq1$.  Then
$\ell\leq c_0n$, $B\leq c_0n$, and the factors in the avoidance probabilities for $A$, $B'$, and
$A\cup B'$ are bounded away from zero.  For the standardized avoidance indicators associated with $A$ and $B'$, their covariance is
\[
 \frac{q(A\cup B')}{q(A)q(B')}-1.
\]
We bound the log of the ratio, separating the effect of overlap from
the dependence due to sampling without replacement. For $u=n-r$, first put $t=0$.  The ratio of the three factors is
\[
 \frac{1-(a+b)/u}{(1-a/u)(1-b/u)}
 =1-\frac{ab}{(u-a)(u-b)}.
\]
This term can be nonzero even for disjoint sets. Its logarithm has absolute value at most $Cab/u^2$.  Replacing $a+b$ by
$a+b-t$ changes this logarithm by at most $Ct/u$.  Summing over
$0\leq r<\ell$ bounds the absolute log-ratio by
$C(\ell t/n+\ell B^2/n^2)$.  This is bounded by a universal constant under
$\ell B\leq c_0n$, so $|e^x-1|\leq C|x|$ gives
\begin{equation}
  \left|
  \frac{q(A\cup B')}{q(A)q(B')}-1
  \right|
  \leq
  C\left\{\frac{\ell t}{n}+\frac{\ell B^2}{n^2}\right\}.
  \label{eq:no-hit-pair-bound}
\end{equation}
This is exactly the absolute covariance of the standardized indicators
associated with \(A\) and \(B'\).  Applying
equation~\eqref{eq:no-hit-pair-bound} to \(A_i,A_j\) and summing,
\[
  \sum_{i,j}|\Cov(W_i,W_j)|
  \leq
  C\left\{
    \frac{\ell}{n}\sum_{i,j}|A_i\cap A_j|
    + \ell B^2
  \right\}.
\]
The assumed intersection bound \(\sum_{i,j}|A_i\cap A_j|\leq nB^2\) proves
equation~\eqref{eq:lit-ball-covariance-control}.  In the lit-ball application,
\(A_i=B_2(i)\) in an undirected graph, and the intersection bound follows from symmetry:
\[
  \sum_{i,j}|A_i\cap A_j|
  =
  \sum_v |\{i:v\in B_2(i)\}|^2
  =
  \sum_v |B_2(v)|^2
  \leq nB^2 .
\]
The rate statement follows from \(\ell B\leq c_0n\), which implies
\(\ell^2B^2\leq c_0^2n^2\).
\end{proof}

\begin{resultrestatement}{Theorem~\ref{thm:hard-budget-regular-tte}}{Budgeted TTE under restricted growth}
Let $G$ be $d$-regular with restricted-growth coefficient $\kappa$. There is
$c_{\rm rg}(\kappa)\in(0,1]$ such that, if $1\leq m\leq n/2$ and
$(d+1)/m\leq c_{\rm rg}(\kappa)$, then
\[
 c(\kappa)\frac{d+1}{m}
 \leq\Rallstar(\cP(\cZ_{\leq m});f,a_{\rm TTE})
 \leq\Runbstar(\cP(\cZ_{\leq m});f,a_{\rm TTE})
 \leq C(\kappa)\frac{d+1}{m}.
\]
\end{resultrestatement}

\begin{proof}[Proof of Theorem~\ref{thm:hard-budget-regular-tte}]
Write $w=d+1$ and $J=\lfloor m/w\rfloor\geq1$.  Choose a maximal set $C$ of
units at pairwise distance at least five, and put $L=|C|$.  The balls
$B_2(c)$, $c\in C$, are disjoint, while maximality implies that the balls
$B_4(c)$ cover the population.  Restricted growth therefore gives
\begin{equation}
 L\geq \frac{n}{\kappa^3w}.
 \label{eq:budget-packing-size}
\end{equation}

Suppose first that $L\geq2J$.  Partition $C$ into $2J$ groups, each containing
at least $L/(4J)$ units.  For each group $H$ of centers, let
\[
 B_H=\{(i,1):i\in\textstyle\bigcup_{c\in H}B_1(c)\}
\]
be its treated-arm perturbation block. These blocks are disjoint and each
has total coefficient magnitude at least $wL/(4Jn)$. If the sampled
edge meets $B_H$ at a vertex $(i,1)$ with $i\in B_1(c)$, then the
assignment treats all $w$ units in $B_1(i)\subseteq B_2(c)$. The
balls $B_2(c)$ are disjoint across centers, so each block met costs
at least $w$ distinct treated units. Thus the sampled edge meets
at most $J$ blocks.  At least $J$ perturbation sets are therefore missing under
every assignment.  The lower bound from disjoint perturbations (Theorem~\ref{thm:fixed-joint-missingness}(i)) gives
\[
 \Rall(P;f,a_{\rm TTE})
 \geq J\left(\frac{wL}{4Jn}\right)^2
 \geq \frac{w}{16\kappa^6m}
\]
for every design $P$ supported on $\cZ_{\leq m}$.

If $L<2J$, then equation~\eqref{eq:budget-packing-size} implies
$m>n/(2\kappa^3)$.  The unrestricted lower bound in
Theorem~\ref{thm:restricted-growth-gcr-tte}(ii), which remains valid for the
smaller budgeted design class, yields
\[
 \Rallstar(\cP(\cZ_{\leq m});f,a_{\rm TTE})
 \geq \frac{w}{\kappa n}
 \geq \frac{w}{2\kappa^4m}.
\]

For the upper bound, restricted growth gives
$B_2^{\max}:=\max_i|B_2(i)|\leq\kappa w$.  Take
$c_{\rm rg}(\kappa)\leq1$ small enough that $w/m\leq c_{\rm rg}(\kappa)$ implies
$B_2^{\max}\leq c_0n$.  Use the lit-ball design with
\[
 \ell=1\vee\left\lfloor
 c_0\min\left\{\frac{m}{w},\frac{n}{B_2^{\max}}\right\}
 \right\rfloor,
\]
where $c_0$ is the constant in
Lemma~\ref{lem:lit-ball-covariance-control}.  The choice of $\ell$ satisfies $\ell w\leq m$ and
$\ell B_2^{\max}\leq c_0n$, including the case $\ell=1$.
\emph{Lit-ball variance calculation.} The treatment-arm term in equation~\eqref{eq:lit-ball-estimator} is the
sample mean of $\ell$ uniformly sampled treatment-arm potential outcomes.
Its variance is at most $C/\ell$. For the control-arm term,
Lemma~\ref{lem:lit-ball-covariance-control} gives
\[
 \Var\left(\frac1n\sum_iW_iY_i(0)\,\middle|\,Y\right)
 \leq\frac1{n^2}\sum_{i,j}|\Cov(W_i,W_j)|\leq\frac C\ell.
\]
The variance of the difference of the two terms is at most twice the sum
of their variances, and is therefore at most $C/\ell$.
Since $m\leq n/2$ and $B_2^{\max}\leq\kappa w$,
\[
 \frac1\ell\leq C\max\left\{\frac wm,\frac{B_2^{\max}}n\right\}
 \leq C(\kappa)\frac wm.
\]
Lemma~\ref{lem:rao-blackwell} transfers unbiasedness and the variance bound
to the observed-data estimator.
\end{proof}

\begin{definition}[Budget expansion]
\label{def:budget-expansion}
A $d$-regular graph is \emph{$(m,c_{\rm exp})$-budget expanding} if
\begin{equation}
 |N[S]|\geq c_{\rm exp}(d+1)|S|
 \qquad\text{whenever } |N[S]|\leq m,
 \label{eq:budget-expansion}
\end{equation}
where $N[S]=\bigcup_{i\in S}B_1(i)$.
\end{definition}

\begin{proposition}[Budgeted TTE on expanding regular graphs]
\label{prop:hard-budget-expanding-tte}
Let $G$ be $d$-regular, $d\geq1$, and
$(m,c_{\rm exp})$-budget expanding, with $d+1\leq m\leq n/2$.  Put
$B=\max_i|B_2(i)|\leq1+d^2$.  Then
\[
 \begin{aligned}
 c\max\left\{\frac{d+1}{m},L_{\rm TTE}(G)\right\}
 &\leq\Rallstar(\cP(\cZ_{\leq m});f,a_{\rm TTE})\\
 &\leq\Runbstar(\cP(\cZ_{\leq m});f,a_{\rm TTE})
 \leq C\max\left\{\frac{d+1}{m},\frac Bn\right\},
 \end{aligned}
\]
where the constants may depend on $c_{\rm exp}$.  In particular, the rate is
$(d+1)/m$ when $mB\leq C_1n(d+1)$ for a fixed $C_1$.
\end{proposition}

\begin{proof}[Proof of Proposition~\ref{prop:hard-budget-expanding-tte}]
For the lower bound, fix an admissible assignment with treated set \(T\), and
let \(A_T=\{i:B_1(i)\subseteq T\}\) be the treated-exposed units.  Since
\(N[A_T]\subseteq T\), equation~\eqref{eq:budget-expansion} gives
\[
  |A_T|\leq \frac{|T|}{c_{\rm exp}(d+1)}
  \leq \frac{m}{c_{\rm exp}(d+1)}.
\]
Define the largest number of treatment-arm potential outcomes observable under the
budget by
\[
 J=\max_{T\subseteq V:\,|T|\leq m}|A_T|.
\]
Treating one closed neighborhood shows that $J\geq1$. The bounds
$|A_T|\leq|T|\leq m$ and $|A_T|\leq m/\{c_{\rm exp}(d+1)\}$ give
\[
 J\leq\min\left\{m,\left\lfloor\frac{m}{c_{\rm exp}(d+1)}\right\rfloor\right\}
 \leq n/2.
\]
Partition the $n$ units into $2J$ blocks, each of size at least
$\lfloor n/(2J)\rfloor$. Put an independent sign on the treatment-arm
potential outcomes in each block and set all other potential outcomes to zero.
Every assignment observes potential outcomes from at most $J$ blocks, so at least
$J$ signs remain hidden. Theorem~\ref{thm:fixed-joint-missingness}(i) gives
\[
 \Rall(P;f,a_{\rm TTE})
 \geq\frac{J}{n^2}\left\lfloor\frac n{2J}\right\rfloor^2
 \geq\frac1{16J}
 \geq\frac{c_{\rm exp}(d+1)}{16m}
\]
for every admissible design $P$. The unrestricted lower bound
$L_{\rm TTE}(G)$ also holds for the smaller budgeted design class.
Taking the maximum of the two lower bounds proves the lower inequality.

For the upper bound, the lit-ball design is hard-budget feasible by
Definition~\ref{def:hard-budget-lit-ball}.  Write
\(B_2^{\max}=\max_i|B_2(i)|\leq 1+d^2\), and choose
\[
  \ell
  =
  1\vee
  \left\lfloor
    c_0\min\left\{\frac{m}{d+1},\,\frac{n}{B_2^{\max}}\right\}
  \right\rfloor
\]
in the case \(B_2^{\max}\leq c_0n\).  Then \(\ell\leq\lfloor m/(d+1)\rfloor\) after reducing
\(c_0\) if necessary, so the design is hard-budget feasible.  The choice also gives $\ell B_2^{\max}\leq c_0n$, so the same lit-ball
variance calculation as in the proof of
Theorem~\ref{thm:hard-budget-regular-tte} gives worst-case variance
at most $C/\ell$ for $\widehat\tau_{{\rm ball},{\rm RB}}$, and
\[
  \frac1\ell
  \leq
  C\max\left\{\frac{d+1}{m},\,\frac{B_2^{\max}}{n}\right\}
  \leq
  C\max\left\{\frac{d+1}{m},\,\frac{d^2}{n}\right\}.
\]
For the remaining case \(B_2^{\max}>c_0n\), augment the lit-ball design as
follows.  Let \(F\) be a fair branch coin.  When \(F=1\), draw a uniformly sampled unit
\(I\) and treat \(B_1(I)\); when \(F=0\), assign all units to control.  The
augmented-data estimator is
\[
 \begin{cases}
 2Y_I^{\rm obs},&F=1,\\
 -2n^{-1}\sum_iY_i^{\rm obs},&F=0.
 \end{cases}
\]
It is unbiased and bounded in absolute value by two, so its variance is at
most four.  Its conditional expectation given the observed data has no larger
variance.  Both branches satisfy the budget because \(m\geq d+1\), and
\(4\leq (4/c_0)B_2^{\max}/n\) proves the claimed upper bound.  Finally,
\(mB\leq C_1n(d+1)\) makes the two upper terms comparable to the budget
lower term and proves the claimed rate.
\end{proof}

\subsection{The direct treatment effect under a treatment budget}
\label{subsec:app-budgeted-dte}

Restrict the assignment space to
$\cZ_{\leq m}=\{z\in\{0,1\}^n:|T(z)|\leq m\}$, as in
Section~\ref{subsec:hard-budget-regular-tte}. The DTE exposure mapping is
unchanged, but the budget limits how many treatment-arm potential outcomes
can be observed together.

\begin{proposition}[Budgeted DTE]
\label{prop:budgeted-dte}
Let $1\leq m\leq n$ be an integer. For the DTE on any graph $G$, the sampling rank of the treated arm over
\(\cZ_{\leq m}\) is
\[
 r_m=\min\{m,\alpha(G)\}.
\]
Consequently, for a universal constant $c>0$,
\begin{equation}
\begin{aligned}
 \Rallstar\bigl(\cP(\cZ_{\leq m});f,a_{\rm DTE}\bigr)
 &\geq \frac{c}{\min\{m,\alpha(G)\}},\\
 \Runbstar\bigl(\cP(\cZ_{\leq m});f,a_{\rm DTE}\bigr)
 &\geq \frac{c}{\min\{m,\alpha(G)\}}.
\end{aligned}
 \label{eq:budgeted-dte-lower}
\end{equation}
If \(G\) is \(d\)-regular, with \(d\) allowed to depend on \(n\), then
\begin{equation}
\begin{aligned}
 \max\left\{\frac1m,\frac1{\alpha(G)}\right\}
 &\lesssim
 \Rallstar\bigl(\cP(\cZ_{\leq m});f,a_{\rm DTE}\bigr)
 \leq
 \Runbstar\bigl(\cP(\cZ_{\leq m});f,a_{\rm DTE}\bigr),\\
 \Runbstar\bigl(\cP(\cZ_{\leq m});f,a_{\rm DTE}\bigr)
 &\lesssim
 \max\left\{\frac1m,\frac{d+1}{n}\right\}.
\end{aligned}
 \label{eq:budgeted-dte-rate}
\end{equation}
In particular, whenever \(m(d+1)\leq c_0n\) for a sufficiently small
universal constant \(c_0>0\), both minimax risks are of order \(1/m\).
\end{proposition}

\begin{proof}
Write $V_+=\{(i,(1,0_{N_i})):i\in[n]\}$ for the treated arm. An assignment observes the treated DTE potential outcome of unit $i$ exactly when
$i$ is treated and all its neighbors are controlled.  Hence the units whose
treated potential outcomes are observed form an independent set and have cardinality
at most $m$.  Conversely, treating any independent set of size at most $m$
observes all of its treated potential outcomes.  The treated-arm sampling rank is
therefore $r(V_+)=r_m=\min\{m,\alpha(G)\}$ over $\cZ_{\le m}$.

If $r_m\leq n/2$, the sampling-rank lower bound (Theorem~\ref{cor:rank-lower}) gives the all-estimator lower
bound $c/r_m$. If $r_m>n/2$, the pairwise conflicts between the treated
and control DTE potential outcomes of each unit give an all-estimator lower bound of order $1/n$,
as in the proof of Proposition~\ref{prop:dte-general-rank}; since $r_m>n/2$,
$1/n\geq1/(2r_m)$.  The same bound holds for
the smaller class of unbiased estimators.

For the upper bound on a \(d\)-regular graph, fix a sufficiently small
universal constant \(c_1>0\) and first suppose
\(d+1\leq c_1n/2\).  Put

\[
 t=\min\left\{m,\left\lfloor\frac{c_1n}{d+1}\right\rfloor\right\},
\]

and use complete randomization of exactly \(t\) treated units, which is
supported on \(\cZ_{\leq m}\).  The two DTE exposure probabilities are
\[
 \pi_{i,+}
 =\frac tn\frac{(n-t)_d}{(n-1)_d}
 \asymp \frac tn,
 \qquad
 \pi_{i,-}
 =\frac{(n-t)_{d+1}}{(n)_{d+1}}
 \asymp1,
\]
because \(t(d+1)/n\leq c_1\); here
\((x)_s=x(x-1)\cdots(x-s+1)\).  It follows that the diagonal term in the HT
bound satisfies \(D_a(P)=O(1/t)\).

We next bound the positive normalized covariance within each arm.
For $k=0,1,2$, write
\[
 \Psi_s^{(k)}=\frac{(n-t)_s}{(n-k)_s}.
\]
When $t\geq k$, this is the probability that $s$ specified units are
controlled, conditional on $k$ other specified units being treated.
For distinct nonadjacent units $i,j$, put $h=|N_i\cap N_j|$. If $t\geq2$,
the treatment-arm joint probability and normalized covariance satisfy
\[
 \pi_{ij,+}=\frac{t(t-1)}{n(n-1)}\Psi_{2d-h}^{(2)},
 \qquad
 1+C_{ij,+}
 =\frac{(t-1)n}{t(n-1)}
   \frac{\Psi_{2d-h}^{(2)}}{(\Psi_d^{(1)})^2}.
\]
For either arm, conditional on $k$ specified treated units, the
remaining $t-k$ treated units form a uniform subset of the other
$n-k$ units. For two disjoint sets of size $s$ in this pool,
conditioning on avoidance of the first set leaves the same number
of draws in a pool smaller by $s$. Thus the conditional probability
of avoiding the second set satisfies
\[
 \prod_{r=0}^{t-k-1}\left(1-\frac{s}{n-k-s-r}\right)
 \leq
 \prod_{r=0}^{t-k-1}\left(1-\frac{s}{n-k-r}\right)
 =\Psi_s^{(k)}.
\]
In the cases used here all factors are nonnegative. Multiplying
by the probability of avoiding the first set gives
$\Psi_{2s}^{(k)}\leq(\Psi_s^{(k)})^2$. In particular,
$\Psi_{2d}^{(2)}\leq(\Psi_d^{(2)})^2$. Moreover,
\[
 \frac{\Psi_d^{(2)}}{\Psi_d^{(1)}}=\frac{n-1}{n-d-1},
 \qquad
 \log\frac{\Psi_{2d-h}^{(2)}}{\Psi_{2d}^{(2)}}
 =\sum_{r=2d-h}^{2d-1}\log\frac{n-2-r}{n-t-r}
 \leq C\frac{th}{n}.
\]
All denominators are a fixed fraction of $n$ when $c_1$ is sufficiently
small. The prefactor $(t-1)n/\{t(n-1)\}$ is at most one, so
\[
 (C_{ij,+})_+\leq C\left(\frac dn+\frac{th}{n}\right).
\]
Here we used $d/n+th/n\leq Cc_1$ to pass from a log-ratio to the ratio.
If $t=1$, or if $i,j$ are adjacent, the two treatment exposures cannot
occur together, so their covariance is nonpositive.

For the control arm, put $w=d+1$ and $h=|B_1(i)\cap B_1(j)|$. Then
\[
 1+C_{ij,-}=\frac{\Psi_{2w-h}^{(0)}}{(\Psi_w^{(0)})^2}
 \leq\frac{\Psi_{2w-h}^{(0)}}{\Psi_{2w}^{(0)}},
 \qquad
 \log\frac{\Psi_{2w-h}^{(0)}}{\Psi_{2w}^{(0)}}
 \leq C\frac{th}{n}.
\]
The inequality again follows from avoidance of disjoint sets under sampling
without replacement. Thus $(C_{ij,-})_+\leq Cth/n$.
Since $d/n\leq c_1/t$, both arms satisfy
\[
 (C_{ij,\pm})_+
 \leq C\left\{\frac1t+\frac tn|B_1(i)\cap B_1(j)|\right\}.
\]
Since regularity gives
\[
 \sum_{i,j}|B_1(i)\cap B_1(j)|=n(d+1)^2,
\]
summing with DTE coefficients of magnitude \(1/n\) gives
\[
 \mathrm{OD}_a^+(P)
 =O\left(\frac1t+\frac{t(d+1)^2}{n^2}\right)
 =O(1/t).
\]
The exposure-wise HT bound in Proposition~\ref{supp:prop:exposure-wise-ht}, with two
target exposures, now gives \(V_{\rm HT}(P)=O(1/t)\).  Since
\[
 \frac1t\asymp
 \max\left\{\frac1m,\frac{d+1}{n}\right\},
\]
this proves the upper side of equation~\eqref{eq:budgeted-dte-rate} when
\(d+1\leq c_1n/2\).  If \(d+1>c_1n/2\), assign probability one half to the
all-control assignment and split the remaining probability uniformly over the
\(n\) assignments that treat one unit.  Exposure HT then has risk \(O(1)\),
which is also \(O((d+1)/n)\).  Finally, if
\(m(d+1)\leq c_0n\), the Caro--Wei bound
\(\alpha(G)\geq n/(d+1)\) gives \(r_m=m\), so the lower and upper bounds both
have order \(1/m\).
\end{proof}

\subsection{Details for the complex-contagion example}
\label{subsec:app-complex-contagion}

\label{subsec:contagion-additional-results}
We prove the rank and packing lemmas used in the lower bound of Theorem~\ref{thm:contagion-regular-lower}, then the Bernoulli upper bound in Theorem~\ref{thm:contagion-regular-upper}. Together these prove Theorem~\ref{thm:contagion-main-bracket}. We also prove Corollary~\ref{cor:contagion-conflict-gap} and the clique and complete-bipartite community rates.

For an integer $s\geq0$, write
$\operatorname{diss}_s(G)=\max\{|S|:\Delta(G[S])\leq s\}$.
Thus $\operatorname{diss}_0(G)=\alpha(G)$.  Let
$L_0(z)=\{i:z_i=1,C_i(z)\leq\theta-1\}$ and
$\mathcal L_\theta(G)=\{S\subseteq V:\Delta(G[S])\leq\theta-1\}$.
The assignment space is the full binary cube, as in the main-text example.

\subsubsection{Arm ranks and the sampling reduction}

\begin{lemma}[Single-arm rank of the treated-high threshold arm]
\label{lem:contagion-high-free}
Let $G$ be $d$-regular with $d\geq\theta$. Then $r(b_{\rm hi})=n$,
and the minimax risk of estimating the high-arm mean on its own is zero.
\end{lemma}

\begin{proof}[Proof of Lemma~\ref{lem:contagion-high-free}]
Under \(z=\1_V\) every unit is treated and has \(C_i(z)=d\geq\theta\), so every
unit attains \(b_{\rm hi}\).  A single deterministic assignment therefore
observes the whole high arm, the high-arm sampling rank is \(n\), and the
degenerate design placing all mass on \(\1_V\) estimates the high-arm mean
without error.
\end{proof}

\begin{lemma}[Sampling rank of the treated-low threshold arm]
\label{lem:contagion-low-capacity}
For every graph \(G\) and every \(\theta\geq1\), the sampling rank of
the treated-low threshold arm \(b_{\rm lo}\) is
\begin{equation}
  r(b_{\rm lo})
  =
  \operatorname{diss}_{\theta-1}(G).
  \label{eq:contagion-low-rank}
\end{equation}
Consequently, when both target exposures are attainable at every unit,
\begin{equation}
  \Rallstar(\cP(\cZ);f,a_\theta)
  \geq
  \frac{c}{\operatorname{diss}_{\theta-1}(G)}
  \geq
  \frac{c}{\theta\,\alpha(G)}
  \label{eq:contagion-capacity-lower}
\end{equation}
for a universal constant \(c>0\).
\end{lemma}

\begin{proof}[Proof of Lemma~\ref{lem:contagion-low-capacity}]
Let \(S\) be a set of units whose treated-low threshold potential outcomes are observed
simultaneously.  Every unit in \(S\) must be treated.  Treating units outside
\(S\) can only increase the numbers of treated neighbors, so the largest feasible
observed set is achieved by assignments of the form \(z=\1_S\).  Under this
assignment, each \(i\in S\) has exactly \(|N_i\cap S|\) treated neighbors.
Thus all units in \(S\) attain the low arm if and only if
\(\Delta(G[S])\leq\theta-1\), proving
equation~\eqref{eq:contagion-low-rank}.

Every graph with maximum degree at most $\theta-1$ contains an independent
set of size at least a $1/\theta$ fraction of its vertices. Hence
$\operatorname{diss}_{\theta-1}(G)\leq\theta\alpha(G)$.
Write $r=\operatorname{diss}_{\theta-1}(G)$. If $r\leq n/2$,
Theorem~\ref{cor:rank-lower} gives $\Rallstar\geq c/r$.
If $r>n/2$, write $\pi_i^{\rm lo},\pi_i^{\rm hi}$ for the two exposure probabilities under any design $P$. Since the two exposures are mutually
exclusive, equation~\eqref{eq:singleton-all} gives
\[
 \Rall(P;f,a_\theta)
 \geq\frac1{n^2}\sum_{i=1}^n(2-\pi_i^{\rm lo}-\pi_i^{\rm hi})
 \geq\frac1n\geq\frac1{2r}.
\]
Taking the infimum over designs and using $r\leq\theta\alpha(G)$ proves
equation~\eqref{eq:contagion-capacity-lower}.
\end{proof}

\paragraph{Reduction to low-arm sampling.}
The lower bound in Theorem~\ref{thm:contagion-regular-lower} is proved through a
sampling problem that keeps only the low arm.  For a law \(Q\) on
\(\mathcal L_\theta(G)\) and an estimator \(\widehat\mu\) of the finite
population mean from the sampled units and their values, define
\begin{equation}
  \mathsf S_\theta(G)
  =
  \inf_{Q}\ \inf_{\widehat\mu}\
  \sup_{x\in[-1,1]^V}\
  \E_{S\sim Q}
  \left[\left\{\widehat\mu(S,x_S)-\frac1n\sum_{i=1}^n x_i\right\}^2\right] .
  \label{eq:low-arm-sampling-risk}
\end{equation}

\begin{lemma}[The contrast dominates low-arm sampling]
\label{lem:contagion-reduction}
For every graph $G$ and threshold $\theta\geq1$ for which both target
exposures are attainable at every unit,
$\Rallstar(\cP(\cZ);f,a_\theta)\geq\mathsf S_\theta(G)$.
\end{lemma}

\begin{proof}
Restrict the outcome class to tables with \(Y_i(1,C_i\geq\theta)=0\) and
\(Y_i(0,\cdot)=0\) for every \(i\), and \(Y_i(1,C_i\leq\theta-1)=x_i\) for some
\(x\in[-1,1]^V\); restricting the class cannot increase the minimax risk.  For
such tables \(\tau_\theta(Y)=-\frac1n\sum_ix_i\), and under an assignment
\(z\) the realized outcome vector is determined by \(z\) together with
\(x_{L_0(z)}\), the remaining observed entries being known constants.  By
Lemma~\ref{lem:contagion-low-capacity}, \(L_0(z)\in\mathcal L_\theta(G)\) for
every \(z\), and \(z=\1_S\) gives \(L_0(z)=S\) for every
\(S\in\mathcal L_\theta(G)\).  A design $P$ therefore induces a law $Q$ on $\mathcal L_\theta(G)$.
Write $S=L_0(Z)$ for the sampled unit set. On the restricted
tables, the observed data are $(Z,Y_{O_f(Z)})=(Z,h_Z(x_S))$,
where $h_z$ places the sampled values in their low-arm coordinates
and fills every other observed coordinate with zero. For a contrast
estimator $T$ and a sample $S$ with $Q(S)>0$, define
\[
 \widehat\mu(S,x_S)
 =-\sum_{\substack{z:L_0(z)=S\\P(Z=z)>0}}
   P(Z=z\mid L_0(Z)=S)\,T(z,h_z(x_S)).
\]
Define $\widehat\mu$ arbitrarily when $Q(S)=0$. The conditional law of $Z$ given
$L_0(Z)=S$ depends on $P$, and not on $x$. Conditional averaging therefore
gives a sampling estimator with no larger squared-error risk. Conversely,
any estimator of the sampled mean gives a contrast estimator by applying
it to $L_0(z)$ and changing its sign. Thus, for fixed $P$, the restricted
minimax risk equals the sampling risk under its induced $Q$.
Taking the infimum over all laws $Q$ gives the claimed lower bound.
\end{proof}

\subsubsection{Packing lower bound on regular graphs}

For the packing lower bound, fix a law \(Q\) on
\(\mathcal L_\theta(G)\), write \(\pi_v=Q(v\in S)\), and let $\mathsf S_\theta(G;Q)$ denote the
infimum over sampling estimators in equation~\eqref{eq:low-arm-sampling-risk}
with $Q$ fixed. We use Theorem~\ref{thm:blocklb-all}(i) in the following form.  If
\(A_1,\ldots,A_m\subseteq V\) are disjoint sets of common size \(\ell\), take
the potential outcome packing whose \(b\)th block is
\(B_b=\{(i,b_{\rm lo}):i\in A_b\}\), with coefficient \(c_b=\ell/n\).  The
observed set meets \(B_b\) exactly when \(A_b\cap S\neq\varnothing\), so
equation~\eqref{eq:blocklb-all} reads
\begin{equation}
  \mathsf S_\theta(G;Q)
  \geq
  \left(\frac{\ell}{n}\right)^2\sum_{b=1}^m Q(A_b\cap S=\varnothing).
  \label{eq:contagion-block-deficit}
\end{equation}
We choose the blocks according to the low-arm inclusion probabilities.

\begin{lemma}[Blocks in a low-inclusion level set]
\label{lem:contagion-occupancy}
Let \(s\in(0,1/4]\), \(V_s=\{v:\pi_v\leq s\}\), and
\(\ell'=\lceil1/(2s)\rceil\).  If \(|V_s|\geq n/2\) and \(n\geq4\ell'\), then
\(\mathsf S_\theta(G;Q)\geq\ell'/(16n)\).
\end{lemma}

\begin{proof}
Split \(V_s\) into \(m=\lfloor|V_s|/\ell'\rfloor\geq n/(4\ell')\) disjoint
blocks of size \(\ell'\).  For any such block,
\(Q(A\cap S\neq\varnothing)\leq\sum_{i\in A}\pi_i\leq\ell's\leq\frac12+s\leq
\frac34\), so \(Q(A\cap S=\varnothing)\geq\frac14\).  Now use
equation~\eqref{eq:contagion-block-deficit}.
\end{proof}

\begin{lemma}[Blocks inside a neighborhood]
\label{lem:contagion-star-block}
Let \(G\) be any graph, let \(v\in V\), let \(F\subseteq N_v\), and let
\(1\leq\ell\leq|F|/(2\theta)\).  Then some \(A\subseteq F\) with \(|A|=\ell\)
satisfies \(Q(A\cap S=\varnothing)\geq\pi_v/2\).
\end{lemma}

\begin{proof}
Draw \(A\subseteq F\) uniformly among \(\ell\)-subsets, independently of \(S\),
and condition on any \(S\in\mathcal L_\theta(G)\) with \(v\in S\).  Then
\(|S\cap N_v|\leq\theta-1\), so at most \(\theta-1\) units of \(F\) lie in
\(S\), and a union bound over them gives
\[
  \Pr\nolimits_A\{A\cap S\neq\varnothing\mid S\}
  \leq\frac{\ell(\theta-1)}{|F|}
  \leq\frac{\theta-1}{2\theta}\leq\frac12 .
\]
Averaging over \(S\) on the event \(\{v\in S\}\) gives
\(\E_A\,Q\{v\in S,\ A\cap S=\varnothing\}\geq\pi_v/2\), so at least one
realization of \(A\) attains the bound.
\end{proof}

The neighborhood-block construction uses the constraint that $v$ and
any $\theta$ of its neighbors cannot all attain the low arm. For
$\theta\geq2$, every pair of potential outcomes within such a set can
be jointly observed.

\begin{lemma}[Packing neighborhood blocks]
\label{lem:contagion-star-packing}
Let $G$ be $d$-regular, $W\subseteq V$, and $1\leq\ell\leq d/2$.  There are units
\(v_1,\ldots,v_m\in W\) and pairwise disjoint sets \(A_{v_b}\) of size \(\ell\)
with \(A_{v_b}\subseteq N_{v_b}\setminus\bigcup_{b'<b}A_{v_{b'}}\), each drawn
from a free part of \(N_{v_b}\) of size at least \(d/2\), and
\(m\geq|W|/(2\ell)\).  Centers may repeat; the selected blocks do not.
At each step any $\ell$-subset of the free neighborhood may be chosen.
\end{lemma}

\begin{proof}
Choose the blocks greedily.  Writing \(B\) for the union of the blocks chosen
so far, call \(v\) available if \(|N_v\cap B|\leq d/2\).  This is the only
step of the argument that uses regularity.  Since \(G\) is
\(d\)-regular, \(\sum_{v\in V}|N_v\cap B|=d|B|\), so at most \(2|B|\) units of
\(V\) fail to be available.  Hence while \(2|B|<|W|\) some \(v\in W\) is
available, and \(N_v\setminus B\) has at least \(d/2\geq\ell\) elements, so
another block can be added.  The process stops only once
\(m\ell=|B|\geq|W|/2\).
\end{proof}

\begin{remark}[The probability-based packing bounds]
\label{rem:contagion-two-certificates}
Fix a design and write \(\pi_v\) for the probability that unit \(v\) attains the low
arm.  We use two block packings: one for units with small
\(\pi_v\) and one for units with large \(\pi_v\).  Where \(\pi_v\) is small, a block of about \(1/\pi_v\)
units is missed with constant probability by the union bound.  Where \(\pi_v\) is large, take the
block inside the neighborhood \(N_v\) instead: attaining the low arm at \(v\)
forces at most \(\theta-1\) treated neighbors, so a block of \(\ell\leq
d/(4\theta)\) units drawn from \(N_v\) misses the sampled low-arm set with probability
at least \(\pi_v/2\).  Blocks of this second kind can be packed disjointly
because \(\ell\) is a small fraction of \(d\).  The low-inclusion packing gives
\(1/(\pi n)\) and the neighborhood packing gives \(d\pi/n\); they cross at \(\pi\asymp d^{-1/2}\),
giving the $\sqrt d/n$ term in equation~\eqref{eq:contagion-regular-lower}. The separate rank bound gives $1/\operatorname{diss}_{\theta-1}(G)$.
\end{remark}

\begin{theorem}[Lower bound on every regular graph]
\label{thm:contagion-regular-lower}
Fix \(\theta\geq1\).  Let \(G\) be \(d\)-regular on \(n\) units with
\(d\geq16\theta\) and \(n\geq4\sqrt d\).  Then
\begin{equation}
  \Rallstar(\cP(\cZ);f,a_\theta)
  \geq
  c_\theta\max\left\{
    \frac1{\operatorname{diss}_{\theta-1}(G)},\
    \frac{\sqrt d}{n}\right\},
  \label{eq:contagion-regular-lower}
\end{equation}
where \(c_\theta\) depends only on \(\theta\).
\end{theorem}

\begin{proof}[Proof of Theorem~\ref{thm:contagion-regular-lower}]
The rank term is Lemma~\ref{lem:contagion-low-capacity}, so it is enough to
prove \(\mathsf S_\theta(G;Q)\geq\sqrt d/(64\theta n)\) for every \(Q\) and
then apply Lemma~\ref{lem:contagion-reduction}.

Set \(s=t=d^{-1/2}\), and split \(V\) into \(V_s=\{v:\pi_v\leq s\}\) and
\(W_t=\{v:\pi_v>t\}\).  One of the two has at least \(n/2\) elements.  Note
\(s\leq1/4\) because \(d\geq16\).

Suppose first \(|V_s|\geq n/2\).  With \(\ell'=\lceil\sqrt d/2\rceil\) the
hypothesis \(n\geq4\sqrt d\) gives \(n\geq4\ell'\), and
Lemma~\ref{lem:contagion-occupancy} gives
\(\mathsf S_\theta(G;Q)\geq\ell'/(16n)\geq\sqrt d/(32n)\).

Otherwise \(|W_t|\geq n/2\).  Put \(\ell=\lfloor d/(4\theta)\rfloor\), so that
\(d/(8\theta)\leq\ell\leq d/2\).  Run the greedy procedure in Lemma~\ref{lem:contagion-star-packing} with
$W=W_t$.  Whenever it chooses an available center $v_b$, apply
Lemma~\ref{lem:contagion-star-block} to that center's currently free
neighborhood.  Choose the resulting block $A_{v_b}$, for which
$Q(A_{v_b}\cap S=\varnothing)\geq\pi_{v_b}/2>t/2$, and only then delete
that block.  The packing lemma permits this choice at every step and yields
$m\geq|W_t|/(2\ell)$ disjoint blocks, so equation~\eqref{eq:contagion-block-deficit} gives
\[
  \mathsf S_\theta(G;Q)
  \geq\left(\frac\ell n\right)^2\frac{|W_t|}{2\ell}\frac t2
  =\frac{\ell\,t\,|W_t|}{4n^2}
  \geq\frac{(d/8\theta)\,d^{-1/2}\,(n/2)}{4n^2}
  =\frac{\sqrt d}{64\,\theta\,n}.
\]
Taking the smaller of the two constants completes the proof.
\end{proof}

\subsubsection{Bernoulli upper bound}

\begin{theorem}[Upper bound on every regular graph]
\label{thm:contagion-regular-upper}
Fix \(\theta\geq1\) and let \(G\) be \(d\)-regular with \(d\geq2\theta\).
Assign treatment by independent coin flips, \(\Pr(Z_i=1)=p=\theta/d\), and
estimate \(\tau_\theta\) by exposure HT with the two constant exposure probabilities
\[
  \pi^{\rm lo}=p\Pr\{\operatorname{Bin}(d,p)\leq\theta-1\},
  \qquad
  \pi^{\rm hi}=p\Pr\{\operatorname{Bin}(d,p)\geq\theta\}.
\]
Writing \(P_p\) for the Bernoulli
design and
\[
  \widehat\tau_\theta
  =\frac1n\sum_{i=1}^n
  \left\{\frac{\1\{f_i(Z)=b_{\rm hi}\}}{\pi^{\rm hi}}
        -\frac{\1\{f_i(Z)=b_{\rm lo}\}}{\pi^{\rm lo}}\right\}Y_i^{\rm obs}
\]
for the exposure HT estimator, which depends only on the observed data,
the strategy \((P_p,\widehat\tau_\theta)\) is unbiased and
\begin{equation}
  \Runbstar(\cP(\cZ);f,a_\theta)\leq C_\theta\frac dn ,
  \label{eq:contagion-regular-upper}
\end{equation}
where \(C_\theta\) depends only on \(\theta\).
\end{theorem}

\begin{proof}[Proof of Theorem~\ref{thm:contagion-regular-upper}]
By regularity the two exposure probabilities do not depend on the unit.
Let $X\sim\operatorname{Bin}(d,\theta/d)$. Since $d\geq2\theta$,
\[
 \Pr(X=0)=(1-\theta/d)^d\geq e^{-2\theta},
 \qquad
 \Pr(X=\theta)
 =\frac{(d)_\theta}{\theta!}\left(\frac\theta d\right)^\theta
   (1-\theta/d)^{d-\theta}
 \geq\frac{(\theta/2)^\theta}{\theta!}e^{-2\theta}.
\]
The first event lies in the low arm and the second in the high arm.
Thus $c_\theta/d\leq\pi^{\rm lo},\pi^{\rm hi}\leq\theta/d$, so both
exposure probabilities are positive and exposure HT is unbiased.

Write \(R_i\) for the observation indicator of one arm at unit \(i\), let
\(\pi_i=\E R_i\) and \(\Delta_{ij}=\Cov(R_i,R_j)\).  Since
\((a-b)^2\leq2a^2+2b^2\), it suffices to bound each arm separately, and for
\(|y_i|\leq1\),
\[
  \Var\left(\frac1n\sum_i\frac{R_iy_i}{\pi_i}\right)
  =\frac1{n^2}\sum_{i,j}\frac{\Delta_{ij}}{\pi_i\pi_j}y_iy_j
  \leq\frac1{n^2}\sum_{i,j}\frac{|\Delta_{ij}|}{\pi_i\pi_j}.
\]

The indicator \(R_i\) is a function of \(Z\) restricted to \(N_i\cup\{i\}\).
Under a product design this implies $\Delta_{ij}=0$ whenever
\({\rm dist}(i,j)\geq3\), so only the diagonal, the edges, and the
distance-two pairs contribute.

The diagonal contributes
\(\frac1{n^2}\sum_i(1-\pi_i)/\pi_i\leq d/(c_\theta n)\).

For each of the $nd$ ordered adjacent pairs, both units must be
treated to attain the chosen arm together. Independence of their
own treatments gives $0\leq\pi_{ij}\leq p^2$, while \(\pi_i\pi_j\asymp_\theta p^2\).  Hence
\(|\Delta_{ij}|/(\pi_i\pi_j)\leq C_\theta\) and the edge contribution is at most
\(C_\theta nd/n^2=C_\theta d/n\).  For \(\theta=1\) and the low arm this ratio
is exactly one, since two adjacent treated units each have a treated neighbor
and \(\pi_{ij}=0\).

For a nonadjacent pair with $k=|N_i\cap N_j|\geq1$, let
$X\sim\operatorname{Bin}(k,p)$ count shared treated neighbors, and put
\[
 g(x)=\Pr\{x+\operatorname{Bin}(d-k,p)\text{ lies in the chosen arm}\}.
\]
Conditional on $X$, the remaining neighbor counts and the two own-treatment
indicators are independent. Hence
\[
 \pi_{ij}=p^2\E[g(X)^2],\qquad
 \pi_i\pi_j=p^2\{\E[g(X)]\}^2,\qquad
 \Delta_{ij}=p^2\Var(g(X)).
\]
Because $0\leq g\leq1$,
\[
 \Var(g(X))\leq\E[\{g(X)-g(0)\}^2]
 \leq\Pr(X\geq1)\leq kp.
\]
It follows that $|\Delta_{ij}|/(\pi_i\pi_j)\leq C_\theta kp$.
Since
\(\sum_{i\neq j}|N_i\cap N_j|=\sum_vd(d-1)\leq nd^2\), the distance-two
contribution is at most \(\frac1{n^2}C_\theta p\,nd^2=C_\theta\theta\,d/n\).

Adding the three pieces gives \(C_\theta d/n\) for each arm, and hence for the
contrast.
\end{proof}

\begin{resultrestatement}{Theorem~\ref{thm:contagion-main-bracket}}{Threshold-contrast bounds on regular graphs}
Fix $\theta\geq1$. Let $G$ be $d$-regular on $n$ units, with
$d\geq16\theta$ and $n\geq4\sqrt d$. Then
\[
 c_\theta\max\left\{\frac1{\operatorname{diss}_{\theta-1}(G)},\frac{\sqrt d}{n}\right\}
 \leq\Rallstar(\cP(\cZ);f,a_\theta)
 \leq\Runbstar(\cP(\cZ);f,a_\theta)
 \leq C_\theta\frac dn.
\]
The constants depend only on $\theta$. The Bernoulli design with $p=\theta/d$
and exposure HT attain the upper bound.
\end{resultrestatement}

\begin{proof}[Proof of Theorem~\ref{thm:contagion-main-bracket}]
Combine Theorems~\ref{thm:contagion-regular-lower}
and~\ref{thm:contagion-regular-upper}.
\end{proof}

\subsubsection{Pairwise conflicts and community examples}

\begin{lemma}[Pairwise projection for threshold contagion]
\label{lem:contagion-conflict-projection}
For the treated-low threshold arm in equation~\eqref{eq:complex-contagion-target}, the
induced conflict subgraph on that arm is \(G\) when \(\theta=1\) and
is edgeless when \(\theta\geq2\).  Moreover, if \(\theta\geq2\) and every
distinct pair \(i,j\) has at least \(\theta\) neighbors of \(i\) outside
\(N_j\cup\{j\}\), then the induced conflict graph on the potential
outcomes in the two arms \(b_{\rm lo},b_{\rm hi}\) has only the same-unit high/low conflicts. A
sufficient condition is girth at least five and \(d\geq\theta+1\) on a
\(d\)-regular graph.
\end{lemma}
\begin{proof}[Proof of Lemma~\ref{lem:contagion-conflict-projection}]
For \(\theta=1\), the low-arm potential outcomes of units \(i\) and \(j\) are jointly
observable if and only if \(i\) and \(j\) are nonadjacent: both units must be
treated, and adjacency would give each a treated neighbor.
Hence the induced conflict graph on the low arm is \(G\).

For \(\theta\geq2\), any two distinct low-arm potential outcomes are jointly
observable by treating exactly those two units; each sees at most one treated
neighbor.  High-arm potential outcomes are jointly observable by the all-treated
assignment.  Finally, for a high-arm potential outcome of unit \(i\) and a low-arm
potential outcome of a distinct unit \(j\), treat \(i\), \(j\), and \(\theta\)
neighbors of \(i\) outside \(N_j\cup\{j\}\).  Then \(i\) has at least
\(\theta\) treated neighbors, while \(j\) has at most one treated neighbor.
The only remaining pairwise conflict is that the same unit cannot occupy both
threshold arms. For the sufficient condition, girth at least five implies
$|N_i\cap N_j|\leq1$ for nonadjacent $i,j$ and no shared neighbors for
adjacent $i,j$. At most one neighbor of $i$ therefore lies in
$N_j\cup\{j\}$, leaving at least $d-1\geq\theta$ available neighbors.
\end{proof}

\begin{remark}[Low-arm pairwise observability at $\theta=1$]
\label{rem:simple-contagion-conflict-exact}
At $\theta=1$, observing the low arm requires the unit to be treated and
all its neighbors to be controlled. A set of low-arm potential outcomes is therefore
jointly observable exactly when its units form an independent set of $G$.
Thus the conflict graph records all low-arm sampling constraints,
and $r(b_{\rm lo})=\alpha(G)$.
\end{remark}

\begin{proof}[Proof of Corollary~\ref{cor:contagion-conflict-gap}]
Girth at least five and $d\geq\theta+1$ give the pairwise-witness condition,
so Lemma~\ref{lem:contagion-conflict-projection} shows that the induced
conflict graph on the two target arms consists of the $n$ same-unit pairs only. The same-unit conflicts give the $1/n$ singleton lower bound,
as in the proof of Lemma~\ref{lem:contagion-low-capacity}.
The neighborhood packing uses an additional joint observation constraint:
a unit and $\theta$ of its neighbors cannot all attain the low arm.
Theorem~\ref{thm:contagion-regular-lower} instead gives the hypergraph-based
lower bound $c_\theta\sqrt d/n$, which is larger by a factor of order
$\sqrt d$.
\end{proof}

\begin{theorem}[Threshold-contrast rate for clique communities]
\label{thm:clique-community-contagion}
Let $G$ be the disjoint union of $G_0$ cliques of size $m=d+1$, so that
$n=G_0m$, and let $1\leq\theta\leq m-2$.  Then
\[
 \Rallstar(\cP(\cZ);f,a_\theta)
 \asymp_\theta \Runbstar(\cP(\cZ);f,a_\theta)
 \asymp_\theta \frac dn .
\]
The independent clique design that assigns each clique either $\theta$ or
$\theta+1$ treated units attains the upper bound.
\end{theorem}

\begin{proof}
For the lower bound, restrict to the low arm.  In a clique, if \(K\) units are
treated, every treated unit has \(K-1\) treated peers.  Hence at most
\(\theta\) low-arm potential outcomes can be observed in one clique, and this is
attained by treating exactly \(\theta\) units.  Across the \(G_0\) cliques,
\[
  r(b_{\rm lo})=\theta G_0=\theta\frac{n}{d+1}.
\]
If \(r(b_{\rm lo})\leq n/2\), the rank lower bound gives
\[
  \Rallstar(\cP(\cZ);f,a_\theta)
  \gtrsim_\theta
  \frac1{G_0}
  \asymp
  \frac{d}{n}.
\]
If \(r(b_{\rm lo})>n/2\), then \(d+1<2\theta\).  The same-unit high- and
low-arm potential outcomes are mutually exclusive, so the all-estimator singleton
lower bound in equation~\eqref{eq:singleton-all} gives \(\Rallstar\geq 1/n\).  Since \(d\) is then bounded by a constant
depending only on \(\theta\), this is again a lower bound of order \(d/n\).

For the upper bound, independently choose the number of treated units in each clique
as $\theta$ or $\theta+1$, each with probability $1/2$, then sample uniformly
with that number of treated units.  A treated unit has the low exposure under the first choice
and the high exposure under the second, so
\[
 q_i^{\rm lo}=\frac{\theta}{2m},\qquad
 q_i^{\rm hi}=\frac{\theta+1}{2m}.
\]
The exposure HT estimator is unbiased.  Its contribution from one clique
is either $-2m/(n\theta)$ times a sum of $\theta$ bounded low-arm potential outcomes,
or $2m/[n(\theta+1)]$ times a sum of $\theta+1$ bounded high-arm potential outcomes.
Its absolute value is at most $2m/n$.  Independence across cliques gives
\[
 \sup_Y\Var(\widehat\tau_\theta\mid Y)
 \leq G_0\frac{4m^2}{n^2}=\frac4{G_0}
 \asymp\frac dn.
\]
Since \(\Rallstar\leq\Runbstar\), the theorem follows.
\end{proof}

\begin{proposition}[Complete bipartite communities]
\label{prop:kdd-community-contagion}
Let \(G\) be the disjoint union of \(n/(2d)\) copies of \(K_{d,d}\) and fix
\(1\leq\theta\leq d-1\).  For the threshold contrast
equation~\eqref{eq:complex-contagion-target},
\[
  \Rallstar(\cP(\cZ);f,a_\theta)
  \asymp_\theta
  \Runbstar(\cP(\cZ);f,a_\theta)
  \asymp_\theta
  \frac dn .
\]
\end{proposition}

\begin{proof}[Proof of Proposition~\ref{prop:kdd-community-contagion}]
Index the components by \(g\) and write \(L_g,R_g\) for the two sides of the
\(g\)th copy of \(K_{d,d}\), so there are \(2n/(2d)=n/d\) sides in all, each of
size \(d\).

For the lower bound, fix a side, say \(L_g\), and pair its two arm blocks
\[
  B_{g}^{\rm lo}=\{(i,b_{\rm lo}):i\in L_g\},
  \qquad
  B_{g}^{\rm hi}=\{(i,b_{\rm hi}):i\in L_g\},
\]
with coefficients \(-d/n\) and \(d/n\).  Every unit of \(L_g\) has neighborhood
exactly \(R_g\), so all treated units of \(L_g\) have the same number of treated
neighbors, namely the number \(X_R\) of treated units in \(R_g\).  If
\(X_R\leq\theta-1\) they all attain \(b_{\rm lo}\); if \(X_R\geq\theta\) they
all attain \(b_{\rm hi}\).  No assignment therefore meets both blocks, so
\(q_P(B_g^{\rm lo})+q_P(B_g^{\rm hi})\leq1\).

Repeat this construction for both sides of every component.
The resulting $m=n/d$ paired perturbation sets are pairwise disjoint,
and each has per-arm weight $w_b=d/n$.  The common number of treated neighbors on each side
shows that no feasible assignment resolves a pair, so their
resolution capacity is $\mathcal M=0$.
Part (i) of Theorem~\ref{thm:arm-paired-certificate} therefore gives
\[
 \Rallstar(\cP(\cZ);f,a_\theta)
 \geq\sum_{b=1}^{n/d}\left(\frac dn\right)^2
 =\frac dn.
\]
The same theorem gives this lower bound for
$\Runbstar(\cP(\cZ);f,a_\theta)$.

For the upper bound, use a design that observes whole side-arm blocks.  In each
component independently pick one of the four blocks
\(B_{g,L}^{\rm lo},B_{g,L}^{\rm hi},B_{g,R}^{\rm lo},B_{g,R}^{\rm hi}\) with
probability \(1/4\) and assign treatment so that it is observed: treat all of \(L_g\)
together with \(\theta-1\) units of \(R_g\) to observe \(B_{g,L}^{\rm lo}\),
with \(\theta\) units of \(R_g\) to observe \(B_{g,L}^{\rm hi}\), and
symmetrically for \(R_g\). Every target potential outcome has exposure probability at least $1/4$, because selecting its own block reveals that
outcome. The outcome may also be revealed when another block is selected.
Use each target outcome's full exposure probability in the exposure HT
estimator. The HT contribution from one component is bounded in absolute
value by $8d/n$, and the component contributions are independent.
Hence the variance of the HT estimator is at most a constant times
\[
 \frac{n}{2d}\left(\frac{d}{n}\right)^2=O(d/n).
\]
Since \(\Rallstar\leq\Runbstar\), the two bounds match.
\end{proof}

\section{Supporting results for Section~\ref{sec:bipartite-interference}}
\label{app:bipartite-proofs}
We prove the shared-pool bounds in Theorem~\ref{thm:bipartite-central-main},
then the conflict-graph comparison in Corollary~\ref{cor:bipartite-full-projection-gap}.

\subsection{Proofs for the shared-pool bipartite example}
\label{subsec:app-bipartite-fixed-count}
The lower bound limits the number of target potential outcomes observed under any
assignment. For the fixed-count upper bound, we first bound the exposure probabilities, then the covariance operator that controls HT variance.

Fix the setting of Theorem~\ref{thm:bipartite-central-main}.  Throughout this
subsection, $0<k'<k<d$, $\Delta=k-k'$,
\[
 s=\sqrt{\min\{k,d-k\}+1},
 \qquad \Delta\leq c_0s,
\]
and $t$ is a nearest integer to $M(k+k')/(2d)$.  Put $p=t/M$,
\[
 b(j;d,p)=\binom djp^j(1-p)^{d-j},\qquad
 q_k=\max_{0\leq x\leq1}b(k;d,x)=b(k;d,k/d),\qquad
 \sigma^2=dp(1-p).
\]
Unless stated otherwise, constants depend only on $c_0$ and on the fixed
comparison constants in $M\asymp n/s$.

\subsubsection{Binomial probabilities near the target doses}
\begin{lemma}[Comparability near the target doses]
\label{lem:bipartite-near-mode}
Assume $d^2\leq M$.  Then $q_k\asymp1/s$, $\sigma^2\asymp s^2$, and
$b(j;d,p)\asymp q_k$ for $j\in\{k,k'\}$.
\end{lemma}

\begin{proof}
For $1\leq j\leq d-1$, the two-sided Stirling bounds give, with $x=j/d$,
\begin{equation}
 b(j;d,p)\asymp\left\{\frac{d}{j(d-j)}\right\}^{1/2}
 \exp\{-d\,D(x\Vert p)\},
 \qquad
 D(x\Vert p)=x\log\frac xp+(1-x)\log\frac{1-x}{1-p},
 \label{eq:bipartite-stirling}
\end{equation}
with universal constants.  The prefactor is of order
$\{\min(j,d-j)\}^{-1/2}$.  Taking $j=k$ and $p=k/d$ gives $q_k\asymp1/s$.

We next compare $\min\{j,d-j\}+1$ with $s^2$ for $j=k'$.  Let
$K_0=16c_0^2$.  If $\min\{k,d-k\}\geq K_0$, then $\Delta\leq c_0s\leq
\min\{k,d-k\}/2$, so $\min\{k',d-k'\}$ lies between $\min\{k,d-k\}/2$ and
$2\min\{k,d-k\}$.  If $\min\{k,d-k\}<K_0$, then $s^2\leq K_0+1$. The function $j\mapsto\min\{j,d-j\}$ changes by at most $\Delta$ between
$k$ and $k'$. Thus
\[
 1\leq\min\{k',d-k'\}\leq K_0+c_0\sqrt{K_0+1},
\]
and both $s^2$ and $\min\{k',d-k'\}+1$ are bounded above and below
by positive constants depending on $c_0$.  In both cases
$\min\{j,d-j\}+1\asymp s^2$ for $j\in\{k,k'\}$.

The mean $dp$ differs from $(k+k')/2$ by at most $d/(2M)\leq1/(2d)$.  If
$k\leq d/2$, then $dp\geq(k+1)/2-1/(2d)$, since $k'\geq1$, and
$1-p\geq1/2-1/(2d)$; thus $\sigma^2\gtrsim k+1\asymp s^2$.  If $k>d/2$,
the same argument applied to $d-dp\geq d-k+1/2-1/(2d)$ and
$p\geq1/5$, which holds because $k>d/2$ and $k'\geq1$ give $(k+k')/(2d)\geq1/4$, gives $\sigma^2\gtrsim d-k+1\asymp s^2$.  In both cases
$\sigma^2\leq\min\{dp,d(1-p)\}\lesssim s^2$.

Finally, $|j-dp|\leq\Delta/2+1/2\leq\Delta$ for $j\in\{k,k'\}$, and the
inequality $D(x\Vert p)\leq(x-p)^2/\{p(1-p)\}$ gives
$d\,D(j/d\Vert p)\leq(j-dp)^2/\sigma^2\leq C\Delta^2/s^2\leq Cc_0^2$.
The exponent in equation~\eqref{eq:bipartite-stirling} is therefore bounded, and
$b(j;d,p)\asymp\{\min(j,d-j)\}^{-1/2}\asymp1/s\asymp q_k$.
\end{proof}

\subsubsection{The capacity lower bound}
\begin{lemma}[Shared-pool lower bound]
\label{lem:bipartite-capacity-lower}
Under the assumptions of Theorem~\ref{thm:bipartite-central-main}, define
\[
 \mathcal J_k=\max_{T\subseteq U}\#\{v:|A_v\cap T|=k\}.
\]
With probability at least $1-n^{-2}$ over the graph draw,
\begin{equation}
 \mathcal J_k\leq C(nq_k+M+\log n)\lesssim n/s,
 \qquad
 \Rallstar(\cP(\{0,1\}^U);f,a_{k,k'})\gtrsim s/n.
 \label{eq:bipartite-capacity-lower}
\end{equation}
\end{lemma}

\begin{proof}
We first bound the dose-$k$ count uniformly over assignments. Fix
$T\subseteq U$ and put $u=|T|$. Probabilities and expectations in this
calculation are over the independent random neighborhoods $A_v$.
For $0<u<M$, the hypergeometric mass satisfies
\[
 \Pr\{|A_v\cap T|=k\}
 =\binom dk\frac{(u)_k(M-u)_{d-k}}{(M)_d}
 \leq\frac{M^d}{(M)_d}b(k;d,u/M)
 \leq e^2q_k,
\]
where $(a)_r=a(a-1)\cdots(a-r+1)$ and $d^2\leq M$ bounds
$M^d/(M)_d\leq\exp(2d^2/M)\leq e^2$. For $u=0$ or $M$, the mass is zero
because $0<k<d$. Writing $X_T=\#\{v:|A_v\cap T|=k\}$, independence gives
\[
 \E e^{X_T}
 =\{1+(e-1)\Pr(|A_v\cap T|=k)\}^n
 \leq\exp(Cnq_k).
\]
Markov's inequality applied to $e^{X_T}$ therefore gives
\[
 \Pr\{X_T>Cnq_k+x\}\leq e^{-x}.
\]
Take $x=M\log2+2\log n$ and sum over the $2^M$ assignments. With
probability at least $1-n^{-2}$,
\[
 \mathcal J_k\leq C(nq_k+M+\log n)\lesssim n/s.
\]
Here Lemma~\ref{lem:bipartite-near-mode} gives $q_k\asymp1/s$.
Also $s\leq\sqrt d\leq M^{1/4}$ and $M\asymp n/s$ imply
$n\lesssim M^{5/4}$, so $\log n\lesssim M\asymp n/s$.

Fix a graph satisfying the bound on $\mathcal J_k$. For $j\in\{k,k'\}$, let
$N_j(T)=\#\{v:|A_v\cap T|=j\}$ and define
\[
 K_0=\max_{T\subseteq U}\min\{N_k(T),N_{k'}(T)\}.
\]
Since the doses are distinct, $K_0\leq n/2$; also $K_0\leq\mathcal J_k$.
If $K_0\geq1$, partition $V$ into $2K_0$ groups $W_b$, each of size at
least $\lfloor n/(2K_0)\rfloor\geq n/(4K_0)$. Pair the perturbation sets
\[
 B_b^+=\{(v,k):v\in W_b\},\qquad
 B_b^-=\{(v,k'):v\in W_b\},
\]
using positive signs and per-arm weights $w_b=|W_b|/n\geq1/(4K_0)$.
An assignment resolves a pair only if its group contains outcomes observed
at both doses. Disjointness of the groups limits the number of resolved
pairs to $\min\{N_k(T),N_{k'}(T)\}\leq K_0$.
Theorem~\ref{thm:arm-paired-certificate}(ii) gives
\[
 \Rallstar(\cP(\{0,1\}^U);f,a_{k,k'})
 \geq(2K_0-K_0)\left(\frac1{4K_0}\right)^2
 =\frac1{16K_0}\gtrsim s/n.
\]
If $K_0=0$, use one group $W_1=V$. No assignment resolves its pair, so
Theorem~\ref{thm:arm-paired-certificate}(i) gives a lower bound of one,
which also bounds $s/n$ up to the fixed comparison constants.
\end{proof}

\subsubsection{The fixed-count upper bound}
Condition on the realized graph and draw a uniform $t$-subset of treatment
units. We first bound the two exposure probabilities. The covariance
calculation then gives a variance bound uniform over bounded potential
outcome tables.

\begin{lemma}[Exposure probabilities under the fixed-count design]
\label{lem:bipartite-fixed-count-propensity}
Suppose $M\geq16$ and $3\leq d\leq\sqrt M$, and let $T$ be a uniform
$t$-subset of $U$.  For every $d$-set $A$,
\[
 \pi_j=\Pr\{|A\cap T|=j\}
 =\frac{\binom tj\binom{M-t}{d-j}}{\binom Md},
 \qquad j\in\{k,k'\},
\]
and $cq_k\leq\pi_j\leq Cq_k$.
\end{lemma}

\begin{proof}
Compare the hypergeometric mass with the binomial mass:
\[
 \frac{\pi_j}{b(j;d,p)}
 =\frac{(t)_j}{t^j}\,
   \frac{(M-t)_{d-j}}{(M-t)^{d-j}}\,
   \frac{M^d}{(M)_d},
\]
The first two factors are at most one, and the last lies between one
and $\exp(2d^2/M)\leq e^2$. To bound the first two factors below, note that
$M\geq d^2$, $k'\geq1$, and $k\leq d-1$ give
\[
 t\geq\frac{M(k+k')}{2d}-\frac12
 \geq\frac{d(k+1)}2-\frac12\geq\frac{dk}{2},
\]
and
\[
 M-t\geq\frac{M(2d-k-k')}{2d}-\frac12
 \geq\frac{d(d-k'+1)}2-\frac12\geq\frac{d(d-k')}{2}.
\]
Thus, for either $(a,r)=(t,j)$ or $(M-t,d-j)$,
\[
 \frac ra\leq\frac2d\leq\frac23,
 \qquad \frac{r^2}{a}\leq2.
\]
Both falling products are positive. Using $\log(1-x)\geq-3x$ for
$0\leq x\leq2/3$ gives
\[
 \log\frac{(a)_r}{a^r}
 =\sum_{h=0}^{r-1}\log(1-h/a)
 \geq-\frac{3r(r-1)}{2a}\geq-3.
\]
Both factors are therefore at least $e^{-3}$, so the ratio
$\pi_j/b(j;d,p)$ is bounded above and below by positive universal constants.
Lemma~\ref{lem:bipartite-near-mode} now gives $\pi_j\asymp q_k$.
\end{proof}

\paragraph*{Covariance operators.}
For a fixed dose $j$, the HT error for its arm mean is an outcome-weighted
sum of centered observation indicators. We regard each centered indicator
as a function of the random treatment set $T$. The next lemma bounds the
covariance operator obtained by averaging over a random neighborhood;
Lemma~\ref{lem:bipartite-sampled-operator} transfers the bound to the
$n$ neighborhoods in the realized graph.

\begin{lemma}[Population covariance operator on a slice]
\label{lem:bipartite-slice-operator}
Let \(M\geq16\), \(2\leq d\leq\sqrt M\), and \(0<t<M\).  Fix a dose \(j\)
with $\pi_j=\Pr_T(|A\cap T|=j)\in(0,1)$ for any fixed $d$-set $A$.
On the space of real functions of uniform $t$-subsets, use the inner product
$\langle h_1,h_2\rangle=\E_T[h_1(T)h_2(T)]$ and define
\[
 g_A(T)=\frac{\1\{|A\cap T|=j\}}{\pi_j}-1,
 \qquad L_j=\|g_A\|_2^2=\pi_j^{-1}-1,
 \qquad \Sigma_j=\E_A[g_A\otimes g_A],
\]
where \(A\) is a uniform \(d\)-subset and
\((g\otimes g)h=\langle g,h\rangle g\).  Then, with \(p=t/M\),
\begin{equation}
 \|\Sigma_j\|_{\rm op}
 \leq\max\left\{
 \frac{(j-dp)^2}{d(M-d)p(1-p)},\
 \frac{L_j}{\binom M2-M}
 \right\}.
 \label{eq:bipartite-slice-operator}
\end{equation}
Under Lemma~\ref{lem:bipartite-fixed-count-propensity}, for
$j\in\{k,k'\}$,
\[
 \|\Sigma_j\|_{\rm op}\leq
 C\left\{\frac{\Delta^2}{s^2M}+\frac{s}{M^2}\right\}\leq\frac CM,
 \qquad L_j\leq Cs.
\]
\end{lemma}

\begin{proof}
Let $F$ map functions of uniform $d$-subsets to functions of uniform
$t$-subsets by
\[
 (Fh)(T)=\E_A[g_A(T)h(A)],\qquad
 (F^*h)(A)=\E_T[g_A(T)h(T)].
\]
Then $FF^*=\Sigma_j$, while $F^*F$ has kernel
$K(A,B)=\langle g_A,g_B\rangle$. These two operators have the same nonzero
eigenvalues. The kernel depends only on \(|A\cap B|\), so it belongs to the
Johnson association scheme.  Its orthogonal degree spaces have dimensions
\[
 D_r=\binom Mr-\binom M{r-1},\qquad 0\leq r\leq d,
\]
and the operator is scalar on each such space
\citep[Theorem~4.1 and Lemma~4.3]{Filmus2016Slice}.  Since \(\E_Ag_A(T)=0\) for every
fixed-size \(T\), the degree-zero eigenvalue is zero.  Writing $\lambda_r$ for the eigenvalue on degree $r$, positivity and the
trace identity give
\[
 \lambda_r\geq0,\qquad
 \sum_{r=1}^d D_r\lambda_r=\E_A\|g_A\|_2^2=L_j.
\]
For \(2\leq r\leq d\leq\sqrt M\), \(D_r\geq D_2=\binom M2-M\).
Indeed
\(D_{r+1}/D_r=(M-2r-1)(M-r+1)/\{(r+1)(M-2r+1)\}\geq1\)
through the needed range.  Thus every eigenvalue of degree at least two is
at most \(L_j/D_2\).

To compute the degree-one eigenvalue, put
\(x_u(A)=\1\{u\in A\}-d/M\).  For fixed $T$, conditional on $|A\cap T|=j$, the inclusion probability
of $u$ in $A$ is $j/t$ when $u\in T$ and $(d-j)/(M-t)$ otherwise. Thus
\[
 \E_A[g_A(T)x_u(A)]
 =\frac{j-dp}{Mp(1-p)}(\1\{u\in T\}-p).
\]
For fixed $A$, conditional on dose $j$, the treatment probability of $u$
is $j/d$ when $u\in A$ and $(t-j)/(M-d)$ otherwise. Hence
\[
 \E_T[g_A(T)(\1\{u\in T\}-p)]
 =\frac{(j-dp)M}{d(M-d)}x_u(A).
\]
Composing these identities gives
\[
 (F^*F)x_u=\lambda_1x_u,\qquad
 \lambda_1=\frac{(j-dp)^2}{d(M-d)p(1-p)}.
\]
The functions $x_u$ span the degree-one space. Combining this eigenvalue
with the bound $\lambda_r\leq L_j/D_2$ for $r\geq2$ proves
equation~\eqref{eq:bipartite-slice-operator}.  For $j\in\{k,k'\}$,
Lemma~\ref{lem:bipartite-near-mode} gives $|j-dp|\leq\Delta$ and
$dp(1-p)=\sigma^2\asymp s^2$, so the degree-one eigenvalue is at most
$C\Delta^2/(s^2M)\leq Cc_0^2/M$.  Lemma~\ref{lem:bipartite-fixed-count-propensity}
gives $L_j\leq C/q_k\leq Cs$, so the remaining eigenvalues are at most
$Cs/M^2\leq C/M$, since $s\leq\sqrt d\leq M$.
\end{proof}

\begin{lemma}[Sampling the covariance operator]
\label{lem:bipartite-sampled-operator}
Let \(g_1,\ldots,g_n\) be independent vectors in a finite-dimensional real
Hilbert space with \(\|g_v\|^2=L>0\) and
a common second-moment operator $\Sigma=\E[g_v\otimes g_v]$ for every $v$,
where $\|\Sigma\|_{\rm op}\leq\gamma$ and $\gamma>0$.
For \(0<\varepsilon<1\), with probability at least \(1-\varepsilon\),
\begin{equation}
 \left\|\frac1n\sum_vg_v\otimes g_v\right\|_{\rm op}
 \leq C\left\{\gamma+
 \frac{L}{n}\log\frac{8(1+L/\gamma)}{\varepsilon}\right\}.
 \label{eq:bipartite-sampled-operator}
\end{equation}
\end{lemma}

\begin{proof}
Write \(X_v=g_v\otimes g_v\).  Then
\(\lambda_{\max}(X_v-\Sigma)\leq L\) and
\(\E(X_v-\Sigma)^2=L\Sigma-\Sigma^2\preceq L\Sigma\).
To apply matrix Bernstein with a variance bound of norm $nL\gamma$, put
\[
 Z_v=\operatorname{diag}(X_v-\Sigma,0),\qquad
 V=\operatorname{diag}(nL\Sigma,nL\gamma).
\]
The $Z_v$ are independent and centered, $\lambda_{\max}(Z_v)\leq L$, and
\[
 \sum_v\E Z_v^2\preceq V,\qquad
 \|V\|_{\rm op}=nL\gamma,\qquad
 \frac{\operatorname{tr}V}{\|V\|_{\rm op}}=1+\frac L\gamma.
\]
Set $\ell=\log\{8(1+L/\gamma)/\varepsilon\}$.
The intrinsic-dimension matrix Bernstein inequality
\citep[Theorem~7.7.1]{Tropp2015Matrix} gives, with probability at least
$1-\varepsilon$,
\[
 \lambda_{\max}\left(\sum_v(X_v-\Sigma)\right)
 \leq C\{\sqrt{nL\gamma\ell}+L\ell\}.
\]
The constant $C$ can be chosen so the deviation threshold also exceeds
$\sqrt{nL\gamma}+L/3$, as required by that theorem.
Divide by $n$, add $\|\Sigma\|_{\rm op}\leq\gamma$, and use
$\sqrt{L\gamma\ell/n}\leq(\gamma+L\ell/n)/2$ to obtain
equation~\eqref{eq:bipartite-sampled-operator}.
\end{proof}

\begin{resultrestatement}{Theorem~\ref{thm:bipartite-central-main}}{Nearby interior doses on a shared pool}
Assume a uniform $d$-out bipartite interference graph, with $n\geq3$,
$M\geq16$, and $3\leq d\leq\sqrt M$. Let $0<k'<k<d$ with
$k-k'\leq c_0s$ for fixed $c_0$, and suppose $M\asymp n/s$.
With probability at least $1-2/n^2$ over the graph draw,
\[
 c\frac sn
 \leq\Rallstar(\cP(\{0,1\}^U);f,a_{k,k'})
 \leq\Runbstar(\cP(\{0,1\}^U);f,a_{k,k'})
 \leq C\frac{s\log n}{n}.
\]
The constants depend only on $c_0$ and the comparison constants in
$M\asymp n/s$.
\end{resultrestatement}

\begin{proof}[Proof of Theorem~\ref{thm:bipartite-central-main}]
Lemma~\ref{lem:bipartite-capacity-lower} proves the lower bound with
probability at least $1-n^{-2}$. We prove the upper bound using the
fixed-count design and exposure HT.

Fix $0<\varepsilon<1$. For each $j\in\{k,k'\}$, apply
Lemma~\ref{lem:bipartite-sampled-operator} to the functions $g_{A_v}$ of
Lemma~\ref{lem:bipartite-slice-operator}. The neighborhoods $A_v$ are
independent, and these functions have common squared norm $L_j\leq Cs$
and second-moment operator $\Sigma_j$ with $\|\Sigma_j\|_{\rm op}\leq C/M$.
Take $\gamma=C/M$ and failure probability $\varepsilon/2$ for each arm.
Since $L_j/\gamma\leq CsM\leq CM^2$, equation~\eqref{eq:bipartite-sampled-operator}
gives, simultaneously for both arms,
\begin{equation}
 \|H_j\|_{\rm op}\leq C\left\{\frac1M+\frac{s\ell_\varepsilon}{n}\right\},
 \qquad
 H_j=\frac1n\sum_v g_{A_v}\otimes g_{A_v},\qquad
 \ell_\varepsilon=\log\frac{8M}{\varepsilon},
 \label{eq:bipartite-empirical-operator}
\end{equation}
with probability at least $1-\varepsilon$ over the graph draw.

Condition on a graph satisfying equation~\eqref{eq:bipartite-empirical-operator}.
For dose $j$, write $y_v=Y_v(j)$ and define
\[
 F_jy=\frac1{\sqrt n}\sum_v y_vg_{A_v}.
\]
Then $F_jF_j^*=H_j$, so $\|F_j\|_{\rm op}^2=\|H_j\|_{\rm op}$.
The HT error for the dose-$j$ mean is the function
$n^{-1}\sum_vy_vg_{A_v}(T)$. Its variance under the fixed-count design is
\begin{equation}
 \Var_T\left(\frac1n\sum_v
       \frac{\1\{|A_v\cap T|=j\}y_v}{\pi_j}\right)
 =\frac1n\|F_jy\|_2^2
 \leq\frac{\|y\|_2^2}{n}\|H_j\|_{\rm op}
 \leq\|H_j\|_{\rm op},
 \label{eq:bipartite-operator-variance}
\end{equation}
since $|y_v|\leq1$. The contrast estimator is
\[
 \widehat\tau_{\rm fix}
 =\frac1n\sum_v\left\{
 \frac{\1\{|A_v\cap T|=k\}}{\pi_k}
 -\frac{\1\{|A_v\cap T|=k'\}}{\pi_{k'}}\right\}Y_v^{\rm obs}.
\]
Lemma~\ref{lem:bipartite-fixed-count-propensity} gives positive exposure probabilities, so the estimator is unbiased. The variance of the difference
is at most twice the sum of the arm variances. Combining
equations~\eqref{eq:bipartite-empirical-operator}
and~\eqref{eq:bipartite-operator-variance} therefore gives
\[
 \sup_{Y\in\cY(f)}\Var_T(\widehat\tau_{\rm fix}\mid Y)
 \leq C\left\{\frac1M+\frac{s\ell_\varepsilon}{n}\right\}.
\]
The operator event depends on the graph, and the bound holds for every
bounded outcome table on that graph.

Take $\varepsilon=n^{-2}$. Since $M\asymp n/s\lesssim n$,
$\ell_\varepsilon=O(\log n)$ and $1/M\asymp s/n$, yielding the upper bound
$Cs\log n/n$. Intersecting the operator-concentration event with the
lower-bound event from Lemma~\ref{lem:bipartite-capacity-lower}
gives probability at least
$1-2/n^2$ and completes the proof.
\end{proof}

\subsubsection{Comparison with the conflict graph}

\begin{proof}[Proof of Corollary~\ref{cor:bipartite-full-projection-gap}]
We first determine the pairwise conflicts, then compare the experiment
that permits all independent sets of the conflict graph with the actual
shared-pool experiment.
Throughout this comparison, the graph is the induced conflict graph
on the two target-dose arms, containing $Y_v(k)$ and $Y_v(k')$ for each
outcome unit $v$.

Let $A,B$ be two $d$-sets with $|A\cap B|=r$. To give $A$ dose $k$ and
$B$ dose $k'\leq k$, an assignment must treat $x$ shared units, $k-x$ units
of $A\setminus B$, and $k'-x$ units of $B\setminus A$. Such an assignment
exists exactly when an integer $x$ satisfies
\begin{equation}
 \max\{0,k-d+r\}\leq x\leq\min\{k',r\}.
 \label{eq:bipartite-pairwise-feasibility}
\end{equation}
Since $0\leq k'\leq k\leq d$, equation~\eqref{eq:bipartite-pairwise-feasibility}
is equivalent to $k-k'\leq d-r$. Thus any same-dose pair of potential
outcomes is jointly observable. For adjacent doses, a cross-dose pair can
be jointly observed exactly when $r<d$, that is, when the neighborhoods are distinct. The two target potential outcomes of
the same unit have a pairwise conflict because they require different doses in
the same neighborhood.

Two independently drawn neighborhoods are identical with probability
$1/\binom Md$. A union bound therefore gives
\begin{equation}
 \Pr\{A_v\ne A_w\text{ for all }v\ne w\}
 \geq1-\frac{n(n-1)}{2\binom Md}.
 \label{eq:bipartite-distinct-neighborhoods}
\end{equation}
On this event, this induced conflict graph is the matching formed by the $n$
same-unit pairs of target potential outcomes.

Consider the sampling experiment that permits every independent set of
this matching. Independently for each unit, choose one of its two target
potential outcomes with probability $1/2$. The corresponding HT estimator is
\[
 \widehat\tau_{\rm pair}
 =\frac2n\sum_v\{\1\{J_v=k\}Y_v(k)-\1\{J_v=k'\}Y_v(k')\},
\]
where $J_v$ is the sampled dose label. The estimator is unbiased and has
variance at most $4/n$, since the unit contributions are independent and
bounded in absolute value by $2/n$. For the lower bound, put independent
Rademacher signs on all target potential outcomes. Every independent set leaves at
least one sign unobserved per unit. Each missing sign contributes $1/n^2$
to the posterior variance of the contrast, so the Bayes risk is at least
$1/n$. Hence both minimax risks in this sampling experiment are of order
$1/n$.

For the shared-pool experiment, apply
Theorem~\ref{thm:bipartite-central-main} with $\Delta=1$ and $c_0=1$.
With probability at least
\[
 1-\frac2{n^2}-\frac{n(n-1)}{2\binom Md},
\]
the induced conflict graph on the two target-dose arms is a matching and the minimax risks satisfy
\[
 c\frac sn\leq\Rallstar(\cP(\{0,1\}^U);f,a_{k,k'})
 \leq\Runbstar(\cP(\{0,1\}^U);f,a_{k,k'})
 \leq C\frac{s\log n}{n}.
\]
The condition $\binom Md/n^2\to\infty$ makes this probability tend to one.
\end{proof}

\end{document}